

\ifx\macrosloaded\relax    \else \let\macrosloaded\relax \fi


\ifx\formloaded\undefined 

\hsize=16cm
\vsize=24.5 cm \advance \vsize -24pt 

\fi

\catcode`@=11


\def\letnext{\let\next= }
\def\defnext{\def\next}   

\newif\ifamstex {\defnext{AmS-TeX} \ifx\fmtname\next \global\amstextrue \fi}

\def\reset#1{\global#1=0 }
\def\useouter #1 {\begingroup
                  \expandafter\let\csname#1\endcsname=\relax 
                  \afterassignment\endgroup\global}
\def\?#1\then#2\else#3\fi{#1\defnext{#2}\else\defnext{#3}\fi \next}
\newif\ifnot
\newif\ifnomessage \nomessagetrue 
\let\domessage=\message \def\message#1{{\nomessagefalse\domessage{#1}}}



\font\boldit=cmbxti10 
\font\ssq=cmssq8 \font\ssqi=cmssqi8 
\font\tenss=cmss10
\font\tenssit=cmssi10
\font\tenssbf=cmssbx10

\def\sans
{\def\rm{\fam0\tenss}\def\it{\fam\itfam\tenssit}\def\bf{\fam\bffam\tenssbf}\rm
}


\font\ninerm=cmr9 \font\nineit=cmti9 \font\ninebf=cmbx9 \font\ninesl=cmsl9
\font\ninett=cmtt9 \font\ninei=cmmi9 \font\ninesy=cmsy9
\font\niness=cmss9 \font\ninessit=cmssi9 \font\ninessbf=cmssbx10 at 9 pt
\skewchar\ninei='177 \skewchar\ninesy='60

\def\ninept
{\def\rm{\fam0\ninerm}\def\it{\fam\itfam\nineit}\def\bf{\fam\bffam\ninebf}%
 \def\sl{\fam\slfam\ninesl}\def\tt{\fam\ttfam\ninett}\rm
 \def\sans{\def\rm{\fam0\niness}\def\it{\fam\itfam\ninessit}%
           \def\sl{\fam\slfam\ninessit}\def\bf{\fam\bffam\ninessbf}\rm}%
 \textfont0\ninerm \textfont1\ninei \textfont2=\ninesy
 \textfont\itfam\nineit \textfont\slfam\ninesl
 \textfont\bffam\ninebf \textfont\ttfam\ninett
 \normalbaselineskip=11pt \normalbaselines
 \setbox\strutbox=\hbox{\vrule height 8 pt depth 3 pt width 0 pt}%
}


\newif\iftoc 
\def\writetoc{\write\tocfile}\def\tocfile{-1 }
\useouter newwrite
\def\opentoc{\newwrite\tocfile\openout\tocfile=\jobname.toc }
\edef\tocheader{\string\NotocSection\space Table of contents.
  \string\noindent\string\medskip}
\def\maketoc{\opentoc\writetoc{\tocheader}\toctrue}

\def\str#1{\string#1 }
{\catcode`@=3 
 \gdef\passord#1$#2\end{#1\ifx@#2@\else\getmath#2\end\fi}
 \gdef\getmath#1${$\noex#1\endmath$\passord}
 \gdef\noex#1\endmath{\ifx@#1@\else\donoex#1\endmath\fi}
 \gdef\donoex#1{\noexpand#1\noex}
 \gdef\nomathexp#1{\passord#1$\end}
}

\newtoks\marktoks

\def\putmark#1{\mark
  {\count1=\the\secno \count2=\the\subsecno \count3=\the\proclno
   \marktoks={\ignorespaces\setchap\space#1}}}
\def\setchap{\ifnum\secno<0 Appendix \lastlabel.\else\lastlabel\fi}
\def\puttoc#1{\toks0={\nomathexp{#1} \str\onpage\folio.}\edef\next
                      {{\str\tocitem \lastlabel=\the\toks0}}\writetoc\next}
\def\setmark#1{\putmark{#1}\puttoc{#1}}

\output=
{\everyhbox{}
 \firstmark \headline={\tenit \ifodd\pageno\hfill\fi \the\marktoks \hfil }
 \let~\relax 
 \plainoutput
}

\newwrite\labfile
\def\tolab{\write\labfile}
\def\writelab{\openout\labfile=\jobname.lab
              \global\let\writelab=\tolab \tolab}
\def\chklab#1\lastlabel#2.{\def\next{#2}\ifx\next\empty}
\def\deflab#1{\chklab#1\lastlabel.
                \immediate\writelab{\string\def\string#1{#1}}
              \fi} 
\def\savecount#1{\writelab{#1=\number#1}}
\toks0={\savecount\secno\savecount\subsecno\savecount\proclno
        \savecount\eqnumber
        \savecount\pageno\writelab{\string\advancepageno}}
\toks2=\expandafter{\bye}
\edef\bye{\the\toks0 \the\toks2}

\newwrite\indexfile
\def\toindex{\write\indexfile}
\def\writeindex{\openout\indexfile=\jobname.inx
                \global\let\writeindex=\toindex \toindex}
\def\inx#1{\writeindex{\nomathexp{#1} @\folio.}}
\def\index#1{\inx{#1}#1}
\def\closeindex{\vfil\eject 
                \ifx\writeindex\toindex \immediate\closeout\indexfile \fi}


\def\vfootnote#1%
{\insert\footins\bgroup \ninept
 \interlinepenalty=\interfootnotelinepenalty
 \splittopskip=\ht\strutbox 
 \splitmaxdepth=\dp\strutbox 
 \floatingpenalty=\@MM 
 \leftskip=1pc \rightskip=1pc \spaceskip=0pt \xspaceskip=\z@skip
 \parindent=\z@ \parskip=3pt minus 1pt \parfillskip=\z@ plus 1 fil
 \textindent{#1}\footstrut\futurelet\next\f@@t
}

\newcount\secno \newcount\subsecno \newcount\proclno \newcount\subsubsecno
\newcount\itemno
\reset\secno \reset\subsecno \reset\proclno \reset\itemno
\newdimen\itemindent \itemindent=\parindent 

\global\subsubsecno=-1 

\newif\ifnewsec 

\let\Sec\S 
\useouter beginsection
 \outer\def\newsection #1.
  {\ifnewsec \endgroup \fi 
   \global\advance\secno by 1
   \ifnum-1<\subsecno \reset\subsecno \fi 
   \ifnum-1<\subsubsecno \reset\subsubsecno \fi 
   \reset\proclno
   \edef\lastlabel{\number\secno}%
   {\parskip0pt\beginsection\ifnomessage\Sec\fi\lastlabel. #1.\par
    \par}
   \setmark{#1}
   \startsec 
  }
\def\subsection #1.
  {\secbigbreak
   \global\advance\subsecno by 1
   \ifnum-1<\subsubsecno \reset\subsubsecno \fi 
   \reset\proclno
   \edef\lastlabel{\number\secno.\number\subsecno}%
    \leftline{\bf\lastlabel.\enspace
              \defnext{#1}\ifx\next\empty\else \boldit#1.\fi}\setmark{#1}%
    \nobreak\smallskip\startsec}
\def\subsubsection #1.
  {\secmedbreak
   \global\advance\subsubsecno by 1 \reset\proclno
   \edef\lastlabel{\number\secno.\number\subsecno.\number\subsubsecno}%
    \leftline{\bf\lastlabel.\enspace
              \defnext{#1}\ifx\next\empty\else \it#1.\fi}\setmark{#1}%
    \nobreak\smallskip\startsec}
\useouter beginsection
 \outer\def\Appendix #1.
  {\ifnum\secno>0 \startapx\fi
   \global\advance\secno by -1
   \ifnum-1<\subsecno \reset\subsecno \fi 
   \ifnum-1<\subsubsecno \reset\subsubsecno \fi 
   \reset\proclno
   \edef\lastlabel
    {\ifcase -\secno \or A\or B\or C\or D\or E\or F\or G\or H\or I\or J\fi}%
   {\parskip0pt\beginsection Appendix \lastlabel. #1.\par\par}%
   \setmark{#1}\startsec}
\def\startapx{\ifx\writetoc\totoc\totoc{\str\Appendices}\fi \secno=0 }
\useouter beginsection
 \outer\def\Supplement #1.
  {\let\lastlabel=\empty
   \global\ifnum\secno<0 \secno=1000 \else\advance\secno by 1 \fi
   \beginsection #1.\par\par\setmark{#1}\startsec}
\useouter beginsection
 \outer\def\NotocSection #1.
  {\let\lastlabel=\empty \reset\secno
   \beginsection #1.\par\par\putmark{#1}\startsec}

\def\sticksection 
{\vbox\bgroup\parskip=0pt \everypar{\egroup\noindent}}

\def\startsec
  {\begingroup \newsectrue\parskip0pt\parindent0pt \everypar{\endgroup}}
\def\secbigbreak
  {\ifnewsec \endgroup
     \ifdim\lastskip<\bigskipamount \removelastskip\bigskip \fi 
   \else \bigbreak
   \fi}
\def\secmedbreak
  {\ifnewsec \endgroup
     \ifdim\lastskip<\medskipamount \removelastskip\medskip \fi 
   \else \medbreak
   \fi}
\def\labelsec#1{\global\let#1=\lastlabel \deflab#1\nobreak}


{\let\x=\expandafter
 \def~#1{\defnext##1{\gdef#1{\ifnomessage\leavevmode##1\fi}}\x\next\x{#1}}
 ~\copyright
 \def~#1{\defnext##1{\gdef#1####1{\ifnomessage\leavevmode##1\else####1\fi}}%
         \x\next\x{#1{##1}}}
 ~\` ~\' ~\v ~\^^_ ~\u ~\^^S ~\^ ~\^^D ~\. ~\H ~\~ ~\" 
 \ifamstex ~\B \else ~\= ~\t \fi
}

\def\quotation
{\begingroup\everypar{}\parfillskip=0pt
 \leftskip=0pt plus \hsize \def\par{\ifhmode\/\endgraf\nobreak\fi}\obeylines
 \let\it=\ssq\let\rm=\ssqi\rm 
 \baselineskip=10pt\parskip=0pt
 \def\author##1, ##2.{\vskip2pt\it ##1, \rm ##2\par\smallskip\endgroup}
}

\newcount\eqnumber \eqnumber=0
\def\neqn(#1#2){\global\advance\eqnumber by 1
                (\ifcat a\noexpand#1%
                        \errmessage{Incorrect label #1#2}\number\eqnumber
                 \else\xdef#1{\number\eqnumber}\deflab#1%
                      #1\fam\z@#2
                 \fi)}
\def\label{\ifinner\else \eqno \fi \neqn}
\def\nn{\label(\lastlabel)}
\newcount\acc
\def\add#1#2{\acc=#1\advance\acc#2\relax\number\acc}

\def\heading "#1"{\smallbreak\noindent{\sl #1}.\quad\ignorespaces}

\let\proclaim=\relax
\ifamstex
\let\amsproclaim=\innerproclaim@
\def\innerproclaim@#1. #2#3\par
  {\secmedbreak \global\advance\proclno by 1 \reset\itemno
   \edef\lastlabel{\ifnum-1<\secno    \number\secno.\fi
                   \ifnum-1<\subsecno \number\subsecno.\fi
                   \ifnum-1<\subsubsecno \number\subsubsecno.\fi
                   \number\proclno}
   {\interlinepenalty=300 \parskip=0pt\def\par{\endgraf\penalty200 }
    \noindent 
    \nottrue 
    \if\let\noexpand#2\ifx\undefined#2\notfalse\fi\fi
    \ifnot\amsproclaim{\lastlabel.\enspace#1}#2#3\innerendproclaim@
    \else\global\let#2\lastlabel \deflab#2%
      \amsproclaim{\lastlabel.\enspace#1}#3\innerendproclaim@
    \fi
   }}
\else
\outer\def\proclaim#1. #2#3\par
  {\secmedbreak \global\advance\proclno by 1 \reset\itemno
   \edef\lastlabel{\ifnum-1<\secno    \number\secno.\fi
                   \ifnum-1<\subsecno \number\subsecno.\fi
                   \ifnum-1<\subsubsecno \number\subsubsecno.\fi
                   \number\proclno}
   {\interlinepenalty=300 \parskip=0pt\def\par{\endgraf\penalty200 }
    \noindent 
    \if\let\noexpand#2
     \ifx\undefined#2
      \global\let#2\lastlabel \deflab#2\letnext\ignorespaces
     \else\defnext{#2}\fi
    \else\defnext{#2}\fi
    \bf\lastlabel.\enspace#1.\enspace \sl\next#3\endgraf
    }%
   \ifdim \lastskip<\medskipamount \removelastskip\penalty55\medskip \fi
   }
\fi

\def\dueto#1.{\unskip~[#1].}

\def\newpar{\ifhmode\par\fi{\parskip0pt\noindent}}
\def\enditem{\newpar\reset\itemno}
\def\item{\doindent\itemindent}

\def\doindent#1#2{\newpar
  \hangindent#1\hangafter0 \llap{#2\enspace}\ignorespaces}

\def\statitemnr#1{(\number#1)}
\def\statitem{\global\advance\itemno by 1 \item{\statitemnr\itemno}}
\def\eqitemnr#1{(\romannumeral#1)}
\def\eqitem{\global\advance\itemno by 1 \item{\eqitemnr\itemno}}
\def\defitemnr#1{\count@=#1\advance\count@96\char\count@}
\def\defitem{\global\advance\itemno by 1 \item{\defitemnr\itemno.}}

\def\proof{\heading "Proof"} \def\proofof#1{\heading "Proof of #1"}
\def\QED{\dimen0=.5em\advance\dimen0\wd\QEDbox\advance\dimen0.4pt
         \unskip\kern\dimen0\vadjust{\vskip-\baselineskip}\unhcopy\tightbox
         \copy\QEDbox{\pretolerance=1000\parfillskip=0.4pt\smallbreak}}
\def\QEDQED{\unskip\nobreak\enskip\kern\wd\QEDbox
            {\setbox\QEDbox=\hbox{\copy\QEDbox\enspace\copy\QEDbox}\QED}}

\def\emph#1{{\it #1}\futurelet\next\doitcor}
\def\newterm#1{{\sl\index{#1}}\futurelet\next\doitcor}
\def\foreign#1{{\sl #1}\futurelet\next\doitcor}
\def\doitcor{\if.\next\else\if,\next\else\/\fi\fi}
\def\needfil{\penalty-100\vfil\penalty1500\vfilneg}
\def\bigdisplay{\ifmmode \errmessage{Display begun in math mode.}
                \else\ifhmode\unskip\vadjust{\needfil}%
                     \else\needfil\noindent\fi$$\fi}
\def\frame#1{\leavevmode\hbox{\vrule\vtop{\vbox{\hrule\kern2pt
     \hbox{\strut\kern2pt#1\kern2pt}}\kern2pt\hrule}\vrule}} 
\def\clap#1{\hbox to 0pt {\hss#1\hss}}
\long\def\centerpar#1
{\begingroup
  \advance\leftskip by 0pt plus 1 fil minus 3pt \rightskip=\leftskip
	\parskip=0pt \parindent=0pt \parfillskip=0pt \linepenalty=100
  \noindent\ignorespaces#1\par
 \endgroup
}


\def\setbaselineskip{\par\afterassignment\setbskip\skip0 =}
\def\setbskip{\null\baselineskip=\skip0 \vskip-\baselineskip}
\def\forceeven{\eject\ifodd\pageno \null\mark{\marktoks{}}\vfil\eject \fi}
\def\forceodd{\eject\ifodd\pageno\else\null\mark{\marktoks{}}\vfil\eject\fi}
\def\topglue{\eject\afterassignment\topglu\skip0 }
\def\topglu{\hrule height0pt\vskip\skip0\vskip-\topskip
     \dimen0=\baselineskip\advance\dimen0-\topskip \prevdepth=\dimen0 }


\def\qdef#1{\?\ifx #1\undefined \then \def #1%
              \else \toks0={#1}\afterassignment\chkdef\def\next
              \fi}
\def\qqdef#1{%
\?\ifx #1\undefined
  \then \message{Patching definition:\string#1.}\let\Qdef=\qdef\def#1%
  \else \toks0={#1}\afterassignment\chkdef\def\next
  \fi}
\let\Qdef=\qqdef
\def\chkdef{\expandafter\let\expandafter\curcs\the\toks0 
  \ifx\next\curcs\else
  \message{Altered definition:\the\toks0= \meaning\next,
                              should have been \meaning\curcs.}\fi}
\ifamstex

\def\mathtex#{\quad\hbox}
\def\mathttex#{\qquad\hbox}
\else
\def\text#1{\quad\hbox{#1}\quad} \def\ttext#1{\qquad\hbox{#1}\qquad}
\def\tex#{\quad\hbox}            \def\ttex#{\qquad\hbox}
\fi


\newcount\t@
\def\mathdef #1=#2 #3
{\count@=\z@
 {\def\t{#2}\nottrue
  \loop \xdef\next{\ifcase \count@
    Ord\or Op\or Bin\or Rel\or Open\or Close\or Punct\or VarFam\else #2\fi}%
    \ifx \next\t \notfalse \fi
  \ifnot \global\advance\count@ by 1 \repeat
 }\multiply\count@"1000
 \t@=\expandafter
 \ifx\expandafter\delimiter#3\divide\t@"1000 \nottrue
 \else
   \ifcat \alpha#3#3\nottrue 
   \else \mathcode`#3 \ifnum\t@<"7000 \nottrue \else \notfalse \fi
   \fi
 \fi
 \advance\count@\t@
 \ifnot \divide\t@"1000 \multiply\t@"1000 \advance\count@-\t@
   \mathchardef#1=\count@ 
 \else \divide\t@"100 \multiply\t@"100 \advance\count@-\t@
   \def\next{#1}\afterassignment\mathd@f \t@=
 \fi
}
\def\mathd@f
{\multiply\t@"100 \advance\count@\t@ \expandafter\mathchardef\next=\count@}
\def\mathorddef #1={\mathdef #1=Ord }

\ifamstex  \let\romath=\roman
\else  \def\romath#1{{\fam\z@#1}}
\fi

\mathorddef\N=N\bffam
\mathorddef\Z=Z\bffam
\mathorddef\Q=Q\bffam
\mathorddef\R=R\bffam
\mathorddef\C=C\bffam
\mathorddef\P=P\bffam
\def\Np{\N_{>0}}
\mathorddef\1=1\bffam
\mathorddef\2=2\bffam
\mathorddef\ii=i\bffam 
\mathorddef\\=\lambda
\mathdef\:=Punct :
\let\Id=\1
\mathorddef\Sym=S\bffam 

\def\opdef#1{\def#1{\mathop\romath{\escapechar\m@ne\string#1}\nolimits}}
\def\altopdef#1#2{\def#1{\mathop\romath{#2}\nolimits}}
\def\li(#1#2..#3){%
  #1_{#2},\penalty\binoppenalty\ldots,\penalty\binoppenalty#1_{#3}}
\def\lix(#1#2.#3.#4){%
  #1_{#2},\penalty\binoppenalty\ldots,\widehat{#1_{#3}},\penalty\binoppenalty
  \ldots,#1_{#4}}
\def\fold#1(#2#3..#4){#2_{#3}#1\cdots#1#2_{#4}}
\def\foldx#1(#2#3.#4.#5){#2_{#3}#1\cdots#1\widehat{#2_{#4}}#1\cdots#1#2_{#5}}
\def\foldwith#1#2#3(#4..#5)
{\defnext##1{#1}\next{{#4}}#2#3#2\next{{#5}}}
\def\({\bigl(}\def\){\bigr)}
\def\inv{^{-1}}
\def\<#1>{{\langle #1\rangle}}

\mathdef\orbit=Bin \backslash
\let\union=\cup 
                      \let\Union=\bigcup
\let\thru=\cap

\let\after=\circ 
\let\eps=\varepsilon
\def\setof#1:#2\endset{{\{\,#1\mid#2\,\}}}
\def\Setof#1:#2\endset{\left\{\,#1\,\,\vrule\,\,#2\,\right\}}

\let\ssubset=\subset \let\subset=\subseteq
 \let\supset=\supseteq
\let\baraccent=\= \def\={\relax\?\ifmmode\then\approx\else\baraccent\fi}

\opdef\sg \opdef\Ker\opdef\Im \opdef\Aut \opdef\End \opdef\rk
\altopdef\kar{char}
\def\Card#1{\##1}

\def\multiset#1{{\{\!\!\{#1\}\!\!\}}}
\def\multisetof#1:#2\endset{\multiset{\,#1\mid#2\,}}
\def\Multisetof#1:#2\endset{\left\{\!\!\left\{\,#1\,\,\vrule
                           \,\,#2\,\right\}\!\!\right\}}
\def\multichoose(#1,#2){\mathchoice
  {\left(\!{#1\choose#2}\!\right)}%
  {\bigl(\!{#1\choose#2}\!\bigr)}%
  {\left(\!{#1\choose#2}\!\right)}%
  {\left(\!{#1\choose#2}\!\right)}%
}




\def\topsmash#1{{\def\finsm@sh{\ht\z@\z@ \box\z@}\smash{#1}}}
\def\botsmash#1{{\def\finsm@sh{\dp\z@\z@ \box\z@}\smash{#1}}}

\newbox\tightbox 
\setbox\tightbox=
\hbox{\penalty -200               \skip0 = 0pt plus 0.05 \hsize
      \hskip \skip0 \penalty    0 
      \hskip \skip0 \penalty  200 \skip0 = .1\hsize
      \hskip \skip0 \penalty  283 
      \hskip \skip0 \penalty  346
      \hskip \skip0 \penalty  400
      \hskip \skip0 \penalty  447 
      \hskip \skip0 \penalty  490
      \hskip \skip0 \penalty  529
      \hskip \skip0 \penalty  567 
      \hskip \skip0 \penalty  600
      \hskip \skip0 \penalty  632
      \hskip 2\hsize 
      \vadjust{}\hfil} 
\newbox\QEDbox
\setbox\QEDbox=\hbox
{\kern-0.2pt\vrule height 6pt \kern-0.4pt 
 \vbox to 6pt
  {\kern -0.2pt
   \hbox to 6.4pt{\hrulefill\kern-0.2pt\vrule}\vfil
   \hrule width 6.2pt
   \kern -0.2pt}%
 \kern-0.4pt \vbox to 6pt{\leaders\vrule\vfil\kern-0.2pt}\kern-0.2pt
}


\ifx\ninept\undefined 

\font\ninerm=cmr9 \font\nineit=cmti9 \font\ninebf=cmbx9 \font\ninesl=cmsl9
\font\ninett=cmtt9 \font\ninei=cmmi9 \font\ninesy=cmsy9
\skewchar\ninei='177 \skewchar\ninesy='60

\def\ninept
{\def\rm{\fam0\ninerm}\def\it{\fam\itfam\nineit}\def\bf{\fam\bffam\ninebf}%
 \def\sl{\fam\slfam\ninesl}\def\tt{\fam\ttfam\ninett}\rm
 \textfont0\ninerm \textfont1\ninei \textfont2=\ninesy
 \textfont\itfam\nineit \textfont\slfam\ninesl
 \textfont\bffam\ninebf \textfont\ttfam\ninett
 \normalbaselineskip=11pt \normalbaselines
 \setbox\strutbox=\hbox{\vrule height 8 pt depth 3 pt width 0 pt}%
}

\fi

\def\caps#1{{\ninept #1}}

\ifamstex 
\let\amsfootnote=\footnote
\def\footnote#1{\amsfootnote"#1"}
\fi


\catcode`_=11 

\newif\ifamsrefs 

\ifamstex
\let\amsref=\ref
\let\amsendref=\endref
\let\amspubl=\publ
\let\amsinbook=\inbook
\let\amskey=\key
\let\amsed=\ed 
\let\amsvol=\vol
\let\amspages=\pages
\fi

\def\label_#1#2{\expandafter\xdef \csname#1.\endcsname{#2}}
\def\ref#1{\ref_#1,,,.}
\def\ref_#1,#2,,#3.{\expandafter\ref_lab\csname#1.\endcsname{#3#2}}
\def\ref_lab#1#2{\ifx #1\relax \errmessage{Reference to undefined label}%
                 \else [#1#2]\fi
                }

\newdimen\refindent
\refindent=36pt

\def\init_refs 
{\begingroup\setbox0=\vbox\bgroup \hbadness 10000 \hfuzz=\maxdimen
 \def\chapter ##1\par{}
 \def\endrefs{\egroup\egroup\endgroup }
 \ifamsrefs
 \else
  \def\endref{}\ifx\the_key\undefined\else
     \label_{\expandafter\strip_space\the_key _}%
	    {\expandafter\strip_space\the_mark _}\fi
     \endgroup}%
 \fi
}
\def\strip_space #1 _{#1}

\def\readrefs{\init_refs \input \jobname.ref }

\useouter amsref \def\ams_start\key #1 \ident#2%
{\label_{#2}{#1}\amsref\amskey{#1}%
}
\outer\def\beginrefs{\noindent\par 
 \ifamsrefs
  \bgroup \let\ref=\ams_start
  \let\publ=\amspubl \let\inbook=\amsinbook \let\key=\amskey
  \let\ed=\amsed \let\vol=\amsvol \let\pages=\amspages
 \else
  \bgroup \ninept
  \frenchspacing\tolerance=900
  \let\ref=\_ref \let\author=\_auth \let\title=\_tit
  \let\vol=\_vol \let\year=\_year
 \fi
}

\outer\def\endrefs{\egroup\smallbreak}

\def\ref_item#1{}
  \ifx#1\undefined \else
    \errmessage{Multiple definition of \string#1 within reference}\fi
  \edef#1{\iffalse}\fi 
}

\def\_ref{\begingroup
  \def\par{\endref\par}\def\ref{\endref\ref}\edef\the_mark{\iffalse}\fi}
\def\_auth{\ref_item\the_author}
\def\_tit{\ref_item\the_title}
\def\key{\ref_item\the_key}

\def\publ{\ref_item\the_publ}

\def\journal{\ref_item\the_journal}
\def\_vol{\ref_item\the_vol}
\def\_year{\ref_item\the_year}
\def\pages{\ref_item\the_pages\range_}

\def\inbook{\ref_item\the_inbook}
\def\ed{\ref_item\the_ed}
\def\report{\ref_item\the_report}

\def\range_{\noexpand\range_}   
\def\amount_{\noexpand\amount_}
\def\do_amount{{\aftergroup\do_pp\aftergroup}}
\def\do_pp{~pp}
\def\expl_pp{\def\range_{pp.~}\let\amount_=\do_amount} 
\def\impl_pp{\let\range_\empty\let\amount_=\do_amount} 

\def\q_plur#1{\expandafter\q_s#1,_}\def\q_s#1,#2_{\ifx _#2_\else s\fi}

\def\_pages{\expandafter\pages_\the_pages\unskip---_}
\def\pages_#1-#2#3-#4_{#1\ifx _#4_\else\if -#2--#3\else --#2#3\fi\fi}

\ifamstex

\def\_opt#1#2{\ifx#1\undefined\else#2\fi}
\def\opt_pages{\_opt\the_pages{\amspages\expandafter\pages_\the_pages---_}}

\def\endref{\ifamsrefs \amsendref \else \ams_convert \fi}

\useouter amsref \def\ams_convert
{}\let\par\endgraf
 \amsref \amskey\the_mark \by\the_author
 \ifx\the_inbook\undefined
   \ifx\the_journal\undefined
     \ifx\the_publ\undefined
       \ifx\the_report\undefined 
	 \errhelp{A reference must at least have a \journal \inbook \publ^^J
		  or \report field defined.}%
	 \errmessage{Unrecognised reference}%
       \else \make_report \fi
     \else \make_book \fi
   \else \make_article \fi
 \else\make_proceedings \fi
 \ifx\the_endnote\undefined\else \finalinfo\the_endnote \fi
 \amsendref
 \endgraf\endgroup
}

\def\make_report
{\paper\the_title \paperinfo\the_report\unskip \opt_pages
 \_opt\the_year{\yr\the_year}%
}
\def\make_book
{\book\the_title \opt_pages
 \_opt\the_series
 {\bookinfo\the_series
  \_opt\the_ed{\csname amsed\q_plur\the_ed\endcsname\the_ed}%
  \_opt\the_vol{\amsvol\the_vol\_opt\the_idno{\unskip, \the_idno}}%
 }%
 \amspubl\the_publ \yr\the_year \_opt\the_ISBN{\bookinfo\the_ISBN}%
}
\def\make_article
{\paper\the_title \jour\the_journal
 \_opt\the_vol{\amsvol\the_vol\_opt\the_idno{\unskip, \the_idno}}%
 \_opt\the_year{\yr\the_year}\opt_pages
}
\def\make_proceedings
{\paper\the_title \opt_pages \amsinbook\the_inbook
 \_opt\the_ed{\csname amsed\q_plur\the_ed\endcsname\the_ed}%
 \_opt\the_series{\bookinfo\the_series
   \_opt\the_vol{\amsvol\the_vol\_opt\the_idno{\unskip, \the_idno}}}%
 \_opt\the_publ{\amspubl\the_publ}\yr\the_year
}

\else 

\def\_opt#1#2{\ifx#1\undefined\else{\def\__{#1\unskip}#2}\fi}
\def\opt_pages#1{\ifx\the_pages\undefined\else{\let\__\_pages#1}\fi}

\def\endref
{}\let\par\endgraf\smallbreak\hangindent\refindent
 \setbox0=\hbox{[\the_mark\unskip]\enspace}\noindent
 \ifdim\wd0<\refindent \hbox to \refindent{\unhbox0\hfil}\else\unhbox0 \fi
 \format_author{\the_author\unskip},
 \ifx\the_inbook\undefined
   \ifx\the_journal\undefined
     \ifx\the_publ\undefined
       \ifx\the_report\undefined 
	 \errhelp{A reference must at least have a \journal \inbook \publ^^J
        	  or \report field defined.}%
         \errmessage{Unrecognised reference}%
       \else \format_report \fi
     \else \format_book \fi
   \else \format_article \fi
 \else\format_proceedings \fi
 \unskip \_opt\the_endnote{. \__}.
 \endgraf\endgroup
}

\def\format_author#1{{\rm#1}}

\def\format_report
{{\sl \the_title\unskip}, \the_report\unskip\expl_pp\opt_pages{, \__}%
 \_opt\the_year{, (\__)}}
\def\format_book
{{\sl\the_title\unskip},\expl_pp\opt_pages{ \__,}%
 \_opt\the_series{ \__\_opt\the_ed{ (\__, ed\q_plur\the_ed.)}%
 \_opt\the_vol{ \__}\_opt\the_idno{, \__},} \the_publ\unskip,
 \the_year\unskip\_opt\the_ISBN{, ISBN \__}}
\def\format_article
{``\the_title\unskip'', {\sl\the_journal\unskip\_opt\the_vol{\/ \bf\__}}%
 \_opt\the_idno{, \__}\_opt\the_year{, (\__)}\impl_pp\opt_pages{, \__}}
\def\format_proceedings
{``\the_title\unskip'',\expl_pp\opt_pages{ \__} in {\sl\the_inbook\unskip}%
 \_opt\the_ed{ (\__, ed\q_plur\the_ed.)},
 \_opt\the_series{\__\_opt\the_vol{ \__}\_opt\the_idno{, \__}, }%
 \_opt\the_publ{\__, }\the_year}

\fi 

\catcode`_=8 
\catcode`@=12 

\catcode`_=11 

\def\let_next{\let\next=}

\let\_authors\empty 
\newtoks\this_author
\newtoks\this_addr

\newtoks\title \newtoks\MSCcodes \newtoks\CRcodes \newtoks\keywords
\newtoks\titlenote

\newtoks\DateSubmitted \newtoks\DateAccepted \newtoks\DatePublished

\newif\if_started

\def\good_break{\unskip\penalty-10\ \ignorespaces}

%

\def\Author#1{\this_author={#1}%
  \ifx =#1\let_next\Author 
  \else \let_next\let_next \afterassignment\_author 
  \fi \next}
\def\_author
 {\this_addr=\expandafter{\default_address}%
  \ifx\next\address \let_next\get_address
  \else \ifx\next\email \let_next\get_email
  \else \ifx\next\URL \let_next\get_URL
  \else \add_author
  \fi\fi\fi \next
}

\def\address#1{\errmessage
 {\string\address\space should only follow \string\Author}}
\def\email#1{\errmessage
 {\string\email\space should only follow \string\Author}}
\def\URL#1{\errmessage
 {\string\URL\space should only follow \string\Author}}

\def\get_address{\begingroup\obeylines\get_addr}
\def\get_addr#1%
{\endgroup\this_addr={#1}%
  \ifx =#1\let_next\get_address 
  \else \let_next\let_next \afterassignment\get_ad
  \fi \next
}
\def\get_ad
{\ifx\next\email \let_next\get_email
 \else \ifx\next\URL \let_next\get_URL
 \else \add_author
 \fi\fi \next
}

\def\get_email#1%
{\ifx =#1\let_next\get_email
 \else \this_addr=\expandafter{\the\this_addr\endgraf{\tt #1}}%
   \let_next\let_next \afterassignment\after_em
 \fi \next
}
\def\after_em{\ifx\next\URL \let_next\get_URL \else \add_author \fi \next}
\def\get_URL#1%
{\ifx =#1\let_next\get_URL
 \else
  \this_addr=\expandafter{\the\this_addr\endgraf
                          {\def~{\char "7E }\tt http://#1}}%
  \let_next\add_author
 \fi \next
}

\def\add_author 
{\toks0=\expandafter{\_authors\\}%
 \xdef\_authors{\the\toks0 {\the\this_author
			    \noexpand\_addr{\the\this_addr}}}%
}


\def\start_art
{\vfil\eject \ifodd\pageno\else \shipout\vbox to \vsize{}\advancepageno\fi
 \topglue 10mm plus 3 mm minus 2mm
 \centerpar{\let\\=\good_break\titlefont \baselineskip=22pt \the\title}%
 \ifx \_authors\empty \message{Warning: no authors specified.}\else
   \bigskip\centerpar
    {\subtitlefont
     \def\par{\ifhmode\endgraf\else\medskip\fi}\obeylines
     \def\_addr##1{\medskip{\addressfont\let\tt=\emailfont\let\\=\par
                   \noindent##1\endgraf}}
     \let\\=\bigskip
     \_authors
    }\bigskip
 \fi
 \dedication
 \aftertitlematter
 \_startedtrue
}


\outer\def\abstract
{\if_started \errmessage{Abstract must be placed at start of article.}%
 \else \start_art \begingroup 
   \initabstract
	 \everypar{\everypar{}\setbox0=\lastbox}
 \fi
}

\outer\def\maintext
{\if_started\medskip\endgroup\else\start_art\fi
 \afterabstractmatter
}


{\obeylines\gdef\default_address%
 {Universit\'e de Poitiers,  D\'epartement de Math\'ematiques,
  UFR Sciences SP2MI, T\'el\'eport~2, BP~30179, %
\ 86962 Futuroscope Chasseneuil Cedex, \ France}}

\font\titlefont=cmssbx10 scaled \magstep 3
\font\subtitlefont=cmss12
\font\addressfont=cmssi9
\font\emailfont=cmtt10

\def\initabstract
{\bigskip \ninept \sans \narrower 
 \leftline{\hskip\leftskip ABSTRACT}\nobreak\smallskip 
}


\def\dedication{}
\def\aftertitlematter{}

\def\afterabstractmatter
{\if &\the\MSCcodes\the\CRcodes\the\keywords\the\titlenote&\else
  {\narrower\ninept\sans
   \if &\the\MSCcodes&\else
     \hangindent=4em \noindent
     {\it Mathematics Subject Classification:\/} \the\MSCcodes.
     \par\fi
   \if &\the\CRcodes&\else
     \hangindent=4em \noindent
     {\it Computing Reviews Subject Classification:\/} \the\CRcodes.
     \par\fi
   \if &\the\keywords&\else
     \hangindent=2em \noindent{\it Keywords and Phrases:\/} \the\keywords.
     \par\fi
   \if &\the\titlenote&\else
     {\def\\{\unskip\hfil\break\ignorespaces}%
      \hangindent=2em \noindent\the\titlenote.
      \par}\fi
  }\bigskip
 \fi
}

\catcode`_=8 

\Author={Marc A. A. van Leeuwen}
\address={Universit\'e de Poitiers,  D\'epartement de Math\'ematiques,
  BP~30179, 86962 Futuroscope Chasseneuil Cedex, \ France}
\email={Marc.van-Leeuwen@math.univ-poitiers.fr}
\URL={www-math.univ-poitiers.fr/~maavl/}

\title={Double crystals of binary and integral matrices}
\MSCcodes={05E10, 05E15}
\keywords={crystals, RSK correspondence, jeu de taquin, growth diagrams,
           Littlewood-Richardson rule}

\DateSubmitted={}
\DateAccepted={}
\DatePublished={}


\font\titlefont=cmbx12
\let\subtitlefont=\ninerm
\let\addressfont=\ninerm
\let\emailfont=\ninett

\catcode`_=11 
\def\start_art
{\mark{\marktoks={}}
 \topglue \smallskipamount
 \begingroup \parskip 0pt
 \centerpar{\let\\=\good_break\titlefont \the\title}%
 \ifx \_authors\empty \message{Warning: no authors specified.}\else
   \medskip\centerpar
    {\subtitlefont
     \def\par{\ifhmode\endgraf\else\medskip\fi}\obeylines
     \def\_addr##1{\par{\addressfont\let\tt=\emailfont\let\\=\par
                   \noindent##1\endgraf}}
     \let\\=\relax
     \_authors
    }\medskip
 \fi
 \endgroup
 \dedication
 \aftertitlematter
 \_startedtrue
}

\def\initabstract
{\narrower \ninept 
 \centerline{\bf Abstract}\nobreak\noindent 
}

\outer\def\maintext
{\if_started\endgroup\else\start_art\fi 
 \afterabstractmatter
}

\let\afterabstractmatter=\empty 

\refindent=3pc

\catcode`_=8 


\ifx\Tableauxloaded\relax   \else \let\Tableauxloaded\relax \fi


\newcount\cols
{\catcode`,=\active\catcode`|=\active
 \gdef\Young(#1){\hbox{$\vcenter
 {\mathcode`,="8000\mathcode`|="8000
  \def,{\global\advance\cols by 1 &}%
  \def|{\cr
        \multispan{\the\cols}\hrulefill\cr
        &\global\cols=2 }%
  \offinterlineskip\everycr{}\tabskip=0pt
  \dimen0=\ht\strutbox \advance\dimen0 by \dp\strutbox
  \halign
   {\vrule height \ht\strutbox depth \dp\strutbox##
    &&\hbox to \dimen0{\hss$##$\hss}\vrule\cr
    \noalign{\hrule}&\global\cols=2 #1\crcr
    \multispan{\the\cols}\hrulefill\cr%
   }
 }$}}
 \gdef\Skew(#1:#2){\hbox{$\vcenter
 {\mathcode`,="8000\mathcode`|="8000
  \dimen0=\ht\strutbox \advance\dimen0 by \dp\strutbox
  \def\boxbeg{\vbox
    \bgroup\hrule\kern-0.4pt\hbox to\dimen0\bgroup\strut\vrule\hss$}%
  \def\boxend{$\hss\egroup\hrule\egroup}%
  \def,{\boxend\boxbeg}%
  \def|##1:{\boxend\vrule\egroup\nointerlineskip\kern-0.4pt
    \moveright##1\dimen0\hbox\bgroup\boxbeg}%
  \def\\##1\\##2:{\boxend\vrule\egroup\nointerlineskip\kern-0.4pt
    \kern ##1\dimen0\moveright##2\dimen0\hbox\bgroup\boxbeg}%
  \moveright#1\dimen0\hbox\bgroup\boxbeg#2\boxend\vrule\egroup
 }$}}
}


\font\sevenit=cmti7 \font\fiveit=cmti7 at 5 pt
\scriptfont\itfam=\sevenit \scriptscriptfont\itfam=\fiveit


\def\smallsquares
{\textfont0=\scriptfont0 \scriptfont0=\scriptscriptfont0
 \textfont1=\scriptfont1 \scriptfont1=\scriptscriptfont1
 \textfont2=\scriptfont2 \scriptfont2=\scriptscriptfont2
 \textfont\itfam=\scriptfont\itfam \scriptfont\itfam=\scriptscriptfont\itfam
 \textfont\bffam=\scriptfont\bffam \scriptfont\bffam=\scriptscriptfont\bffam
 \setbox0=\hbox{\raise0.4pt\hbox{$($}}
 \setbox\strutbox=\hbox{\vrule width 0pt height\ht0 depth\dp0 }
}

\def\bigsquares
{\setbox\strutbox=\hbox
 {\vrule width 0pt height1.3\ht\strutbox depth1.3\dp\strutbox}
}

\newbox\hpair \newbox\vpair
\def\installboxes
{\dimen0=\ht\strutbox \advance\dimen0 by \dp\strutbox
 \setbox0=\vbox{\hrule width \dimen0}%
 \setbox1=\rlap{\strut\vrule}%
 \setbox\hpair=\rlap{\raise\ht\strutbox\copy0 \lower\dp\strutbox\copy0}%
 \setbox\vpair=\rlap{\raise\dimen0\copy1 \kern\dimen0 \copy1}%
}

\newcount\ribbonlength
\newcount\joff \newcount\jmax \newcount\imid \newcount\jmid


\def\ribbon #1,#2:#3;#4 
{\setbox2=\hbox
 {\count0=0 \count1=0
  \copy1 \rlap{\lower\dp\strutbox\copy0}\if|#4|\else\ribs#4 \fi
  \raise\count1\dimen0\rlap
    {\kern\count0\dimen0 \raise\ht\strutbox\copy0 \copy1\raise0.4pt\copy1}%
  \global\joff=\count0
  \divide\dimen0 by 2 \raise\imid\dimen0\rlap{\kern\jmid\dimen0
    \hbox to 2\dimen0{\hfil$#3$\hfil}}%
 }%
 \lower#1\dimen0\rlap{\kern#2\dimen0\box2}%
 \advance\joff by #2\relax \ifnum \joff>\jmax \global\jmax=\joff \fi
}
\def\ribs#1#2
{\raise\count1\dimen0\rlap
   {\kern\count0\dimen0 \copy\ifnum #1=0 \hpair \else \vpair \fi}%
 \advance\ribbonlength by -2
 \ifnum \ribbonlength=0 \jmid=\count0 \imid=\count1 \fi
 \advance\count#1 by 1
 \ifnum \ribbonlength=1 \jmid=\count0 \imid=\count1 \fi
 \ifnum \ribbonlength>-1 \ifnum \ribbonlength<2
   \advance \jmid by \count0 \advance \imid by \count1
 \fi\fi
 \if |#2|\else \ribs#2 \fi 
}

\def\rtab#1
{\installboxes \ribbonlength=#1 \global\jmax=0
 \setbox0=\hbox{\aftergroup\endrtab \aftergroup}}
\def\endrtab
{\advance\dimen0 by \jmax\dimen0 
 \wd0=\dimen0 \vcenter{\box0}}

\def\edgeseq#1,#2;#3 
{\setbox2=\hbox{\count0=0 \count1=0 \edges#3
  \global\joff=\count0
  }%
 \rlap{\kern#2\dimen0
       \dimen1=#1\dimen0 \advance\dimen1 by -\dimen0 \lower\dimen1\box2
      }%
 \advance\joff by #2\advance\joff by -1 
 \ifnum \joff>\jmax \global\jmax=\joff \fi
}
\def\edges#1#2
{\rlap
 {\kern\count0\dimen0
  \ifnum #1=0
     \dimen1=\count1\dimen0 \advance\dimen1 by -\dp\strutbox
     \raise\dimen1\copy0
  \else \raise\count1\dimen0\copy1
  \fi}%
 \advance\count#1 by 1
 \if |#2|\else \edges#2 \fi 
}


\def\labeledges#1,#2:#3;#4 
{\setbox2=\hbox{\count0=0 \count1=0 \divide\dimen0 by 2 \edgelabs#3;#4 }%
 \rlap{\kern#2\dimen0
       \dimen1=-\ht\strutbox 
       \advance\dimen1 by \fontdimen22\textfont2 
       \advance\dimen1 by #1\dimen0 \lower\dimen1\box2}%
}
\def\edgelabs#1#2;#3#4
{\advance\count#3 by 1 
 \rlap {\kern\count0\dimen0 \raise\count1\dimen0 \clap{$#1$}}%
 \advance\count#3 by 1 
 \if |#4|\else \edgelabs#2;#4 \fi 
}

\def\arrow#1,#2:#3
{\smash{\lower#1\dimen0\rlap{$\kern#2.5\dimen0
 \if r#3 \,\rightarrow \else
 \if l#3 \llap{$\leftarrow\,$} \else
 \if u#3 \setbox0=\hbox to 0pt{\hss$\uparrow$\hss}\ht0=0pt 
         \raise.5\dimen0 \box0 \else
 \if d#3 \setbox0=\hbox to 0pt{\hss$\downarrow$\hss}\dp0=-2pt
         \lower.5\dimen0 \box0 \else
 \if i#3 \setbox0=\hbox to 0pt{\hss$\nwarrow\,$}\ht0=0pt
         \raise.5\dimen0 \box0 \else
 \if o#3 \setbox0=\hbox to 0pt{$\,\searrow$\hss}\dp0=-2pt
         \lower.5\dimen0 \box0 \else
 \fi\fi\fi\fi\fi\fi
 $}}%
}


\def\ifabsent#1\intoks#2%
{\defnext##1#1##2##3\end {\ifx\ifabsent##2}%
 \expandafter\next\the#2#1\ifabsent\end
}

\let\intoks=\if 

\def\ifremoved#1\intoks#2{%
    \defnext##1#1##2#1##3\end{\ifx#1##3\global#2{##1##2}}%
    \expandafter\next\the#2#1#1\end}

\def\addtotoks#1#2{\ifabsent#2\intoks#1%
    \edef\next{\global#1{\the#1\noexpand#2}}\next\fi}

\newwrite\frwfile
\def\tofrw{\immediate\write\frwfile}
\def\writefrw{\immediate\openout\frwfile=\jobname.frw
              \global\let\writefrw=\tofrw \tofrw}

\newtoks\oldforward 
\newtoks\newforward 

\begingroup
  \openin15=\jobname.frw
  \ifeof15 \closein15
  \else
   \closein15
   \def\forward#1=#2{\gdef#1{#2}\addtotoks\oldforward#1}
   \input\jobname.frw
  \fi
\endgroup

\newif\ifforward
\def\deadref{{\bf?}}


\def\forward#1
{\ifx#1\undefined \global\let#1\deadref \forwardtrue
   \immediate\write16{Warning: unresolved forward reference to \string#1.}%
 \else \ifabsent#1\intoks\oldforward \forwardfalse \else\forwardtrue \fi
 \fi
 \ifforward \addtotoks\newforward#1\fi
 #1%
}

\def\definefrw#1
{\ifremoved#1\intoks\newforward 
  \ifx#1\lastlabel\else 
    \ifx#1\deadref\else 
     \immediate\write16{Warning:
       forward reference to \string#1 is incorrect, run again!
       Changed from #1 to \lastlabel.}%
  \fi\fi
  \writefrw{\string\forward\string#1={\lastlabel}}
  \global\let#1\lastlabel \deflab#1
  \letnext\ignorespaces 
 \else 
  \ifremoved#1\intoks\oldforward 
   \global\let#1\lastlabel \deflab#1\letnext\ignorespaces 
  \else\defnext{#1}
  \fi
 \fi
}

\def\labelsec#1%
{\ifx\undefined#1
 \global\let#1\lastlabel \deflab#1%
 \else\definefrw{#1}%
 \fi
 \nobreak 
}

\outer\def\proclaim#1. #2#3\par
  {\secmedbreak \global\advance\proclno by 1 \reset\itemno
   \edef\lastlabel{\ifnum-1<\secno    \number\secno.\fi
                   \ifnum-1<\subsecno \number\subsecno.\fi
                   \ifnum-1<\subsubsecno \number\subsubsecno.\fi
                   \number\proclno}
   {\interlinepenalty=300 \parskip=0pt\def\par{\endgraf\penalty200 }
    \noindent 
    \if\let\noexpand#2
     \ifx\undefined#2
      \global\let#2\lastlabel \deflab#2\letnext\ignorespaces
     \else \definefrw{#2}
     \fi
    \else\defnext{#2}
    \fi
    \bf\lastlabel.\enspace#1.\enspace \sl\next#3\endgraf
    }%
   \ifdim \lastskip<\medskipamount \removelastskip\penalty55\medskip \fi
   }

\def\neqn(#1#2)%
{\global\advance\eqnumber by 1 \edef\lastlabel{\number\eqnumber}%
  (\ifcat a\noexpand#1\errmessage{Incorrect label #1#2}\lastlabel
   \else
    \ifx\undefined#1\global\let#1\lastlabel \deflab#1%
    \else\definefrw{#1}%
    \fi
    \lastlabel\fam0 #2
   \fi
  )}


\font\fivesl=cmsl5
\font\sevensl=cmsl7
\scriptfont\slfam=\sevensl \scriptscriptfont\slfam=\fivesl

\mathorddef\comp=C2 
\mathorddef\Part=P2 
\def\bcomp{\comp^{\scriptscriptstyle\set2}}
\mathorddef\Mat=M2 
\def\bMat{\Mat^{\scriptscriptstyle\set2}}

\def\tr{^{\sl t}}
\def\repex^#1{^{(#1)}}

\def\Kr#1{\mathchar"465B\,#1\,\mathchar"565D\,} 

\def\set#1{{[#1]}}
\def\sing#1{\{#1\}} 
\def\edgprm#1#2{\eps(#1,#2)}
\altopdef\sym{Sym}

\def\[#1,#2]{{#1[\if!#2!\,\else#2\fi]}}
\def\d#1#2{{\delta_{#1}\if|#2|\else+#2\fi}}

\def\scalar<#1,#2>{{\bigl<\,#1\bigm|#2\,\bigr>}}
\def\Scalar<#1,#2>{\left<\,#1\,\,\vrule\,\,#2\,\right>}

\def\Mtt(#1,#2,#3,#4){\({#1\atop#3}~{#2\atop#4}\)} 

\let\leh=\leftharpoondown
\let\geh=\rightharpoonup
\let\lev=\leftharpoonup
\let\gev=\rightharpoondown
\def\nleh{\,\not\nobreak\!\leh}

\opdef\Tabl
\def\mat#1#2{\Mat_{#1,#2}}
\def\bmat#1#2{\bMat_{#1,#2}}
\def\IE#1{\Tabl(#1)}
\def\BE#1{\Tabl^{\scriptscriptstyle\set2}(#1)}
\opdef\LR
\def\IL#1{\LR(#1)}
\def\BL#1{\LR^{\scriptscriptstyle\set2}(#1)}

\opdef\row \opdef\col
\def\rsum#1{\row(#1)}
\def\csum#1{\col(#1)}

\opdef\wt
\def\weight#1{\wt(#1)}
\opdef\cwt \opdef\rwt

\opdef\SST

\def\SSTx#1{\SST(#1)}

\altopdef\height{ht}

\def\commas#1#2 {#1\if|#2|\else,\commas#2 \fi}
\def\s#1 {s_{(\commas#1 )}}
\def\m#1 {m_{(\commas#1 )}}
\def\h#1 {h_{(\commas#1 )}}
\def\e#1 {e_{(\commas#1 )}}

\newdimen\diaght

\def\xleft{^{\scriptscriptstyle\leftarrow}}
\def\xright{^{\scriptscriptstyle\rightarrow}}
\def\xraise{^{\scriptscriptstyle\uparrow}}
\def\xlower{^{\scriptscriptstyle\downarrow}}

\def\ltm{e\xleft}  
\def\rtm{f\xright} 
\def\upm{e\xraise} 
\def\dnm{f\xlower} 

\def\nlm{n\xleft}  
\def\nrm{n\xright} 
\def\num{n\xraise} 
\def\ndm{n\xlower} 

\def\hr{\sigma^{\scriptscriptstyle\leftrightarrow}} 
\def\vr{\sigma^{\scriptscriptstyle\updownarrow}}    

\def\horact{^{\scriptscriptstyle\leftrightarrow}} 
\def\vertact{^{\scriptscriptstyle\updownarrow}} 

\opdef\Diag 
\def\diagram#1{\Diag^{\scriptscriptstyle\set2}(#1)}
\def\diagon#1{\Diag(#1)}
\opdef\NF 
\opdef\SF 
\def\frdiag#1#2{\NF(#1,#2)}
\opdef\rev
\def\revert#1#2{\rev_{#2}(#1)}
\def\sliced#1#2#3{\SF(#1,#2,#3)}

\def\e#1#2{e_{#2/#1}}

\def\lenw{\leq_{\scriptscriptstyle\nwarrow}}

\def\lene{\leq_{\scriptscriptstyle\nearrow}}

\opdef\Pic \opdef\Int \opdef\Bin

\def\is#1{\pi(#1)} 

\def\matrixcross#1#2#3
{\multispan#3\strut&\kern-0.5em\smash
 {\hbox{\setbox0=\hbox{0}\dimen0=\ht0\dimen1=\dp0
        \advance\dimen0 by #1\baselineskip\advance\dimen1 by #2\baselineskip
        \vrule height \dimen0 depth\dimen1}}\hfill\cr
 \noalign{\kern-\baselineskip\kern-0.2pt\hrule\kern-0.2pt}\strut
}
\def\matrixh
{\strut\cr
 \noalign{\kern-\baselineskip\kern-0.2pt\hrule\kern-0.2pt}\strut
}
\def\matrixv#1#2#3%
{\multispan#3\strut&\kern-0.5em\smash
 {\hbox{\setbox0=\hbox{0}\dimen0=\ht0\dimen1=\dp0
        \advance\dimen0 by #1\baselineskip\advance\dimen1 by #2\baselineskip
        \vrule height \dimen0 depth\dimen1}}\hfill\cr
 \noalign{\kern-\baselineskip}\strut
}

\input epsf    

\refindent=3pc

\readrefs

\maketoc

\abstract
 We introduce a set of operations that we call crystal operations on matrices
 with entries either in $\sing{0,1}$ or in~$\N$. There are horizontal and
 vertical forms of these operations, which commute with each other, and they
 give rise to two different structures of a crystal graph of type~$A$ on these
 sets of matrices. They provide a new perspective on many aspects of the
 RSK~correspondence and its dual, and related constructions. Under a
 straightforward encoding of semistandard tableaux by matrices, the operations
 in one direction correspond to crystal operations applied to tableaux, while
 the operations in the other direction correspond to individual moves
 occurring during a jeu de taquin slide. For the (dual) RSK~correspondence, or
 its variant the Burge correspondence, a matrix~$M$ can be transformed by
 horizontal respectively vertical crystal operations into each of the matrices
 encoding the tableaux of the pair associated to~$M$, and the inverse of this
 decomposition can be computed using crystal operations too. This
 decomposition can also be interpreted as computing Robinson's correspondence,
 as well as the Robinson-Schensted correspondence for pictures. Crystal
 operations shed new light on the method of growth diagrams for describing the
 RSK and related correspondences: organising the crystal operations in a
 particular way to determine the decomposition of matrices, one finds growth
 diagrams as a method of computation, and their local rules can be deduced
 from the definition of crystal operations. The Sch\"utzenberger involution
 and its relation to the other correspondences arise naturally in this
 context. Finally we define a version of Greene's poset invariant for both of
 the types of matrices considered, and show directly that crystal operations
 leave it unchanged, so that for such questions in the setting of matrices
 they can take play the role that elementary Knuth transformations play for
 words.

\maintext

\def\LRgensymeq{50}
\def\LRthmone{5.1}
\def\LRthmtwo{5.2}
\def\LRgeninteq{54}
\def\LRgenbineq{55}

\subsecno=0  

\secno=-1    

\newsection Introduction.

The Robinson-Schensted correspondence between permutations and pairs of
standard Young tableaux, and its generalisation by Knuth to matrices and
semistandard Young tableaux (the \caps{RSK}~correspondence) are intriguing not
only because of their many surprising combinatorial properties, but also by
the great variety in ways in which they can be defined. The oldest
construction by Robinson was (rather cryptically) defined in terms of
transformations of words by ``raising'' operations. The construction by
Schensted uses the well known insertion of the letters of a word into a Young
tableau. While keeping this insertion procedure, Knuth generalised the input
data to matrices with entries in~$\N$ or in~$\sing{0,1}$. He also introduced a
``graph theoretical viewpoint'' (which could also be called poset theoretical,
as the graph in question is the Hasse diagram of a finite partially ordered
set) as an alternative construction to explain the symmetry of the
correspondence; a different visualisation of this construction is presented in
the ``geometrical form'' of the correspondence by Viennot, and in Fulton's
``matrix-ball'' construction. A very different method of describing the
correspondence can be given using the game of ``\foreign{jeu de taquin}''
introduced by Lascoux and Sch\"utzenberger. Finally a construction for the
\caps{RSK}~correspondence using ``growth diagrams'' was given by Fomin; it
gives a description of the correspondence along the same lines as Knuth's
graph theoretical viewpoint and its variants, but it has the great advantage
of avoiding all iterative modification of data, and computes the tableaux
directly by a double induction along the rows and columns of the matrix.

The fact that these very diverse constructions all define essentially the same
correspondence (or at least correspondences that are can be expressed in terms
of each other in precise ways) can be shown using several notions that
establish bridges between them. For instance, to show that the
``rectification'' process using jeu de taquin gives a well defined map that
coincides with the one defined by Schensted insertion, requires (in the
original approach) the additional consideration of an invariant quantity for
jeu de taquin that is a special case of Greene's invariant for finite posets,
and of a set of elementary transformations of words introduced by Knuth. A
generalisation of the \caps{RSK}~correspondence from matrices to ``pictures''
was defined by Zelevinsky, which reinforces the link with
Littlewood-Richardson tableaux already present in the work of Robinson; it
allows the correspondences considered by Robinson, Schensted, and Knuth to be
viewed as derived from a single general correspondence. The current author has
shown that that correspondence can alternatively be described using (two forms
of) jeu de taquin for pictures instead of an insertion process, and that in
this approach the use of Greene's invariant and elementary Knuth operations
can be avoided. A drawback of this point of view is that the complications of
the Littlewood-Richardson rule are built into the notion of pictures itself;
for understanding that rule we have also given a description that is simpler
(at the price of loosing some symmetry), where semistandard tableaux replace
pictures, and ``crystal'' raising and lowering operations replace one form of
jeu de taquin, so that Robinson's correspondence is completely described in
terms of jeu de taquin and crystal operations.

In this paper we shall introduce a new construction, which gives rise to
correspondences that may be considered as forms of the
\caps{RSK}~correspondence (and variants of it). Its interest lies not so much
in providing yet another computational method for that correspondence, as in
giving a very simple set of rules that implicitly define it, and which can be
applied in sufficiently flexible ways to establish links with nearly all known
constructions and notions related to it. Simpler even than semistandard
tableaux, our basic objects are matrices with entries in~$\N$ or
in~$\sing{0,1}$, and the basic operations considered just move units from one
entry to an adjacent entry. As the definition of those operations is inspired
by crystal operations on words or tableaux, we call them crystal operations on
matrices.

Focussing on small transformations, in terms of which the main correspondences
arise only implicitly and by a non-deterministic procedure, our approach is
similar to that of jeu de taquin, and to some extent that of elementary Knuth
transformations. By comparison our moves are even smaller, they reflect the
symmetry of the \caps{RSK}~correspondence, and they can be more easily related
to the constructions of that correspondence by Schensted insertion or by
growth diagrams. Since the objects acted upon are just matrices, which by
themselves hardly impose any constraints at all, the structure of our
construction comes entirely from the rules that determine when the transfer of
a unit between two entries is allowed. Those rules, given in definitions
\forward\binmovedef\ and~\forward\intmovedef\ below, may seem somewhat strange
and arbitrary; however, we propose to show in this paper is that in many
different settings they do precisely the right thing to allow interesting
constructions. One important motivating principle is to view matrices as
encoding semistandard tableaux, by recording the weights of their individual
rows or columns; this interpretation will reappear all the time. All the same
it is important that we are not forced to take this point of view: sometimes
it is clearest to consider matrices just as matrices.

While the above might suggest that we introduced crystal operations in an
attempt to find a unifying view to the numerous approaches to the
\caps{RSK}~correspondence, this paper in fact originated as a sequel
to~\ref{alternating Schur}, at the end of which paper we showed how two
combinatorial expressions for the scalar product of two skew Schur functions,
both equivalent to Zelevinsky's generalisation of the Littlewood-Richardson
rule, can be derived by applying cancellations to two corresponding
expressions for these scalar products as huge alternating sums. We were led to
define crystal operations in an attempt to organise those cancellations in
such a way that they would respect the symmetry of those expressions with
respect to rows and columns. We shall remain faithful to this original
motivation, by taking that result as a starting point for our paper; it will
serve as motivation for the precise form of the definition crystal operations
on matrices. That result, and the Littlewood-Richardson rule, do not however
play any role in our further development, so the reader may prefer to take the
definitions of crystal operations as a starting point, and pick up our
discussion from there.

This paper is a rather long one, even by the author's standards, but the
reason is not that our constructions are complicated or that we require
difficult proofs in order to justify them. Rather, it is the large number of
known correspondences and constructions for which we wish to illustrate the
relation with crystal operations that accounts for much of the length of the
paper, and the fact that we wish to describe those relations precisely rather
than in vague terms. For instance, we arrive fairly easily at our central
theorem~\forward\decompthm, which establishes the existence of bijective
correspondences with the characteristics of the \caps{RSK}~correspondence and
its dual; however, a considerable additional analysis is required to identify
these bijections precisely in terms of known correspondences, and to prove the
relation found. Such detailed analyses require some hard work, but there are
rewards as well, since quite often the results have some surprising aspects;
for instance the correspondences of the mentioned theorem turn out to be more
naturally described using column insertion than using row insertion, and in
particular we find for integral matrices the Burge correspondence rather than
the \caps{RSK}~correspondence. We do eventually find a way in which the
\caps{RSK}~correspondence arises directly from crystal operations, in
proposition~\forward\RSKprop, but this is only after exploring various
different possibilities of constructing growth diagrams.

Our paper is organised as follows. We begin directly below by recalling from
\ref{alternating Schur} some basic notations that will be used throughout.
In \Sec\forward\copsec\ we introduce, first for matrices with entries
in~$\sing{0,1}$ and then for those with entries in~$\N$, crystal operations
and the fundamental notions related to them, and we prove the commutation of
horizontal and vertical operations, which will be our main technical tool. In
\Sec\forward\crystalsec\ we mention a number of properties of crystal graphs,
which is the structure one obtains by considering only vertical or only
horizontal operations; in this section we also detail the relation between
crystal operations and jeu de taquin. In \Sec\forward\Knuthsec\ we start
considering double crystals, the structure obtained by considering both
vertical and horizontal crystal operations. Here we construct our central
bijective correspondence, which amounts to a decomposition of every double
crystal as a Cartesian product of two simple crystals determined by one same
partition, and we determine how this decomposition is related to known Knuth
correspondences. In \Sec\forward\growthsec\ we present the most novel aspect
of the theory of crystal operations on matrices, namely the way in which the
rules for such operations lead to a method of computing the decomposition of
theorem~\forward\decompthm\ using growth diagrams. The use of growth diagrams
to compute Knuth correspondences is well known of course, but here the fact
that such a growth diagram exists, and the local rule that this involves, both
follow just from elementary properties of crystal operations, without even
requiring enumerative results about partitions. In \Sec\forward\chainsec\ we
study the relation between crystal operations and the equivalents in terms of
matrices of increasing and decreasing subsequences, and more generally of
Greene's partition-valued invariant for finite posets. Finally, in
\Sec\forward\proofsec\ we provide the proofs of some results, which were
omitted in text of the preceding sections to avoid distracting too much from
the main discussion. (However for all our central results the proofs are quite
simple and direct, and we have deemed it more instructive to give them
right away in the main text.)

\subsection Notations.

We briefly review those notations from \ref{alternating Schur} that will be
used in the current paper. We shall use the Iverson symbol, the notation
$\Kr{{\it condition}}$ designating the value~$1$ if the Boolean {\it
condition} holds and $0$~otherwise. For $n\in\N$ we denote by $\set{n}$ the
set $\setof{i\in\N}:i<n\endset=\sing{0,\ldots,n-1}$. The set $\comp$ of
compositions consists of the sequences $(\alpha_i)_{i\in\N}$ with
$\alpha_i\in\N$ for all~$i$, and $\alpha_i=0$ for all but finitely many~$i$;
it may be thought of as $\Union_{n\in\N}\N^n$ where each $\N^n$ is considered
as a subset of~$\N^{n+1}$ by extension of its vectors by an entry~$0$. Any
$\alpha\in\comp$ is said to be a composition of the number
$|\alpha|=\sum_{i\in\N}\alpha_i$. We put
$\bcomp=\setof\alpha\in\comp:\forall{i}\:\alpha_i\in\set2\endset$; its
elements are called binary compositions. The set $\Part\ssubset\comp$ of
partitions consists of compositions that are weakly decreasing sequences. The
operators `$+$' and~`$-$', when applied to compositions or partitions, denote
componentwise addition respectively subtraction.

The diagram of~$\\\in\Part$, which is a finite order ideal of~$\N^2$, is
denoted by~$\set\\$, and the conjugate partition of~$\\$ by~$\\\tr$. For
$\kappa,\\\in\Part$ the symbol $\\/\kappa$ is only used when
$\set\kappa\subset\set\\$ and is called a skew shape; its diagram
$\set{\\/\kappa}$ is the set theoretic difference $\set\\-\set\kappa$, and we
write $|\\/\mu|=|\\|-|\mu|$. For $\alpha,\beta\in\comp$ the relation
$\alpha\leh\beta$ is defined to hold whenever one has
$\beta_{i+1}\leq\alpha_i\leq\beta_i$ for all~$i\in\N$; this means that
$\alpha,\beta\in\Part$, that $\set\alpha\subset\set\beta$, and that
$\set{\beta/\alpha}$ has at most one square in any column. When $\mu\leh\\$,
the skew shape $\\/\mu$ is called a horizontal strip. If $\mu\tr\leh\\\tr$
holds, we call $\\/\mu$ a vertical strip and write $\mu\lev\\$; this condition
amounts to $\mu,\\\in\Part$ and $\\-\mu\in\bcomp$ (so $\set{\\/\mu}$ has at
most one square in any row).

A semistandard tableau~$T$ of shape $\\/\kappa$ (written
$T\in\SSTx{\\/\kappa}$) is a sequence $(\\^{(i)})_{i\in\N}$ of partitions
starting at~$\kappa$ and ultimately stabilising at~$\\$, of which successive
members differ by horizontal strips: $\\^{(i)}\leh\\^{(i+1)}$ for
all~$i\in\N$. The weight $\weight{T}$ of~$T$ is the composition
$(|\\^{(i+1)}/\\^{(i)}|)_{i\in\N}$. Although we shall work mostly with such
tableaux, there will be situations where it is more natural to consider
sequences in which the relation between successive members is reversed
($\\^{(i)}\geh\\^{(i+1)}$) or transposed ($\\^{(i)}\lev\\^{(i+1)}$), or both
($\\^{(i)}\gev\\^{(i+1)}$); such sequences will be called reverse and/or
transpose semistandard tableaux. The weight of transpose semistandard tableaux
is then defined by the same expression as that of ordinary ones, while for
their reverse counterparts it is the composition
$(|\\^{(i)}/\\^{(i+1)}|)_{i\in\N}$.

The set $\Mat$ is the matrix counterpart of~$\comp$: it consists of
matrices~$M$ indexed by pairs $(i,j)\in\N^2$, with entries in~$\N$ of which
only finitely many are nonzero (note that rows and columns are indexed
starting from~$0$). It may be thought of as the union of all sets
of finite matrices with entries in~$\N$, where smaller matrices are identified
with larger ones obtained by extending them with entries~$0$. The set of such
matrices with entries restricted to $\set2=\sing{0,1}$ will be denoted
by~$\bMat$; these are called binary matrices. For matrices $M\in\Mat$, we
shall denote by $M_i$ its row~$i$, which is $(M_{i,j})_{j\in\N}\in\comp$,
while $M\tr_j=(M_{i,j})_{i\in\N}\in\comp$ denotes its column~$j$. We denote by
$\rsum{M}=(|M_i|)_{i\in\N}$ the composition formed by the row sums of~$M$, and
by $\csum{M}=(|M\tr_j|)_{j\in\N}$ the composition formed by its column sums,
and we define $\mat\alpha\beta =\setof{M\in\Mat}:
\rsum{M}=\alpha,\csum{M}=\beta\endset$ and $\bmat\alpha\beta
=\mat\alpha\beta\thru\bMat$.

In the remainder of our paper we shall treat $\bMat$ and $\Mat$ as analogous
but separate universes, in other words we shall never consider a binary matrix
as an integral matrix whose entries happen to be~$\leq1$ or vice versa; this
will allow us to use the same notation for analogous constructions in the
binary and integral cases, even though their definition for the integral case
is not an extension of the binary one.


\subsubsecno=-1 

\newsection Crystal operations on matrices.
\labelsec\copsec

The motivation and starting point for this paper are formed by a number of
expressions for the scalar product between two Schur functions in terms over
enumerations of matrices, which were described in~\ref{alternating Schur}.
To present them, we first recall the way tableaux were encoded by matrices in
that paper.

\subsection Encoding of tableaux by matrices.
\labelsec\encodingsec

A semistandard tableau~$T=(\\^{(i)})_{i\in\N}$ of shape~$\\/\kappa$ can be
displayed by drawing the diagram~$\set{\\/\kappa}$ in which the squares of
each strip $\set{\\^{(i+1)}/\\^{(i)}}$ are filled with entries~$i$. Since the
columns of such a display are strictly increasing and the rows weakly
increasing, such a display is uniquely determined by its shape plus one of the
following two informations: (1)~for each column~$C_j$ the set of entries
of~$C_j$, or (2)~for each row~$R_i$ the multiset of entries of~$R_i$. Each of
those informations can be recorded in a matrix: the binary matrix~$M\in\bMat$
in which $M_{i,j}\in\{0,1\}$ indicates the absence or presence of an
entry~$i$ in column~$C_j$ of the display of~$T$ will be called the
\newterm{binary encoding} of~$T$, while the integral matrix~$N\in\Mat$ in
which $N_{i,j}$ gives the number of entries~$j$ in row~$R_i$ of the display
of~$T$ will be called the \newterm{integral encoding} of~$T$. In terms of the
shapes~$\\^{(i)}$ these matrices can be given directly by
$M_i=(\\^{(i+1)})\tr-(\\^{(i)})\tr$ for all~$i$ and
$N\tr_j=(\\^{(j+1)})-(\\^{(j)})$ for all~$j$, cf.~\ref{alternating Schur,
definition~1.2.3}. Note that the columns~$M\tr_j$ of the binary encoding
correspond to the columns~$C_j$, and the rows $N_i$ of the integral encoding
to the rows~$R_i$. While this facilitates the interpretation of the matrices,
it will often lead to an interchange of rows and columns between the binary
and integral cases; for instance from the binary encoding~$M$ the weight
$\weight{T}$ can be read off as $\rsum{M}$, while in terms of the integral
encoding~$N$ it is~$\csum{N}$. Here is an example of the display of a
semistandard tableau~$T$ of shape $(9,8,5,5,3)/(4,1)$ and
weight~$(2,3,3,2,4,4,7)$, with its binary and integral encodings $M$ and~$N$,
which will be used in examples throughout this paper:
\bigdisplay
T\:
\Skew(4:0,2,4,5,5|1:0,1,3,4,6,6,6|0:1,1,2,4,5|0:2,3,4,6,6|0:5,6,6)
\kern-1em,\quad
M\:
\pmatrix{0&1&0&0&1&0&0&0&0\cr
         1&1&1&0&0&0&0&0&0\cr
         1&0&1&0&0&1&0&0&0\cr
	 0&1&0&1&0&0&0&0&0\cr
	 0&0&1&1&1&0&1&0&0\cr
	 1&0&0&0&1&0&0&1&1\cr
	 0&1&1&1&1&1&1&1&0\cr}
,\quad
N\:
\pmatrix{1&0&1&0&1&2&0\cr
         1&1&0&1&1&0&3\cr
         0&2&1&0&1&1&0\cr
	 0&0&1&1&1&0&2\cr
	 0&0&0&0&0&1&2\cr}
.
$$

To reconstruct $T$ from its binary or integral encoding, one needs to know the
shape $\\/\kappa$ of~$T$, which is not recorded in the encoding; since
$\\$~and~$\kappa$ are related by $\\\tr-\kappa\tr=\csum{M}$ in the binary case
and by $\\-\kappa=\rsum{N}$ in the integral case, it suffices to know one of
them. Within the sets $\bMat$ and $\Mat$ of all binary respectively integral
matrices, each shape $\\/\kappa$ defines a subset of matrices that occur as
encodings of tableaux of that shape: we denote by $\BE{\\/\kappa}\subset\bMat$
the set of binary encodings of tableaux $T\in\SSTx{\\/\kappa}$, and by
$\IE{\\/\kappa}\subset\Mat$ the set of integral encodings of such tableaux.
The conditions that define such subsets, which we shall call ``tableau
conditions'', can be stated explicitly as follows.

\proclaim Proposition. \Tabprop
Let $\\/\kappa$ be a skew shape. For $M\in\bMat$ one has
$M\in\BE{\\/\kappa}$ if and only if\/ $\csum{M}=\\\tr-\kappa\tr$, and
$\kappa\tr+\sum_{i<k}M_i\in\Part$ for all~$k\in\N$. For $M\in\Mat$ one has
$M\in\IE{\\/\kappa}$ if and only if\/ $\rsum{M}=\\-\kappa$, and
$(\kappa+\sum_{j<l}M\tr_j)\leh(\kappa+\sum_{j\leq{l}}M\tr_j)$ for all~$l\in\N$.

\proof
This is just a verification that an appropriate tableau encoded by the matrix
can be reconstructed if and only if the given conditions are satisfied. We
have seen that if $M\in\bMat$ is the binary encoding of some
$(\\^{(i)})_{i\in\N}\in\SSTx{\\/\kappa}$, then
$M_i=(\\^{(i+1)})\tr-(\\^{(i)})\tr$ for all~$i$, which together with
$\\^{(0)}=\kappa$ implies $(\\^{(k)})\tr=\kappa\tr+\sum_{i<k}M_i$ for
$k\in\N$. A sequence of partitions $\\^{(i)}$ satisfying this condition exists
if and only if each value $\kappa\tr+\sum_{i<k}M_i$ is a partition. If so,
each condition $\\^{(i)}\leh\\^{(i+1)}$ will be automatically satisfied, since
it is equivalent to $(\\^{(i)})\tr\lev(\\^{(i+1)})\tr$, while by construction
$(\\^{(i+1)})\tr-(\\^{(i)})\tr=M_i\in\bcomp$; therefore $(\\^{(i)})_{i\in\N}$
will be a semistandard tableau. Moreover $\csum{M}=\\\tr-\kappa\tr$ means that
$\kappa\tr+\sum_{i<k}M_i=\\\tr$ for sufficiently large~$k$, and therefore that
the shape of the semistandard tableau found will be $\\/\kappa$.

Similarly if $M\in\Mat$ is the integral encoding of some
$(\\^{(i)})_{i\in\N}\in\SSTx{\\/\kappa}$, then we have seen that
$M\tr_j=(\\^{(j+1)})-(\\^{(j)})$ for all~$j$, which together with
$\\^{(0)}=\kappa$ implies $\\^{(l)}=\kappa+\sum_{j<l}M\tr_j$ for $l\in\N$. By
definition the sequence $\smash{(\\^{(i)})_{i\in\N}}$ so defined for a given
$\kappa$ and $M\in\Mat$ is a semistandard tableau if and only if
$\\^{(l)}\leh\\^{(l+1)}$ for all~$l\in\N$ (which implies that all $\\^{(l)}$
are partitions), in other words if and only if
$(\kappa+\sum_{j<l}M\tr_j)\leh(\kappa+\sum_{j\leq{l}}M\tr_j)$ for
all~$l\in\N$. The value of $\\^{(l)}$ ultimately becomes $\kappa+\rsum{M}$, so
the semistandard tableau found will have shape $\\/\kappa$ if and only if
$\rsum{M}=\\-\kappa$.
\QED

Littlewood-Richardson tableaux are semistandard tableaux satisfying some
additional conditions, and the Littlewood-Richardson rule expresses certain
decomposition multiplicities by counting such tableaux (details, which are not
essential for the current discussion, can be found in~\ref{L-R intro}).
In~\ref{alternating Schur, theorems \LRthmone~and~\LRthmtwo}, a generalised
version of that rule is twice stated in terms of matrices, using respectively
binary and integral encodings. A remarkable aspect of these formulations is
that the additional conditions are independent of the tableau conditions that
these matrices must also satisfy, and notably of the shape~$\\/\kappa$ for
which they do so; moreover, the form of those additional conditions is quite
similar to the tableau conditions, but with the roles of rows and columns
interchanged. We shall therefore consider these conditions separately, and
call them ``Littlewood-Richardson conditions''.

\proclaim Definition. \LRdef
Let $\nu/\mu$ be a skew shape. The set $\BL{\nu/\mu}\subset\bMat$ is
defined by $M\in\BL{\nu/\mu}$ if and only if\/ $\rsum{M}=\nu-\mu$, and
$\mu+\sum_{j\geq{l}}M\tr_j\in\Part$ for all~$l\in\N$, and the set
$\IL{\nu/\mu}\subset\Mat$ is defined by $M\in\IL{\nu/\mu}$ if and only if\/
$\csum{M}=\nu-\mu$, and $(\mu+\sum_{i<k}M_i)\leh(\mu+\sum_{i\leq{k}}M_i)$ for
all~$k\in\N$.

Thus for integral matrices, the Littlewood-Richardson conditions for a given
skew shape are just the tableau conditions for the same shape, but applied to
the transpose matrix. For binary matrices, the relation is as follows: if $M$
is a finite rectangular binary matrix and $M'$ is obtained from~$M$ by a
quarter turn counterclockwise, then viewing $M$~and~$M'$ as elements
of~$\bMat$ by extension with zeroes, one has $M\in\BL{\\/\kappa}$ if and only
if $M'\in\BE{\\/\kappa}$. Note that rotation by a quarter turn is not a well
defined operation on~$\bMat$, but the matrices resulting from the rotation of
different finite rectangles that contain all nonzero entries of~$M$ are all
related by the insertion or removal of some initial null rows, and such
changes do not affect membership of any set~$\BE{\\/\kappa}$ (they just give a
shift in the weight of the tableaux encoded by the matrices).

\subsection Commuting cancellations.

We can now concisely state the expressions mentioned above for the scalar
product between two skew Schur functions, which were given in \ref{alternating
Schur}. What interests us here is not so much what these expressions compute,
as the fact that one has different expressions for the same quantity. We shall
therefore not recall the definition of this scalar product
$\scalar<s_{\\/\kappa},s_{\nu/\mu}>$, but just note that the theorems
mentioned above express that value as
$\Card{\(\BE{\\/\kappa}\thru\BL{\nu/\mu}\)}$ and as
$\Card{\(\IE{\\/\kappa}\thru\IL{\nu/\mu}\)}$, respectively (the two sets
counted encode the same set of tableaux). Those theorems were derived via
cancellation from equation~\ref{alternating Schur, (\LRgensymeq)}, which
expresses the scalar product as an alternating sum over tableaux. That
equation involves a symbol~$\edgprm\alpha\\$, combinatorially defined for
$\alpha\in\comp$ and $\\\in\Part$ with values in~$\sing{-1,0,1}$. For our
current purposes the following characterisation of this symbol will suffice:
in case $\alpha$ is a partition one has $\edgprm\alpha\\=\Kr{\alpha=\\}$, and
in general if $\alpha,\alpha'\in\comp$ are related by
$(\alpha'_i,\alpha'_{i+1})=(\alpha_{i+1}-1,\alpha_i+1)$ for some~$i\in\N$, and
$\alpha'_j=\alpha_j$ for all~$j\notin\sing{i,i+1}$, then
$\edgprm\alpha\\+\edgprm{\alpha'}\\=0$ for any~$\\$. Another pair of
equations~\ref{alternating Schur, (\LRgenbineq,~\LRgeninteq)} has an opposite
relation to equation~\ref{alternating Schur, (\LRgensymeq)}, as they contain
an additional factor of the form~$\edgprm\alpha\\$ in their summand, but they
involve neither tableau conditions nor Littlewood-Richardson conditions. These
different expressions, stated in the form of summations over all binary or
integral matrices but whose range is effectively restricted by the use of the
Iverson symbol, and ordered from the largest to the smallest
effective range, are as follows. For the binary case they are
\bigdisplay
\eqalignno
{\scalar<s_{\\/\kappa},s_{\nu/\mu}>
&=\sum_{M\in\bMat}
   \edgprm{\kappa\tr+\csum{M}}{\\\tr}\edgprm{\mu+\rsum{M}}\nu
&\label(\binbruteeq)
\cr
&=\sum_{M\in\bMat}  \Kr{M\in\BE{\\/\kappa}}  \edgprm{\mu+\rsum{M}}\nu
&\label(\LRbintabeq)
\cr
&=\sum_{M\in\bMat}  \Kr{M\in\BE{\\/\kappa}}  \Kr{M\in\BL{\nu/\mu}}
, & \label(\LRbineq)
\cr
\noalign{\vskip\belowdisplayskip
\vbox{\noindent and for the integral case}
\vskip\abovedisplayskip}
 \scalar<s_{\\/\kappa},s_{\nu/\mu}>
&=\sum_{M\in\Mat} \edgprm{\kappa+\rsum{M}}\\\edgprm{\mu+\csum{M}}\nu
&\label(\intbruteeq)
\cr&=\sum_{M\in\Mat}  \Kr{M\in\IE{\\/\kappa}}  \edgprm{\mu+\csum{M}}\nu
&\label(\LRinttabeq)
\cr
&=\sum_{M\in\Mat}  \Kr{M\in\IE{\\/\kappa}}  \Kr{M\in\IL{\nu/\mu}}
.& \label(\LRinteq)
\cr
}
$$

The expressions in (\binbruteeq)--(\LRbintabeq)--(\LRbineq) as well as those
in (\intbruteeq)--(\LRinttabeq)--(\LRinteq) are related to one another by
successive cancellations: in each step one of the factors $\edgprm\alpha\\$ is
replaced by an Iverson symbol that selects only terms for which the mentioned
factor~$\edgprm\alpha\\$ already had a value~$1$; this means that all
contributions from terms for which that factor~$\edgprm\alpha\\$ has been
replaced by~$0$ cancel out against each other.

The symmetry between the tableau- and Littlewood-Richardson conditions allows
us to achieve the cancellations form (\binbruteeq)~to~(\LRbineq) and from
(\intbruteeq)~to~(\LRinteq) in an alternative way, handling the second factor
of the summand first, so that halfway those cancellations one has
\bigdisplay
\eqalignno
{
\scalar<s_{\\/\kappa},s_{\nu/\mu}>
&=\sum_{M\in\bMat}
   \edgprm{\kappa\tr+\csum{M}}{\\\tr} \Kr{M\in\BL{\nu/\mu}}
&\label(\LRbinalteq)
\cr
\noalign{\vskip\belowdisplayskip
\vbox{\noindent in the binary case, and in the integral case}
\vskip\abovedisplayskip}
 \scalar<s_{\\/\kappa},s_{\nu/\mu}>
&=\sum_{M\in\Mat} \edgprm{\kappa+\rsum{M}}\\ \Kr{M\in\IL{\nu/\mu}}
.&\label(\LRintalteq)
\cr}
$$
Indeed, the cancellation form (\binbruteeq) to~(\LRbinalteq) is performed just
like the one from (\binbruteeq) to~(\LRbintabeq) would be for matrices rotated
a quarter turn (and for $\\/\kappa$ in place of $\nu/\mu$), while the
cancellation form (\intbruteeq) to~(\LRintalteq) is performed just like the
one from (\intbruteeq) to~(\LRinttabeq) would be for the transpose matrices.
Slightly more care is needed to justify the second cancellation phase, since
the Littlewood-Richardson condition in the second factor of the summand does
not depend merely on row or column sums, as the unchanging first factor did in
the first phase. In the integral case, the second cancellation phase can be
seen to proceed like the cancellation from (\LRinttabeq)~to~(\LRinteq) with
matrices transposed, but in the binary case the argument is analogous, but not
quite symmetrical to the one used to go from (\LRbintabeq)~to~(\LRbineq). Of
course, we already knew independently of this argument that the right hand
sides of (\LRbinalteq)~and~(\LRintalteq) describe the same values as those of
(\LRbineq)~and~(\LRinteq).

Although for the two factors of the summand of (\binbruteeq)~or~(\intbruteeq)
we can thus apply cancellations to the summation in either order, and when
doing so each factor is in both cases replaced by the same Iverson symbol, the
actual method as indicated in~\ref{alternating Schur} by which terms would be
cancelled is not the same in both cases. This is so because in the double
cancellations leading from (\binbruteeq) to (\LRbineq) or from (\intbruteeq)
to~(\LRinteq), whether passing via (\LRbintabeq) respectively (\LRinttabeq) or
via (\LRbinalteq) respectively (\LRintalteq), the first phase of cancellation
has rather different characteristics than the second phase. The first phase is
a Gessel-Viennot type cancellation; it is general (in that it operates on all
terms of the initial summation) and relatively simple (it just needs to make
sure that a matrix cancels against one with the same row- or column sums). By
contrast the second phase is a Bender-Knuth type cancellation that only
operates on terms that have survived the first phase (for matrices satisfying
the pertinent tableau condition), and it has to be more careful, in order to
assure that whenever such a term is cancelled it does so against a term that
also survived the first phase.

The question that motivated the current paper is whether it is possible to
find an alternative way of defining the cancellations that has the same effect
on the summations (so we only want to change the manner in which cancelling
terms are paired up), but which has the property that the cancellation of
terms failing one of the (tableau- or Littlewood-Richardson-) conditions
proceeds in the same way, whether it is applied as the first or as the second
phase. This requires the definition of each cancellation to be general (in
case it is applied first), but also to respect the survival status for the
other cancellation (in case it is applied second). The notion of crystal
operations on matrices described below will allow us to achieve this goal. We
shall in fact see that for instance the cancellation that cancels terms not
satisfying the Littlewood-Richardson condition for $\nu/\mu$ is defined
independently of the shape $\\/\kappa$ occurring in the tableau condition; in
fact it respects the validity of the tableau condition for all skew shapes at
once.

\subsection Crystal operations for binary matrices.
\labelsec \binmatsubsec

Consider the cancellation of terms that fail the Littlewood-Richardson
condition, either going from~(\LRbintabeq) to~(\LRbineq), or
from~(\binbruteeq) to~(\LRbinalteq). Since the condition $M\in\BL{\nu/\mu}$
involves partial sums of all columns of~$M$ to the right of a given one, this
condition can be tested using a right-to-left pass over the columns of~$M$,
adding each column successively to composition that is initialised as~$\mu$,
and checking whether that composition remains a partition. If it does
throughout the entire pass, then there is nothing to do, since in particular
the final value $\mu+\rsum{M}$ will be a partition, so that
$\edgprm{\mu+\rsum{M}}\nu =\Kr{M\in\BL{\nu/\mu}} =\Kr{\mu+\rsum{M}=\nu}$. If on
the other hand the composition fails to be a partition at some point, then one
can conclude immediately that $M\notin\BL{\nu/\mu}$, so the term for~$M$
cancels. Up to this point there is no difference between a Gessel-Viennot type
cancellation and a Bender-Knuth type cancellation.

Having determined that the term for~$M$ cancels, one must find a matrix~$M'$
whose term cancels against it. The following properties will hold in all
cases. Firstly $M'$ will be obtained from~$M$ by moving entries within
individual columns, so that $\csum{M}=\csum{M'}$. Secondly, the columns of~$M$
that had been inspected at the point where cancellation was detected will be
unchanged in~$M'$, so that the term for~$M'$ is sure to cancel for the same
reason as the one for~$M$. Thirdly, a pair of adjacent rows is selected that
is responsible for the cancellation; all moves take place between these rows
and in columns that had not been inspected, with the effect of interchanging
the sums of the entries in those columns between those two rows. In more
detail, suppose $\beta$ is the first composition that failed the test to be a
partition, formed after including column~$l$ (so
$\beta=\mu+\sum_{j\geq{l}}M\tr_j$), then there is at least one index~$i$ for
which $\beta_{i+1}=\beta_i+1$; one such~$i$ is chosen in a systematic way (for
instance the minimal one) and all exchanges applied in forming~$M'$ will be
between pairs of entries $M_{i,j},M_{i+1,j}$ with $j<l$. As a result the
partial row sums $\alpha=\sum_{j<l}M\tr_j$ and $\alpha'=\sum_{j<l}(M')\tr_j$
will be related by $(\alpha'_i,\alpha'_{i+1})= (\alpha_{i+1},\alpha_i)$ (the
other parts are obviously unchanged), so that $\mu+\rsum{M}=\alpha+\beta$ and
$\mu+\rsum{M'}=\alpha'+\beta$ are related in a way that ensures
$\edgprm{\mu+\rsum{M}}\nu+\edgprm{\mu+\rsum{M'}}\nu=0$, so that the terms
of~$M$ and~$M'$ may cancel out.

Within this framework, there remains some freedom in constructing~$M'$, and
here the Gessel-Viennot and Bender-Knuth types of cancellation differ. If our
current cancellation occurs as the first phase, in other words if we are
considering the cancellation from~(\binbruteeq) to~(\LRbinalteq), then the
fact that we have ensured
$\edgprm{\mu+\rsum{M}}\nu+\edgprm{\mu+\rsum{M'}}\nu=0$ suffices for the
cancellation of the terms of~$M$ and~$M'$, and $M'$ can simply be constructed
by interchanging all pairs of bits $(M_{i,j},M_{i+1,j})$ with $j<l$, which is
what the Gessel-Viennot type cancellation does (of course such exchanges only
make any difference if the bits involved are unequal). If however our current
cancellation occurs as the second phase (so we are considering the
cancellation from~(\LRbintabeq) to~(\LRbineq)), then we must in addition make
sure that $M\in\BE{\\/\kappa}$ holds if and only if $M'\in\BE{\\/\kappa}$
does. This will not in general be the case for the exchange just described,
which is why the Bender-Knuth type of cancellation limits the number of pairs
of bits interchanged, taking into account the shape $\\/\kappa$ for which the
tableau condition must be preserved. The (easy) details of how this is done do
not concern us here, but we note that among the pairs of unequal bits whose
interchange is avoided, there are as many with their bit~`$1$' in row~$i$ as
there are with their bit~`$1$' in row~$i+1$, so that the relation between
$\beta$~and~$\beta'$ above is unaffected. The alternative construction given
below similarly ensures that $M\in\BE{\\/\kappa}$ holds if and only if
$M'\in\BE{\\/\kappa}$ does, but since it is defined independently
of~$\\/\kappa$, it works for all shapes at once, and it can be applied to any
matrix, unlike the Bender-Knuth cancellation which is defined only for
(encodings of) tableaux of shape~$\\/\kappa$.

Our fundamental definition will concern the interchange of a single pair of
distinct adjacent bits in a binary matrix; this will be vertical adjacent pair
in the discussion above, but for the cancellation of terms failing the tableau
condition we shall also use the interchange of horizontally adjacent bits. Our
definition gives a condition for allowing such an interchange, which is
sufficiently strict that at most one interchange at a time can be authorised
between a given pair of adjacent rows or columns and in a given direction
(like moving a bit~`$1$' upwards, which of course also involves a bit~`$0$'
moving downwards). Multiple moves (up to some limit) of a bit in the same
direction between the same rows or columns can be performed sequentially,
because the matrix obtained after an interchange may permit the interchange of
a pair of bits that was not allowed in the original matrix.

\proclaim Definition. \binmovedef
Let $M\in\bMat$ be a binary matrix.
\nobreak\smallskip
\defitem
A vertically adjacent pair of bits $(M_{k,l},M_{k+1,l})$ is called
interchangeable in~$M$ if $M_{k,l}\neq{M}_{k+1,l}$, and if\/
$\sum_{j=l'}^{l-1}M_{k,j}\geq\sum_{j=l'}^{l-1}M_{k+1,j}$ for all $l'\leq{l}$,
while $\sum_{j=l+1}^{l'}M_{k,j}\leq\sum_{j=l+1}^{l'}M_{k+1,j}$ for all
$l'\geq{l}$.
\nobreak\smallskip
\defitem
A horizontally adjacent pair of bits
$(M_{k,l},M_{k,l+1})$ is called interchangeable in~$M$ if
$M_{k,l}\neq{M}_{k,l+1}$ and if\/
$\sum_{i=k'}^{k-1}M_{i,l}\leq\sum_{i=k'}^{k-1}M_{i,l+1}$ for all $k'\leq{k}$,
while $\sum_{i=k+1}^{k'}M_{i,l}\geq\sum_{i=k'}^{k-1}M_{i,l+1}$ for all
$k'\geq{k}$.
\nobreak\smallskip
\enditem
Applying a upward, downward, leftward, or rightward move to~$M$ means
interchanging an interchangeable pair of bits, which is respectively of the
form $(M_{k,l},M_{k+1,l})=(0,1)$, \ $(M_{k,l},M_{k+1,l})=(1,0)$, \
$(M_{k,l},M_{k,l+1})=(0,1)$, \ or $(M_{k,l},M_{k,l+1})=(1,0)$.


These operations are inspired by crystal (or coplactic) operations, and we
shall call them \newterm{crystal operations} on binary matrices. Indeed,
horizontal moves correspond to coplactic operations (as defined in~\ref{L-R
intro, \Sec3}) applied to the concatenation of the increasing words with
weights given by the (nonzero) rows of~$M$, from top to bottom; vertical moves
correspond to coplactic operations on the concatenation of increasing words
with weights given by the columns of~$M$, taken from right to left. Applied to
the binary encoding of a semistandard tableau~$T$, vertical moves correspond
to coplactic operations on~$T$.

This definition has a symmetry with respect to rotation of matrices: if a pair
of bits in a finite binary matrix is interchangeable, then the corresponding
pair of bits in the matrix rotated a quarter turn will also be
interchangeable. However the definition does not have a similar symmetry with
respect to transposition of matrices, and this makes it a bit hard to
memorise. As a mnemonic we draw the matrices $1~0\choose0~1$ and
$0~1\choose1~0$ with a line between the pairs of bits that are not
interchangeable (and they will not be interchangeable whenever they occur in a
$2\times2$ submatrix of this form, since the conditions allowing interchange
can only get more strict when a matrix is embedded in a larger one); the pairs
not separated by a line are in fact interchangeable in the given $2\times2$
matrices:
$$
  \left( {1\atop 0} \kern4pt \vrule \kern4pt {0\atop 1} \right)
, \kern0.4\hsize
  {0\quad1\overwithdelims()1\quad0}
.\nn
$$
As a somewhat larger example, consider vertical moves in the binary matrix
$$M=
\pmatrix{1&0&0&1&0&1&1&0&0&0&0&0&1\cr
         0&1&1&1&1&0&0&1&1&0&1&1&1\cr
         0&1&0&1&0&0&1&1&1&0&1&0&1\cr}
.\nn
$$
The pair $1\choose0$ at the top right is interchangeable, because in every
initial part of the remainder of rows $0$~and~$1$, the pairs~$0\choose1$ are
at least as numerous as the pairs $1\choose0$. Since they are in fact always
strictly more numerous, the pair $0\choose1$ in column~$1$ is also
interchangeable (the closest one comes to violating the second inequality in
\binmovedef\defitemnr1 is the equality
$\sum_{j=2}^6M_{0,j}=3=\sum_{j=2}^6M_{1,j}$, and the first inequality poses no
problems). None of the remaining pairs in rows $0$~and~$1$ are interchangeable
however; for the pair in column~$2$ the first inequality
in~\binmovedef\defitemnr1 fails for~$l'=1$ since $M_{0,1}=0\not\geq
M_{1,1}=1$, and in fact this inequality continues to fail for~$l'=1$ and all
further columns (often there are other inequalities that fail as well, but one
may check that for column~$7$ the mentioned inequality is the only one that
fails). In rows $1$~and~$2$, only the pair $1\choose0$ in column~$4$ is
interchangeable (while all inequalities are also satisfied for columns~$5$
and~$12$, these columns contain pairs of equal bits $0\choose0$ and
$1\choose1$, which are never interchangeable). As an example of successive
moves in the same direction, one may check that, in rows $0$~and~$1$, after
interchanging the pair~$0\choose1$ in column~$1$, one may subsequently
interchange similar pairs in columns $7$, $8$, $10$, and~$11$, in that order.

Let us now show our claim that at most one move at a time is possible between
any given pair of rows or columns and in any given direction. Consider the
case of adjacent rows, $i,i+1$ and suppose they contain two interchangeable
vertically adjacent pairs of bits in columns $j_0<j_1$. Then one has two
opposite inequalities for the range of intermediate columns, which implies
that $\sum_{j=j_0+1}^{j_1-1}M_{i,j} =\sum_{j=j_0+1}^{j_1-1}M_{i+1,j}$. One can
also see that the interchangeable pair in column~$j_0$ is~$1\choose0$ and the
one in column~$j_1$ is~$0\choose1$, since any other values would contradict
definition~\binmovedef. So there can be at most one downward move and at most
one upward move that can be applied between rows $i$ and~$i+1$, with the
downward move being to the left of the upward move if both occur. Similarly,
at most one leftward move and at most one rightward move can be applied
between a given pair of adjacent columns, and if both can, the leftward move
is in a row above of the rightward move.

These uniqueness statements justify the following crucially important
definition. In keeping with the usual notation for crystal operations, we use
the letter~$e$ for raising operations and the letter~$f$ for lowering
operations, but since we have a horizontal and a vertical variant of either
one, we attach an arrow pointing in the direction in which the bit~`$1$'
moves.

\proclaim Definition.
(binary raising and lowering operations) Let $M\in\bMat$.
\defitem
    If $M$ contains an interchangeable pair of bits in rows~$i$~and~$i+1$,
    then the matrix resulting from the interchange of these bits is denoted
    by~$\upm_i(M)$ if the interchange is an upward move, or by~$\dnm_i(M)$ if
    the interchange is a downward move. If for a given $i\in\N$ the matrix~$M$
    admits no~upward or no~downward move interchanging any pair bits in
    rows~$i$~and~$i+1$, then the expression $\upm_i(M)$ respectively
    $\dnm_i(M)$ is undefined.
\defitem
    If $M$ contains an interchangeable pair of bits in columns~$j$~and~$j+1$,
    then the matrix resulting from the interchange of these bits is denoted
    by~$\ltm_l(M)$ if the interchange is a leftward move, or by~$\rtm_l(M)$ if
    the interchange is a rightward move. If for a given $j\in\N$ the
    matrix~$M$ admits no~leftward or no~rightward move interchanging any pair
    bits in columns~$j$~and~$j+1$, then the expression $\ltm_l(M)$
    respectively $\rtm_l(M)$ is undefined.

Since an interchangeable pair of bits remains so after it has been
interchanged, it follows that whenever $\upm_i(M)$ is defined then so is
$\dnm_i(\upm_i(M))$, and it is equal to $M$. Similarly each of the expressions
$\upm_i(\dnm_i(M))$, \ ,$\ltm_j(\rtm_j(M))$ and $\rtm_j(\ltm_j(M))$ is defined
as soon as its inner application is, in which case it designates~$M$. Our next
concern will be characterising when expressions such as $\upm_i(M)$ are
defined, and more generally determining the number of times each of the
operations $\upm_i$, $\dnm_i$, $\ltm_j$ and $\rtm_j$ can be successively
applied to a given matrix~$M$, which we shall call the potential of~$M$ for
these operations.

\proclaim Definition. \nmdefs
For $M\in\bMat$ and $i,j\in\N$, the numbers
$\num_i(M),\ndm_i(M),\nlm_j(M),\nrm_j(M)\in\N$ are defined by
$$\eqalignno
{ \num_i(M)
&=\max\setof\textstyle\sum_{j\geq{l}}(M_{i+1,j}-M_{i,j}):l\in\N\endset
,& \label(\numdefeq)
\cr
  \ndm_i(M)
&=\max\setof\textstyle\sum_{j<l}(M_{i,j}-M_{i+1,j}):l\in\N\endset
,& \nn
\cr
  \nlm_j(M)
&=\max\setof\textstyle\sum_{i<k}(M_{i,j+1}-M_{i,j}):k\in\N\endset
,& \label(\nlmdefeq)
\cr
  \nrm_j(M)
&=\max\setof\textstyle\sum_{i\geq{k}}(M_{i,j}-M_{i,j+1}):k\in\N\endset
.& \label(\nrmdefeq)
\cr
}
$$

\proclaim Proposition. \nbinmvprop
For $M\in\bMat$ and $i,j\in\N$, the numbers of times each of $\upm_i$,
$\dnm_i$, $\ltm_j$, and $\rtm_j$ can be successively applied to~$M$ are
respectively given by the numbers $\num_i(M)$, $\ndm_i(M)$, $\nlm_j(M)$, and
$\nrm_j(M)$. Moreover $\ndm_i(M)-\num_i(M)=\rsum{M}_i-\rsum{M}_{i+1}$ and
$\nrm_j(M)-\nlm_j(M)=\csum{M}_j-\csum{M}_{j+1}$.

\proof
Suppose $M$ is an $n\times{m}$ binary matrix (so all entries outside that
rectangle are zero) that admits an upward move interchanging a pair of bits
$0\choose1$ in column~$l$ of rows $i,i+1$. Then it follows from
$\sum_{j=l+1}^{m-1}(M_{i+1,j}-M_{i,j})\geq0$ that
$\num_i(M)\geq\sum_{j\geq{l}}(M_{i+1,j}-M_{i,j})>0$. Conversely if
$\num_i(M)>0$, then let~$l<m$ be the maximal index for which the maximal value
of\/ $\sum_{j\geq{l}}(M_{i+1,j}-M_{i,j})$ is attained. One then verifies that
$M$ admits an upward move in column~$l$ of rows $i,i+1$: the fact that the
pair in that position is $0\choose1$ follows from the maximality of~$l$, and
failure of one of the inequalities in \binmovedef\defitemnr1 would
respectively give a value $l'<l$ for which a strictly larger sum is obtained,
or a value $l'+1>l$ for which a weakly larger sum is obtained, either of which
contradicts the choice of~$l$.

The statement concerning $\upm_i$ can now be proved by induction
on~$\num_i(M)$. For $\num_i(M)=0$ we have just established that no upward
moves in rows $i,i+1$ are possible. So suppose $\num_i(M)>0$ and let $M'$ be
obtained from~$M$ by an upward move in column~$l_0$. Then replacing $M$ by
$M'$ decreases the sums $\sum_{j\geq{l}}(M_{i+1,j}-M_{i,j})$ by~$2$ for all
$l\leq{l_0}$, while those sums are unchanged for $l>l_0$. The sums for
$l\leq{l_0}$ therefore become at most $\num_i(M)-2$, while the sums for
$l>l_0$ remain at most $\num_i(M)-1$ (since $l_0$ was the maximal index for
which the value~$\num_i(M)$ is attained for~$M$, as we have seen). Therefore
the maximal sum for~$M'$ is attained for the index $l_0+1$, and its value is
$\num_i(M') =\sum_{j\geq{l_0+1}}(M_{i+1,j}-M_{i,j})=\num_i(M)-1$; by induction
$\upm_i$ can be applied precisely that many times to~$M'$, and so it can be
applied $\num_i(M)$ times to~$M$ as claimed. The statements for $\ltm_j$,
$\dnm_i$, and $\rtm_j$ follow from the statement we just proved by considering
the (finite) matrices obtained from~$M$ by turning it one, two, or three
quarter turns. The statements in the final sentence of the proposition are
clear if one realises that for instance $\sum_{j<l}(M_{i,j}-M_{i+1,j})$ and
$\sum_{j\geq{l}}(M_{i+1,j}-M_{i,j})$ differ by $\rsum{M}_i-\rsum{M}_{i+1}$
independently of~$l$, so that their maxima $\ndm_i(M)$ and $\num_i(M)$ are
attained for the same (set of) values of~$l$, and also differ by
$\rsum{M}_i-\rsum{M}_{i+1}$.
\QED

With respect to the possibility of successive moves between a pair of adjacent
rows or columns, we can make a distinction between pairs whose interchange is
forbidden in~$M$ but can be made possible after some other exchanges between
those rows or columns, and pairs whose interchange will remain forbidden
regardless of such exchanges. We have seen that when a move is followed by a
move in the opposite direction, the latter undoes the effect of the former; it
follows that if a given move can be made possible by first performing one or
more moves between the same pair of rows or columns, then one may assume that
all those moves are in the same direction. Moreover we have seen for instance
that successive upward moves between two rows always occur from left to right;
this implies that if a pair $0\choose1$ in column~$l$ of rows~$i,i+1$ is not
interchangeable due to a failure of some instance of the second inequality
in~\binmovedef\defitemnr1 (which only involves columns~$j>l$), then this
circumstance will not be altered by any preceding upward moves between the
same rows, and the move will therefore remain forbidden. On the other hand if
the second inequality in~\binmovedef\defitemnr1 is satisfied for all~$l'>l$,
then the value of $\sum_{j\geq{l}}(M_{i+1,j}-M_{i,j})$ is larger than the one
obtained by replacing the bound~$l$ by any~$l'>l$; it may still be less than
that overall maximum~$\num_i(M)$, but that value can be lowered by successive
upward moves between rows~$i,i+1$, which must necessarily occur in
columns~$j<l$, until the pair $0\choose1$ considered becomes interchangeable.

We may therefore conclude that, in the sense of repeated moves between two
adjacent rows, failure of an instance of the first inequality
in~\binmovedef\defitemnr1 gives a temporary obstruction for a candidate upward
move, while failure of an instance of the second inequality gives a permanent
obstruction. For candidate downward moves the situation is reversed. The
following alternative description may be more intuitive. If one represents
each pair $0\choose1$ by ``$($'', each pair $1\choose0$ by ``$)$'', and all
remaining pairs by ``$-$'' (or any non-parenthesis symbol), then for all
parentheses that match another one in the usual sense, the pairs in the
corresponding columns are permanently blocked. The remaining unmatched
parentheses have the structure ``$)\cdots)(\cdots($'' of a sequence of right
parentheses followed by a sequence of left parentheses (either of which might
be an empty sequence). An upward move between these rows is possible in the
column corresponding to the leftmost unmatched~``('' if it exists, an a
downward move between these rows is possible in the column corresponding to
the rightmost unmatched~``)'' if it exists. In either case the move replaces
the parenthesis by an opposite one, and since it remains unmatched, we can
continue with the same description for considering subsequent moves. In this
view it is clear that all unmatched parentheses can be ultimately inverted,
and that upward moves are forced to occur from left to right, and downward
moves from right to left. For instance, in the $3\times13$ matrix given as an
example after definition~\binmovedef, the sequence of symbols for the two
topmost rows is ``$)\,(\,({-}(\,)\,)\,(\,({-}(\,({-}$'', and from this it is
clear that one downward move is possible in column~$0$, or at most $5$
successive upward moves in columns $1$, $7$, $8$, $10$, and~$11$; for the
bottommost rows we have the sequence ``${-}{-}){-}){-}({-}{-}{-}{-}){-}$'' and
only successive downward moves are possible, in columns $4$~and~$2$. For moves
between adjacent columns the whole picture described here must be rotated a
quarter turn (clockwise or counterclockwise, this makes no difference).

We now consider the relation of the definitions above to the tableau- and
Littlewood-Richardson conditions on matrices. The first observation is that
these conditions can be stated in terms of the potentials for raising (or for
lowering) operations.

\proclaim Proposition. \bincondcharprop
Let $M\in\bMat$ and let $\\/\kappa$ and $\mu/\nu$ be skew shapes.
\statitem $M\in\BE{\\/\kappa}$ if and only if $\csum{M}=\\\tr-\kappa\tr$ and
          $\nlm_j(M)\leq\kappa\tr_j-\kappa\tr_{j+1}$ for all~$j\in\N$.
\statitem $M\in\BL{\nu/\mu}$ if and only if $\rsum{M}=\nu-\mu$ and
          $\num_i(M)\leq\mu_i-\mu_{i+1}$ for all~$i\in\N$.
\enditem
The second parts of these conditions can also be stated in terms of the
potentials of~$M$ for lowering operations, as
$\nrm_j(M)\leq\\\tr_j-\\\tr_{j+1}$ for all~$j\in\N$, respectively as
$\ndm_i(M)\leq\nu_i-\nu_{i+1}$ for all~$i\in\N$.

\proof
In view of the expressions in definition~\nmdefs, these statements are just
reformulations of the parts of proposition~\Tabprop\ and definition~\LRdef\
that apply to binary matrices.
\QED

The next proposition shows that vertical and horizontal crystal operations on
matrices respect the tableau conditions respectively the Littlewood-Richardson
conditions for all skew shapes at once.

\proclaim Proposition. \bininvarprop
If binary matrices $M,M'\in\bMat$ are related by $M'=\upm_i(M)$ for some
$i\in\N$, then $\nlm_j(M)=\nlm_j(M')$ and $\nrm_j(M)=\nrm_j(M')$ for
all~$j\in\N$. Consequently, the conditions $M\in\BE{\\/\kappa}$ and
$M'\in\BE{\\/\kappa}$ are equivalent for any skew shape~$\\/\kappa$. Similarly
if $M$~and~$M'$ are related by $M'=\rtm_j(M)$ for some $j\in\N$, then
$\num_i(M)=\num_i(M')$ and $\ndm_i(M)=\ndm_i(M')$ for all~$i\in\N$, and
$M\in\BL{\nu/\mu}\iff M'\in\BL{\nu/\mu}$ for any skew shape~$\nu/\mu$.

\proof
It suffices to prove the statements about $M'=\upm_i(M)$, since those
concerning $M'=\ltm_j(M)$ will then follow by applying the former to matrices
obtained by rotating $M$~and~$M'$ a quarter turn counterclockwise. For the
case considered it will moreover suffice to prove $\nlm_j(M)=\nlm_j(M')$ for
any~$j\in\N$, since $\nrm_j(M)=\nrm_j(M')$ will then follow from
$\csum{M}=\csum{M'}$, and the equivalence of $M\in\BE{\\/\kappa}$ and
$M'\in\BE{\\/\kappa}$ will be a consequence of proposition~\bincondcharprop.
One may suppose that the pair of bits being interchanged to obtain $M'$
from~$M$ is in column $j$~or~$j+1$, since otherwise $\nlm_j(M)=\nlm_j(M')$ is
obvious from~(\nlmdefeq). Let $(p_k)_{k\in\N}$ be the sequence of partial sums
for~$M$ of which $\nlm_j(M)$ is the maximum, in other words
$p_k=\sum_{i'<k}(M_{i',j+1}-M_{i',j})$, and let $(p'_k)_{k\in\N}$ be the
corresponding sequence for~$M'$. Then the only index~$k$ for which $p_k\neq
p'_k$ is $k=i+1$: one has $p'_{i+1}=p_{i+1}-1$ if the move occurred in
column~$j$, or $p'_{i+1}=p_{i+1}+1$ if it occurred in column~$j+1$. The only
way in which this change could make $\nlm_j(M)=\max_kp_k$ differ from
$\nlm_j(M')=\max_kp'_k$ is if $k=i+1$ were the unique index for which
$p_k=\nlm_j(M)$ (in the former case) or for which $p'_k=\nlm_j(M')$ (in the
latter case). That would in particular require that the indicated value be
strictly larger than $p_k=p'_k$ and than $p_{k+2}=p'_{k+2}$, so $M$~or~$M'$
would have to contain a submatrix $0~1\choose1~0$ at the intersection of rows
$i,i+1$ and columns $j,j+1$, while the other matrix would differ by the
interchange of one of those two vertically adjacent pairs of bits. But we have
seen that in such a submatrix neither of those two pairs of bits can be
interchangeable, which excludes this possibility, and one therefore has
$\nlm_j(M)=\nlm_j(M')$ in all cases.
\QED

One can summarise the last two propositions as follows: Littlewood-Richardson
conditions can be stated in terms of the potentials for vertical moves, which
moves preserve tableau conditions, while tableau conditions can be stated in
terms of the potentials for horizontal moves, which moves preserve
Littlewood-Richardson conditions.

We shall now outline the way in which crystal operations can be used to define
cancellations either of terms for matrices not in $\BE{\\/\kappa}$ or of those
not in $\BL{\nu/\mu}$, in the summations of~(\binbruteeq), (\LRbintabeq), or
(\LRbinalteq). One starts by traversing each matrix~$M$ as before, searching
for a violation of the condition in question, and of an index that witnesses
it; this amounts to finding a raising operation~$e$ (i.e., some
$\ltm_j$~or~$\upm_i$) for which the potential of~$M$ is larger than allowed by
the second part of proposition~\bincondcharprop\ \statitemnr1~or~\statitemnr2.

Now consider the set of matrices obtainable from~$M$ by a sequence of
applications of~$e$ or of its inverse lowering operation~$f$; these form a
finite ``ladder'' in which the operation~$e$ moves up, and $f$ moves down.
Note that the potential for~$e$ increases as one \emph{descends} the ladder.
The condition of having a too large a potential for~$e$ determines a lower
portion containing~$M$ of the ladder, for which all corresponding terms must
be cancelled, and the witness chosen for such a cancellation will be the same
one as chosen for~$M$ (there may also be terms cancelled in the remaining
upper portion of the ladder, but their witnesses will be different). Now the
negation of the expression of the form~$\edgprm\alpha\\$, which is required
for cancellation, can be obtained by reversing the indicated lower part of the
ladder. Since a pair of matrices whose terms cancel are thus linked by a
sequence of horizontal or vertical moves, their status for any
Littlewood-Richardson respectively tableau condition (the kind for which one
is not cancelling) will be the same, which allows this cancellation to be used
as a second phase (starting from (\LRbinalteq)~or~(\LRbintabeq)).

Let us fill in the details of the description above, for the cancellation of
terms for matrices not in $\BL{\nu/\mu}$, in other words leading from
(\binbruteeq) to~(\LRbinalteq) of from (\LRbintabeq) to~(\LRbineq). As
described at the beginning of this subsection, we start by finding the maximal
index~$l$ such that the composition $\beta=\mu+\sum_{j\geq{l}}M\tr_j$ is not
a partition, and choosing an index~$i$ for which $\beta_{i+1}=\beta_i+1$;
this implies that $\num_i(M)>\mu_i-\mu_{i+1}$, so the potential of~$M$
for~$e=\upm_i$ exceeds the limit given in \bincondcharprop~\statitemnr2. For
convenience let us use the notation $e^d$ to stand for $(\upm_i)^d$ when
$d>0$, for the identity when $d=0$, and for $(\dnm_i)^{-d}$ when $d<0$; then
the ladder mentioned above is $\setof{e^d(M)}:
-\ndm_i(M)\leq{d}\leq\num_i(M)\endset$, and its lower part that must be
cancelled because of a too large potential for~$\upm_i$ is $\setof{e^d(M)}:
-\ndm_i(M)\leq{d}<\num_i(M)-(\mu_i-\mu_{i+1})\endset$.

From the maximality of~$l$ it follows that~$M$ contains a pair $0\choose1$ in
column~$l$ of rows~$i,i+1$, and that this pair is not permanently blocked for
upward moves in rows~$i,i+1$ (in other words, one has $\sum_{j=l+1}^m
M_{i,j}\leq\sum_{j=l+1}^m M_{i+1,j}$ for all $m>l$); indeed the pair will be
interchanged in $e^d(M)$ when $d$ satisfies
$\num_i(M)-(\mu_i-\mu_{i+1})\leq{d}\leq\num_i(M)$, i.e., in the mentioned
upper part of the ladder. So the lower part of the ladder is precisely the
part in which that pair is not interchanged, and the matrices in this part
will give rise to the same indices $l$~and~$i$ as~$M$ to witness their
cancellation. The expression for~$d$ such that the term for~$M$ cancels
against the one for~$e^d(M)$ can be found as follows. If $M$ is at the bottom
of the ladder ($\ndm_i(M)=0$) then $d$ has the value
$\num_i(M)-(\mu_i-\mu_{i+1})-1$ that gives the topmost value of the bottom
part of the ladder, and $d$ decreases with the level $\ndm_i(M)$ of~$M$, so
the expression is $d=\num_i(M)-(\mu_i-\mu_{i+1})-1-\ndm_i(M)$. Putting
$\alpha=\mu+\rsum{M}$, this can also be written as
$d=\alpha_{i+1}-\alpha_i-1$, by proposition~\nbinmvprop. Since each
application of~$\upm_i$ increases the sum of entries in row~$i$ while
decreasing the sum in row~$i+1$, the value $\alpha'=\mu+\rsum{M'}$ for the
matrix $M'=e^d(M)$ satisfies
$(\alpha'_i,\alpha'_{i+1})=(\alpha_{i+1}-1,\alpha_i+1)$ while its remaining
components are unchanged from~$\alpha$, which ensures that
$\edgprm\alpha\nu+\edgprm{\alpha'}\nu=0$. The fact that $M$ and~$M'$ are
related by vertical moves implies that $M\in\BE{\\/\kappa}\iff
M'\in\BE{\\/\kappa}$ for any skew shape $\\/\kappa$, so the terms for
$M$~and~$M'$ do indeed cancel, whether we are considering the passage
from~(\binbruteeq) to~(\LRbinalteq) or the one from (\LRbintabeq)
to~(\LRbineq).

For the cancellations involved in passing from~(\binbruteeq) to~(\LRbintabeq)
and from~(\LRbinalteq) to~(\LRbineq) the description is similar, but rotated a
quarter turn counterclockwise: the initial scan of the matrix is by rows from
top to bottom, and the raising operations~$\num_i$ are replaced by raising
operations~$\nlm_j$.

The reader may have been wondering whether we have been going through all
these details just to obtain more aesthetically pleasing descriptions of the
reductions (\binbruteeq)$\to$(\LRbintabeq)$\to$(\LRbineq) and
(\binbruteeq)$\to$(\LRbinalteq)$\to$(\LRbineq) (and maybe the reader even
doubts whether that goal was actually obtained). But crystal operations turn
out to be useful in other ways than just to define cancellations, and several
such applications will be given below; those applications alone largely
justify the definition of crystal operations. We have nevertheless chosen to
introduce them by considering cancellations, because that provides a
motivation for the precise form of their definition and for treating separate
cases for binary and integral matrices; such motivation might otherwise not be
evident. For our further applications it is of crucial importance that
horizontal and vertical moves are compatible in a stronger sense than
expressed in proposition~\bininvarprop. Not only do moves in one direction
leave invariant the potentials for all moves in perpendicular directions, they
actually commute with those moves, as stated in the following lemma.

\proclaim Lemma. \bclemma (binary commutation lemma)
Let $M,M',M''\in\bMat$ be related by $M'=\upm_i(M)$ and $M''=\ltm_j(M)$ for
some $i,j\in\N$; then $\ltm_j(M')=\upm_i(M'')$. The same holds when $\upm_i$
is replaced both times by $\dnm_i$ and/or $\ltm_j$ is replaced both times
by~$\rtm_j$.

\proof
Note that the expressions $\ltm_j(M')$ and $\upm_i(M'')$ are defined since
$\nlm_j(M')=\nlm_j(M)>0$ and $\num_i(M'')=\num_i(M)>0$ by
proposition~\bininvarprop. The variants given in the second part of the lemma
can be deduced from the initial statement either by rotation symmetry or by a
suitable change of roles between the four matrices involved. So we shall focus
on proving the initial statement.

Suppose first that the pairs of bits interchanged in the moves
$\upm_i\:M\mapsto{M'}$ and $\ltm_j\:M\mapsto{M''}$ are disjoint. In this case
we shall argue that these pairs of bits are in the same position as the pairs
of bits interchanged in the moves $\upm_i\:M''\mapsto\upm_i(M'')$ and
$\ltm_j\:M'\mapsto\ltm_j(M')$, respectively; then it will be obvious that
$\ltm_j(\upm_i(M))=\upm_i(\ltm_j(M))$. To this end must show that the
conditions in definition~\binmovedef, which are satisfied in~$M$ for each of
the two pairs of bits considered, remain valid after the other pair is
interchanged. Since the values of one pair of bits is not affected by the
interchange of the other pair, we only need to worry about the four
inequalities in that definition. Depending on the relative positions of the
two pairs, at most one of those inequalities can have an instance for which
the values being compared change, but since we do not know which one, this
does not help us much; nevertheless the four cases are quite similar, so we
shall treat only the first one explictly. Each inequality, with its
quantification, can be reformulated as stating that some maximum of partial
sums does not exceed~$0$ (actually it equals~$0$); for instance the first
inequality is equivalent to `$\max\setof\sum_{j=l'}^{l-1}(M_{k+1,j}-M_{k,j}):
0\leq{l'}\leq{l}\endset\leq0$' (this condition applies for $k=i$ if the move
of $\upm_i$ occurs in column~$l$). That maximum of partial sums is of the same
type as the one in one of the equations (\numdefeq)--(\nrmdefeq), but for a
truncated matrix; in the cited case they are the partial sums of (\numdefeq)
but computed for~$M$ truncated to its columns~$j<l$. Therefore the same
reasoning as in the proof of proposition~\bininvarprop\ shows that although
one of the partial sums may change, their maximum remains the same, so that
the pair of bits considered remains interchangeable.

Now suppose that to the contrary the pairs of bits $0\choose1$ and
$(0\enspace1)$ being interchanged in~$M$ do overlap. Then after performing one
interchange, the pair of bits in the position of the other pair can no longer
be interchangeable, as its bits will have become equal. There is a unique
$2\times2$ submatrix of~$M$ that contains the two overlapping pairs, and since
it contains both a vertical and a horizontal interchangeable pair of bits, its
value can be neither $0~1\choose1~0$ nor $1~0\choose0~1$. Therefore it will
contain either $0~1\choose1~1$ if the two pairs overlap in their bit~`$0$' (at
the top left), or $0~0\choose0~1$ if the two pairs overlap in their bit~`$1$'
(at the bottom right). In either case it is not hard to see that the
overlapping bit, after having been interchanged horizontally or vertically, is
again (in its new position) part of an interchangeable pair within the
$2\times2$ submatrix, in the direction perpendicular to the move performed;
the other bit of that pair is the one in the corner diametrically opposite to
the old position of the overlapping bit in the submatrix considered (the
bottom right corner in the former case and the top left corner in the latter
case). This is so because comparing that new pair with the interchangeable
pair that used to be in the remaining two squares of the $2\times2$ submatrix,
the only difference for each of the pertinent inequalities of
definition~\binmovedef\ is the insertion or removal of a bit with the same
value in each of the two sums being compared, which does not affect the result
of the comparison. Therefore the succession of two raising operations, applied
in either order, will transform the submatrix $0~1\choose1~1$ into
$1~1\choose1~0$, or the submatrix $0~0\choose0~1$ into $1~0\choose0~0$, as
illustrated below.
\QED

\bigdisplay
\abovedisplayshortskip=\bigskipamount
\def\mapleft#1{\smash
       {\mathop{\hbox to 2.5em{\leftarrowfill}}\limits^{\textstyle#1}}}
  \matrix
  { \pmatrix{1&1\cr1&0\cr} & \mapleft{\ltm_j} & \pmatrix{1&1\cr0&1\cr} \cr
    \llap{$\upm_i$}\bigg\uparrow & & \llap{$\upm_i$}\bigg\uparrow \cr
    \pmatrix{1&0\cr1&1\cr} & \mapleft{\ltm_j} & \pmatrix{0&1\cr1&1\cr} \cr
  }
\kern 0.15 \hsize
  \matrix
  { \pmatrix{1&0\cr0&0\cr} & \mapleft{\ltm_j} & \pmatrix{0&1\cr0&0\cr} \cr
    \llap{$\upm_i$}\bigg\uparrow & & \llap{$\upm_i$}\bigg\uparrow \cr
    \pmatrix{0&0\cr1&0\cr} & \mapleft{\ltm_j} & \pmatrix{0&0\cr0&1\cr} \cr
  }
\nn
$$

\subsection Crystal operations for integral matrices.
\labelsec \intmatsubsec

Motivated by the existence of cancellations
(\intbruteeq)$\to$(\LRinttabeq)$\to$(\LRinteq) and
(\intbruteeq)$\to$(\LRintalteq)$\to$(\LRinteq), we shall now define operations
like those defined in the previous subsection, but for integral instead of
binary matrices. Much of what will be done in this subsection is similar to
what was done before, so we shall focus mainly on the differences with the
situation for binary matrices.

A first difference is the fact that for integral matrices the operation of
interchanging adjacent entries is too restrictive to achieve the desired kind
of modifications. We shall therefore regard each matrix entry~$m$ as if it
were a pile of $m$ units, and the basic type of operation will consist of
moving some units from some entry $m>0$ to a neighbouring entry, which amounts
to decreasing the entry~$m$ and increasing the neighbouring entry by the same
amount. We shall call this a transfer between the two entries; as in the
binary case we shall impose conditions for such a transfer to be allowed.
Another difference is the kind of symmetry implicitly present in equation
(\LRinteq) compared with (\LRbineq), which in fact stems from the difference
between the cases of binary and integral matrices in the relation of
definition~\LRdef\ to proposition~\Tabprop, which we already observed
following that definition. As a consequence the rules for allowing transfers
will not be symmetric with respect to rotation by a quarter turn, but instead
they will be symmetric with respect to transposition of the integral matrices
and with respect to rotation by a half turn.

This new type of symmetry somewhat simplifies the situation, but there is also
a complicating factor, due to the fact that the tableau- and
Littlewood-Richardson conditions are more involved for integral matrices than
for binary ones. In the binary case it sufficed to construct a sequence of
compositions by cumulating rows or columns, and to test each one individually
for being a partition. But in the integral case one must test for each pair
$\alpha,\beta$ of successive terms in the sequence whether $\alpha\leh\beta$,
in other words whether $\beta/\alpha$ is a horizontal strip. That test amounts
to verifying $\beta_{i+1}\leq\alpha_i$ for all~$i$, since
$\alpha_i\leq\beta_i$ already follows from the circumstance that $\beta_i$ is
obtained by adding a matrix entry to~$\alpha_i$. Thus if we focus on the
inequalities involving the parts $i$~and~$i+1$ of the compositions in the
sequence, then instead of just checking that part~$i+1$ never exceeds part~$i$
of the same composition, one must test the stronger requirement that
part~$i+1$ of the \emph{next} partition in the sequence still does not exceed
that (old) part~$i$.

This will mean for the analogues of definitions \binmovedef~and~\nmdefs, that
the final entries in partial sums in two adjacent rows or columns will not be
in the same column~or~row, but in a diagonal position with respect to each
another (always in the direction of the main diagonal). This also means that
the conditions required to allow a transfer must take into account some of the
units that are present in the matrix entries between which the transfer takes
place, but which are not being transferred themselves (in the binary case no
such units exist). Although the precise form of the following definition could
be more or less deduced from the properties we seek, we shall just state it,
and observe afterwards that it works.

\proclaim Definition. \intmovedef
Let $M\in\Mat$, $k,l\in\N$, and $a\in\Z-\{0\}$.
\defitem Suppose that $\sum_{j=l'}^{l-1}(M_{k+1,j+1}-M_{k,j})\geq\max(a,0)$
    for all $l'<l$, or if $l=0$ that $M_{k+1,0}\geq{a}$, and suppose
    $\sum_{j=l}^{\smash{l'}-1}(M_{k,j}-M_{k+1,j+1})\geq\max(-a,0)$ for all
    $l'>l$. Then we allow the entries $(M_{k,l},M_{k+1,l})$ to be replaced by
    $(M_{k,l}+a,M_{k+1,l}-a)$; this is called an upward transfer of $a$ units
    between rows $k$~and~$k+1$ if $a>0$, or a downward transfer of $-a$
    units between those rows if $a<0$.
\defitem Suppose that $\sum_{i=k'}^{k-1}(M_{i+1,l+1}-M_{i,l})\geq\max(a,0)$
    for all $k'<k$, or if $k=0$ that $M_{0,l+1}\geq{a}$, and suppose
    $\sum_{i=k}^{\smash{k'}-1}(M_{i,l}-M_{i+1,l+1})\geq\max(-a,0)$ for all
    $k'>k$. Then we allow the entries $(M_{k,l},M_{k,l+1})$ to be replaced by
    $(M_{k,l}+a,M_{k,l+1}-a)$; this is called a leftward transfer of $a$ units
    between columns $l$~and~$l+1$ if $a>0$, or a rightward
    transfer of $-a$ units between those columns if $a<0$.

\heading "Remarks"
(1)~The occurrence of the quantity~$a$ in the inequalities has the effect of
   cancelling its contribution to the entry from which it would be
   transferred. It follows that the transfer can always be followed by a
   transfer of~$a$ units in the opposite sense between the same entries, which
   reconstructs the original matrix.
(2)~The exceptional conditions $M_{0,l+1}\geq{a}$ and $M_{k+1,0}\geq{a}$
   compensate for the absence of any inequality where $a$ occurs in the way
   just mentioned. They serve to exclude the introduction of negative entries
   by a transfer; note that for instance this is taken care of for upward
   moves with $l>0$ by the condition $M_{k+1,l}-M_{k,l-1}\geq{a}$, and for
   downward moves by $M_{k,l}-M_{k+1,l+1}\geq-a$. Hence the cases $k=0$ and
   $l=0$ are not treated any differently from the others.
(3)~We could have restricted the definition to the cases $a=1$ and $a=-1$,
   since transfers with $|a|>1$ can be realised by repeated transfers with
   $|a|=1$. The current formulation was chosen for readability, and because it
   permits a straightforward and interesting generalisation to matrices with
   non-negative real coefficients.

We shall call these transfers crystal operations on integral matrices. It can
be seen that horizontal transfers correspond to $|a|$ successive coplactic
operations on the word formed by concatenating weakly decreasing words whose
weights are given by the (nonzero) rows of~$M$, taken from top to bottom;
vertical transfers correspond $|a|$ successive to coplactic operations on the
word similarly formed by concatenating weakly decreasing words with weights
given by the columns of~$M$, taken from left to right. Horizontal transfers in
the integral encoding of a semistandard tableau~$T$ correspond to coplactic
operations on~$T$.

Here is a small example of vertical transfers; for an example of horizontal
transfers one can transpose the matrix. Consider vertical moves between the
two nonzero rows of the integral matrix
$$
  M=\pmatrix
{ 1&2&1&3&3&1&2&4&0\cr
  2&1&1&4&2&0&5&2&0\cr
}.\nn
$$
For the three leftmost columns no transfer is possible because
$\sum_{j=l}^2(M_{0,j}-M_{1,j+1})<0$ for $l\in\sing{0,1,2}$. For column~$3$
however an upward transfer of at most $2$ units is possible, where the limit
comes from the value $\sum_{j=1}^2(M_{1,j+1}-M_{0,j})=2$ (and starting the
sum at~$j=0$ would give the same value). No downward transfer is possible in
that column because $\smash{\sum_{j=3}^5}(M_{0,j}-M_{1,j+1})=0$; however a
downward move of at most~$4$ units is possible in column~$7$. If all~$4$ units
are  transferred downwards in that column it is replaced by $0\choose6$, and
no further downwards between these rows will be possible; however if $2$~units
are transferred upwards in column~$3$ one obtains the matrix
$$
  M'=\pmatrix
{ 1&2&1&5&3&1&2&4&0\cr
  2&1&1&2&2&0&5&2&0\cr
},\nn
$$
and an upwards transfer in~$M'$ is still possible, namely in column~$0$ where
both units from row~$1$ can be transferred to row~$0$.

We now consider the ordering of transfers between a given pair of adjacent
rows or columns. Since any simultaneous transfer of more than one unit can be
broken up into smaller transfers, we cannot claim that at most one transfer in
a given direction between a given pair of rows or columns is possible, but it
is important that the only choice concerns the amount transferred, not the
pair of entries where the transfer takes places. Now if vertical transfers of
$a$ and $a'$ units are possible between rows $i$~and~$i+1$, respectively in
columns $j_0$~and~$j_1$ with $j_0<j_1$, then the inequalities
in~\intmovedef\defitemnr1 give the inequalities
$\min(a,0)\geq\sum_{j=j_0}^{j_1-1}(M_{j+1,i+1}-M_{j,i})\geq\max(a',0)$ of
which the middle member must be~$0$, and which together with $a,a'\neq0$
therefore imply $a>0>a'$. This shows that the transfer in column~$j_0$ is
upwards and the one in column~$j_1$ is downwards. Thus successive upward
transfers between the same pair of rows take place in columns whose index
decreases weakly, and successive downward transfers between these rows take
place in in columns whose index increases weakly; the same is true with the
words ``rows'' and ``columns'' interchanged. Note that these directions are
opposite to the ones in the binary case for transfers between rows, but they
are the same as the ones in the binary case for transfers between columns. The
following definition is now justified.

\proclaim Definition. \intcopdef (integral raising and lowering operations)
Let $M\in\Mat$ and $i,j\in\N$.
\defitem
    If $M$ admits an upward or downward transfer of a single unit between rows
    $i$~and~$i+1$, then the resulting matrix is denoted by~$\upm_i(M)$
    respectively by~$\dnm_i(M)$. If~$M$ admits no~upward or no~downward
    transfer between rows $i$~and~$i+1$, then the expression $\upm_i(M)$
    respectively $\dnm_i(M)$ is undefined.
\defitem
    If $M$ admits a leftward or rightward transfer of a single unit between
    columns $j$~and~$j+1$, then the resulting matrix is denoted by~$\ltm_j(M)$
    respectively by~$\rtm_j(M)$. If~$M$ admits no~leftward or no~rightward
    transfer between columns $j$~and~$j+1$, then the expression $\ltm_j(M)$
    respectively $\rtm_j(M)$ is undefined.

If an upward, downward, leftward, or rightward transfer of $a>1$ units between
a given pair of rows~$i,i+1$ or columns~$j,j+1$ is defined, the the resulting
matrix can also be obtained by a succession of $a$ raising or lowering
operations, as $(\upm_i)^a(M)$, $(\dnm_i)^a(M)$, $(\ltm_j)^a(M)$ or
$(\rtm_j)^a(M)$, respectively. A succession of even more transfers between the
same rows or columns and in the same direction may be possible, and the
potentials for such transfers are given by the following expressions.

\proclaim Definition. \nmidefs
For $M\in\Mat$ and $i,j\in\N$, the numbers
$\num_i(M),\ndm_i(M),\nlm_j(M),\nrm_j(M)\in\N$ are defined by
$$\eqalignno
{ \num_i(M)
&=\max\setof\textstyle{M}_{i+1,0}+\sum_{j<l}(M_{i+1,j+1}-M_{i,j})
      :l\in\N\endset
,& \label(\intnumeq)
\cr
  \ndm_i(M)
&=\max\setof\textstyle\sum_{j\geq{l}}(M_{i,j}-M_{i+1,j+1}):l\in\N\endset
,& \label(\intndmeq)\cr
  \nlm_j(M)
&=\max\setof\textstyle{M}_{0,j+1}+\sum_{i<k}(M_{i+1,j+1}-M_{i,j})
      :k\in\N\endset
,& \label(\intnlmeq)
\cr
  \nrm_j(M)
&=\max\setof\textstyle\sum_{i\geq{k}}(M_{i,j}-M_{i+1,j+1}):k\in\N\endset
.& \label(\intnrmeq)
\cr}
$$

\proclaim Proposition. \nintmvprop
For $M\in\Mat$ and $i,j\in\N$, the numbers of times each of $\upm_i$,
$\dnm_i$, , $\ltm_j$ and $\rtm_j$ can be successively applied to~$M$ are given
respectively by the numbers $\num_i(M)$, $\ndm_i(M)$, $\nlm_j(M)$, and
$\nrm_j(M)$. Moreover $\ndm_i(M)-\num_i(M)=\rsum{M}_i-\rsum{M}_{i+1}$ and
$\nrm_j(M)-\nlm_j(M)=\csum{M}_j-\csum{M}_{j+1}$.

\proof
If $M$ admits an upward transfer of $a$ units between rows $i$~and~$i+1$ and
in column~$l$, then $\num_i(M)\geq
M_{i+1,0}+\sum_{j<l}(M_{i+1,j+1}-M_{i,j})\geq M_{i+1,0}+a>0$ if $l>0$, while
one has $\num_i(M)\geq{M}_{i+1,0}\geq{a}>0$ in case $l=0$, so $\num_i(M)$ is
nonzero either way. Similarly if $M$ admits a downward transfer of $a$ units
between those rows then $\ndm_i(M)\geq{a}>0$. Now suppose conversely that
$\num_i(M)>0$ or that $\ndm_i(M)>0$. In the former case let $l_0$ be the
minimal index~$l$ for which the maximum in~(\intnumeq) is attained, and in the
latter case let $l_1$ be the maximal index~$l$ for which the maximum
in~(\intndmeq) is attained; in either case let the maximum exceed all values
at smaller respectively at larger indices by~$a$. Then it is easily verified
that an upward transfer of $a$ units in column~$l_0$, respectively a downward
transfer of $a$ units in column~$l_1$, is possible between rows~$i,i+1$. We
know that this is the only column in which a transfer in that direction
between rows~$i,i+1$ is possible. Moreover the expressions in (\intnumeq)
and~(\intndmeq) have a constant difference $\rsum{M}_{i+1}-\rsum{M}_i$ for
every~$l$, so we may conclude that, whenever a transfer in either direction
between $M_{i,l}$ and~$M_{i+1,l}$ is possible, the index $l$ realises the
maxima defining $\num_i(M)$ and~$\ndm_i(M)$ in both expressions. Then the fact
that the transfer can be followed by an inverse transfer shows that an upward
transfer of $a$ units decreases $\num_i(M)$ by~$a$, and that a downward
transfer of $a$ units decreases $\ndm_i(M)$ by~$a$; a straightforward
induction on $\num_i(M)$ or $\ndm_i(M)$ will complete the proof. The
statements concerning $\nlm_j(M)$ and $\nrm_j(M)$ follow by transposition of
the matrices involved.
\QED

Like for binary case, we can describe the possibilities for transfers between
a given pair of rows or columns using parentheses. We shall describe the case
of rows~$i,i+1$, where the notion of matching parentheses is most easily
visualised; for the case of adjacent columns one must apply transposition of
the matrix (unlike the binary case, where rotation must be used). Since
successive transfers between two rows progress in the sense opposite to the
binary case, namely from right to left for successive upward moves, we must
invert the correspondence between parentheses and moves: an upward transfer of
a unit will correspond to the transformation of ``)'' into ``(''. Each unit in
row~$i$ is a potential candidate for an downward transfer, and therefore
represented by~``('', and each unit in row~$i+1$ is represented by~``)''. All
these symbols are gathered basically from left to right to form a string of
parentheses, but the crucial point is how to order the symbols coming from a
same column~$j$. The form of the summations in the definitions above makes
clear that the rule must be to place the $M_{i,j}$ symbols~``(''to
the~\emph{right} of the $M_{i+1,j}$ symbols~``)'', so that these cannot match
each other; rather the symbols from $M_{i,j}$ may match those
from~$M_{i+1,j+1}$. Now, as for the binary case, the units that may be
transferred correspond to the unmatched parentheses, and the order in which
they are transferred is such that the corresponding parentheses remain
unmatched: upward moves transform unmatched symbols~``)'' into ``('' from
right to left, and downward moves transform unmatched symbols~``('' into ``)''
from left to right. In the example given the string of parentheses is
$\let~=\mid )) \underline{(~)} \underline{((~) (~))})) \underline{(((~)) (((~
(~)))))} \underline{((~))} (((( $, where we have inserted bars to separate the
contributions from different columns, and underlined the maximal substrings
with balanced parentheses; this makes clear that $4$~successive upward
transfers are possible in columns $3,3,0,0$, or $4$~successive downward
transfers, all in (the final nonzero) column~$7$.

Note that any common number of units present in entries $M_{i,j}$ and
$M_{i+1,j+1}$ will correspond to matching parentheses both for transfers
between rows $i,i+1$ and between columns~$j,j+1$. Therefore no transfer
between those rows or those columns will alter the value
$\min(M_{i,j},M_{i+1,j+1})$ (but it can be altered by transfers in other pairs
of rows or columns). In fact one may check that those instances of the
inequalities in definition~\intmovedef\ whose summation is reduced to a single
term forbid any such transfer involving $M_{i,j}$ or~$M_{i+1,j+1}$ when that
entry is strictly less than the other, and in case it is initially greater or
equal than the other they forbid transfers that would make it become strictly
less.

The relation of the above definitions to tableau- and Littlewood-Richardson
conditions for integral matrices is given by the following proposition, which
like the one for the binary case is a direct translation of the pertinent
parts of proposition~\Tabprop\ and definition~\LRdef.

\proclaim Proposition. \intcondcharprop
Let $M\in\Mat$ and let $\\/\kappa$ and $\mu/\nu$ be skew shapes.
\statitem $M\in\IE{\\/\kappa}$ if and only if $\rsum{M}=\\-\kappa$ and
          $\num_i(M)\leq\kappa_i-\kappa_{i+1}$ for all~$i\in\N$.
\statitem $M\in\IL{\nu/\mu}$ if and only if $\csum{M}=\nu-\mu$ and
          $\nlm_j(M)\leq\mu_j-\mu_{j+1}$ for all~$j\in\N$.
\enditem
The second parts of these conditions can also be stated in terms of the
potentials of~$M$ for lowering operations, as
$\ndm_i(M)\leq\\_i-\\_{i+1}$ for all~$i\in\N$, respectively as
$\nrm_j(M)\leq\nu_j-\nu_{j+1}$ for all~$j\in\N$.
\QED

And as is the binary case, the potentials for transfers in a perpendicular
direction are unchanged.

\proclaim Proposition. \intinvarprop
If integral matrices $M,M'\in\Mat$ are related by $M'=\upm_i(M)$ for some
$i\in\N$, then $\nlm_j(M)=\nlm_j(M')$ and $\nrm_j(M)=\nrm_j(M')$ for
all~$j\in\N$. Consequently, the conditions $M\in\IL{\nu/\mu}$ and
$M'\in\IL{\nu/\mu}$ are equivalent for any skew shape~$\nu/\mu$. Similarly if
$M$~and~$M'$ are related by $M'=\rtm_j(M)$ for some $j\in\N$, then
$\num_i(M)=\num_i(M')$ and $\ndm_i(M)=\ndm_i(M')$ for all~$i\in\N$, and one
has $M\in\IE{\\/\kappa}\iff M'\in\IE{\\/\kappa}$ for any skew
shape~$\\/\kappa$.

\proof
Like for the binary case we may focus on proving $\nlm_j(M)=\nlm_j(M')$ when
$M'=\upm_i(M)$, and we may assume that the upward transfer involved in passing
from~$M$ to~$M'$ occurs in column $j$~or~$j+1$. Suppose first that it occurs
in column~$j$. Then the only partial sum in (\intnlmeq) that differs between
$M$~and~$M'$ is the one for~$k=i+1$, which has $M_{i+1,j+1}-M_{i,j}$ as final
term; it will decrease by~$1$ in~$M'$. But this decrease will not affect the
maximum taken in that equation unless it was strictly greater than the
previous partial sum, for~$k=i$, which means that $M_{i,j}<M_{i+1,j+1}$;
however that would contradict the supposition that $\upm_i$ involves
$M_{i,j}$, so it does not happen. Suppose next that the upward transfer occurs
in column~$j+1$. Then the only one of the values of which the maximum is taken
in (\intnlmeq) that differs between $M$~and~$M'$ is the one for~$k=i$: it is
$M_{i,j+1}$ if $i=0$, and otherwise contains a partial sum with final term
$M_{i,j+1}-M_{i-1,j}$; it will increase by~$1$ in~$M'$ in either case. But
that increase will not affect the maximum unless the value for~$k=i$ is at
least as great as the one for~$k=i+1$, which means $M_{i,j}\geq{M}_{i+1,j+1}$,
but this would contradict the supposition that $\upm_i$ decreases
$M_{i+1,j+1}$, so this does not happen either. Therefore
$\nlm_j(M)=\nlm_j(M')$ holds in all cases.
\QED

Using crystal operations for integral matrices, it is possible to describe
cancellations that will realise the transitions
(\intbruteeq)$\to$(\LRinttabeq)$\to$(\LRinteq) and
(\intbruteeq)$\to$(\LRintalteq)$\to$(\LRinteq), as we have done for binary
matrices. We shall describe the cancellation of the terms not satisfying
$M\in\IL{\nu/\mu}$, which realises (\LRinttabeq)$\to$(\LRinteq) and
(\intbruteeq)$\to$(\LRintalteq); the cancellation of the terms not satisfying
$M\in\IE{\kappa/\\}$ is similar, transposing all operations. One starts
traversing~$M$ by rows from top to bottom, searching for the first index~$k$
(if any) for which one has $(\mu+\sum_{i<k}M_i)\nleh(\mu+\sum_{i\leq{k}}M_i)$,
and then choosing the maximal witness~$j$ for this, i.e., with
$(\mu+\sum_{i<k}M_i)_j<(\mu+\sum_{i\leq{k}}M_i)_{j+1}$. This implies that
$\nlm_j(M)>\mu_j-\mu_{j+1}$, and similarly to the binary case the cancellation
is along the ``ladder'' of matrices obtainable from~$M$ by operations $\ltm_j$
and~$\rtm_j$, reversing its lower portion where the potential for~$\ltm_j$
exceeds~$\mu_j-\mu_{j+1}$: the term for~$M$ is cancelled against the term for
$M'=(\ltm_j)^d(M)$, where $d=\nlm_j(M)-(\mu_j-\mu_{j+1}+1)-\nrm_j(M)$, and as
before $(\ltm_j)^d$ means $(\rtm_j)^{-d}$ if $d<0$.

One can also write $d=\alpha_{j+1}-\alpha_j-1$ where $\alpha=\mu+\csum{M}$,
from which it easily follows that the terms for $M$~and~$M'$ do indeed cancel
each other. The fact that the same values $k$~and~$j$ will be found for~$M'$
(so that we have actually defined an involution on the set of cancelling
terms) is most easily seen as follows, using the parenthesis description given
above. Throughout the lower portion of the ladder, the unmatched
parentheses~``)'' from the left, up to and including the one corresponding to
the unit of~$M_{k,j+1}$ that ``causes'' the inequality $(\mu+\sum_{i<k}M_i)_j
<(\mu+\sum_{i\leq{k}}M_i)_{j+1}$, are unchanged, so $M'_i=M_i$ for all $i<k$,
and the entry $M'_{k,j+1}$ remains large enough that the mentioned inequality
still holds for~$M'$. Since no coefficients have changed that could cause any
index larger than~$j$ to become a witness for
$(\mu+\sum_{i<k}M'_i)\nleh(\mu+\sum_{i\leq{k}}M'_i)$, the index~$j$ will still
be the maximal witness of that relation for~$M'$. (The change to the entry
$M_{k,j}$ could cause the index $j-1$ to become, or cease to be, another
witness of this relation for~$M'$, so it is important here that $j$ be chosen
as a \emph{maximal} witness, unlike in the binary case where any systematic
choice works equally well.)

Again it will be important in applications to strengthen
proposition~\intinvarprop, in a way that is entirely analogous to the binary
case. The proof in the integral case will be slightly more complicated
however.

\proclaim Lemma. \iclemma (integral commutation lemma)
Let $M,M',M''\in\Mat$,  be related by $M'=\upm_i(M)$ and $M''=\ltm_j(M)$ for
some $i,j\in\N$; then $\ltm_j(M')=\upm_i(M'')$. The same holds when $\upm_i$
is replaced both times by $\dnm_i$ and/or $\ltm_j$ is replaced both times
by~$\rtm_j$.

\proof
As in the binary case, it suffices to prove the initial statement, and both
members of the equation $\ltm_j(M')=\upm_i(M'')$ are well defined by
proposition~\intinvarprop. That equation obviously holds in those cases where
the applications of $\ltm_j$~and~$\upm_i$ occurring in it are realised by
unit transfers between the same pairs of entries as in the application of these
operations to~$M$. This will first of all be the case when the pairs of
entries involved in the transfers $\upm_i\:M\mapsto{M'}$ and
$\ltm_j\:M\mapsto{M''}$ are disjoint: in that case, for each inequality in
definition~\intmovedef\ required for allowing one transfer, an argument
similar to the one given in the proof of proposition~\intinvarprop\ shows that
the minimum over all values $l'$~or~$k'$ of its left hand member is unaffected
by the other transfer.

In the remaining case where the pairs of entries involved in the two transfers
overlap, they lie inside the $2\times2$ square at the intersection of rows
$i,i+1$ and columns~$j,j+1$. Since we are considering upward and leftward
transfers, there are only two possibilities: if $M_{i,j}\geq{M_{i+1,j+1}}$,
then both transfers are towards~$M_{i,j}$, while if $M_{i,j}<{M_{i+1,j+1}}$
then both transfers are from~$M_{i+1,j+1}$. In either case there is for each
transfer only one of the inequalities in definition~\intmovedef\ whose left
hand member is affected by the other transfer, namely the one involving just
$M_{i,j}$ and~${M_{i+1,j+1}}$. In case both transfers are towards~$M_{i,j}$,
the inequalities in question both read $M_{i,j}-M_{i+1,j+1}\geq\max(-1,0)=0$,
and they continue to hold when $M_{i,j}$ is increased in $M'$ and in~$M''$. In
case both transfers are from~$M_{i+1,j+1}$, the inequalities both read
$M_{i+1,j+1}-M_{i,j}\geq\max(1,0)=1$, and here the other transfer \emph{may}
invalidate the inequality: indeed the inequality will continue to hold, when
$M_{i+1,j+1}$ is decreased by~$1$ in $M'$ and in~$M''$, if and only if
$M_{i+1,j+1}\geq{M_{i,j}}+2$. Therefore the only case where the unit transfers
involved in $\ltm_j(M')$ and in~$\upm_i(M'')$ are \emph{not} in the same
squares as in $\ltm_j(M)$ and in~$\upm_i(M)$, respectively, is when
$M_{i+1,j+1}=M_{i,j}+1$. In that case an elementary verification, comparing
the required inequalities in~$M'$ and~$M''$ to those that are known to hold
in~$M$, shows that those transfers are from $M'_{i,j+1}$ to~$M'_{i,j}$
respectively from $M''_{i+1,j}$ to~$M''_{i,j}$; then $\ltm_j(M')=\upm_i(M'')$
still holds, as illustrated by the display below of the relevant $2\times2$
squares. Instead of that verification, one may also argue that the transfer
involved in $\ltm_j(M')$ or in~$\upm_i(M'')$, which as we saw cannot occur in
the same place as the corresponding transfer in~$M$, cannot occur in a pair of
squares disjoint from those of the perpendicular transfer preceding it either
(such a transfer in disjoint squares cannot validate a relevant inequality any
more than it could invalidate one in the first case considered above); a
transfer towards the entry at~$(i,j)$ is therefore the only possibility that
remains.
\QED

\bigdisplay \advance \abovedisplayshortskip by \medskipamount
\def\mapleft#1{\smash
       {\mathop{\hbox to 2.5em{\leftarrowfill}}\limits^{\textstyle#1}}}
\matrix
{ \pmatrix{x+1&y\cr{z}&x\cr} & \mapleft{\ltm_j} & \pmatrix{x&y+1\cr{z}&x\cr}
    \rlap{$\displaystyle{}=M''$}
\cr
  \llap{$\upm_i$}\bigg\uparrow & & \llap{$\upm_i$}\bigg\uparrow
\cr
  \llap{$\displaystyle M'={}$}
    \pmatrix{x&y\cr{z}+1&x\cr} & \mapleft{\ltm_j} & \pmatrix{x&y\cr{z}&x+1\cr}
  \rlap{$\displaystyle{}=M$}
\cr}
\nn
$$

By a simple induction, the commutation lemma implies that
$(\ltm_j)^m((\upm_i)^n(M))=(\upm_i)^n((\ltm_j)^m(M))$ whenever $(\upm_i)^n(M)$
and $(\ltm_j)^m(M)$ are defined (this is true in the binary case as well). It
is interesting to consider how this commutation is realised if one assumes
that the operations applied to~$M$ can each be performed by a single transfer.
Much of the proof above can be applied to such transfers as well, but the
final case is more complicated. One can see that transfers remain in the same
position and are by the same amount, unless those in $(\upm_i)^n(M)$ and
$(\ltm_j)^m(M)$ are both from the entry $M_{k+1,l+1}$ and one has
$M_{k,l}<M_{k+1,l+1}<M_{k,l}+n+m$; then in fact $M_{k,l}+\max(n,m)\leq
M_{k+1,l+1}$ by the assumption that $(\upm_i)^n(M)$ and $(\ltm_j)^m(M)$ are
defined by single transfers. In this case one can put $d=M_{k+1,l+1}-M_{k,l}$
and decompose $(\ltm_j)^m=(\ltm_j)^a\after(\ltm_j)^b$ and
$(\upm_i)^n=(\upm_i)^a\after(\upm_i)^c$ where $a=m+n-d$ and $b=m-a$, $c=n-a$;
then the commutation takes place as in the diagram below, which shows the
evolution of the $2\times2$ square at the intersection of rows~$i,i+1$ and
columns~$j,j+1$. One sees that only in the top-left commuting square the
transfers realising the same operation are between different pairs of entries.
\bigdisplay
\def\mapleft#1{\smash
       {\mathop{\hbox to 2.5em{\leftarrowfill}}\limits^{\textstyle#1}}}
\matrix
{ \pmatrix{x+a&y+c\cr{}z+b&x\cr} & \mapleft{(\ltm_j)^a}
 & \pmatrix{x&y+a+c\cr{}z+b&x\cr} & \mapleft{(\ltm_j)^b}
 & \pmatrix{x&y+a+c\cr{}z&x+b\cr}
    \rlap{$\displaystyle{}=M''$}
\cr
  \llap{$(\upm_i)^a$}\bigg\uparrow & &
  \llap{$(\upm_i)^a$}\bigg\uparrow & &
  \llap{$(\upm_i)^a$}\bigg\uparrow
\cr
 \pmatrix{x&y+c\cr{}z+a+b&x\cr} & \mapleft{(\ltm_j)^a}
 & \pmatrix{x&y+c\cr{}z+b&x+a\cr} & \mapleft{(\ltm_j)^b}
 & \pmatrix{x&y+c\cr{}z&x+a+b\cr}
\cr
  \llap{$(\upm_i)^c$}\bigg\uparrow & &
  \llap{$(\upm_i)^c$}\bigg\uparrow & &
  \llap{$(\upm_i)^c$}\bigg\uparrow
\cr
  \llap{$\displaystyle M'={}$}
 \pmatrix{x&y\cr{}z+a+b&x+c\cr} & \mapleft{(\ltm_j)^a}
 & \pmatrix{x&y\cr{}z+b&x+a+c\cr} & \mapleft{(\ltm_j)^b}
 & \pmatrix{x&y\cr{}z&x+a+b+c\cr}
   \rlap{$\displaystyle{}=M$}
\cr}
\nn
$$
We see that when considering non-unit transfers between pairs of adjacent
columns and rows, a transfer that can be realised in a single step can become
a transfer broken into several steps after the commutation. The opposite can
also happen, and it suffices to consider the inverse operations in the diagram
to obtain an example of that situation.

 \par

\subsecno=-1 

\newsection Crystal graphs.
\labelsec\crystalsec

We shall now discuss some basic properties of crystal operations on matrices,
starting with those that involve only horizontal or only vertical operations.
So for the moment we assume that one of the binary and the integral cases is
chosen, as well as one of the two orientations horizontal and vertical; then
raising and lowering operations can be indicated simply by $e_i$~and~$f_i$,
respectively. Doing so we essentially find the structure of crystal graphs of
type~$A_n$, so the results in this section are generally known in some form.
Our goal here is to give an brief overview of facts, formulated for matrices,
and of terminology used.

A \newterm{crystal graph} of matrices is a set closed under the application
of operations $e_i$~and~$f_i$ for $i\in\N$ whenever these are defined, with
matrices considered as vertices, linked by directed edges defined and labelled
by those operations. A connected crystal graph will be called a
\newterm{crystal}. Since for any nonzero matrix there is always some lowering
operation that can be applied, such graphs are necessarily infinite. It will
in some cases be more convenient to fix some $n>0$ and to restrict our
attention to matrices with at most $n$ rows (in case of vertical operations)
or at most $n$ columns (in case of horizontal operations), in which case the
set of crystal operations is limited to $e_i$~and~$f_i$ for $i\in\set{n-1}$;
then crystals will be finite.

Crystals have a complicated structure, but a number of properties can be
established directly. Whether considering finite or infinite sets of crystal
operations, it is easy to see that the crystals that occur do not depend, up
to isomorphism, on whether one uses horizontal or vertical operations. In the
case of integral matrices, transposition defines an isomorphism between the
two structures on the set of all such matrices. In the binary case, rotation
of matrices by a quarter turn counterclockwise defines an isomorphism from
crystal graphs defined by vertical moves to crystal graphs defined by
horizontal moves, but such a rotation can only be defined for a subset of
binary matrices at a time. For any crystal for vertical moves, there is
some~$m\in\N$ so that all its nonzero matrix entries are contained in the
leftmost~$m$ columns, and an appropriate rotation maps these columns to the
topmost~$m$ rows in reverse order.

Using suitable rotations by a half turn, one can similarly show that when
crystal operations $e_i$~and~$f_i$ are limited to $i\in\set{n-1}$, the
structure obtained from a crystal graph by interchanging the labels $e_i$ and
$f_{n-2-i}$ is again a crystal graph of the same type. (It can be shown using
other known properties that each crystal thus gives rise to a crystal
isomorphic to itself, but that is not obvious at this point.)

It is also true that crystals defined using binary matrices are isomorphic to
those defined using integral matrices. Here the correspondence is rather less
obvious than in the above cases, so we shall state it more formally. The
statement basically follows from the fact that coplactic operations on
semistandard tableaux can be defined either using the Kanji reading order
(down columns, from right to left), or using the Semitic reading order
(backwards by rows, from top to bottom); see~\ref{L-R intro,
proposition~3.2.1}. Translating the operations on tableaux in terms of
matrices, using binary respectively integral encodings, gives crystal
operations on matrices, and leads to the following formulation.

\proclaim Proposition. \bininteqprop
Let $T\in\SSTx{\\/\kappa}$ be a semistandard tableau, and let its binary and
integral encodings be $M\in\bMat$ and $N\in\Mat$, respectively. Then for
$m\in\N$, \ $\upm_m(M)$ is defined if and only if $\ltm_m(N)$ is defined, and
if so, $\upm_m(M)$ and $\ltm_m(N)$ are respectively the binary and integral
encodings of one same tableau $T'\in\SSTx{\\/\kappa}$. The same is true when
$\upm_m$ and $\ltm_m$ are replaced by $\dnm_m$ and $\rtm_m$, respectively.

This proposition immediately implies that the correspondence of $M$~and~$N$ as
encodings of the same semistandard tableau of shape~$\\/\kappa$ can be
extended to a correspondence between the crystal for vertical moves that
contains~$M$ and the crystal for horizontal transfers that contains~$N$, so
that these two structures are isomorphic. Note that a tableau may satisfy the
tableau condition for many different shapes (indeed any tableau satisfies
infinitely many tableau conditions, though not all of them are equally
interesting), so for a crystal of one type, one obtains isomorphisms with many
different crystals of the other type. This indicates that it is quite common
to have distinct but isomorphic crystals of matrices; that observation could
also have been made after commutation lemmas of the previous section.

When we consider in a crystal graph only the edges labelled $e_i$~and~$f_i$
for a single index~$i$, then as we have seen, each vertex~$M$ is part of a
finite ladder. The cancellations we considered in the previous section were
defined, after selecting an appropriate ladder, by reversal of a lower portion
of a ladder. It is natural to consider also the simpler operation of reversing
each complete ladder for a fixed index~$i$. This defines an involution on the
set of all matrices, which in the case of vertical operations preserves all
column sums while interchanging the row sums $\rsum{M}_i$
and~$\rsum{M}_{i+1}$, and in case of horizontal operations preserves all row
sums while interchanging $\csum{M}_i$ and~$\csum{M}_{i+1}$. Now if one takes
any crystal graph for which the multivariate generating series of its matrices
by row sums (for vertical operations) or by column sums (for horizontal
operations) exists (i.e., its coefficients remain finite), then that series
will be invariant under all permutations of the variables; if the monomials of
the series are of bounded degree, this means that the series is a symmetric
function. (While the crystal graph of all binary or integral matrices is too
large to have such a series, each individual crystal does indeed define a
symmetric function.)

This is similar to the situation for Bender-Knuth involutions, which can be
used not only to define cancellations, but also to show the symmetry of the
weight generating series of skew semistandard tableaux, in other words of skew
Schur functions. Indeed, the crystal graphs on the sets $\BE{\\/\kappa}$ or
$\IE{\\/\kappa}$ for some skew shape~$\\/\kappa$ are simple examples of graphs
to which the symmetry result obtained by ladder reversal can be applied, and
here the generating series is the skew Schur function $s_{\\/\kappa}[X_\N]$ in
both cases. Although this result can also be obtained by using the simpler
Bender-Knuth involutions, there is at least one sense in which it is
``better'' to use the ladder reversing involutions, namely that they actually
define a symmetric group action of the underlying set of any crystal graph:
the ladder reversing involutions satisfy the Coxeter relations of type~$A$,
whereas the Bender-Knuth involutions do not. (On the other hand ladder
reversal does not respect the potentials for crystal operations at other
indices, so the action is not by graph isomorphisms of any kind.) This
symmetric group action means that one not only knows, for any crystal graph
for vertical operations and any pair of compositions~$\alpha,\alpha'$ whose
parts are related by some permutation, that there are as many matrices~$M$ in
the crystal graph with $\rsum{M}=\alpha$ as with $\rsum{M}=\alpha'$, but that
one also has a natural bijection between those sets of matrices.

We shall now state the above matters more formally. We first define the
involutions generating the action, which come in four flavours: binary or
integral, and horizontal or vertical. However, like for the crystal
operations, there is no explicit distinction in the notation between the
binary and the integral case, so only two symbols are used, $\vr_i$ for the
``vertical'' involutions and $\hr_j$ for the ``horizontal'' ones.

\proclaim Definition.
The symbols $\vr_i$ and $\hr_j$ for $i,j\in\N$ denote involutions both on
$\bMat$ and on $\Mat$, given in either case by the following expressions:
$$\eqalignno
{\vr_i(M)& =(\upm_i)^\delta(M)
         \ttex{where $\delta=\num_i(M)-\ndm_i(M)=\rsum{M}_{i+1}-\rsum{M}_i$,}
& \label(\vreq)
\cr \noalign{\nobreak\smallskip}
  \hr_j(M)&=(\ltm_j)^\delta(M)
         \ttex{where $\delta=\nlm_j(M)-\nrm_j(M)=\csum{M}_{j+1}-\csum{M}_j$.}
& \label(\hreq)
\cr}
$$
and where $(\upm_i)^\delta$ and $(\ltm_j)^\delta$ are to be interpreted as
$(\dnm_i)^{-\delta}$ and $(\rtm_j)^{-\delta}$, respectively, when $\delta<0$,
and as the identity when $\delta=0$.

It is immediately obvious that one has the following identities, where $s_i$
denotes the transposition of $i$~and~$i+1$, acting as usual on compositions by
permuting their parts:
$$
  \eqalign{\rsum{\vr_i(M)}&=s_i\cdot\rsum{M}\cr
           \csum{\vr_i(M)}&=        \csum{M}\cr
          }
\ttext{and}
  \eqalign{\rsum{\hr_j(M)}&=        \rsum{M}\cr
           \csum{\hr_j(M)}&=s_j\cdot\csum{M}\cr
          }
.\nn
$$
Moreover, when $\rsum{M}_i=\rsum{M}_{i+1}$ one has $\vr_i(M)=M$, and similarly
one has $\hr_j(M)=M$ when $\csum{M}_j=\csum{M}_{j+1}$; this contrasts with the
Bender-Knuth involutions, which may send one tableau whose weight already has
equal values at the positions being interchanged to another such tableau. It
follows from the commutation lemmas that $\vr_i$ and $\hr_j$ commute for all
$i,j\in\N$. 

\proclaim Lemma. \actionlem
The operations $\vr_i$ satisfy the Coxeter relations $(\vr_i)^2=1$, \
$\vr_i\vr_k=\vr_k\vr_i$ when $|k-i|\geq2$, and
$\vr_i\vr_{i+1}\vr_i=\vr_{i+1}\vr_i\vr_{i+1}$. Similar relations hold for
the operations~$\hr_j$.

Via the correspondence of crystal operations to coplactic operations on
semistandard tableaux, this lemma is essentially equivalent to
\ref{Lothaire, theorem~5.6.3}, but we shall give a proof at the end of our
paper. The following theorem summarises the main conclusions that follow from
the relations given.

\proclaim Theorem.
There are two commuting actions of the group $\Sym_\infty$, denoted by
$\pi\vertact(M)$ and $\pi\horact(M)$, on each of $\bMat$ and~$\Mat$; they are
defined in both cases by $s_i\vertact(M)=\vr_i(M)$ and
$s_i\horact(M)=\hr_i(M)$ for a generator~$s_i$ of~$\Sym_\infty$ (the
transposition $(i,i+1)$).
Any operation $\pi\vertact$ permutes the parts of~$\rsum{M}$ by~$\pi$ while
leaving $\csum{M}$ invariant, and if $\pi$ fixes $\rsum{M}$ then $\pi\vertact$
fixes~$M$; similar statements hold for~$\pi\horact$.

\proof
The only point that is not immediate is the final claim that $\pi\vertact$
fixes~$M$, which only follows directly from the observed properties of
$\vr_i$~and~$\hr_j$ in case the stabiliser of $\rsum{M}$ is generated by a
subset of the generators~$s_i$ (i.e., any equalities among the parts of
$\rsum{M}$ occur in consecutive ranges). But that stabiliser is always
\emph{conjugate} to such a (parabolic) subgroup, which allows reduction to
the given case.
\QED

We continue with some simpler but important observations. While every matrix
admits some lowering operation, both vertically and horizontally, there are
matrices that permit no raising operations in a given direction at all. Indeed
it follows from propositions \bincondcharprop~and~\intcondcharprop\ that a
$M\in\bMat$ or $N\in\Mat$ admits no upwards operations if and only if
$M\in\BL{\nu/(0)}$ where $\nu=\csum{M}\in\Part$ respectively $N\in\IE{\\/(0)}$
where $\\=\csum{M}\in\Part$, while the matrix admits no leftwards operations
if and only if $M\in\BE{\\/(0)}$ respectively $N\in\IL{\nu/(0)}$, where
$\\=\rsum{M}$ and $\nu=\rsum{N}$. The given values of $\\$~and~$\nu$ are the
only ones for which the respective conditions could hold, so these properties
state that the fact that no raising operations are possible means that the
matrix encodes, according to the case considered, a semistandard tableau of
straight shape, or some (skew) Littlewood-Richardson tableau.

Fixing one direction as before, if we start with a given matrix, and repeat
the step of applying some raising operation that can be applied, until no such
operation exists, then we end up with a matrix of the type just described. The
process cannot go on indefinitely since at every step we move a unit to a row
or column with a lower number, which decreases an obvious statistic (the sum
over all units of their row respectively column number). The matrix found at
the end of the process will be said to be obtained from the original matrix by
\newterm{exhausting} raising operations (of the type considered). We shall
show later that the final matrix found does not depend on the choices made
during the exhaustion process (in other words, the process is confluent). This
is an important and quite nontrivial statement; it implies that every crystal
contains a unique vertex in which no raising operations are possible. For the
moment however it suffices to know that some final matrix can always be
reached in this way.

Similar properties are known to hold for jeu de taquin on skew semistandard
tableaux, which can be used to ``rectify'' a skew tableau to one of straight
shape by a nondeterministic procedure, but whose final result is independent
of the choices made. In fact there is a close relation between that game and
crystal operations on matrices encoding the tableaux. In fact it is known, see
\ref{L-R intro, theorem~3.4.2}, that jeu de taquin slides applied to a
tableau~$T$ correspond to sequences of coplactic operations applied to a its
``companion tableau'', which is a tableau whose integral encoding is the
transpose of the integral encoding of~$T$. Since coplactic operations on
tableaux correspond to horizontal crystal operations on their integral
encodings, this means that jeu de taquin slides correspond to sequences of
vertical crystal operations. A correspondence exists also for binary
encodings. More precisely, one has the following.

\proclaim Proposition. \jdtrelprop
Let $T\in\SSTx{\\/\kappa}$ and let $T'\in\SSTx{\\'/\kappa'}$ be obtained
from~$T$ by an inward jeu de taquin slide starting at the corner $(i_0,j_0)$
of~$\kappa$ and ending in the corner $(i_1,j_1)$ of~$\\$. The their respective
binary encodings $M,M'$ are related by
$M'=\ltm_{j_1-1}(\cdots(\ltm_{j_0+1}(\ltm_{j_0}(M))\cdots))$, and their
respective integral encodings $N,N'$ are related by
$N'=\upm_{i_1-1}(\cdots(\upm_{i_0+1}(\upm_{i_0}(N))\cdots))$.

\proof
In the description of a jeu de taquin slide in terms of a skew tableau with an
empty square into which entries as successively slid, one can at each point
associate both a binary and an integral encoding matrix, which count the
entries in columns respectively rows in the usual fashion (ignoring the empty
square). It is a straightforward verification that under this correspondence
each leftward slide of a tableau entry results in applying a leftward move in
corresponding pair of columns of the binary matrix while leaving the integral
matrix unchanged, and that each upward slide of a tableau entry results in
applying an upward unit transfer in corresponding pair of rows of the integral
matrix while leaving the binary matrix unchanged. Indeed, the inequalities in
definition~\binmovedef\defitemnr2 reflect the fact that the entries above and
below a horizontally sliding entry, in the same pair of columns, satisfy weak
increase along rows, while the inequalities in
definition~\intmovedef\defitemnr1 reflect the fact that the entries to the
left and right of a vertically sliding entry, and in the same rows, satisfy
strict increase down columns. The proposition immediately follows.
\QED

Thus one can see that if a skew tableau~$T$ has rectification~$P$, then the
binary encoding of~$P$ can be obtained from that of~$T$ by exhausting leftward
moves, while the integral encoding of~$P$ can be obtained from that of~$T$ by
exhausting upward transfers. However this reasoning does not yet imply that
these are the only possible results of such exhaustion of raising operations,
even if one admits the uniqueness of the rectification of~$T$, since jeu de
taquin slides correspond only to specific composition of raising operations,
and there might exist transitions that can be realised by raising operations
but not by jeu de taquin.

Since the raising operations on matrices in the direction perpendicular to the
one just considered correspond to coplactic operations on tableaux, one easily
recognises a relation between the commutation lemmas of the previous section
and the commutation of coplactic operations with jeu de taquin, as stated for
instance in \ref{L-R intro, theorem~3.3.1}. In the next section we shall
elaborate on this relation, developing a theory for matrices that is similar
(but more symmetric) than the one for semistandard tableaux.

\subsecno=0 

 \par

\newsection Double crystals.
\labelsec\Knuthsec

We shall now return to considering at the same time vertical and horizontal
operations. A set of all matrices that can be obtained from a given one by
using both vertical and horizontal crystal operations, provided with the
structure of labelled directed graph defined by those operations, will be
called a \newterm{double crystal}. Due to the commutation lemmas, double
crystals will be much easier to study than ordinary crystals. We shall see
that they can be used to define in a very natural way correspondences quite
similar to the \caps{RSK}~correspondence and to its dual correspondence that
was (also) defined in~\ref{Knuth PJM}, respectively. In fact, the
correspondence that is most naturally obtained for integral matrices is the
Burge variant of the \caps{RSK}~correspondence, described in~\ref{Fulton,
A.4.1}.

\subsection Knuth correspondences for binary and integral matrices.

We have just seen that by exhausting either vertical or horizontal raising
operations, one can transform any integral or binary matrix into one that
encodes, according to the case considered, a straight semistandard tableau or
a Littlewood-Richardson tableau. Now by exhausting both vertical and
horizontal raising operations, one can obviously transform any integral or
binary matrix into one that encodes a Littlewood-Richardson tableau of
straight shape. It is well known that the possibilities for such tableaux are
extremely limited: in their display, any row~$i$ can contain only entries~$i$,
whence such tableaux are completely determined by either their shape or their
weight, which moreover are the same partition~$\\$. One easily sees that the
binary and integral encodings of such tableaux are the matrices described by
the following definition.

\proclaim Definition. \diagdef
Let $\\\in\Part$. We denote by $\diagram\\\in\bMat$ the binary matrix
$\(\Kr{(i,j)\in\set\\}\){}_{i,j\in\N}$ whose bits~`$1$' occur precisely in the
shape of the Young diagram~$\set\\$. We denote by $\diagon\\\in\Mat$ the
integral matrix~$=\(\Kr{i=j}\\_i\){}_{i,j\in\N}$, which is diagonal with the
parts of~$\\$ as diagonal entries.

One immediately sees that this binary matrix satisfies $\rsum{\diagram\\}=\\$
and $\csum{\diagram\\}=\\\tr$, while this integral matrix satisfies
$\rsum{\diagon\\}=\csum{\diagon\\}=\\$. From the consideration of tableaux one
may deduce that these matrices are precisely the ones for which all vertical
and horizontal raising operations are exhausted, but it is easy and worth
while to prove this directly.

\proclaim Lemma. \normalformlemma
A binary matrix $M\in\bMat$ satisfies $\forall{i}\:\num_i(M)=0$ and
$\forall{j}\:\nlm_j(M)=0$ if and only if $M=\diagram\\$ for some $\\\in\Part$.
Similarly an integral matrix $N\in\bMat$ satisfies $\forall{i}\:\num_i(N)=0$
and $\forall{j}\:\nlm_j(N)=0$ if and only if $N=\diagon\\$ for some
$\\\in\Part$.

\proof
One easily sees that no raising operations are possible for $\diagram\\$ or
for~$\diagon\\$, since for $\diagram\\$ there are not even candidates for
upward of leftward moves, while for $N=\diagon\\$ and all~$i$ the units
of~$N_{i+1,i+1}$ are blocked by those of $N_{i,i}$ for upward or leftward
transfers. This takes care of the ``if'' parts. For the ``only if'' parts we
shall apply induction on the sum of all matrix entries, the case of null
matrices being obvious (with $\\=(0)$).

Suppose that $\forall{i}\:\num_i(M)=0$ and $\forall{j}\:\nlm_j(M)=0$, and let
$(k,l)$ be the indices of the bit $M_{k,l}=1$ for which $l$ is maximal, and
$k$ is minimal for that~$l$. If $k>0$ one would have $\num_{k-1}(M)>0$, so
$k=0$. But then also $M_{0,j}=1$ for all $j<l$, since otherwise the largest
$j$ violating this would give $\nlm_j(M)>0$. Let $M'=(M_{i+1,j})_{i,j\in\N}$
be the matrix obtained by removing the topmost row~$M_0$. The hypothesis
$\nlm_j(M)=0$ implies $\nlm_j(M')=0$ for $j<l$ since
$(M_{0,j},M_{0,j+1})=(1,1)$; one also has $\nlm_j(M')=0$ for $j\geq{l}$, for
lack of any bits $M'_{i,j+1}=1$. Therefore $M'$ satisfies the hypotheses of
the lemma, and by induction one has $M'=\diagram{\\'}$ for some $\\'\in\Part$,
with obviously $\\'_0\leq{l}$. Then $M=\diagram\\$, where $\\$ is the
partition obtained by prefixing~$l$ to the parts of~$\\'$ (formally: $\\_0=l$
and $\\_{i+1}=\\'_i$ for~$i\in\N$).

Now suppose that $\forall{i}\:\num_i(N)=0$ and $\forall{j}\:\nlm_j(N)=0$ (the
integral case). The conditions $\num_i(N)=0$ and $\nlm_j(N)=0$ respectively
imply that $N_{i+1,0}=0$ and $N_{0,j+1}=0$. But one cannot have
$\rsum{N}_0=0$, since then for the smallest index~$i$ such that
$\rsum{N}_{i+1}>0$ one would have $\num_i(N)=\rsum{N}_{i+1}>0$. Therefore
$N_{0,0}>0$ and we can apply induction to the matrix
$N'=(N_{i+1,j+1})_{i,j\in\N}$ obtained by removing the initial row and column
from~$N$, since clearly $\num_i(N')=\num_{i+1}(N)=0$ and
$\nlm_j(N')=\nlm_{j+1}(N)=0$. Then $N'=\diagon{\\'}$ for some $\\'\in\Part$,
and $N_{0,0}\geq\\'_0$ follows either from $\num_0(N)=0$ or from
$\nlm_0(N)=0$. So one has $N=\diagon\\$, where $\\\in\Part$ is obtained by
prefixing~$N_{0,0}$ to the parts of~$\\'$.
\QED

Now that we have characterised the matrices at which the process of applying
raising operations can stop, and with the commutation lemmas at our disposal,
the general structure of double crystals can be easily analysed, as follows.

Fix a matrix $M\in\bMat$ or $M\in\Mat$. Suppose that $R\xraise$ is some
sequence of raising operations $\upm_i$ for varying $i\in\N$ that exhausts
such operations when applied to~$M$, in other words such that the resulting
matrix $R\xraise(M)$ satisfies $\forall{i}\:\num_i(R\xraise(M))=0$. Similarly
let $R\xleft$ be some sequence of raising operations $\ltm_j$ for varying
$j\in\N$ that exhausts such operations when applied to~$M$; then the resulting
matrix $R\xleft(M)$ satisfies $\forall{j}\:\nlm_j(R\xleft(M))=0$. By
lemma~\bclemma\ or lemma~\iclemma, the expressions $R\xleft(R\xraise(M))$ and
$R\xraise(R\xleft(M))$ are both defined and they designate the same matrix,
which we shall call~$N$. By proposition \bininvarprop\ or~\intinvarprop, one
has $\forall{i}\:\num_i(N)=0$ and $\forall{j}\:\nlm_j(N)=0$, and so
lemma~\normalformlemma\ states that $N=\diagram\\$ or $N=\diagon\\$ (according
to the case considered) for some~$\\\in\Part$. Since $R\xleft$ preserves row
sums one has $\rsum{R\xraise(M)}=\rsum{N}$, which composition equals~$\\$ both
in the binary and in the integral case; similarly $R\xraise$ preserves column
sums, so one has $\csum{R\xleft(M)}=\csum{N}$, which equals~$\\$ in the
integral case and $\\\tr$ in the binary case. Thus either of the values
$\rsum{R\xraise(M)}$ or $\csum{R\xleft(M)}$ completely determines~$\\$.

Now let $S\xraise$ and $S\xleft$ be other sequences of upward respectively
leftward crystal operations that can be applied to~$M$ and that exhaust such
operations. Then the matrices $S\xleft(R\xraise(M))$ and
$S\xraise(R\xleft(M))$ satisfy the conditions of lemma~\normalformlemma, and
since their row sums respectively their column sums are the same as those
of~$N$, they must both be equal to~$N$. From $R\xleft$ we can form the
``inverse'' operation $R\xright$ by composing the lowering operations $\rtm_j$
corresponding to the operations $\ltm_j$ used, in the opposite order. Applying
$R\xright$ to the equation $R\xleft(R\xraise(M))=N =S\xraise(R\xleft(M))
=R\xleft(S\xraise(M))$, one finds $R\xraise(M)=S\xraise(M)$, and one similarly
finds $R\xleft(M)=S\xleft(M)$. So the matrices obtained from~$M$ by exhausting
vertical or horizontal crystal operations are uniquely defined. The matrix~$N$
obtained by exhausting both types of operations is also uniquely determined
by~$M$, and we shall call $N$ the \newterm{normal form} of~$M$.

Moreover from the knowledge of $R\xraise(M)$ and $R\xleft(M)$ one can uniquely
reconstruct~$M$. To state this more precisely, let $P$~and~$Q$ be matrices of
the same type (binary or integral) that satisfy $\forall{i}\:\num_i(P)=0$ and
$\forall{j}\:\nlm_j(Q)=0$ (which implies that $\rsum{P}$ and~$\csum{Q}$ are
partitions), and that satisfy moreover $\rsum{P}=\csum{Q}\tr$ in the binary
case, or $\rsum{P}=\csum{Q}$ in the integral case. Let $R\xleft(P)$ and
$R\xraise(Q)$ be obtained by from~$P$ by exhausting leftward crystal
operations, respectively from~$Q$ by exhausting upwards crystal operations (as
before $R\xleft$ and $R\xraise$ denote particular sequences of operations that
realise these transformations). Both matrices satisfy the conditions of
lemma~\normalformlemma, and by the hypothesis on column and row sums, the
partition~$\\$ in the conclusion of that lemma is the same in both cases;
therefore $R\xleft(P)=R\xraise(Q)$. As before we can form inverse operations
$R\xright$ and $R\xlower$ of $R\xleft$ and of $R\xraise$, respectively. We
successively apply these inverse operations to the equation
$R\xleft(P)=R\xraise(Q)$, and by commutation of $R\xright$ and $R\xlower$
(which follows that of $R\xleft$ and $R\xraise$) one finds
$R\xlower(P)=R\xright(Q)$. In the members of this final equation we have found
a matrix~$M$ for which $R\xraise(M)=P$ and $R\xleft(M)=Q$.

To show uniqueness of~$M$, suppose that $M'$ is another matrix, and that
$S\xraise$ and $S\xleft$ are other sequences of upwards respectively leftwards
crystal operations for which $P=S\xraise(M')$ and $Q=S\xleft(M')$. Then
$S\xleft(P)$ is the normal form of $P$, and therefore equal to~$R\xleft(P)$.
From this equality and $P=S\xraise(M')$ one deduces by commutation that
$S\xraise(S\xleft(M')) =S\xraise(R\xleft(M'))$, and since
$S\xleft(M')=Q=R\xleft(M)$ this gives $S\xraise(R\xleft(M))
=S\xraise(R\xleft(M'))$. But $S\xraise$ and $R\xleft$ are invertible
operations, so it follows that $M=M'$.

The above reasoning shows the extraordinary usefulness of the commutation
lemmas. But despite the relative ease with which the conclusions are reached,
one should not forget that each application of a crystal operation has to be
justified; for instance in the final part we could talk about $S\xleft(P)$
only because $P=S\xraise(M')$ and by assumption $S\xleft(M')$ is defined. Also
it should not be inferred from $S\xleft(P)=R\xleft(P)$ that $S\xleft$
and~$R\xleft$ are equivalent sequences of crystal operations, or even that
they can be applied to the same set of matrices. We now state more formally
the statement that follows from the argument given.

\proclaim Theorem. \decompthm (binary and integral decomposition theorem)
\statitem
There is a bijective correspondence between binary matrices $M\in\bMat$ on one
hand, and pairs $P,Q\in\bMat$ of such matrices satisfying
$\forall{i}\:\num_i(P)=0$ and $\forall{j}\:\nlm_j(Q)=0$ and
$\rsum{P}=\csum{Q}\tr$ on the other hand, which is determined by the condition
that $P$ can be obtained from~$M$ by a finite sequence of operations chosen
from $\setof\upm_i:i\in\N\endset$ and $Q$ can be obtained from~$M$ by applying
a finite sequence of operations chosen from $\setof\ltm_j:j\in\N\endset$. In
particular $\csum{P}=\csum{M}$ and $\rsum{Q}=\rsum{M}$.
\statitem
There is a bijective correspondence between integral matrices $M\in\Mat$ on
one hand, and pairs $P,Q\in\Mat$ of such matrices satisfying
$\forall{i}\:\num_i(P)=0$ and $\forall{j}\:\nlm_j(Q)=0$ and
$\rsum{P}=\csum{Q}$ on the other hand, which is determined by the condition
that $P$ can be obtained from~$M$ by a finite sequence of operations chosen
from $\setof\upm_i:i\in\N\endset$ and $Q$ can be obtained from~$M$ by applying
a finite sequence of operations chosen from $\setof\ltm_j:j\in\N\endset$.
In particular $\csum{P}=\csum{M}$ and $\rsum{Q}=\rsum{M}$.
\enditem
These correspondences can be explicitly computed by the constructions given
above.
\QED

\proclaim Corollary.
The crystal defined by any binary or integral matrix, and by one type
(vertical or horizontal) of crystal operations, contains a unique vertex at
which the potentials for all raising operations are zero. Moreover, the crystal
is determined up to isomorphism by the potentials for the lowering operations
at that vertex.

\proof
We have seen that exhausting upward or leftward crystal operations in all
possible ways from a given initial matrix always leads to the same final
matrix (the matrix $P$ respectively~$Q$ of theorem~\decompthm), where the
potentials for all raising operations are zero. Then whenever two matrices
$M,M'$ are linked by a crystal operation in the chosen direction, this final
matrix will be the same whether $M$~or~$M'$ is taken as point of departure.
Consequently, every matrix in the crystal considered leads to the same final
matrix, which is the unique matrix in the crystal where the potentials for all
raising operations are zero.

For the final part of the claim, consider the common normal form
$N=\diagram\\$ or $N=\diagon\\$ of each of $M$,~$P$, and~$Q$ in
theorem~\decompthm. Since horizontal crystal operations induce isomorphisms of
crystals defined by vertical crystal operations (due to their commutation) and
vice versa, the crystal for vertical operations containing~$P$ is isomorphic
to the one for the same operations containing~$N$, and similarly for the
crystal for horizontal operations containing~$Q$; it will suffice to show that
$\\$ can be deduced from the sequences of potentials for the appropriate
lowering operations at $P$ respectively from those at~$Q$. Now
$\ndm_i(P)=\ndm_i(N)=\\_i-\\_{i+1}$ for all~$i$ in both the binary and the
integral case, while $\nrm_j(Q)=\nrm_j(N)$ for all~$j$, which equals
$\\\tr_j-\\\tr_{j+1}$ in the binary case, or $\\_j-\\_{j+1}$ in the integral
case. Since $\\$ can be reconstructed from the sequence of differences its
consecutive parts or of those of~$\\\tr$ (for instance
$\\_k=\sum_{i\geq{k}}d_i$ for all~$k\in\N$, where $d_i=\\_i-\\_{i+1}$), this
completes the proof.
\QED

Part~\statitemnr2 of the theorem is very similar to the statement of
bijectivity of the \caps{RSK}~correspondence (see for instance \ref{Stanley
EC2, theorem~7.11.5}), since the matrices~$P$ and~$Q\tr$ both lie in
$\IE{\\/(0)}$ where $\\$ is the partition $\rsum{P}=\csum{Q}$, whence they
encode straight semistandard tableaux of equal shape~$\\$; the statements
$\csum{P}=\csum{M}$ and $\rsum{Q}=\rsum{M}$ give the usual relation between
the weights of these tableaux and the column and row sums of~$M$. However, the
construction of the bijection is quite different, and as we shall see, the
bijection itself is does not correspond directly to the~\caps{RSK}
correspondence either. Yet it is appropriate to call the bijection of
part~\statitemnr2 a ``Knuth correspondence''. This is a generic term for a
large class of generalisations of the \caps{RSK}~correspondence that was used
in~\ref{Fomin Schur} and in~\ref{spin Knuth}, without being formally defined.
It actually refers to the result of the type of construction defined
in~\ref{Fomin Schur} using growth diagrams; that construction is rather
different from the one above, but as we shall see below, the bijection
considered does allow an alternative construction using growth diagrams.

As for the bijection of Part~1, it has some resemblance to dual \caps{RSK}
correspondence, the second bijection that was defined together with the usual
\caps{RSK}~correspondence in~\ref{Knuth PJM}. That correspondence is a
bijection between binary matrices and pairs of straight tableaux that are less
symmetric than those associated to integral matrices: one can define them to
be both semistandard but of transpose shapes (as in the original paper), or
one can define one of them to be semistandard and the other transpose
semistandard, in which case their shapes will be equal (as in \ref{Stanley
EC2, \Sec7.14}); in fact at the end of the original construction one of the
tableaux is transposed to make both semistandard. In our bijection however,
while $Q$ is the binary encoding of a straight semistandard tableau of
shape~$\\=\csum{Q}\tr$, it is not entirely natural to associate a straight
semistandard or transpose semistandard tableau to~$P$. One could rotate~$P$ a
quarter turn counterclockwise to get the binary encoding of a semistandard
tableau of shape~$\rsum{P}\tr=\\\tr$, but that matrix, and the weight of the
tableau it encodes, depend on the number of columns of the matrix that is
rotated, which should be large enough to contain all bits~`$1$' of~$P$. The
type of straight tableau most naturally associated to~$P$ is a reverse
transpose semistandard tableau of shape~$\rsum{P}=\\$, namely the tableau
given by the sequence of partitions $(\\^{(l)})_{l\in\N}$ where
$\\^{(l)}=\sum_{j\geq{l}}P\tr_j$.

\subsection Relation with jeu de taquin and with Robinson's correspondence.
\labelsec\Robsec

Let us give an example of the bijections of theorem~\decompthm. We could take
arbitrary matrices~$M$ in both cases, but it will be instructive to use the
binary and integral encodings of a same semistandard tableau~$T$. We shall
take
\bigdisplay
  T=\Skew(4:0,2,4,5,5|1:0,1,3,4,6,6,6|0:1,1,2,4,5|0:2,3,4,6,6|0:5,6,6)
,\label(\Tdef)
$$
that was also used as illustration in~\Sec\encodingsec. Taking its binary
encoding~$M$ in theorem~\decompthm\statitemnr1, one has
\bigdisplay
\def\mapleft#1{\smash
       {\mathop{\hbox to 2.5em{\leftarrowfill}}\limits^{\textstyle#1}}}
\matrix
{\llap{$N={}$}
 \pmatrix{1&1&1&1&1&1&1&1&0\cr
          1&1&1&1&1&1&1&1&0\cr
	  1&1&1&1&1&0&0&0&0\cr
	  1&1&1&0&0&0&0&0&0\cr
	  1&0&0&0&0&0&0&0&0\cr
	  0&0&0&0&0&0&0&0&0\cr
	  0&0&0&0&0&0&0&0&0\cr}
&\mapleft{R\xleft}
&\pmatrix{1&1&1&0&1&1&1&1&1\cr
          1&1&1&1&1&1&1&1&0\cr
	  1&1&1&1&1&0&0&0&0\cr
	  0&1&0&1&1&0&0&0&0\cr
	  0&0&1&0&0&0&0&0&0\cr
	  0&0&0&0&0&0&0&0&0\cr
	  0&0&0&0&0&0&0&0&0\cr}
 \rlap{${}=P$}
\cr\noalign{\nobreak\smallskip}
  \llap{$R\xraise$}\bigg\uparrow & & \llap{$R\xraise$}\bigg\uparrow
\cr\noalign{\nobreak\smallskip}
 \llap{$Q={}$}
 \pmatrix{1&1&0&0&0&0&0&0&0\cr
          1&1&1&0&0&0&0&0&0\cr
	  1&1&0&1&0&0&0&0&0\cr
	  1&0&1&0&0&0&0&0&0\cr
	  0&0&1&1&1&1&0&0&0\cr
	  1&0&0&0&1&0&1&1&0\cr
	  0&1&1&1&1&1&1&1&0\cr}
&\mapleft{R\xleft}
&\pmatrix{0&1&0&0&1&0&0&0&0\cr
          1&1&1&0&0&0&0&0&0\cr
	  1&0&1&0&0&1&0&0&0\cr
	  0&1&0&1&0&0&0&0&0\cr
	  0&0&1&1&1&0&1&0&0\cr
	  1&0&0&0&1&0&0&1&1\cr
	  0&1&1&1&1&1&1&1&0\cr}
 \rlap{${}=M$,}
\cr
}
\label(\PQbindiagram)
$$
where the composite raising operations are for instance $R\xleft=
\def~#1{\>\ltm_{#1}} \ltm_7~6~5~2~3~4~0~1~2~3~1~0$ and
$R\xraise=
\def~#1#2{\penalty\binoppenalty\>\ifnum#2=1 \upm_#1\else(\upm_#1)^#2\fi}
(\upm_1)^3~25~37~47~57~02~13~23~33~44~02~13~24~34~11~21~01~11~01$; as usual
operations are applied from right to left, and we have systematically chosen
to climb the ladder with the lowest possible index first, and to climb it
completely (so without switching ladders when a lower index ladder becomes
available halfway up).
As an illustration of theorem~\decompthm\statitemnr2, taking for~$M$ the
integral encoding of~$T$, one has the diagram
\bigdisplay
\def\mapleft#1{\smash
       {\mathop{\hbox to 2.5em{\leftarrowfill}}\limits^{\textstyle#1}}}
\matrix
{\llap{$N={}$}
 \pmatrix{8&0&0&0&0&0&0\cr
          0&8&0&0&0&0&0\cr
	  0&0&5&0&0&0&0\cr
	  0&0&0&3&0&0&0\cr
	  0&0&0&0&1&0&0\cr}
&\mapleft{R\xleft}
&\pmatrix{2&1&1&0&2&2&0\cr
          0&2&0&1&1&1&3\cr
	  0&0&2&0&1&0&2\cr
	  0&0&0&1&0&0&2\cr
	  0&0&0&0&0&1&0\cr}
 \rlap{${}=P$}
\cr\noalign{\nobreak\smallskip}
  \llap{$R\xraise$}\bigg\uparrow & & \llap{$R\xraise$}\bigg\uparrow
\cr\noalign{\nobreak\smallskip}
 \llap{$Q={}$}
 \pmatrix{5&0&0&0&0&0&0\cr
          2&5&0&0&0&0&0\cr
	  1&2&2&0&0&0&0\cr
	  0&1&2&2&0&0&0\cr
	  0&0&1&1&1&0&0\cr}
&\mapleft{R\xleft}
&\pmatrix{1&0&1&0&1&2&0\cr
          1&1&0&1&1&0&3\cr
	  0&2&1&0&1&1&0\cr
	  0&0&1&1&1&0&2\cr
	  0&0&0&0&0&1&2\cr}
 \rlap{${}=M$,}
\cr
}
\label(\PQintdiagram)
$$
where the composite raising operations are for instance
(using the same rule to choose ladders) given by $R\xleft
\def~#1#2{\penalty\binoppenalty\>\ifnum#2=1 \ltm_#1\else(\ltm_#1)^#2\fi}
=(\ltm_1)^3~25~37~47~57~02~13~23~33~44~02~13~24~34~11~21~01~11~01
$
and $R\xraise= \def~#1#2{\>\ifnum#2=1 \upm_#1\else(\upm_#1)^#2\fi}
\upm_2~32~11~23~01~13~02$.

Comparing (\PQbindiagram) and (\PQintdiagram), one may observe that in both
cases the normal form $N$ is parametrised by the same partition
$\\=(8,8,5,3,1)$. This can be explained using proposition~\bininteqprop, which
implies that the same composite raising operations can be taken for $R\xraise$
in the binary case as for $R\xleft$ in the integral case, up to the
replacement of $\upm_i$ by~$\ltm_i$; therefore the binary matrix~$P$ and the
integral matrix~$Q$ both encode the same tableau~$L$, of the shape
$(9,8,5,5,3)/(4,1)$ of~$M$. The raising operations in the given direction are
exhausted, so this is a Littlewood-Richardson tableau, namely
\bigdisplay
  L=\Skew(4:0,0,0,0,0|1:0,0,1,1,1,1,1|0:0,1,1,2,2|0:1,2,2,3,3|0:2,3,4)
.\label(\Ldef)
$$
The partition $\weight{L}=(8,8,5,3,1)$ is equal to $\rsum{P}$ in the binary
case and to $\csum{Q}$ in the integral case, and in both cases it is therefore
also equal to the partition~$\\=\rsum{N}$ that parametrises~$N$. One can also
read off $\\$ from the matrix~$Q$ in the binary case, or from the matrix~$P$
in the integral case, namely as $\\=\csum{Q}\tr$ respectively as
$\\=\rsum{P}$. Those matrices are in fact the binary and integral encodings of
the rectification~$S$ of~$T$, which has shape~$\\$, since exhausting inward
jeu de taquin slides applied to~$T$ translates into particular sequences that
exhaust leftward respectively upward crystal operations when applied to~$M$,
by proposition~\jdtrelprop. That rectification is easily computed to be
\bigdisplay
  S=\Young(0,0,1,2,4,4,5,5|1,1,3,4,5,6,6,6|2,2,4,6,6|3,6,6|5)
\label(\Sdef)
$$
and its binary and integral encodings are indeed the binary matrix~$Q$ above
and the integral matrix~$P$, respectively. The matrices $L$ and~$S$ are the
ones associated to~$T$ under Robinson's bijection, as is described
in~\ref{L-R intro, corollary~3.3.4}; in the notation of that reference one has
$(L,S)={\cal R}(T)$.

Thus in a certain sense both bijections of theorem~\decompthm\ model
Robinson's bijection. However that point of view depends on interpreting
matrices as encodings of tableaux, which might not always be natural.
In any case one should be aware that one matrix can simultaneously encode
tableaux of many different shapes. This point of view also gives somewhat
incomplete information about double crystals, because even though the
exhaustion of appropriate raising operations can be realised using jeu de
taquin, it is not easy to do the same with individual raising operations.

\subsection Identifying the factors in the decomposition.

Since any matrix can easily be seen to encode a skew semistandard tableau of
an appropriately chosen shape, the relation established above gives a complete
characterisation of the bijections in theorem~\decompthm\ that is independent
of crystal operations. However, it describes the matrices $P$~and~$Q$ in
rather different terms, which seems unnatural given the characterisation of
theorem itself. As we have observed above, the theorem in fact suggests a
close relationship with the well known \caps{RSK}~correspondence (and its
dual), in which the two straight tableaux associated to a matrix play a much
more comparable role. We shall now proceed to characterise the bijections of
our theorem in the language of those correspondences.

It is well known that the ``$P$-tableau'' of the pair that is associated under
the \caps{RSK}~correspondence to an integral matrix~$M$ can be computed using
rectification by jeu de taquin; since the matrix~$P$ of
\decompthm\statitemnr2 can also be determined using rectification, this
provides a useful starting point to relate the two correspondences. In more
detail, the mentioned $P$-tableau can be found by rectifying a semistandard
tableau~$H$ whose shape is a horizontal strip, and whose ``rows'', taken from
left to right, successively contain words (each of which is necessarily weakly
increasing) whose weights given by the rows of~$M$, from top to bottom. This
observation was first made, for the case of Schensted correspondence,
in~\ref{Schutzenberger CRGS}, see also \ref{Stanley EC2, A1.2.6}; the
generalisation to the \caps{RSK}~correspondence follows directly from its
definition. For instance, the initial skew tableau corresponding to the matrix
\bigdisplay
  M=\pmatrix{0&2&3&1\cr2&2&0&1\cr1&0&4&0\cr}
\ttext{would be}
  H=\Skew(11:0,2,2,2,2|6:0,0,1,1,3|0:1,1,2,2,2,3)
.
$$
Note that by filling the horizontal strip from left to right, the vertical
order of rows is interchanged: the topmost row of~$M$ determines the bottom
row of~$H$. Therefore the relation between $M$ and~$H$ is not that of a matrix
encoding a tableau; rather the integral encoding of~$H$ is the matrix formed
by reversing the order of the rows of~$M$. As for the matrix~$P$
of~\decompthm\statitemnr2, we have seen above that it is the integral encoding
of the rectification~$S$ of any tableau~$T$ of which $M$ is the integral
encoding. There are many possibilities for~$T$, all of which differ only by
horizontal shifts of their rows and which therefore have the same
rectification. It is in particular always possible to take for~$T$ a tableau
whose shape is a horizontal strip. Comparing with what we said for the
\caps{RSK}~correspondence, and using the fact that both for that
correspondence and for the one of theorem~\decompthm\statitemnr2,
transposing~$M$ results in interchanging $P$ and~$Q$ (or actually $P$ and
$Q\tr$ for our theorem), we arrive at the following description.

\proclaim Proposition. \RSKcompareprop
Let $M\in\Mat$, which we view as a some sufficiently large finite rectangular
matrix. The corresponding matrices $P,Q\in\Mat$ of
theorem~\decompthm\statitemnr2 can be found as follows: $P$ is the integral
encoding of the $P$-symbol under the \caps{RSK}~correspondence of the matrix
obtained from~$M$ by reversing the order of its rows, and $Q\tr$ is the
integral encoding of the $Q$-symbol under the \caps{RSK}~correspondence of the
matrix obtained from~$M$ by reversing the order of its columns.
\QED

Note that reversal of the order of rows or columns is not a well defined
operation on~$\Mat$, but the resulting ambiguity about the matrices of which
the $P$-symbol and the $Q$-symbol are to be taken is harmless, since it only
affects the other symbol associated to each matrix, the one that is not used in
the proposition. For instance the matrix obtained by the reversal of the order
of rows is only determined up to a vertical shift, but this only induces a
shift in the entries of the associated $Q$-symbol, which is unused.

The description of the proposition can be reformulated as follows: $P$ and
$Q\tr$ are the integral encodings of the insertion tableau for the traversal
of~$M$ in reverse Semitic reading order (by rows, read from left to right,
ordered from bottom to top), row-inserting at each entry $M_{i,j}$ copies of
the number~$j$, respectively of the recording tableau for a similar insertion,
but of negated values and traversing $M$ in the opposite (Semitic) order. The
latter tableau can also be described either as the Sch\"utzenberger dual of
the recording tableau for the initial insertion process, or as the insertion
tableau for a traversal of $M$ in Kanji reading order, where at each entry
$M_{i,j}$ copies of the number~$i$ are row-inserted, but in all these cases
this involves two separate applications of an \caps{RSK} or Sch\"utzenberger
algorithm.

However the pair of tableaux of which $P$~and~$Q\tr$ are encodings is in fact
the pair associated to~$M$ under the ``Burge correspondence''. This is a
variation of the \caps{RSK}~correspondence that was first defined (for a
somewhat specialised case) in~\ref{Burge}, and which is discussed in detail
in~\ref{Fulton, A.4.1} and in~\ref{spin Knuth, \Sec3.2}. To be more specific,
this gives us the following description.

\proclaim Proposition. \Burgecorporp
The matrices $P,Q\in\Mat$ corresponding to $M\in\Mat$ under the bijection of
theorem~\decompthm\statitemnr2 are respectively the integral encoding of the
$P$-symbol and the transpose of the integral encoding of the $Q$-symbol that
are associated to $M$ under the Burge correspondence. This $P$-symbol can be
obtained by traversing~$M$ in the Semitic reading order (by rows from right to
left) and successively \emph{column} inserting, into an initially empty
tableau, for each position $(i,j)$ encountered $M_{i,j}$ copies of the
number~$j$; the $Q$-symbol is the recording tableau $(\\^{(i)})_{i\in\N}$ for
this procedure, where $\\^{(i)}$ is the shape of the tableau under
construction before the traversal of row~$i$ of~$M$.

The sequence of numbers that is being inserted is also the Semitic reading of
any tableau of which $M$ is the integral encoding; for instance for the
matrix~$M$ in~(\PQintdiagram) this is the sequence $\def\;{;\penalty4000
}5,6,6\;2,3,4,6,6\;1,1,2,4,5\;0,1,3,4,6,6,6\;0,2,4,5,5$, read \emph{from right
to left}. (For column insertion it is best to imagine the entries being
inserted as coming from the right, so writing this sequence backwards is
actually quite natural; for instance, in the notation of
\ref{Fulton} the insertion tableau may be written as
$5\to(6\to(6\to(2\to\cdots\to(2\to(4\to(5\to(5))))\cdots)))$, where we
added redundant parentheses to stress the right-associativity of the
operator~`$\to$'.) Since column insertion can be simulated by jeu de taquin
just like row insertion can, with as only difference that the entries inserted
are initially arranged into a tableau according to the Semitic reading order,
it is clear that the tableau obtained by column inserting the sequence given
is equal to the rectification~$S$ of the tableau~$L$ given in~(\Sdef).
This argument justifies the ``$P$''-part of the proposition above; to justify
its ``$Q$''-part, the easiest argument is to recall that the Burge
correspondence enjoys the same symmetry property as the
\caps{RSK}~correspondence, and that the bijection of
theorem~\decompthm\statitemnr2 has a similar symmetry.

We illustrate this computation for the integral matrix~$M$ of~(\PQintdiagram).
We are performing column insertion, and the display below should be read from
right to left; it shows the state of the insertion process after each row of
$M$ is processed. The entries being inserted are written above the arrows
between the tableaux; they too should be read from right to left to get their
order of insertion.
\bigdisplay
\hbox
{$
 \smallsquares
 \diaght=\ht\strutbox \advance\diaght by \dp\strutbox
 \multiply \diaght by 5 \advance\diaght by -\dp\strutbox
 \def\tableau(#1){\vcenter to\diaght{\Young(#1)\vss}}
 \def\ins(#1){\buildrel #1\over\longleftarrow}
 \tableau(0,0,1,2,4,4,5,5|1,1,3,4,5,6,6,6|2,2,4,6,6|3,6,6|5)\ins(5,6,6)
 \tableau(0,0,1,2,4,4,5,5|1,3,4,5,6,6,6|2,2,4,6|3,6)\ins(2,3,4,6,6)
 \tableau(0,0,1,2,4,4,5,5|1,3,5,6,6,6|2,4)\ins(1,1,2,4,5)
 \tableau(0,0,2,4,4,5,5|1,3,6,6,6)\ins(0,1,3,4,6,6,6)
 \tableau(0,2,4,5,5)\ins(0,2,4,5,5)
$}
\label(\Burgeexamp)
\vadjust{\smallskip}
$$
Indeed the final tableau coincides with $S$ of~(\Sdef) that is encoded by $P$
of~(\PQintdiagram), and the sequence of shapes
$(0)\leh(5)\leh(7,5)\leh(8,7,2)\leh(8,8,4,2)\leh(8,8,5,3,1)$ defines the
semistandard Young tableau
\bigdisplay
  \bar{L}=\Young(0,0,0,0,0,1,1,2|1,1,1,1,1,2,2,3|2,2,3,3,4|3,3,4|4)
\label(\Lbardef)
$$
encoded by~$Q\tr$; it is the companion Young tableau of the
Littlewood-Richardson tableau~$L$ of~(\Ldef).

For the bijection of~\decompthm\statitemnr1, one can use a similar argument to
identify the matrix~$Q$, which is the binary encoding of the rectification of
any tableau~$T$ encoded by the binary tableau~$M$. There always exists such a
tableau~$T$ whose shape is a vertical strip, and as in the case of a
horizontal strip, the rectification can then be computed as the insertion
tableau either for the column insertion of the entries taken from top right to
bottom left, of for the row insertion of the entries taken in the opposite
order. In either case the shape of the tableau grows by a vertical strip after
the insertion of the entries for each complete column of~$T$, which defines a
transpose semistandard recording tableau; since its vertical strips are most
naturally indexed by the numbers of the columns of~$T$ traversed, the
recording tableau for the right-to-left column insertion is actually a reverse
tableau. As we have observed above, the binary matrix~$P$ of
\decompthm\statitemnr1 corresponds to a reverse transpose semistandard
tableau, and this is in fact the recording tableau of the top right to bottom
left column insertion of the entries of~$T$ (for the row insertion of the
entries in the opposite order, it is the Sch\"utzenberger dual of the
recording tableau). One thus has the following description of the bijection
between $M$ and $(P,Q)$ of \decompthm\statitemnr1.

\proclaim Proposition. \Knuthbinaryprop
The matrices $P,Q\in\bMat$ corresponding to $M\in\bMat$ under the bijection of
theorem~\decompthm\statitemnr1 can be computed as follows. The matrix~$Q$ is
the binary encoding of the semistandard Young tableau~$S$ obtained by
traversing the matrix~$M$ in the Kanji reading order (down columns, from right
to left), and successively \emph{column} inserting, into an initially empty
tableau, for each bit $M_{i,j}$ equal to~$1$ an entry~$i$. The matrix~$P$
determines the reverse transpose semistandard Young tableau
$R=(\\^{(j)})_{j\in\N}$ that is the recording tableau for this insertion
process: $\\^{(j)}=\sum_{l\geq{j}}P\tr_l$ is the shape of the tableau under
construction after the the traversal of column~$j$ of~$M$.

We illustrate this computation for the binary matrix~$M$ of our running
example. As before we display successive stages of column insertion from right
to left; the strictly increasing sequence of numbers inserted at each
rightward arrow is written from top to bottom, as it is a column of the
tableau~$T$ of (\Tdef).
\bigdisplay
\advance \belowdisplayskip by 5pt
\vcenter
{\advance\displaywidth by -4em 
 \smallsquares
 \diaght=\ht\strutbox \advance\diaght by \dp\strutbox
 \multiply \diaght by 5 \advance\diaght by -\dp\strutbox
 \def\tableau(#1){\vcenter to\diaght{\Young(#1)\vss}}
 \def\ins(#1){\;\tower\longleftarrow#1,.\;}
 \def\tower#1#2,#3.{\setbox0=\hbox{$\buildrel#2\over{#1}$}%
     \if|#3|\box0
     \else \tower{\box0 }#3.\fi}
 \displaylines{\hfill
 \ins(6,4,3)
 \tableau(0,2,4,5,5|4,6,6,6|5|6)
 \ins(6,5,4,0)
 \tableau(2,4,5,5|6,6,6)
 \ins(4,2)
 \tableau(4,5,5|6,6)
 \ins(6,4)
 \tableau(5,5|6)
 \ins(6,5)
 \tableau(5)
 \ins(5)
 \qquad\cr\qquad
 \tableau(0,0,1,2,4,4,5,5|1,1,3,4,5,6,6,6|2,2,4,6,6|3,6,6|5)
 \ins(5,2,1)
 \tableau(0,0,2,4,4,5,5|1,1,3,5,6,6,6|2,4,4,6|3,6,6|6)
 \ins(6,3,1,0)
 \tableau(0,2,4,4,5,5|1,3,5,6,6,6|2,4,6|4,6|6)
 \ins(6,4,2,1)
 \tableau(0,2,4,5,5|3,4,6,6,6|4,5|6,6)
 \;\longleftarrow
 \hfill\cr}
}\label(\Knuthbinex)
$$
Indeed the final tableau coincides with the tableau $S$ of~(\Sdef) of which
the matrix~$Q$ of~(\PQbindiagram) is the binary encoding, and the reverse
transpose semistandard Young tableau~$R$ given by the sequence of shapes
$(8,8,5,3,1)\gev (7,7,4,3,1)\gev (6,6,3,2,1)\gev (5,5,2,2)\gev (5,4,1,1)\gev
(4,3)\gev (3,2)\gev (2,1)\gev (1)\gev (0)$ corresponds to the binary
matrix~$P$ of~(\PQbindiagram). This tableau can be displayed by filling for
each $i$ the squares of the vertical strip $\set{\\^{(i)}/\\^{(i+1)}}$ with
entries~$i$, giving
\bigdisplay R=\Young(8,7,6,5,4,2,1,0|7,6,5,4,3,2,1,0|4,3,2,1,0|4,3,1|2)
.\label(\Rdef)
$$

We shall finally compare the bijection of \decompthm\statitemnr1 with the dual
\caps{RSK}~correspondence. One should first of all realise that, unlike
for the ordinary \caps{RSK} and Burge correspondences, there is no universal
agreement of which pair of tableaux corresponds to a given binary matrix
(versions differ by transposition of one or both tableaux). We shall use the
form in \ref{Stanley EC2}, where the entries coming from a conventional reading
of~$M$ are row-inserted into a transpose semistandard Young tableau, and the
recording tableau is semistandard. Note that the following statement directly
relates to the tableaux coming from the dual \caps{RSK}~correspondence applied
to~$M$ (even if a Sch\"utzenberger dual is involved), whereas for the integer
correspondence we could not use the tableaux corresponding to $M$ itself under
the \caps{RSK}~correspondence.

\proclaim Proposition. \dualRSKrowprop
The matrices $P,Q\in\bMat$ corresponding to $M\in\bMat$ under the bijection of
theorem~\decompthm\statitemnr1 are related as follows to the pair $(R^*,S)$
that is associated to~$M$ under the dual \caps{RSK} correspondence as
described in~\ref{Stanley EC2, theorem~7.14.2}: the transpose semistandard
Young tableau $R^*$ is the Sch\"utzenberger dual of the reverse transpose
semistandard Young tableau~$R$ determined by~$P$, and $S$ is the semistandard
Young tableau encoded by~$Q$.

\subsection Relation with decomposition of pictures.
\labelsec\picturesec

We have established relations between the bijections of theorem~\decompthm\
and known correspondences defined by various insertion procedures, but the
descriptions could use some clarification. We have justified the
correspondences with insertion tableaux via their relation to jeu de taquin,
but the correspondences with the recording tableau are only justified using
symmetry properties for the integral case, and only by an example in the
binary case. We shall now present a correspondence that explains both factors
in a more symmetric way, and allows to transparently justify the descriptions
of the previous subsection.

If an integral matrix~$M$ satisfies $M\in\IE{\\/\kappa}$, it encodes a
semistandard tableau~$T$ of shape~$\\/\kappa$, and if it also satisfies
$M\in\IL{\nu/\mu}$, then $M\tr$ encodes a semistandard tableau~$\bar{T}$ of
shape~$\nu/\mu$, which is a companion tableau of~$T$. The pair $(T,\bar{T})$
can be represented in a symmetric manner by the notion of a \newterm{picture}
$\\/\kappa\to\nu/\mu$, which designates a particular kind of bijection between
skew diagrams $\set{\\/\kappa}\to\set{\nu/\mu}$, in a way that we shall now
describe.

In fact two different flavours of pictures can be defined, both of which have
the property that the set of pictures $\\/\kappa\to\nu/\mu$ is in bijection
with $\IE{\\/\kappa}\thru\IL{\nu/\mu}$, and both of which give rise to
essentially the same theory. Our reference for the theory of pictures will
be~\ref{pictures}, but we shall define the opposite flavour of pictures than
was used there, since that simplifies the correspondences that are of interest
here. This means that we must adapt the results we cite to our situation; in
fact the form of these results changes little between the two flavours of
pictures, except for explicit references to the Schensted row insertion
procedure (for instance in theorem~3.2.1 and lemma~3.3.2 of~\ref{pictures}),
which must be replaced by column insertion. The pictures used in the present
paper are defined as follows (the difference with the other flavour of
pictures is the use of `$\lene$' rather than its opposite partial ordering).

\proclaim Definitions. \picturedef
\defitem The partial orderings `$\lenw$' and `$\lene$' on $\N^2$ are
         defined by $(i,j)\lenw(k,l)$ if and only if $i\leq{k}$ and
         $j\leq{l}$, respectively by $(i,j)\lene(k,l)$ if and only if
         $i\leq{k}$ and $j\geq{l}$.
\defitem For skew shapes $\\/\kappa$ and~$\nu/\mu$, a picture
         $\\/\kappa\to\nu/\mu$ is a bijection
	 $f\:\set{\\/\kappa}\to\set{\nu/\mu}$ between their diagrams, such
	 that $s\lenw{t}$ implies $f(s)\lene{f(t)}$ and $f(s)\lenw{f(t)}$
	 implies~$s\lene{t}$ for all $s,t\in\set{\\/\kappa}$. The set of
	 pictures $\\/\kappa\to\nu/\mu$ is denoted by
	 $\Pic(\\/\kappa,\nu/\mu)$.
\defitem For a picture $f\in\Pic(\\/\kappa,\nu/\mu)$, the matrices
         $\Int(f)\in\mat{\\-\kappa}{\nu-\mu}$ and
         $\Bin(f)\in\bmat{\nu-\mu}{\\\tr-\kappa\tr}$ are respectively given by
         $\Int(f)_{i,k} =\Card\setof(j,l)\in\N^2:(i,j)\in\set{\\/\kappa},
         (k,l)\in\set{\nu/\mu}, f(i,j)=(k,l)\endset$ and by
         $\Bin(f)_{k,j}=\Card\setof(i,l)\in\N^2:(i,j)\in\set{\\/\kappa},
         (k,l)\in\set{\nu/\mu}, f(i,j)=(k,l)\endset$.

The expressions in part~\defitemnr3 can be described in words as follows:
$\Int(f)_{i,k}$ counts the number of squares of row~$i$ that the picture~$f$
maps to row~$k$, while $\Bin(f)_{k,j}$ counts the number of squares of
column~$j$ that~$f$ maps to row~$k$.

From the form of part~\defitemnr2 it is clear that the inverse of a picture
$\\/\kappa\to\nu/\mu$ is a picture $\nu/\mu\to\\/\kappa$. The point of
part~\defitemnr3 of the definition is the following fact.

\proclaim Proposition.
For any pair of skew shapes $\\/\kappa$ and~$\nu/\mu$, the operations
$\Int$~and~$\Bin$ define bijections from $\Pic(\\/\kappa,\nu/\mu)$
respectively to $\IE{\\/\kappa}\thru\IL{\nu/\mu}$ and to
$\BE{\\/\kappa}\thru\BL{\nu/\mu}$.

\proof
For a given matrix~$M$, a picture~$f$ with $\Int(f)=M$ or $\Bin(f)=M$ is
straightforwardly constructed, defining images of the squares
of~$\set{\\/\kappa}$ in the Semitic (in the integral case) or the Kanji (in
the binary case) order. The conditions of \picturedef\defitemnr2 allow only
one candidate image for each square, and one checks that the tableau and
Littlewood-Richardson conditions for~$M$ assure that this is a valid choice.
\QED

When $\nu/\mu$ is a horizontal strip, the conditions $M\in\IL{\nu/\mu}$ and
$M\in\BL{\nu/\mu}$ just require $\csum{M}=\nu-\mu$ respectively
$\rsum{M}=\nu-\mu$, so in this case $\Pic(\\/\kappa,\nu/\mu)$ is in bijection
with the set of tableaux in $\SSTx{\\/\kappa}$ of weight~$\nu-\mu$. Moreover,
if $f$ so corresponds to a tableau~$T$, then $\Int(f)$ and $\Bin(f)$ are
respectively the integral and binary encodings of~$T$. Thus we find our
encoding of semistandard tableaux by integral or binary matrices as a special
case of representing pictures by matrices using $\Int(f)$ or $\Bin(f)$.
Similarly, if $\nu/\mu$ is arbitrary but $\\/\kappa$ is a horizontal strip,
then the pictures in $\Pic(\\/\kappa,\nu/\mu)$ are in bijection, via their
inverse picture, with tableaux in $\SSTx{\nu/\mu}$ of weight~$\\-\kappa$. In
the latter case $\Int(f)=\Int(f\inv)\tr\in\IL{\nu/\mu}$ is the transpose of
the integral encoding of the tableau thus corresponding to~$f$; however there
is in this case no similar interpretation for~$\Bin(f)$ (which is in fact a
rather sparse matrix, having column sums of at most~$1$). On the other hand,
if $\\/\kappa$ is a vertical strip rather than a horizontal one, then it is
most practical to consider $\Bin(f)$ instead of $\Int(f)$: in this case
$\Pic(\\/\kappa,\nu/\mu)$ is in bijection with the set of reverse transpose
semistandard tableaux of shape $\nu/\mu$ and weight $\\\tr-\kappa\tr$, and
their correspondence to the matrix $\Bin(f)\in\BL{\nu/\mu}$ is exactly the one
that we used before.

A central result concerning pictures is the decomposition theorem
\ref{pictures, theorem~3.2.1}, which generalises the Robinson-Schensted and
\caps{RSK}~correspondences; it gives a bijection between
$\Pic(\\/\kappa,\nu/\mu)$ and the union over all~$\pi\in\Part$ of Cartesian
products $\Pic(\pi,\nu/\mu)\times\Pic(\\/\kappa,\pi)$. This bijection can be
defined by a version of the Schensted algorithm using (for our notion of
pictures) column insertion, translated into the language of pictures, as
described in~\ref{pictures, lemma~3.3.2}. It follows that for specific shapes
$\\/\kappa$ and~$\nu/\mu$, this decomposition specialises under the
appropriate correspondences to the column insertion Robinson-Schensted
correspondence, to the Burge correspondence, and to the version of the
dual~\caps{RSK} correspondence using column insertion described in our
proposition~\Knuthbinaryprop. For the Robinson-Schensted correspondence one
takes $\\/\kappa$ and $\nu/\mu$ to be shapes that are both horizontal and
vertical strips, so that $\Pic(\\/\kappa,\nu/\mu)$ consists of all bijections
$\set{\\/\kappa}\to\set{\nu/\mu}$, and that $\Pic(\pi,\nu/\mu)$ and
$\Pic(\\/\kappa,\pi)$ are both in bijection with the set of standard Young
tableaux of shape~$\pi$. For the other two cases, one takes $\nu/\mu$ to be a
horizontal strip, and $\\/\kappa$ to be a horizontal respectively a vertical
strip. Then the operations $\Int$ and $\Bin$ define bijections
$\Pic(\\/\kappa,\nu/\mu)\to\mat{\\-\kappa}{\nu-\mu}$ respectively
$\Pic(\\/\kappa,\nu/\mu)\to\bmat{\nu-\mu}{\\\tr-\kappa\tr}$, and the sets
$\Pic(\pi,\nu/\mu)$ and $\Pic(\\/\kappa,\pi)$ are in bijection with sets of
semistandard Young tableaux of shape~$\pi$ as described above. Thus these
latter two special cases of the decomposition of pictures give rise, under
application of $\Bin$ and $\Int$, precisely to the bijections that occur in
our theorem~\decompthm.

The description of the decomposition of pictures using the column insertion
Schensted algorithm is rather different from the definition of our bijections
in theorem~\decompthm. However, the decomposition of pictures has an
alternative description that is given in \ref{pictures, theorem~5.4.1}. This
description is in terms of two forms of jeu de taquin that are defined for
pictures, one that changes the shape $\\/\kappa$ of the domain, and another
that changes the shape~$\nu/\mu$ of the image of the picture. The two types of
operations commute \ref{pictures, theorem~5.3.1}, and a decomposition of
pictures based on these operations follows in a manner quite similar to the
way we obtained theorem~\decompthm. From the fact that for pictures Schensted
column insertion can be emulated by jeu de taquin (as it can for tableaux), it
follows that the two descriptions define the same decomposition of pictures.
Under the operation $\Int$, each domain slide of a picture corresponds to a
sequence of upward crystal operations on integer matrices, just like in
proposition~\jdtrelprop, and each image slide corresponds to a sequence of
leftward crystal operations; hence the decomposition of pictures translates
directly into the decomposition of integral matrices of
theorem~\decompthm\statitemnr2. For the operation $\Bin$, the situation is
similar, except that domain slides correspond to leftward crystal operations
and image slides correspond to upward crystal operations. Thus the
decomposition of pictures also translates into the decomposition of binary
matrices of theorem~\decompthm\statitemnr2.

With these links between the decomposition of matrices and that of pictures,
and the fact that the latter can be computed using Schensted column insertion,
the descriptions of the previous subsection that involve column insertion are
easy to understand. For instance, the facts that in proposition~\Burgecorporp\
the integral matrix~$M$ is traversed in Semitic order while in
proposition~\Knuthbinaryprop\ the binary matrix~$M$ is traversed in Kanji
ordering both reflect the fact that the column insertions used to compute the
first factor of the decomposition of a picture traverse its domain from
top-right to bottom-left (in fact any order will do that never makes a
decreasing step for~`$\lene$'). This is so because $\mat\alpha\beta$ is
modelled by $\Pic(\\/\kappa,\nu/\mu)$ where $\\/\kappa$ and $\nu/\mu$ are
horizontal strips with $\\-\kappa=\alpha$ and $\nu-\mu=\beta$, and with the
rows of the matrix corresponding to those of $\\/\kappa$, while
$\bmat\alpha\beta$ is modelled by $\Pic(\\/\kappa,\nu/\mu)$ where $\\/\kappa$
is a vertical strip and $\nu/\mu$ a horizontal strip, with
$\\\tr-\kappa\tr=\beta$ and $\nu-\mu=\alpha$, and with the columns of the
matrix corresponding to those of $\\/\kappa$. In both cases top-right to
bottom-left traversal of $\set{\\/\kappa}$ corresponds to the indicated
traversal of the matrix. The description of the recording tableaux is also a
direct translation of the way recording tableaux are defined for pictures. For
propositions \RSKcompareprop~and~\dualRSKrowprop, which use row insertion, the
more complicated formulations reflect a less direct correspondence to the
decomposition of pictures; they can however be understood by using the
symmetry properties of pictures that are discussed in \ref{pictures, \Sec4};
we shall not elaborate the details here.

 \par

\newsection Growth diagrams.
\labelsec\growthsec

After what was said in the last subsection, one might be tempted to conclude
that crystal operations on matrices and the correspondences of
theorem~\decompthm\ are merely shadows of similar structures for pictures
under the maps `$\Int$' and~`$\Bin$', and that as such they do not add much
insight, although they are somewhat simpler to present than pictures and the
operations of jeu de taquin on them. In this section we shall however
establish a link with the descriptions as those in \ref{Roby} and~\ref{Fomin
Schur} of Knuth correspondences using growth diagrams, which link does not
exist in the theory of pictures. In doing so we build a bridge between the
previously unrelated approaches to Knuth correspondences via jeu de taquin on
one hand, and via growth diagrams on the other hand (the approach via jeu de
taquin often also uses invariants, an in \ref{Schutzenberger CRGS}
or~\ref{Fulton}; we shall discuss these below in \Sec\forward\chainsec).

\subsection Submatrices and implicit shapes.

For a discussion of the notion of growth diagrams we refer to \ref{spin Knuth,
\Sec\Sec2,3}. In general, a growth diagram associates ``shapes'' to the
elements of the grid~$\N\times\N$, and it also associates ``entries'' to
squares of the grid; in our case the shapes will be partitions, and the
entries are our usual matrix entries, lying either in~$\set2=\{0,1\}$ or
in~$\N$. Formally a square of the grid is an element of~$\N\times\N$, just
like a grid point is, but the ``square''~$(i,j)$ is supposed to have as
corners the ``grid points'' in $\{i,i+1\}\times\{j,j+1\}$. We shall build a
growth diagram from a given matrix~$M$ by associating to the grid point
$(k,l)$ a shape determined by the submatrix
$(M_{i,j})_{i\in\set{k},j\in\set{l}}$, of which that grid point is the lower
right hand corner. The shape will in fact be the one that parametrises the
normal form of that submatrix, which by our usual interpretation of finite
matrices is identified with the matrix obtained from~$M$ by clearing the
entries outside $\set{k}\times\set{l}$; we shall call this the implicit shape
of the submatrix. The formal definitions are as follows.

\proclaim Definitions.
Let $M\in\bMat$ or~$M\in\Mat$, and let $N$ be the normal form of~$M$,
obtainable from~it by exhausting the operations chosen from
$\setof\upm_i:i\in\N\endset\union\setof\ltm_j:j\in\N\endset$. We define
\defitem $\is{M}=\rsum{N}$, which we call the implicit shape of~$M$;
\defitem $M_{I\times{J}}=(\Kr{(i,j)\in{I}\times{J}}M_{i,j})_{i,j\in\N}$  for
         intervals $I,J\subset\N$, which is the matrix of the same type (binary
         or integral) as~$M$ and obtained from if by replacing all entries
         $M_{i,j}$ with $(i,j)\notin{I}\times{J}$ by~$0$.

The term implicit shape refers to the fact that it cannot be easily read off
directly from~$M$, although in \Sec\forward\chainsec\ we shall give a
characterisation that does not require modifying~$M$. Once either of the
factors $P,Q$ associated to~$M$ under the decomposition of theorem~\decompthm\
has been determined, the implicit shape of~$M$ can however be read off, for
instance as $\is{M}=\rsum{P}$; in other words it is the shape of the tableaux
corresponding to~$M$ under the Burge correspondence or the column insertion
dual \caps{RSK}~correspondence. If $M$ is the encoding of a skew semistandard
tableau~$T$, then $\is{M}$ is the shape of the rectification of~$T$.

Growth diagrams should satisfy a local rule, which we called ``shape datum''
in~\ref{spin Knuth}, that relates the shapes at the four corners of any grid
square~$(i,j)$ and the matrix entry in that square. This rule is such that
when the shapes at $(i+1,j)$ and $(i,j+1)$ are arbitrarily fixed, it defines a
bijection between the possible values for the pair consisting of the shape
at~$(i,j)$ and the matrix entry on one hand, and the possible values for the
shape at~$(i+1,j+1)$ on the other hand. When the shape datum is prescribed,
this permits deducing the entire growth diagram from a partial specification
of it, namely from the shapes along any lattice path in the grid from bottom
left to top right, together with all matrix entries to the bottom right of
that path. In our current situation this will eventually provide alternative
methods to compute implicit shapes, but at this point we just have the matrix
of shapes $(\is{M_{\set{i}\times\set{j}}})_{i,j\in\N}$, which is defined
without any reference to a shape datum. So we must first prove that it indeed
defines a growth diagram, and determine the shape datum for which it does so.

\subsection Growth diagrams for integral matrices.

The key to proving that $(\is{M_{\set{i}\times\set{j}}})_{i,j\in\N}$ defines a
growth diagram will be the study of a particular way of reducing~$M$ to its
normal form. Consider for the moment integral matrices, and take the
matrix~$M$ of~(\PQintdiagram). We start by exhausting the raising operations
in $\setof\upm_i:i\in\set2\endset$ and then exhaust those in
$\setof\ltm_j:j\in\set5\endset$; this gives the following transformations
\bigdisplay
  \pmatrix{1&0&1&0&1&2&0\cr
           1&1&0&1&1&0&3\cr
	   0&2&1&0&1&1&0\cr
	   0&0&1&1&1&0&2\cr
	   0&0&0&0&0&1&2\cr}
\buildrel (\upm_\set2)^* \over
\longrightarrow
\pmatrix{2&1&1&0&2&2&0\cr
         0&2&0&1&0&1&3\cr
	 0&0&1&0&1&0&0\cr
\matrixh
	 0&0&1&1&1&0&2\cr
	 0&0&0&0&0&1&2\cr}
\buildrel (\ltm_\set5)* \over
\longrightarrow
\pmatrix{8&0&0&0&0&0&0\cr
         0&4&0&0&0&0&3\cr
	 0&0&2&0&0&0&0\cr
\matrixcross316
	 0&1&1&1&0&0&2\cr
	 0&0&0&0&1&0&2\cr}
.\label(\imgrowthex)
$$
One sees that the submatrix indexed by $\set3\times\set6$ has been reduced to
the normal form $\diagon{(8,4,2)}$. We shall first argue that this is
indeed the normal form of the original submatrix $M_{\set3\times\set6}$, and
that this reduction is no accident (a similar reduction can be achieved for
any integral matrix).

\proclaim Proposition. \submatprop
Fix $M\in\Mat$, and $k,l\in\Np$. Let $M'\in\Mat$ be the matrix obtained
from~$M$ by exhausting the operations in $\setof\upm_i:i\in\set{k-1}\endset
\union\setof\ltm_j:j\in\set{l-1}\endset$, and let $N\in\Mat$ be the normal form
of~$M_{\set{k}\times\set{l}}$, which is obtained from it by exhausting the
same set of crystal operations. Then $M'_{\set{k}\times\set{l}}=N$.

\proof
It is clear that the normalisation of~$M_{\set{k}\times\set{l}}$ cannot
involve any raising operations other than the kind being applied to~$M$, as
$M_{\set{k}\times\set{l}}$ contains no units that could be transferred by such
raising operations. It is also clear that the operations applied to~$M$ cannot
transfer units across the boundaries of the rectangle $\set{k}\times\set{l}$:
they either transfer units within that rectangle or within its complement. If
we extract, from a sequence that exhausts these operations applied to~$M$,
those that involve a transfer within $\set{k}\times\set{l}$, then the
extracted sequence can be applied to $M_{\set{k}\times\set{l}}$, and each
operation will involve the same transfer as it did in the original sequence;
this is so because any valid transfer in a matrix remains valid in any
submatrix containing the squares involved in the transfer. Therefore the
extracted sequence reduces $M_{\set{k}\times\set{l}}$ to
$M'_{\set{k}\times\set{l}}$, and it remains to show that the latter matrix is
a normal from, which then necessarily equals~$N$. This means that we have to
prove that $\num_i(M'_{\set{k}\times\set{l}})
=\nlm_j(M'_{\set{k}\times\set{l}}) =0$ for any $i\in\set{k-1}$ and
$j\in\set{l-1}$, while we know that $\num_i(M') =\nlm_j(M') =0$. Now if $U|V$
denotes an integral matrix vertically split into submatrices $U,V$, then it
follows from (\intnumeq) that $\num_i(U|V)\geq\num_i(U)$, so $\num_i(U|V)=0$
implies $\num_i(U)=0$. For a horizontal division into submatrices, (\intnlmeq)
similarly gives $\nlm_j\({S\over T}\)\geq\nlm_j(S)$, so that $\nlm_j\({S\over
T}\) =0$ implies $\nlm_j(S)=0$. These relations allow us to conclude
$\num_i(M'_{\set{k}\times\set{l}}) =\nlm_j(M'_{\set{k}\times\set{l}}) =0$ as
desired.
\QED

Note that while the final part of this argument amounts to showing that if any
raising operation can be applied to a submatrix of the form
$M_{\set{k}\times\set{l}}$, then the same operation could also be applied to
the full matrix~$M$, this does not mean that these operations would
necessarily involve the same transfer. Indeed, if one has
$\num_i(U|V)>\num_i(U)>0$, then the first application(s) of~$\upm_i$ to $U|V$
will be between entries of~$V$, and only once $\num_i(U|V)$ and~$\num_i(U)$
have become equal will the remaining applications of~$\upm_i$ transfer entries
within~$U$. The same is true for leftward moves with respect to a horizontal
split. In the example given one sees that the exhaustion of operations
in~$\ltm_\set5$ involves leftward transfers below as well as above the
horizontal line; the former precede the latter in case of transfers between
the same pair of columns.

Let us now consider the relation between the four shapes
$\is{M_{\set{i}\times\set{j}}}$ for $(i,j)\in\sing{k,k+1}\times\sing{l,l+1}$
and the matrix entry~$M_{k,l}$. All of these values only depend in the
submatrix $M_{\set{k+1}\times\set{l+1}}$, and none of them changes when one
applies operations in $\setof\upm_i:i\in\set{k-1}\endset
\union\setof\ltm_j:j\in\set{l-1}\endset$. Therefore let $A$ be the matrix
obtained from~$M_{\set{k+1}\times\set{l+1}}$ by exhausting these operations,
and to simplify the notation put $\\=\is{M_{\set{k}\times\set{l}}}$, \
$\mu=\is{M_{\set{k}\times\set{l+1}}}$, \
$\nu=\is{M_{\set{k+1}\times\set{l}}}$,
and~$\kappa=\is{M_{\set{k+1}\times\set{l+1}}}$, which are also the implicit
shapes of the corresponding submatrices of~$A$ (here we abandon our previous
conventions for the use of $\\,\ldots,\kappa$, and adopt those of \ref{spin
Knuth}). By proposition~\submatprop\ we have
$A_{\set{k}\times\set{l}}=\diagon\\$.

The normalisation of $A_{\set{k+1}\times\set{l}}$ will only involve upward
transfers, since the leftward operations in
$\setof\ltm_j:j\in\set{l-1}\endset$ remain exhausted throughout. This
normalisation transfers any units present in (the bottom) row~$k$ of
$A_{\set{k+1}\times\set{l}}$ upwards, until they arrive on the main diagonal,
where together with the units that were already there in
$A_{\set{k}\times\set{l}}=\diagon\\$ they form the matrix $\diagon\nu$. It
follows that row~$k$ of $A_{\set{k+1}\times\set{l}}$ must have been equal to
the composition $\nu-\\$ (and in particular that it contained no units
strictly above the main diagonal). A similar argument shows that the final
column~$l$ of $A_{\set{k}\times\set{l+1}}$ is equal to $\mu-\\$. The matrix
$A$ is therefore completely determined by the values of $\\$, $\mu$, $\nu$,
and~$M_{k,l}$: one has $A_{\set{k}\times\set{l}}=\diagon\\$, \
$A_{i,l}=\mu_i-\\_i$ for $i\in\set{k}$, \ $A_{k,j}=\nu_j-\\_j$ for
$j\in\set{l}$, and $A_{k,l}=M_{k,l}$. This establishes one part of requirement
for families of shapes $(\is{M_{\set{i}\times\set{j}}})_{i,j\in\N}$ to define
a growth diagram, namely that $\\$, $\mu$, $\nu$, and~$M_{k,l}$ together
determine~$\kappa$; it also provides a procedure to compute this direction of
the shape datum (namely construct the matrix $A$ and find its implicit shape).
It is easy to see that one has $|\kappa|=|\mu|+|\nu|-|\\|+M_{i,j}$ as required
by the definition of a shape datum \ref{spin Knuth, definition~2.2.1}.

Let us apply the procedure for computing~$\kappa$ in the case of the
matrix~$M$ of~(\PQintdiagram) that was used in the example above, for $k=3$
and $l=6$; so our goal is to compute $\kappa=\is{M_{\set4\times\set7}}$. From
the last matrix in~(\imgrowthex) one reads off $\\=(8,4,2)$, \
$\mu-\\=(0,3,0)$ so that $\mu=(8,7,2)$, \ $\nu-\\=(0,1,1,1,0,0)$ so that
$\nu=(8,5,3,1)$, and $M_{3,6}=2$. We may now determine the following diagram,
from which it can be read off that $\kappa=(8,8,4,2)$ (unlike the previously
drawn diagrams we now let the arrows between the matrices point in the
direction opposite to the crystal operations that are being applied; this is
done so that the rectangle $\set{i}\times\set{j}$ for which the submatrix
has been normalised to $\is{A_{\set{i}\times\set{j}}}$ grows from left to
right and from top to bottom).
\bigdisplay
\def\mapright#1{\smash
       {\mathop{\hbox to 2.5em{\rightarrowfill}}\limits^{\textstyle#1}}}
\matrix
{\llap{$A={}$}
 \pmatrix{8&0&0&0&0&0&0\cr
          0&4&0&0&0&0&3\cr
	  0&0&2&0&0&0&0\cr
\matrixcross306
	  0&1&1&1&0&0&2\cr}
&\mapright{(\ltm_\set6)^*}
&\pmatrix{8&0&0&0&0&0&0\cr
          0&7&0&0&0&0&0\cr
	  0&0&2&0&0&0&0\cr
\matrixh
	  0&1&2&2&0&0&0\cr}
\cr\noalign{\nobreak\smallskip}
  \llap{$(\upm_\set3)^*$}\bigg\downarrow &
& \llap{$(\upm_\set3)^*$}\bigg\downarrow
\cr\noalign{\nobreak\smallskip}
  \pmatrix{8&0&0&0&0&0&0\cr
          0&5&0&0&0&0&3\cr
	  0&0&3&0&0&0&1\cr
\matrixv306
	  0&0&0&1&0&0&1\cr}
&\mapright{(\ltm_\set6)^*}
&\pmatrix{8&0&0&0&0&0&0\cr
          0&8&0&0&0&0&0\cr
	  0&0&4&0&0&0&0\cr
	  0&0&0&2&0&0&0\cr}
 \rlap{${}=\diagon\kappa$.}
\cr
}\nn
$$

To complete our argument we shall also have to show that conversely $\mu$,
$\nu$, and $\kappa$ determine $\\$ and~$M_{k,l}$, but before we do so we must
specify which relations should hold among the shapes for the above
construction to apply. First of all, each implicit shape
$\is{M_{\set{i}\times\set{j}}}$ can have at most $\min(i,j)$ nonzero parts.
Apart from that and $\\\subset\mu$ and $\\\subset\nu$ (which must hold since
neither $\nu-\\$ nor $\mu-\\$ have negative entries), there is only one
requirement for $\\$, $\mu$, and~$\nu$, namely that the matrix $A$ constructed
from the pieces $\diagon\\$, $\mu-\\$, $\mu-\\$ and~$M_{k,l}$ as indicated
must satisfy the hypothesis of initial exhaustion, in other words that
$\num_i(A)=\nlm_j(A)=0$ for all $i\in\set{k-1}$ and $j\in\set{l-1}$. This
amounts to the requirement $\\_{i+1}-\\_i+(\mu-\\)_{i+1}\leq0$ or equivalently
$\mu_{i+1}\leq\\_i$ for $i\in\set{k-1}$, and similarly $\nu_{j+1}\leq\\_j$ for
$j\in\set{l-1}$. Together with the limit on the number of nonzero parts, this
implies that $\\$, $\mu$, and~$\nu$ satisfy the relations $\\\leh\mu$ and
$\\\leh\nu$. Conversely, for any $\\,\mu,\nu\in\Part$ satisfying
$\mu_{\min(k,l+1)} =\nu_{\min(k+1,l)} =0$, and the relations $\\\leh\mu$, and
$\\\leh\nu$, the construction of~$A$ (with an arbitrary entry $A_{k,l}\in\N$)
provides an example where $\\$, $\mu$, and~$\nu$ occur as the implicit shapes
of appropriate submatrices.

Now assume that $\mu,\nu,\kappa\in\Part$ are given with $\mu_{\min(k,l+1)}
=\nu_{\min(k+1,l)} =0$, \ $\mu\leh\kappa$, and $\nu\leh\kappa$. We want to
find~$A$ as above, giving rise to these shapes. We know its normal form
$\diagon\kappa$, as well as the normal forms $\diagon\mu$ and $\diagon\nu$ of
$A_{\set{k}\times\set{l+1}}$ and $A_{\set{k+1}\times\set{l}}$, respectively.
Exhausting the operations $\setof\ltm_j:j\in\set{l}\endset$ applied to~$A$
must give the $(k+1)\times(l+1)$ matrix~$B$ whose restriction to
$\set{k}\times\set{l+1}$ equals $\diagon\mu$ and whose row~$k$ equals
$\kappa-\mu$. The normal form $\diagon\kappa$ of~$B$ can be written as
$R\xraise(B)$ for some sequence~$R\xraise$ of upward crystal operations in
$\setof\upm_i:i\in\set{k}\endset$. Then $R\xraise$ can also be applied to~$A$,
exhausting upward transfers and resulting in the $(k+1)\times(l+1)$ matrix~$C$
whose restriction to $\set{k+1}\times\set{l}$ equals $\diagon\mu$ and whose
column~$l$ equals $\kappa-\nu$. Then inverting the operations in~$R\xraise$,
it is clear that $A=R\xlower(C)$, and from this matrix one can extract the
values $\\=\is{A_{\set{k}\times\set{l}}}$ and $M_{k,l}=A_{k,l}$.

As an example of the computation in the inverse direction, let partitions
$\kappa=(13,9,9,5)$, $\mu=(11,9,8)$, and $\nu=(10,9,8,2)$ be given, and take
$k=3$ and $l=4$ (the minimal possible values, given $\mu$ and~$\nu$).
We form a commutative diagram like the one drawn above, but starting with the
matrices $B$, $C$, and~$\diagon\kappa$ below, which correspond to the
initial data. One easily determines composite raising operations transforming
$B$ or~$C$ into $\diagon\kappa$, for instance
$R\xraise=(\upm_0)^2\>(\upm_1)^2\>(\upm_2)^3$ and
$R\xleft=(\ltm_0)^3\>(\ltm_1)^3\>(\ltm_2)^4\>(\ltm_3)^7$ (the exponents are
just initial partial sums of the parts of the compositions
$\kappa-\mu=(2,0,1,5)$ and $\kappa-\nu=(3,0,1,3)$). Then one computes~$A$ by
inverting one of these composite operations, for instance $A=R\xlower(C)$
where $R\xlower=(\dnm_2)^3\>(\dnm_1)^2\>(\dnm_0)^2$, and from~$A$ one then
reads off $\\=(9,9,6)$ and $M_{3,4}=3$.
\bigdisplay
\def\mapright#1{\smash
       {\mathop{\hbox to 2.5em{\rightarrowfill}}\limits^{\textstyle#1}}}
\matrix
{\llap{$A={}$}
 \pmatrix{9&0&0&0&2\cr
          0&9&0&0&0\cr
	  0&0&6&0&2\cr
\matrixcross304
	  1&0&2&2&3\cr}
&\mapright{R\xleft}
&\pmatrix{11&0&0&0&0\cr
           0&9&0&0&0\cr
	   0&0&8&0&0\cr
\matrixh
	   2&0&1&5&0\cr}
 \rlap{${}=B$}
\cr\noalign{\nobreak\smallskip}
  \llap{$R\xraise$}\bigg\downarrow & & \llap{$R\xraise$}\bigg\downarrow
\cr\noalign{\nobreak\smallskip}
 \llap{$C={}$}
 \pmatrix{10&0&0&0&3\cr
           0&9&0&0&0\cr
	   0&0&8&0&1\cr
\matrixv304
	   0&0&0&2&3\cr}
&\mapright{R\xleft}
&\pmatrix{13&0&0&0&0\cr
           0&9&0&0&0\cr
	   0&0&9&0&0\cr
	   0&0&0&5&0\cr}
 \rlap{${}=\diagon\kappa$.}
\cr
}\label(\Burgeshapedatumex)
$$

As is clear from the examples, the only aspect of the reduction
$A\to\diagon\kappa$ that is not completely straightforward is the
transformation during the first reduction phase of the row or column that is
\emph{outside} the submatrix being normalised in that phase. Consider for
instance applying upward transfers to reduce the submatrix
$A_{\set{k+1}\times\set{l}}$ to its normal form $\diagon\nu$. Here column~$l$,
which initially describes the composition $\alpha=\mu-\\$ and the
entry~$M_{k,l}$, is transformed into a composition~$\alpha'$ that will
determine $\kappa=\nu+\alpha'$. The operations $\upm_i$ can be applied by
weakly decreasing indices~$i$. For fixed~$i$, the effect on rows $i,i+1$ of
exhausting the operation~$\upm_i$ is a transformation of the following form
(with $\beta=\nu-\\$):
$$
  \left({0\atop\beta_0}~{\cdots\atop\cdots}~{0\atop\beta_{i-1}}
 ~{\\_i\atop\beta_i} ~{0\atop\nu_{i+1}} ~{0~\cdots~0\atop0~\cdots~0}
 ~{\alpha_i\atop{m}}\right)
\buildrel (\upm_i)^* \over \longrightarrow
 \left({\beta_0\atop0}~{\cdots\atop\dots}~{\beta_{i-1}\atop0}
 ~{\nu_i\atop0} ~{0\atop\nu_{i+1}} ~{0~\cdots~0\atop0~\cdots~0}
 ~{m'\atop\alpha'_{i+1}}\right)
,\label(\Burgestep)
$$
where $\alpha'_{i+1}=\min(m,\\_i-\nu_{i+1})$, and
$m'=\alpha_i+m-\alpha'_{i+1}$. The number~$m$ is determined by the previous
raising operations, except for the initial case $i=k-1$, where $m=M_{k,l}$.
Since the diagonal entry $\nu_i$ is unchanged by further applications
of~$\upm_{i-1}$, and $\nu_i\geq\\_i\geq\nu_{i+1}+\alpha'_{i+1}$, there is no
need to reapply~$\upm_i$ later. Thus we have established a simple iterative
procedure for determining the transformation of the final column of the
matrix, which will then determine~$\kappa$.

We rewrite this procedure into an algorithm for directly computing~$\kappa$
from $\\$, $\mu$, $\nu$, and~$M_{k,l}$, by taking into account $\alpha=\mu-\\$
and $\kappa=\nu+\alpha'$; in this procure the roles of $\mu$ and of~$\nu$ will
be symmetrical. Adding $\nu_{i+1}$ to the equation for $\alpha'_{i+1}$ we
obtain $\kappa_{i+1}=\min(m+\nu_{i+1},\\_i)$. The amount $m'-\alpha_i$ added
to~$\alpha_i$ is equal to $m-\alpha'_{i+1}=m+\nu_{i+1}-\kappa_{i+1}$; we shall
use this ``carry''~$c$ in our algorithm rather than~$m$. To that end we
replace~$m$ in our description by the sum of the original entry
$\alpha_{i+1}=\mu_{i+1}-\\_{i+1}$ and the value~$c$ left by the previous step;
initially we take that expression with $c=M_{k,l}$. Then we get
$\kappa_{i+1}=\min(d,\\_i)$ where $d=\mu_{i+1}-\\_{i+1}+c+\nu_{i+1}$, and
$d-\kappa_{i+1}$ will replace~$c$. The index $i+1$ now occurs much more often
than~$i$, and we shall therefore shift the index one up, so that $i+1$
becomes~$i$. Finally, the last part of~$\kappa$ that can be nonzero is
$\kappa_{\min(k,l)}$, so we can start at $i=\min(k,l)$ instead of at $i=k$ and
gain in efficiency and symmetry. After these transformations we get the
following algorithm:
\bigdisplay\vcenter
{\settabs\+{\bf do }& {\bf do }&\cr
\+  $c:=M_{k,l}$ ; \cr
\+  {\bf for }$i$ {\bf from }$\min(k,l)$ {\bf down to }1 \cr
\+  & {\bf do }$d:=\mu_i+\nu_i-\\_i+c$
     ;~ $\kappa_i:=\min(d,\\_{i-1})$ ;~ $c:=d-\kappa_i$
     {\bf od} ;\cr
\+  $\kappa_0:=\mu_0-\\_0+c+\nu_0$.\strut \cr
}\label(\kappacompalg)
$$
If we execute the algorithm for the values of the example above $\\=(8,4,2)$,
\ $\mu=(8,7,2)$, \ $\nu=(8,5,3,1)$, $k=3$, \ $l=6$, and $M_{3,6}=2$ then we
find for $i=2$: $d:=3$, $\kappa_3=2$, $c:=1$; for $i=1$: $d:=4$, $\kappa_2=4$,
$c:=0$; for $i=0$: $d:=8$, $\kappa_1=8$, $c:=0$, and finally $\kappa_0=8$, so
that indeed one finds $\kappa=(8,8,4,2)$.

Note that the algorithm ensures that $\kappa_{i+1}\leq\\_i$ for all~$i$, which
means that $\\\leh\kappa$ will always hold; this is a stronger condition than
follows just from the known relations $\\\leh\mu\leh\kappa$ and
$\\\leh\nu\leh\kappa$. In fact the correspondence between the four partitions
and the entry $M_{k,l}$ can be characterised by the following properties:
(1)~there exist numbers $c_i\in\N$ for $0\leq{i}\leq\min(k,l)$ such that
putting $c_{-1}=0$ one has $c_i+\mu_i+\nu_i=\\_i+\kappa_i+c_{i-1}$ for all
those~$i$; (2)~$\kappa_{i+1}\leq\\_i$ for all~$i$, and $c_i>0$ implies
$\kappa_{i+1}=\\_i$; (3)~$c_{\min(k,l)}=M_{k,l}$. For $i>0$, $c_i$ will be the
value of~$c$ in the algorithm at the start of the loop indexed by~$i$, while
$c_0$ is the final value of~$c$. One may eliminate the dependence of this
statement on $(k,l)$ by having $i$ range over~$\N$ in property~(1) and
replacing property~(3) by $\lim_{i\to\infty}c_i=M_{k,l}$ (in this case the
sequence $c_i$ becomes stationary for $i\geq\min(k,l)$). This characterisation
is useful for deducing the following algorithm inverse to the one above, which
will compute $\\$ and $M_{k,l}$ from $\mu$, $\nu$, and~$\kappa$.
\bigdisplay\vcenter
{\settabs\+{\bf do }& {\bf do }&\cr
\+  $c:=0$ ; \cr
\+  {\bf for }$i$ {\bf from }0 {\bf to }$\min(k,l)$ \cr
\+  & {\bf do }$d:=\mu_i+\nu_i-\kappa_i-c$
     ;~ $\\_i:=\max(d,\kappa_{i+1})$ ;~ $c:=\\_i-d$
     {\bf od} ;\cr
\+  $M_{k,l}:=c$.\strut \cr
}\label(\lambdacompalg)
$$
Executing this with as above $\kappa=(13,9,9,5)$, $\mu=(11,9,8)$,
$\nu=(10,9,8,2)$, $k=3$ and $l=4$ gives for $i=0$: $d:=8$, $\\_0=9$, $c:=1$;
for $i=1$: $d:=8$, $\\_1=9$, $c:=1$; for $i=2$: $d:=6$, $\\_2=6$, $c:=0$; for
$i=3$: $d:=-3$, $\\_3=0$, $c:=3$, and finally $M_{3,4}=3$, which together with
$\\=(9,9,6)$ is the result found before.

It should come as no surprise that the shape datum we have just described is
in fact the shape datum for the Burge correspondence, which was given
in~\ref{spin Knuth, \Sec3.2}; after all we know from
proposition~\Burgecorporp\ that the implicit shape computed for~$M\in\Mat$ by
the growth diagram for our shape datum is the shape of the tableaux associated
to~$M$ under the Burge correspondence. Although none of the descriptions given
in that paper coincide exactly with the one given in~(\kappacompalg), it is
not hard to see their equivalence.

Those descriptions compute~$\set\kappa$ by making a bottom-left to top-right
pass over the squares that could be added to~$\set\\$ (namely the topmost
squares of the columns of the complement of that diagram), using a ``stock of
squares'' whose size may grow and shrink. Initially the stock contains
$M_{k,l}$ squares, and at each square~$s$ visited, a square is added to the
stock if $s\in\set{\mu/\\}$ and if $s\in\set{\nu/\\}$ (if both conditions
hold, two squares are added), and then a square from the stock is contributed
to~$\set\kappa$ unless the stock is empty (so if $s$ lies either in
$\set{\mu/\\}$ or $\set{\nu/\\}$ it will certainly be in~$\set\kappa$, and if
it lies in both there will in addition be a net increase of one square in the
stock; if it lies in neither and the stock was not empty then $s$ will still
be in~$\set\kappa$ while the size of the stock is decremented).
Algorithm~(\kappacompalg) can be seen to compute this same diagram
$\set\kappa$ a whole row at the time. The value
$d=\mu_i+\nu_i-\\_i+c=\\_i+c+(\mu-\\)_i+(\nu-\\)_i$ represents the maximal
length row~$i$ of $\set\kappa$ can have, obtained by adding to $\\_i$ the size
of the stock when arriving at this row and the sizes of rows~$i$ of
$\set{\mu/\\}$ and of~$\set{\nu/\\}$. Either that length does not exceed the
limit $\\_{i-1}$, in which case the one sets $\kappa_i=d$ and the stock
becomes empty, or $d$ does exceed $\\_{i-1}$, in which case $\kappa_i$ is set
to that limit and the size of the stack becomes $d-\\_{i-1}=d-\kappa_i$.

It is interesting to see how this description of the shape datum for the Burge
correspondence arises naturally from the rules governing the normalisation
process of integral matrices. Summarising, we have:

\proclaim Theorem. \integralgrowththm
If for any $M\in\Mat$ one associates to each grid point $(i,j)\in\N^2$ the
shape $\is{M_{\set{i}\times\set{j}}}$, then one obtains a growth diagram for
the shape datum of the Burge correspondence.
\QED

Having determined the shape datum involved, we can now determine the family of
shapes $\is{M_{\set{i}\times\set{j}}}$ directly from the matrix~$M$, without
actually performing any crystal operations: it suffices to use the fact that
$\is{M_{\set{i}\times\set{j}}}=(0)$ when $i=0$ or $j=0$, and to construct the
growth diagram for the shape datum. Diagram~1 illustrates this for the
integral matrix~$M$ of our running example.

\midinsert
$$\abovedisplayskip=0pt
\vcenter{
\catcode`@=11
\def~{\unhcopy\QEDbox}
\setbox\z@=\hbox{~}\dimen@=.5\ht\z@ \advance\dimen@.5\dp\z@
\let\f=\hfil
\def\d#1 {\vtop\bgroup\smallskip
          \vbox to\z@{\vss\llap{$#1$\kern6pt}\kern 4pt}%
          \futurelet\n\E}
\def\E{\?\ifx\n\f \then\smallskip\egroup \else\afterassignment\e\count@ \fi}
\def\e{\hbox{\loop~\advance\count@\m@ne \ifnum\count@>\z@ \repeat}%
       \nointerlineskip\futurelet\n\E}
\def\g#1{\omit#1\hfill\ignorespaces}
\def\h#1{\lower8pt\hbox{#1} \vrule height 9pt }
\def\z{\omit\vtop{\smallskip\hbox{$\circ$}}\hfil}
\offinterlineskip
\everycr{\noalign{\nobreak}}
\tabskip=\z@
\halign{\h# \tabskip=17 pt &&\d#\f\cr
\omit\hfil\strut\vrule
    & \g0 & \g1 & \g2 & \g3 & \g4 &  \g5  &  \g6  &  \g7 \cr
\noalign{\hrule}
0 &\z& \z &\z   &\z     & \z      & \z      &   \z      &   \z      \cr
1 &\z& 1 1&0 1  &1 2    &0 2      &1 3      &2 5        &0 5        \cr
2 &\z& 1 2&1 2 1&0 3 1  &1 3 2    &1 5 2    &0 7 2      &3 7 5      \cr
3 &\z& 0 2&2 3 2&1 4 2 1&0 4 3 1  &1 6 3 2  &1 8 4 2    &0 8 7 2    \cr
4 &\z& 0 2&0 3 2&1 4 2 2&1 4 3 2 1&1 6 4 3 1&0 8 5 3 1  &2 8 8 4 2  \cr
5 &\z& 0 2&0 3 2&0 4 2 2&0 4 3 2 1&0 6 4 3 1&1 8 5 3 1 1&2 8 8 5 3 1\cr
}
}
$$
\centerpar{{\bf Diagram~1.} Growth diagram with integral matrix and
shapes $(\is{M_{\set{i}\times\set{j}}})_{i,j\in\N}$.}
\endinsert

From this growth diagram one not only can read off $\is{M}$ in the bottom
right corner, but also the matrices $P,Q$ associated to~$M$ under the
bijection of~\decompthm\statitemnr2. To see that, observe that the operation
$M\mapsto{P}$ of exhausting upward transfers commutes with the operation
$M\mapsto{M_{\N\times\set{j}}}$ of discarding all columns beyond a certain
point; this is just the fact observed above that $\num_i(U|V)=0$ implies
$\num_i(U)=0$. Therefore the matrix $P_{\N\times\set{j}}$ can be obtained
from~$M_{\N\times\set{j}}$ by exhausting upward transfers, and one has
$\rsum{P_{\N\times\set{j}}} =\is{P_{\N\times\set{j}}}
=\is{M_{\N\times\set{j}}}$. Putting $\\^{(j)} =\rsum{P_{\N\times\set{j}}}$ for
all~$j$ defines the semistandard Young tableau $S =(\\^{(j)})_{j\in\N}$ of
which~$P$ is the integral encoding, and we see that $S$~coincides with the
sequence $(\is{M_{\N\times\set{j}}})_{j\in\N}$ that can be read off along the
lower border of the growth diagram. Similarly $Q\tr$ is the integral encoding
of the semistandard Young tableau $\bar{L}
=(\is{M_{\set{i}\times\N}})_{i\in\N}$ that can be read off along the right
hand border. Indeed in the example of growth diagram~$1$, one finds $S$
of~(\Sdef) and $\bar{L}$ of~(\Lbardef), whose integral encodings are the
matrices $P$ and~$Q\tr$ of (\PQintdiagram).

This argument confirms (and provides an alternative proof of)
proposition~\Burgecorporp, which can be interpreted as giving one particular
way in which the growth diagram can be computed: column insertion from right
to left of the letters of the weakly increasing word~$w$ whose weight is
given by the row~$M_i$ gives the transition from the semistandard Young
tableau $(\is{M_{\set{i}\times\set{j}}})_{j\in\N}$, formed by the shapes along
row~$i$ of the grid, to the tableau
$(\is{M_{\set{i+1}\times\set{j}}})_{j\in\N}$ along the next row of the grid.
Indeed one may check that the successive stages of the insertion process in
(\Burgeexamp) correspond to the shapes along the successive rows of the growth
diagram above. (It may be a bit confusing that whereas the letters of~$w$
are inserted from right to left, the growth diagram is extended from left to
right; however, completing the square~$(i,j)$ in the growth diagram does not
correspond to the insertion of the letters~$j$ of~$w$, but rather to
determining the effect on entries~$\leq{j}$ of the insertion of~$w$,
relative to the effect of that insertion on entries~$<j$.)

\proclaim Summary.
Given $M\in\Mat$, the shapes $\is{M_{\set{i}\times\set{j}}}$ can be determined
as follows for all $i,j\in\N$ (as there are only finitely many distinct
matrices $M_{\set{i}\times\set{j}}$, the computation is effectively finite).
Starting with $\is{M_{\set{i}\times\set{j}}}=(0)$ whenever $i=0$ or~$j=0$, one
completes a growth diagram by repeatedly computing shapes
$\kappa=\is{M_{\set{i+1}\times\set{j+1}}}$ from previously computed shapes
$\\=\is{M_{\set{k}\times\set{j}}}$, \ $\mu=\is{M_{\set{k}\times\set{j+1}}}$,
and $\nu=\is{M_{\set{k+1}\times\set{j}}}$, and the entry $M_{k,l}$, using the
algorithm~(\kappacompalg). From this growth diagram one can read off the
semistandard Young tableaux $S=(\is{M_{\N\times\set{j}}})_{j\in\N}$ and
$\bar{L}=(\is{M_{\set{i}\times\N}})_{i\in\N}$ associated to~$M$ under the
Burge correspondence, as well as the matrices $P,Q$ associated to~$M$ under
the bijection of theorem~\decompthm\statitemnr2: $P$ is the integral encoding
of~$S$ and $Q\tr$ is the integral encoding of~$\bar{L}$. One can reconstruct
$M$ from~$(P,Q)$ or from $(S,\bar{L})$ by reversing the computation, using
algorithm~(\lambdacompalg) to compute $\\$ and~$M_{k,l}$ given $\mu$, $\nu$,
and~$\kappa$.

\subsection Growth diagrams for binary matrices.

We shall now consider the family of shapes $\is{M_{\set{i}\times\set{j}}}$ for
binary matrices. With respect to the integral case there is a complication,
because one cannot in general reduce a submatrix $M_{\set{k}\times\set{l}}$ to
its normal form using raising operations acting on the entire matrix. For
instance after exhausting $\ltm_0$, $\ltm_1$, and $\upm_0$ for the binary
matrix $M$ of~(\PQbindiagram), its submatrix~$M_{\set2\times\set3}$ is reduced
to $1~1~0\choose1~0~1$ which is not a normal form; in isolation the bottom
right entry could be moved up by~$\upm_0$, but in the whole matrix it is
blocked by the parts of rows $0$~and~$1$ further to the right.

The difficulty comes from the reversal of left and right in the rule
$\num_i(U|V)\geq\num_i(U)$ that we used for integral matrices: for binary
matrices one instead has $\num_i(U|V)\geq\num_i(V)$ (for horizontal crystal
operations there is no such difference with the integral case: the relation
$\nlm_j({S\over{T}})\geq\nlm_j(S)$ still holds for binary matrices). So we
cannot exhaust upward moves in the submatrix~$U$; however we have the rule
$\ndm_i(U|V)\geq\ndm_i(U)$ that allows us to eventually exhaust any
\emph{downward} crystal operation~$\dnm_i$ on~$U$ (after a finite number of
applications of~$\dnm_i$ that move bits in~$V$). One cannot exhaust all
downward moves on a matrix, but if $k>0$ is such that $M=M_{\set{k}\times\N}$,
then it will suffice to exhaust the operations $\dnm_i$ for~$i\in\set{k-1}$.
If moreover leftward moves are also exhausted, then the matrix obtained
from~$M$ by a quarter turn clockwise, mapping the topmost~$k$ rows to the
leftmost $k$ columns, will satisfy the conditions of lemma~\normalformlemma;
it follows that $M$ is obtained by the inverse rotation from
$\diagram{\\\tr}$, where $\\\tr=\csum{M}$, and in particular that $M$ is
determined by its implicit shape $\\$ together with~$k$. It is in fact an
upside-down version of $\diagram\\$, in other words the region occupied by its
bits~`$1$' is the French style display\footnote* {One might have imagined that
the occurrence of Young diagrams inside matrices, as in
lemma~\normalformlemma, would give a
\foreign{coup de gr\^ace} to the French style of displaying diagrams.
\foreign{Mais non}.} of the diagram of~$\\$, in which the parts of~$\\$
correspond to the rows taken from bottom to top.

\proclaim Definition. \frdiagdef
Let $\\\in\Part$ and $k\in\N$ with $\\_k=0$. We denote by $\frdiag\\k$
the binary matrix $\(\Kr{(k-1-i,j)\in\set\\}\){}_{i,j\in\N}$, which is called
the French normal form (relative to~$k$) parametrised by~$\\$.

Note that $\csum{\frdiag\\k}=\\\tr$, while
$\rsum{\frdiag\\k}=(\li(\\k-1..0))$; we shall write this composition obtained
by reversing the $k$ initial parts of~$\\$ as $\revert\\k$. Now by exhausting
the crystal operations in $\setof\dnm_i:i\in\set{k-1}\endset
\union\setof\ltm_j:j\in\set{l-1}\endset$, the submatrix
$M_{\set{k}\times\set{l}}$ of any binary matrix can be reduced to the
form~$\frdiag\\k$, and this will be for $\\=\is{M_{\set{k}\times\set{l}}}$. We
call this process French normalisation of the submatrix
$M_{\set{k}\times\set{l}}$. For instance for the binary matrix $M$
of~(\PQbindiagram) and $(k,l)=(6,4)$, this exhaustion of crystal operations
gives:
\bigdisplay
 \pmatrix{0&1&0&0&1&0&0&0&0\cr
          1&1&1&0&0&0&0&0&0\cr
	  1&0&1&0&0&1&0&0&0\cr
	  0&1&0&1&0&0&0&0&0\cr
	  0&0&1&1&1&0&1&0&0\cr
	  1&0&0&0&1&0&0&1&1\cr
	  0&1&1&1&1&1&1&1&0\cr}
\buildrel (\dnm_\set5)^*
,(\ltm_\set3)*
\over\longrightarrow
 \pmatrix{0&0&0&0&0&0&0&0&0\cr
          0&0&0&0&1&0&0&0&0\cr
	  1&0&0&0&0&0&0&0&0\cr
	  1&1&0&0&0&1&0&0&0\cr
	  1&1&1&1&1&0&0&0&0\cr
	  1&1&1&1&1&0&1&1&1\cr
\matrixcross604
	  0&1&1&1&1&1&1&1&0\cr}
.\label(\nfeq)
$$
We see the submatrix $M_{\set6\times\set4}$ has been transformed into
$\frdiag{(4,4,2,1)}6$, whence we may conclude that
$\is{M_{\set6\times\set4}}=(4,4,2,1)$. As in the integral case, we shall now
compare the results of the French normalisation of submatrices
$M_{\set{i}\times\set{j}}$ for $(i,j)\in\sing{k,k+1}\times\sing{l,l+1}$. Let
$A$ be obtained from $M_{\set{k+1}\times\set{l+1}}$ by exhausting the
operations in $\setof\dnm_i:i\in\set{k-1}\endset
\union\setof\ltm_j:j\in\set{l-1}\endset$, and as before put
$\\=\is{M_{\set{k}\times\set{l}}}$, \ $\mu=\is{M_{\set{k}\times\set{l+1}}}$, \
$\nu=\is{M_{\set{k+1}\times\set{l}}}$,
and~$\kappa=\is{M_{\set{k+1}\times\set{l+1}}}$. Then
$A_{\set{k}\times\set{l}}=\frdiag\\k$, and since French normalisation of
$A_{\set{k}\times\set{l+1}}$ and of $A_{\set{k+1}\times\set{l}}$ only involves
horizontal respectively vertical moves, and it transforms these submatrices
into $\frdiag\mu{k}$ respectively into $\frdiag\nu{k+1}$, column~$l$ of
$A_{\set{k}\times\set{l+1}}$ is equal to the composition $\revert{\mu-\\}k$,
while row~$k$ of $A_{\set{k+1}\times\set{l}}$ is equal to~$\nu\tr-\\\tr$. In
particular one always has $\\\lev\mu$ and $\\\leh\nu$. In the example
displayed above, where $\\=(4,4,2,1)$ and the $7\times5$ matrix~$A$ can be
read off from the matrix on the right, one finds
$\mu=\\+(1,1,0,0,1,0)=(5,5,2,1,1)$ and $\nu=(\\\tr+(0,1,1,1))\tr=(4,4,4,2)$.

As in the case of integral matrices we see that $A$, and therefore $\kappa
=\is{A}$, is determined by the shapes $\\$, $\mu$, $\nu$, and the bit
$M_{k,l}$; moreover for any such data the construction provides an
example~$M=A$ of a matrix with the specified implicit shapes of submatrices
and bit~$M_{k,l}$, provided that one has $\\\lev\mu$, $\\\leh\nu$, and
$\mu_k=\nu\tr_l=0$. By the same procedure as for integral matrices one can
also reconstruct $\\$ and~$M_{k,l}$ from $\mu$, $\nu$, and~$\kappa$. Thus for
any $\mu,\nu\in\Part$ with $\mu_k=\nu\tr_l=0$ one has a bijective
correspondence between $(\\,M_{k,l})$ on one hand and $\kappa$ on the other,
where these values are restricted by the conditions $\\\lev\mu$, $\\\leh\nu$,
$M_{k,l}\in\set2$, $\nu\lev\kappa$ and~$\mu\leh\kappa$. Following the
terminology of \ref{spin Knuth, definition~6.1}, these bijective
correspondences for all~$\mu,\nu\in\Part$ together define an (asymmetric)
shape datum for $(\Part,\lev,\leh,\set2)$. By assigning
$\is{M_{\set{i}\times\set{j}}}$ to each grid point~$(i,j)$, one obtains a
growth diagram for this asymmetric shape datum. Note that the shapes
associated to grid points that are horizontal neighbours differ by a vertical
strip, while those associated to vertical neighbours differ by a horizontal
strip (we have no choice); in \ref{spin Knuth} an opposite convention was used
(its shape datum was for $(\Part,\leh,\lev,\set2)$, not for
$(\Part,\lev,\leh,\set2)$).

Before considering this shape datum more closely, let us compute an example.
The example of (\nfeq) for $(k,l)=(6,4)$ shows that during the French
normalisation of $A_{\set{k+1}\times\set{l}}$ all bits~`$1$' either remain in
place or move one place down, while during the French normalisation of
$A_{\set{k}\times\set{l+1}}$ some bits~`$1$' may move several places to the
left (the bit at $(1,4)$ moves all the way to~$(1,0)$); however this
example is too simple to illustrate well what can happen in general.
Therefore consider instead $\\=(5,3,2)$, $\mu=(6,3,3,1,1)$, $\nu=(6,3,2,2)$,
and $M_{k,l}=1$ with $(k,l)=(5,6)$. For the corresponding matrix~$A$ one has:
\bigdisplay
\def\mapright#1{\smash
       {\mathop{\hbox to 2.5em{\rightarrowfill}}\limits^{\textstyle#1}}}
\matrix
{\llap{$A={}$}
 \pmatrix{0&0&0&0&0&0&1\cr
          0&0&0&0&0&0&1\cr
	  1&1&0&0&0&0&1\cr
	  1&1&1&0&0&0&0\cr
	  1&1&1&1&1&0&1\cr
 \matrixcross506
	  1&1&0&0&0&1&1\cr}
&\mapright{(\ltm_\set6)^*}
&\pmatrix{1&0&0&0&0&0&0\cr
          1&0&0&0&0&0&0\cr
	  1&1&1&0&0&0&0\cr
	  1&1&1&0&0&0&0\cr
	  1&1&1&1&1&1&0\cr
 \matrixh
	  1&1&0&1&0&0&1\cr}
\cr\noalign{\nobreak\smallskip}
  \llap{$(\dnm_\set5)^*$}\bigg\downarrow &
& \llap{$(\dnm_\set5)^*$}\bigg\downarrow
\cr\noalign{\nobreak\smallskip}
 \pmatrix{0&0&0&0&0&0&1\cr
          0&0&0&0&0&0&1\cr
	  1&1&0&0&0&0&0\cr
	  1&1&0&0&0&0&1\cr
	  1&1&1&0&0&0&1\cr
 \matrixv506
	  1&1&1&1&1&1&1\cr}
&\mapright{(\ltm_\set6)^*}
&\pmatrix{1&0&0&0&0&0&0\cr
          1&0&0&0&0&0&0\cr
	  1&1&0&0&0&0&0\cr
	  1&1&1&0&0&0&0\cr
	  1&1&1&1&0&0&0\cr
	  1&1&1&1&1&1&1\cr}
\rlap{${}=\frdiag\kappa6$,}
\cr
}
\label(\Frenchgrowth)
$$
for $\kappa=(7,4,3,2,1,1)$. The parts of this computation that are not easy to
predict are the following ones. Firstly there is, for the exhaustion of
downward moves applied to~$A$, the transformation of column~$l$ of~$A$ (which
represents $\revert{\mu-\\}k$ followed by~$M_{k,l}$) into
$\revert{\kappa-\nu}{k+1}$; here some bits~`$1$' move down a row (like in the
example~$M_{2,6}$; this corresponds to a square of $\set{\mu/\\}$ for which
$\set{\kappa/\nu}$ contains a corresponding square in the \emph{same} row),
while some bits~`$1$' remain in the same row (in the example bits $M_{0,6}$,
$M_{1,6}$, $M_{4,6}$, and $M_{5,6}$; this corresponds to a square of
$\set{\mu/\\}$ (or $M_{k,l}=1$) for which $\set{\kappa/\nu}$ contains a
corresponding square in the \emph{next} row (respectively in row~$0$)).
Secondly there is, for the exhaustion of leftward moves applied to~$A$, the
transformation of row~$k$ of~$A$ (which represents $\nu\tr-\\\tr$ followed
by~$M_{k,l}$) into $\kappa\tr-\mu\tr$; here the bits~`$1$' may exchange with
zero or more bits~`$0$' initially situated directly to their left (in the
example bit~$M_{5,5}$ interchanges successively with two bits~`$0$' to its
left while $M_{5,0}$, $M_{5,1}$, and $M_{5,6}$ stay in place; such a bit
corresponds to a square of $\set{\nu/\\}$ (or $M_{k,l}=1$), and the
corresponding square of $\set{\kappa/\mu}$ lies in a column weakly to the left
of it (respectively in a column~$\leq{l}$)).

Giving details of either of these transformations will give a description of
how to compute~$\kappa$. However, although these descriptions follow by a
straightforward application of the definitions, they are rather technical, and
both are not stated in the same terms (the former is primarily by rows of the
shapes involved, the latter primarily by columns), which makes it hard to see
that they give the same result for~$\kappa$ (as we know they must). Therefore
we shall instead describe this shape datum using the approach
that was detailed in \ref{spin Knuth, \Sec6}; we have to adapt that
description to the interchange of `$\leh$' and~`$\lev$' in the relations
between shapes at horizontal and vertical neighbours in the grid.

Given $\mu,\nu\in\Part$ for which there exists at least one $\\\in\Part$ with
$\mu\gev\\\leh\nu$, or equivalently at least one $\kappa\in\Part$ with
$\mu\leh\kappa\gev\nu$, the possibilities for such $\\$ and~$\kappa$ have a
rather simple description: in each case the diagrams of these partitions are
the disjoint union of a fixed set of ``obligatory'' squares and an arbitrary
subset of a fixed set of ``optional'' squares. In case of~$\\$, the Young
diagram $\set\mu\thru\set\nu$ is the union of all obligatory and optional
squares (and therefore the maximal choice for~$\set\\$), and its subset~$S$ of
optional squares consists of those that are both at the end of a row
of~$\set\mu$ and at the end of a column of~$\set\nu$, in equation:
$S=\setof(i,j)\in\N^2:(i,j)=(\nu\tr_j-1,\mu_i-1)\endset$. In the case
of~$\kappa$, the set of obligatory squares is the Young diagram
$\set\mu\union\set\nu$, and the set~$T$ of optional squares consists of those
that are both just beyond the end of a column of~$\set\mu$ and just beyond the
end of a row of~$\set\nu$, in other words
$T=\setof(i,j)\in\N^2:(i,j)=(\mu\tr_j,\nu_i)\endset$. It was shown
in~\ref{spin Knuth} that $\Card{T}=\Card{S}+1$, which is in agreement with the
fact that there are twice as many possibilities for~$\kappa$ as for~$\\$;
moreover, there is an injective map matching every $s\in{S}$ with the first
$t\in{T}$ in a row strictly below~$s$ (or equivalently in a column weakly to
the left of~$s$), which leaves the topmost square $t_0$ of~$T$ unmatched. Now
the correspondence between $(\\,M_{k,l})$ and $\kappa$ is determined by the
fact that for every optional square $s\in{S}$ \emph{absent} from~$\set\\$ the
corresponding optional square $t\in{T}$ is \emph{present} in~$\set\kappa$, and
vice versa, while one has $t_0\in\set\kappa$ if and only if $M_{k,l}=1$. We
shall call this shape datum the ``row-insertion shape datum for
$(\Part,\lev,\leh,\set2)$''.

\proclaim Theorem. \binarygrowththm
If for any $M\in\bMat$ one associates to each grid point $(k,l)\in\N^2$ the
shape $\is{M_{\set{k}\times\set{l}}}$, then one obtains a growth diagram for
the row-insertion shape datum for $(\Part,\lev,\leh,\set2)$.

We omit the proof; the assiduous reader may check that the rules governing
French normalisation of submatrices do indeed lead to this shape datum.
Whether one starts normalising $A_{\set{k+1}\times\set{l}}$ or
$A_{\set{k}\times\set{l+1}}$, the description obtained will not directly give
a correspondence between $S$ and~$T$, but rather one between the sets of rows
meeting $\set{\mu/\\}$ and $\set{\kappa/\nu}$, respectively between the sets
of columns meeting $\set{\nu/\\}$ and $\set{\kappa/\mu}$. On one hand this
means that a large number of steps may intervene between encountering the bit
for some $s\in{S}$ and the one for the corresponding $t\in{T}$, and on the
other hand the occurrence of such squares may be hidden by the fact that a
square in~$S$ shares a row or a column with a different square in~$T$. As a
consequence, cases need to be distinguished that play no role in the above
description.

Now that the shape datum is known, it can be used to build up the growth
diagram without actually performing (French) normalisation submatrices; again
one starts by associating the empty partition~$(0)$ to every grid
point~$(i,j)$ with $i=0$ or $j=0$. The construction for the binary matrix~$M$
of~(\PQbindiagram) is given in diagram~2, with the shapes associated to grid
points given by their diagram \foreign{\`a la fran\c{c}aise}.
\midinsert
$$\abovedisplayskip=0pt
\vcenter{
\catcode`@=11
\def~{\unhcopy\QEDbox}
\setbox\z@=\hbox{~}\dimen@=.5\ht\z@ \advance\dimen@.5\dp\z@
\let\f=\hfil
\def\french#1 #2.{\?\if|#2|\then\futurelet\n\E#1 \else\french#2.#1 \fi}
\def\d#1 {\vtop\bgroup\smallskip
          \vbox to\z@{\vss\llap{$#1$\kern6pt}\kern 4pt}%
          \french}
\def\E{\?\ifx\n\f \then\smallskip\egroup \else\afterassignment\e\count@ \fi}
\def\e{\ifnum\count@=0 \phantom~
       \else\hbox{\loop~\advance\count@\m@ne \ifnum\count@>\z@ \repeat}\fi
       \nointerlineskip\futurelet\n\E}
\def\g#1{\omit#1\hfill\ignorespaces}
\def\h#1{\lower8pt\hbox{#1} \vrule height 9pt }
\def\zz{\vtop{\smallskip\hbox{$\circ$}}\hfil}
\def\z{\omit\zz}
\offinterlineskip
\everycr{\noalign{\nobreak}}
\tabskip=\z@
\halign{\h# \tabskip=15 pt &&\d#.\f\cr
\omit\hfil\strut\vrule
    & \g0 & \g1 & \g2 & \g3 & \g4 & \g5 & \g6 & \g7 & \g8 & \g9 \cr
\noalign{\hrule}
0 &\z & \z &\z &\z & \z & \z & \z & \z & \z & \z \cr
1 &\z &\omit \d0 {\f} .\zz
                &1 1       &0 1       &0 1       &1 2         &0 2
        &0 2         &0 2         &0 2         \cr
2 &\z & 1 1 0   &1 2 1     &1 3 1     &0 3 1     &0 3 2       &0 3 2
        &0 3 2       &0 3 2       &0 3 2       \cr
3 &\z & 1 1 1 0 &0 2 1 1   &1 3 2 1   &0 3 2 1   &0 3 2 2     &1 4 2 2
      &0 4 2 2  &0 4 2 2   &0 4 2 2   \cr
4 &\z&0 1 1 0 0 &1 2 2 1 0 &0 3 2 2 0 &1 4 2 2 0 &0 4 2 2 1   &0 4 3 2 1
     &0 4 3 2 1 &0 4 3 2 1 &0 4 3 2 1 \cr
5 &\z&0 1 1 0 0 &0 2 2 1 0 &1 3 3 2 0 &1 4 4 2 0 &1 5 4 2 1   &0 5 4 3 1
     &1 6 4 3 1 &0 6 4 3 1 &0 6 4 3 1   \cr
6 &\z &1 1 1 1 0 0 &0 2 2 1 1 0 &0 3 3 2 1 0 &0 4 4 2 1 0 &1 5 5 2 1 1
      &0 5 5 3 1 1 &0 6 5 3 1 1 &1 7 5 3 1 1 &1 8 5 3 1 1 \cr
7 &\z &0 1 1 1 0 0 &1 2 2 2 1 0 &1 3 3 3 2 0 &1 4 4 4 2 0 &1 5 5 5 2 1
      &1 6 5 5 3 1 &1 7 6 5 3 1 &1 8 7 5 3 1 &0 8 8 5 3 1 \cr
}
}
$$
\centerpar{{\bf Diagram~2.} Growth diagram with binary matrix and
shapes $(\is{M_{\set{i}\times\set{j}}})_{i,j\in\N}$.}
\endinsert

The shapes along the right border define the semistandard tableau~$S$
of~(\Sdef) whose binary encoding is the matrix~$Q$ of~(\PQbindiagram), for
essentially the same reasons as in the integral case: one has
$\csum{Q_{\set{i}\times\N}} =\is{Q_{\set{i}\times\N}}\tr
=\is{M_{\set{i}\times\N}}\tr$ for all~$i\in\N$. There is however no such
direct relationship between matrix~$P$ of~(\PQbindiagram) and the shapes along
the lower border; this is not surprising since those shapes define a transpose
semistandard tableau, whereas $P$ corresponds to the \emph{reverse} transpose
semistandard tableau $R$ of~(\Rdef). Indeed, although like in the integral
case one still has $\is{P_{\N\times\set{j}}}=\is{M_{\N\times\set{j}}}$ for
all~$j$, these shapes are not given by $\rsum{P_{\N\times\set{j}}}$, since
upward moves need not be exhausted in~$P_{\N\times\set{j}}$. To read off~$R$
from a growth diagram, one should instead normalise the complementary
submatrices $M_{\N\times(\N-\set{j})}$, as we shall do in the next subsection.
On the other hand, our growth diagram can be interpreted as recording the
insertion process of proposition~\dualRSKrowprop, which gives as insertion
tableau the Sch\"utzenberger dual $R^*$ of~$R$:
\bigdisplay
\smallsquares
\diaght=\ht\strutbox \advance\diaght by \dp\strutbox
\multiply \diaght by 5 \advance\diaght by -\dp\strutbox
\def\tableau(#1){\vcenter to\diaght{\Young(#1)\vss}}
\def\ins(#1){\buildrel #1\over\longrightarrow}
\ins(1,4)\tableau(1,4)
\ins(0,1,2)\tableau(0,1,2|1,4)
\ins(0,2,5)\tableau(0,1,2,5|0,2|1,4)
\ins(1,3)\tableau(0,1,2,3|0,1,5|1,2|4)
\ins(2,3,4,6)\tableau(0,1,2,3,4,6|0,1,2,3|1,2,5|4)
\ins(0,4,7,8)\tableau(0,1,2,3,4,6,7,8|0,1,2,3|1,2,5|4)
\ins(1,2,3,4,5,6,7)\tableau(0,1,2,3,4,5,6,7|0,1,2,3,4,6,7,8|0,1,2,3,4|1,2,5|4)
.
$$
Indeed, since the shapes in our diagram are governed by the row-insertion
shape datum for $(\Part,\lev,\leh,\set2)$, the passage from the transpose
semistandard Young tableau defined by the shapes along row~$i$ to those along
row~$i+1$ is given by the row insertion (according to the rules for transpose
tableaux) of the letters of the strictly increasing word of weight~$M_i$, from
left to right. This proves:

\proclaim Proposition. \bingdprop
For $M\in\bMat$, the transpose semistandard
tableau~$R^*=(\is{M_{\N\times\set{j}}})_{j\in\N}$ and the semistandard
tableau~$S=(\is{M_{\set{i}\times\N}})_{i\in\N}$ form the pair associated
to~$M$ in~\ref{Stanley EC2, theorem~7.14.2}.
\QED

From~$M$ one can obtain a matrix~$\widetilde{P}$ whose rows give the weights
of those of~$R^*$, but in reverse order, by exhausting \emph{downward} crystal
operations, up to some limit. Exhausting $\setof\dnm_i:i\in\set6\endset$ in
our example gives
$$
  \widetilde{P}=
\pmatrix
{0&0&0&0&0&0&0&0&0\cr
 0&0&0&0&0&0&0&0&0\cr
 0&0&0&0&1&0&0&0&0\cr
 0&1&1&0&0&1&0&0&0\cr
 1&1&1&1&1&0&0&0&0\cr
 1&1&1&1&1&0&1&1&1\cr
 1&1&1&1&1&1&1&1&0\cr
}
\ttext{which matches}
  R^*=\Young(0,1,2,3,4,5,6,7|0,1,2,3,4,6,7,8|0,1,2,3,4|1,2,5|4)
.\label(\Rstardef)
$$

\subsection Other types of growth diagrams.

As was suggested above, one can form a growth diagram from which both tableaux
$R$ and~$S$ encoded by the matrices $P$ and~$Q$ corresponding to $M\in\bMat$
under the bijection of theorem~\decompthm\statitemnr1 can be directly read
off, by computing the family of shapes
$(\is{M_{\set{i}\times(\N-\set{j})}})_{i,j\in\N}$ rather than
$(\is{M_{\set{i}\times\set{j}}})_{i,j\in\N}$. The shapes defining~$R$, a
weakly decreasing sequence in which adjacent shapes differ by vertical strips,
can then be found along the bottom border, while the shapes defining~$S$ are
now found along the left border of the growth diagram. (In the direction from
left to right one might prefer calling this a shrinkage diagram rather than a
growth diagram; we shall call it a horizontally reversed growth diagram.) That
$S$ will be found is obvious, since we are considering the implicit shapes of
the same sequence of matrices; that fact that $R$ can be read off as
indicated, follows from the fact that upward moves are exhausted in the
submatrices $P_{\N\times(\N-\set{j})}$ for all~$j$, whence one has
$\rsum{P_{\N\times(\N-\set{j})}} =\is{P_{\N\times(\N-\set{j})}}
=\is{M_{\N\times(\N-\set{j})}}$. We illustrate this in diagram~3 for the
binary matrix~$M$ of~(\PQbindiagram).
\midinsert
$$\abovedisplayskip=0pt
\vcenter{
\catcode`@=11
\def~{\unhcopy\QEDbox}
\setbox\z@=\hbox{~}\dimen@=.5\ht\z@ \advance\dimen@.5\dp\z@
\let\f=\hfil
\def\d#1 {\vtop\bgroup\smallskip
          \vbox to\z@{\vss\llap{$#1$\kern6pt}\kern 4pt}%
          \futurelet\n\E}
\def\E{\?\ifx\n\f \then\smallskip\egroup \else\afterassignment\e\count@ \fi}
\def\e{\hbox{\loop~\advance\count@\m@ne \ifnum\count@>\z@ \repeat}%
       \nointerlineskip\futurelet\n\E}
\def\g#1{\omit#1\hfill\ignorespaces}
\def\h#1{\lower8pt\hbox{#1} \vrule height 9pt }
\def\zz{\vtop{\smallskip\hbox{$\circ$}}\hfill}
\def\z{\omit\zz}
\offinterlineskip
\everycr{\noalign{\nobreak}}
\tabskip=\z@
\halign{\h# \tabskip=15 pt &\vtop\bgroup\smallskip\futurelet\n\E #\f &&\d#\f\cr
\omit\hfil\strut\vrule
    & \g0 & \g1 & \g2 & \g3 & \g4 & \g5 & \g6 & \g7 & \g8 & \g9 \cr
\noalign{\hrule}
0 &\z & \z &\z &\z & \z & \z & \z & \z & \z & \z \cr
1 & 2         & 0 2          &1 1         &0 1       &0 1       &1 \f\zz
    &0 \f\zz &0 \f\zz &0 \f\zz &0 \f\zz \cr
2 & 3 2       & 1 2 2        &1 1 1       &1 1       &0 1       &0 \f\zz
    &0 \f\zz &0 \f\zz &0 \f\zz &0 \f\zz \cr
3 & 4 2 2     & 1 3 2 1      &0 2 1 1     &1 2       &0 2       &0 1
    &1 \f\zz &0 \f\zz &0 \f\zz &0 \f\zz \cr
4 & 4 3 2 1   & 0 3 3 1 1    &1 2 2 1     &0 2 1     &1 2       &0 1
    &0 \f\zz &0 \f\zz &0 \f\zz &0 \f\zz \cr
5 & 6 4 3 1   & 0 5 3 3 1    &0 4 2 2 1   &1 3 2 1   &1 3 1     &1 2
    &0 1     &1 \f\zz &0 \f\zz &0 \f\zz \cr
6 & 8 5 3 1 1 & 1 7 4 3 1    &0 6 3 2 1   &0 5 2 2   &0 5 1 1   &1 4
    &0 3     &0 2     &1 1     &1 \f\zz \cr
7 & 8 8 5 3 1 & 0 7 7 4 3 1  &1 6 6 3 2 1 &1 5 5 2 2 &1 5 4 1 1 &1 4 3
    &1 3 2   &1 2 1   &1 1     &0 \f\zz \cr
}
}
$$
\centerpar{{\bf Diagram~3.} Horizontally reversed growth diagram with binary
matrix and shapes $(\is{M_{\set{i}\times(\N-\set{j})}})_{i,j\in\N}$.}
\endinsert

Calling this a horizontally reversed growth diagram means that for any pair of
indices $k,l\in\N$, if we put $\\=\is{M_{\set{k}\times(\N-\set{l+1})}}$, \
$\mu=\is{M_{\set{k}\times(\N-\set{l})}}$, \
$\nu=\is{M_{\set{k+1}\times(\N-\set{l+1})}}$, and
$\kappa=\is{M_{\set{k+1}\times(\N-\set{l})}}$, there should be a relation
between these shapes and the matrix entry~$M_{k,l}$ that determines for fixed
$\mu,\nu$ a bijection between the pair $(\\,M_{k,l})$ and the shape~$\kappa$.
Indeed, if we determine the implicit shapes of submatrices by exhausting the
pertinent raising operations (as in computing ordinary rather than French
normal forms), then we have a situation in which everything is rotated a
quarter turn clockwise with respect to that of the previous subsection. Taking
our usual binary matrix~$M$ of~(\PQbindiagram) and $(k,l)=(5,2)$, we can read
off the shapes $\\=\is{M_{\set5\times(\N-\set3)}}=(3,2,1)$, \
$\mu=\is{M_{\set5\times(\N-\set2)}}=(4,2,2,1)$, and
$\nu=\is{M_{\set6\times(\N-\set3)}}=(5,2,2)$ from diagram~3; from this and
$M_{5,2}=0$ we can compute~$\kappa=(6,3,2,1)$ as follows (the two leftmost
columns that are of no interest here are represented by dots):
\bigdisplay
\def\mapleft#1{\smash
       {\mathop{\hbox to 2.5em{\leftarrowfill}}\limits^{\textstyle#1}}}
\matrix
{\pmatrix{\cdot&\cdot&1&1&1&1&0&0&0\cr
          \cdot&\cdot&1&1&0&0&0&0&0\cr
	  \cdot&\cdot&1&1&0&0&0&0&0\cr
	  \cdot&\cdot&1&0&0&0&0&0&0\cr
	  \cdot&\cdot&0&0&0&0&0&0&0\cr
\matrixh 
	  \cdot&\cdot&0&0&1&0&1&1&0\cr}
&\mapleft{(\ltm_{\N-\set2})^*}
&\pmatrix{\cdot&\cdot&1&1&1&1&0&0&0\cr
          \cdot&\cdot&0&1&1&0&0&0&0\cr
	  \cdot&\cdot&1&1&0&0&0&0&0\cr
	  \cdot&\cdot&1&0&0&0&0&0&0\cr
	  \cdot&\cdot&0&0&0&0&0&0&0\cr
 \matrixcross503
	  \cdot&\cdot&0&0&1&0&1&1&0\cr}
\cr\noalign{\nobreak\smallskip}
  \llap{$(\upm_\set5)^*$}\bigg\downarrow &
& \llap{$(\upm_\set5)^*$}\bigg\downarrow
\cr\noalign{\nobreak\smallskip}
 \pmatrix{\cdot&\cdot&1&1&1&1&1&1&0\cr
          \cdot&\cdot&1&1&1&0&0&0&0\cr
	  \cdot&\cdot&1&1&0&0&0&0&0\cr
	  \cdot&\cdot&1&0&0&0&0&0&0\cr
	  \cdot&\cdot&0&0&0&0&0&0&0\cr
 	  \cdot&\cdot&0&0&0&0&0&0&0\cr}
&\mapleft{(\ltm_{\N-\set2})^*}
&\pmatrix{\cdot&\cdot&1&1&1&1&1&1&0\cr
          \cdot&\cdot&1&1&1&0&0&0&0\cr
	  \cdot&\cdot&0&1&1&0&0&0&0\cr
	  \cdot&\cdot&1&0&0&0&0&0&0\cr
	  \cdot&\cdot&0&0&0&0&0&0&0\cr
 \matrixv503
	  \cdot&\cdot&0&0&0&0&0&0&0\cr}
.
\cr
}
$$
The rotation of the matrices induces the transposition of their implicit
shapes and those of their submatrices; therefore the relation here between
$\\\tr$, $\nu\tr$, $\mu\tr$, and $\kappa\tr$ (in that order) will be the same
as was found between $\\$, $\mu$, $\nu$, and~$\kappa$ in the situation of the
previous subsection. It follows that one still has $\\\lev\mu\leh\kappa$ and
$\\\leh\nu\lev\kappa$ (which is also obvious from the computation considered),
but the shape datum that determines~$\kappa$ will be the
\emph{column}-insertion shape datum for $(\Part,\lev,\leh,\set2)$; its
definition differs from the row-insertion shape datum by the fact that another
injective map $S\to{T}$ is used, one that maps each $s\in{S}$ to the first
$t\in{T}$ in a column strictly to its right (and in a row weakly above it).

\proclaim Proposition. \revbingdprop
The bijection of theorem~\decompthm\statitemnr1 between $M\in\bMat$ and the
pair $(P,Q)$, or the pair of tableaux $(R,S)$ encoded by them as in
proposition~\Knuthbinaryprop, can be computed using a horizontally reversed
growth diagram for the column-insertion shape datum for
$(\Part,\lev,\leh,\set2)$, reading off
$R=(\is{M_{\N\times(\N-\set{j}}})_{j\in\N}$ along the lower border and
$S=(\is{M_{\set{i}\times\N}}_{i\in\N}$ along the left border.
\QED

We may interpret proposition~\Knuthbinaryprop\ as a way to determine such
horizontally reversed growth diagrams: the transition from the semistandard
Young tableau formed by the shapes along column~$j+1$ of the grid and the one
formed by those along column~$j$ is given by column insertion in strictly
increasing order of the letters of a word~$w$ whose weight is given by the
column~$M\tr_j$. This may be verified by comparing the intermediate tableaux
in~(\Knuthbinex) with the vertical sequences of shapes in the display above.
Computing the shape~$\kappa$ at the bottom left corner of a square~$(i,j)$
amounts to finding the shape occupied by letters~$\leq{i}$ after
column-insertion of~$w$, given the shape~$\mu$ occupied by letters~$<i$ after
that insertion, the horizontal strip $\nu/\\$ filled with letters~$i$ before
the insertion, and the number $M_{i,j}\in\{0,1\}$ of letters~$i$ in~$w$.

We have seen how $\is{M}$ can be computed using either an ordinary or a
horizontally reversed growth diagram, which were determined by studying the
exhaustion of appropriate crystal operations on submatrices from the top left
respectively from the top right, and the two situations are the same up to a
rotation of everything by a quarter turn. It should by now be clear that two
more variations to this construction can be obtained by further rotations.
Thus we obtain computations using a vertically and a doubly reversed growth
diagram, which are respectively governed by the column-insertion an by the
row-insertion shape datum for $(\Part,\lev,\leh,\set2)$. The relations this
gives between the results of applying the corresponding variations of the dual
\caps{RSK}~algorithm to~$M$ are both straightforward and in principle known
(see for instance \ref{Fulton}), so we shall not give more details here.
However, for integral matrices similarly defined variant growth diagrams are
not all related by simple symmetries, as we shall now see.

For integral matrices $M\in\Mat$ we have seen that any upper left submatrix
$M_{\set{i}\times\set{j}}$ can be reduced inside~$M$ to its normal form by
exhausting the applicable upward and leftward crystal operations. If one views
$M$ as a finite rectangular matrix, then one can similarly bring any submatrix
containing one of the three remaining corners into a standard form by
exhausting an appropriate set of crystal operations, but the result is not the
same in all cases. In the case of a lower right submatrix the situation is
rotated a half turn with respect to that of an upper left submatrix, so
exhaustion of all downward and rightward crystal operations between pairs of
rows respectively columns of the submatrix will end with a rotated version
of~$\diagon\\$: all units are on the diagonal passing through the bottom right
corner, their numbers weakly increasing towards that corner. This rotated
version has the same implicit shape~$\\$ as $\diagon\\$, as can be verified by
any of the many ways we have seen to compute implicit shapes. Therefore the
family of shapes $(M_{(\N-\set{i})\times(\N-\set{j})})_{i,j\in\N}$ defines a
doubly reversed growth diagram in, which for all~$k,l\in\N$ the four shapes
$\\=\is{M_{(\N-\set{k+1})\times(\N-\set{l+1})}}$, \
$\mu=\is{M_{(\N-\set{k+1})\times(\N-\set{l})}}$, \
$\nu=\is{M_{(\N-\set{k})\times(\N-\set{l+1})}}$, and
$\kappa=\is{M_{(\N-\set{k})\times(\N-\set{l})}}$ at the corners of a grid
square and the entry~$M_{k,l}$ of that square are related by the shape datum
of the Burge correspondence, that can be computed by algorithms
(\kappacompalg) and~(\lambdacompalg).

For submatrices containing the top right or bottom left corner of~$M$ however,
the situation is not symmetric to those of the top left and bottom right
corners; on the other hand these two remaining cases are related by a rotation
by a half turn, so we shall focus (arbitrarily) on the bottom left
submatrices. A normalised form for the submatrix
$M_{(\N-\set{k})\times\set{l}}$ that can be obtained by crystal operations
acting on the entire matrix, is the one where upward and rightward crystal
operations between pairs of rows respectively columns of the submatrix are
exhausted, in other words the operations in
$\setof\upm_i:i\in\N-\set{k}\endset\union\setof\rtm_j:j\in\set{l-1}\endset$.
The result is not however that all units of the submatrix end up on one same
diagonal or anti-diagonal, as we illustrate for our usual integral matrix~$M$
and its submatrix $(k,l)=(1,6)$:
\bigdisplay
\pmatrix{1&0&1&0&1&2&0\cr
         1&1&0&1&1&0&3\cr
	 0&2&1&0&1&1&0\cr
	 0&0&1&1&1&0&2\cr
	 0&0&0&0&0&1&2\cr}
\buildrel (\upm_{\N-\set1})^*
,(\rtm_\set5)* \over
\longrightarrow
\pmatrix{0&1&0&1&0&3&0\cr
\matrixcross136
         0&0&1&1&3&0&3\cr
	 0&0&0&1&1&3&0\cr
	 0&0&0&0&1&1&3\cr
	 0&0&0&0&0&1&1\cr}
.\label(\slicetransform)
$$
One sees that entries $M_{i,j}$ with $(i,j)\in(\N-\set{k})\times\set{l}$ have
become constant along diagonals, and the implicit shape of
$M_{(\N-\set{k})\times\set{l}}$ can be read off from the row sums of the
submatrix after the transformation: in~(\slicetransform) on has
$\is{M_{(\N-\set1)\times\set6}}=(5,5,2,1)$. This must be so, for if one would
isolate $M_{(\N-\set{k})\times\set{l}}$ and exhaust leftward crystal
operations on it, then the result would on one hand be a shifted version of
$\diagon\\$ for $\\=\is{M_{(\N-\set{k})\times\set{l}}}$, and on the other hand
it would have unchanged row sums. Whereas the~$\\_i$~units forming such a row
sum are all grouped together in one entry of~$\diagon\\$, they are distributed
in the corresponding row of $M_{(\N-\set{k})\times\set{l}}$ according to the
partitioning of the squares of row~$i$ of~$\set\\$ by the lengths of their
columns: the unit corresponding to a square lies on the diagonal near the top
right of $(\N-\set{k})\times\set{l}$ that has the same length as the column of
the square. This is easily checked, since the entry $M_{k+i,l-1}$ lies at the
end of the row of $M_{(\N-\set{k})\times\set{l}}$ with sum~$\\_i$ and also of
a diagonal of length~$i+1$, and its value is $M_{k+i,l-1}=\\_i-\\_{i+1}$,
which is the number of columns of length~$i+1$ of the diagram~$\set\\$. Thus
the submatrix shows how that diagram can be sliced vertically into parts for
each column length; for instance in the example the diagram $\set{(5,5,2,1)}$
appears as cut up into two single columns, of length $4$~and~$3$, and a group
of the remaining $3$~columns of length~$2$, while there are $0$~columns of
length~$1$. We now define more formally the ``sliced form'' matrix that thus
corresponds to~$\\$, for given $k,l$.

\proclaim Definition.
Let $\\\in\Part$ and $k,l\in\N$ with $\\_l=0$. We denote by
$\sliced\\kl$ the integral matrix~$M$ defined by
$M_{i,j}=\Kr{i\geq{k}}\Kr{j<l}(\\_{m-1}-\\_m)$ for $i,j\in\N$, where
$m=i-k-j+l$.

\proclaim Proposition.
Let $k,l\in\N$ with $l>0$ be fixed; then $M\in\Mat$ satisfies
$M=M_{(\N-\set{k})\times\set{l}}$, \ $\num_i(M)=0$ for $i\geq{k}$, and
$\nrm_j(M)=0$ for $j\in\set{l-1}$, if and only if there exists $\\\in\Part$
with $\\_l=0$ such that $M=\sliced\\kl$.

\proof
It is clear that any $M=\sliced\\kl$ with $\\_l=0$ satisfies the given
conditions: for that case all the contributions to the sums in (\intnumeq)
and~(\intnrmeq) are~$\leq0$. Conversely suppose that $M$ satisfies those
conditions. If its entries were not constant along all diagonals within
$(\N-\set{k})\times\set{l}$, then on the lowest diagonal with non-constant
entries one could find a pair $M_{i,j}\neq{M}_{i+1,j+1}$ with $k\leq{i}$ and
$j+1<l$; if $M_{i,j}<M_{i+1,j+1}$ one would have $\num_i(M)>0$ while if
$M_{i,j}>M_{i+1,j+1}$ one would have $\nrm_j(M)>0$, either of which
contradicts the assumption. Therefore the entries of~$M$ are constant along
such diagonals; moreover $\num_i(M)=0$ implies $M_{i+1,0}=0$ for all
$i\geq{k}$, so the entries are zero on all diagonals below the one containing
$M_{k,0}$. Then putting $\\_i=\sum_{j\in\set{l}}M_{k+i,j}$ one checks that
$\\\in\Part$, \ $\\_l=0$, and $M=\sliced\\kl$.
\QED

Now consider an integral matrix $A=A_{(\N-\set{k})\times\set{l+1}}$ that
satisfies $\num_i(A)=\nrm_j(A)=0$ for all $i\geq{k}+1$ and $j<l-1$, so that in
particular its submatrix $A_{(\N-\set{k+1})\times\set{l}}$ is of the form just
described. Let us try to find the relation between the entry~$A_{k,l}$ just
above and to the right of that submatrix, and the shapes at the grid points
forming the corners of its square: $\\=\is{A_{(\N-\set{k+1})\times\set{l}}}$,
\ $\mu=\is{A_{(\N-\set{k+1})\times\set{l+1}}}$, \
$\nu=\is{A_{(\N-\set{k})\times\set{l}}}$, and
$\kappa=\is{A_{(\N-\set{k})\times\set{l+1}}}$. The matrix after the
transformation in~(\slicetransform) provides an example with $(k,l)=(0,6)$,
$M_{k,l}=0$, and $\\=(5,5,2,1)$; the computation below shows that in this case
$\mu=(8,5,5,2)$, $\nu=(8,5,3,1,1)$, and $\kappa=(8,8,5,3,1)$.
\bigdisplay
\def\mapright#1{\smash
       {\mathop{\hbox to 2.5em{\rightarrowfill}}\limits^{\textstyle#1}}}
\matrix
{\pmatrix{0&1&0&2&2&3&0\cr
\matrixv136
          0&0&1&0&2&2&3\cr
	  0&0&0&1&0&2&2\cr
	  0&0&0&0&1&0&2\cr
	  0&0&0&0&0&1&0\cr}
&\mapright{(\rtm_\set6)^*}
&\pmatrix{0&0&1&2&2&3&0\cr
          0&0&0&1&2&2&3\cr
	  0&0&0&0&1&2&2\cr
	  0&0&0&0&0&1&2\cr
	  0&0&0&0&0&0&1\cr}
 \rlap{${}=\sliced\kappa07$}
\cr\noalign{\nobreak\smallskip}
  \llap{$(\upm_{\N-\set0})^*$}\bigg\uparrow &
& \llap{$(\upm_{\N-\set0})^*$}\bigg\uparrow
\cr\noalign{\nobreak\smallskip}
 \llap{$A={}$}
\pmatrix{0&1&0&1&0&3&0\cr
\matrixcross136
         0&0&1&1&3&0&3\cr
	 0&0&0&1&1&3&0\cr
	 0&0&0&0&1&1&3\cr
	 0&0&0&0&0&1&1\cr}
&\mapright{(\rtm_\set6)^*}
&\pmatrix{0&0&1&1&0&3&0\cr
\matrixh
         0&0&0&2&3&0&3\cr
	 0&0&0&0&2&3&0\cr
	 0&0&0&0&0&2&3\cr
	 0&0&0&0&0&0&2\cr}
\cr
}\nn
$$

We shall now analyse in detail what happens during this procedure. The
description will be heavy going at times, but it eventually leads to the shape
datum for the \caps{RSK}~correspondence; as this forms the most direct link of
our constructions to that correspondence, we insist on presenting it, with its
inevitable technical details. The reader may choose to either bear with us, or
take our word for it and skip to~(\forward\Knuthgrowtheq).

The parts of the transformations above that respectively reduce
$A_{(\N-\set{k})\times\set{l}}$ to $\sliced\nu{k}l$ by upward transfers (the
part to the left of the vertical line for the left upward arrow), and
$A_{(\N-\set{k+1})\times\set{l+1}}$ to $\sliced\mu{k+1}{l+1}$ by rightward
transfers (below the horizontal line for the bottom arrow), are less
transparent that the corresponding parts in~(\Burgeshapedatumex), but we can
describe the former as follows (the latter is quite similar).
From~$\nrm_j(A)=0$ it follows that $A_{k,j}\leq{A}_{k+1,j+1}$
for~$j\in\set{l-1}$. Therefore, if one starts by exhausting $\num_k$ applied
to~$A$, the entry~$A_{k,j}$ will block as many units of $A_{k+1,j+1}$; any
further units of~$A_{k+1,j+1}$ however will not be permanently blocked, and
they will eventually (when all transfers further to the right are completed)
be transferred one row up to the entry~$A_{k,j+1}$. As a result of this
exhaustion, every entry~$A_{k+1,j+1}$ will have become equal to the initial
value of~$A_{k,j}$, and $A_{k+1,0}$ will have become~$0$. We illustrate this
for the example above, showing just the first 6 columns of rows $0,1$:
$$
  \pmatrix{0&1&0&1&0&3\cr 0&0&1&1&3&0\cr}
\buildrel (\upm_0)^* \over \longrightarrow
  \pmatrix{0&1&0&2&2&3\cr 0&0&1&0&1&0\cr}
.\nn
$$
After this exhaustion, the situation in the first~$l$ columns of rows
$k+1,k+2$ is the same as the original situation in rows~$k,k+1$, but shifted
one place to the right so that the final entries disappear, and with
entries~$0$ shifted in at the left. One can therefore continue by exhausting
$\num_{k+1}$, which will involve shifted copies of the same moves (except
those that took place in the final column), and then similarly exhaust
$\li(\num k+2..k+l-1)$; after this, the entries in
$A_{(\N-\set{k})\times\set{l}}$ will have become constant along diagonals, so
that submatrix is then equal to $\sliced\nu{k}l$. Thus we see that row~$k$ of
the submatrix will already have obtained its final value after the initial
exhaustion of $\upm_k$, and $\nu$ can already be read off from this row.

Now for any~$j>0$, the initial value of $A_{k+1,l-j}$ represents the number
$\\_{j-1}-\\_j$ of columns of length~$j$ of~$\set\\$, and the final value of
the entry $A_{k,l-j}$ above it represents the number of columns of that length
of~$\set\nu$. Then the units that are transferred upwards from $A_{k+1,l-j}$
to $A_{k,l-j}$ represent columns whose length is $j$ in~$\set\\$ and remains
so in~$\set\nu$, while the units that already were in $A_{k,l-j}$ initially
(and remain there), and whose number equals the number of units that remain in
$A_{k+1,l-(j-1)}$ if $j>1$, represent columns whose length grows from $j-1$
in~$\set\\$ to~$j$ in~$\set\nu$. The latter columns are the ones that for
which row~$j-1$ of $\set{\nu/\\}$ contains a square. From this one deduces
that $\nu/\\$ is a horizontal strip, and that the initial value of the reverse
$(\foldwith{A_{k,#1}},\ldots(l-1..0))$ of the topmost row of
$A_{(\N-\set{k})\times\set{l}}$ gives the value of $\nu-\\$. In the example
above one has $\nu-\\=(8,5,3,1,1,0)-(5,5,2,1,0,0) =(3,0,1,0,1,0)
=(\foldwith{A_{0,#1}},\ldots(5..0))$.

We might have reached the same conclusion with somewhat less effort by
considering column sums, which are unaffected by upward transfers, but the
details will be needed anyway to understand the full transformation of~$A$ by
upward moves, notably the changes to column~$l$. We now know that this column
initially contains, below the entry~$A_{k,l}$, a vertical copy of the
composition $\mu-\\$, and after the transformation it will contain a similar
copy of $\kappa-\nu$ (but shifted one place less, starting at square~$(k,l)$).
When we consider the exhaustion of some $\upm_{k+i}$, we can focus exclusively
on columns $l-1$~and~$l$, since what we have seen above implies that at that
time one has $A_{k+i,j}\leq{A}_{k+i+1,j+1}$ for all $j\in\set{l-1}$, so that
the presence of columns to the left of column~$l-1$ will not block any upward
transfers in column~$l$.

The first $2\times2$ submatrix to consider is
$A_{\sing{k,k+1}\times\sing{l-1,l}}$ whose entries are initially (i.e., at the
start of the exhaustion of $\upm_{k}$) the following: $A_{k,l-1}$ is the final
entry $\nu_0-\\_0$ of the composition $\revert{\nu-\\}l$ that is in the first
$l$ columns of row~$k$; $A_{k,l}$ is arbitrary; $A_{k+1,l-1}$ is the
top-rightmost entry $\\_0-\\_1$ of the submatrix
$A_{(\N-\set{k+1})\times\set{l}}=\sliced\\{k+1}l$; $A_{k+1,l}$ is the top
entry $\mu_0-\\_0$ of the vertical copy of $\mu-\\$ in column~$l$. At the
start of the exhaustion of $\upm_{k+i}$, the values in the $2\times2$
submatrix $A_{\sing{k+i,k+i+1}\times\sing{l-1,l}}$ are similarly:
$A_{k+i,l-1}=\nu_i-\\_i$ is a copy of the initial value of $A_{k,l-1-i}$;
$A_{k+i,l}$ is some value~$a$ left after the exhaustion of $\upm_{k+i-1}$
(which we shall determine presently); $A_{k+i+1,l-1}$ is the final entry
$\\_i-\\_{i+1}$ of row~$k+1+i$ of $\sliced\\{k+1}l$; and finally
$A_{k+i+1,l}=\mu_i-\\_i$. The upward moves in column~$l$ are then given by
$$
  \pmatrix{\nu_0-\\_0& A_{k,l}\cr \\_0-\\_1 & \mu_0-\\_0\cr}
  \buildrel (\upm_k)^{\max(\mu_0-\nu_0,0)} \over \longrightarrow
  \pmatrix{\nu_0-\\_0& A_{k,l}+\max(\mu_0-\nu_0,0)\cr
           \\_0-\\_1 & \min(\mu_0,\nu_0)-\\_0\cr}
\nn
$$
for the transfer between rows $k$~and~$k+1$ (which will be followed by a
transfer in column~$l-1$ changing it from $\nu_0-\\_0\choose\\_0-\\_1$ into
$\nu_0-\nu_1\choose\nu_1-\\_1$, but that part was already known), and
similarly in rows $k+i$~and~$k+i+1$:
$$
  \pmatrix{\nu_i-\\_i& a\cr \\_i-\\_{i+1} & \mu_i-\\_i\cr}
  \buildrel (\upm_{k+i})^{\max(\mu_i-\nu_i,0)} \over \longrightarrow
  \pmatrix{\nu_i-\\_i& a+\max(\mu_i-\nu_i,0)\cr
           \\_i-\\_{i+1} & \min(\mu_i,\nu_i)-\\_i\cr}
.\label(\bumptransf)
$$
It is now clear that $a=\min(\mu_{i-1},\nu_{i-1})-\\_{i-1}$, and one easily
checks that for all~$i\geq0$, after completion of the exhaustion
of~$\upm_{k+i}$ and of the subsequent exhaustion of~$\upm_{k+i+1}$, the
submatrix $A_{\sing{k+i,k+i+1}\times\sing{l-1,l}}$ will have become
$$ \pmatrix
{\nu_i-\nu_{i+1}    & a+\max(\mu_i-\nu_i,0)\cr
 \nu_{i+1}-\nu_{i+2}&\min(\nu_i,\mu_i)-\\_i
                    +\max(\mu_{i+1}-\nu_{i+1},0)\cr
}.\label(\Knuthstep)
$$
Although $\max(\mu_{i+1}-\nu_{i+1},0)=\max(\mu_{i+1},\nu_{i+1})-\nu_{i+1}$
units have been added to the bottom right coefficient by the exhaustion
of~$\upm_{k+i+1}$, we do not expect that this re-enables applications
of~$\upm_{k+i}$ (since such applications were not required for reducing the
submatrix $A_{(\N-\set{k})\times\set{l}}$ to $\sliced\nu{k}l$ either); indeed
it follows from $\\_i\geq\max(\mu_{i+1},\nu_{i+1})$ that one has
$\nu_i-\nu_{i+1}\geq\min(\nu_i,\mu_i)-\\_i
+\max(\mu_{i+1},\nu_{i+1})-\nu_{i+1}$, whence further upward transfers are
blocked. We may conclude that after the mentioned sequence of operations we
will have determined the new value of column~$l$ representing the composition
$\kappa-\nu$.

In spite of the complications of the reduction, the resulting relation between
$\\$, $\mu$, $\nu$, $\kappa$, and $A_{k,l}$ is surprisingly simple, as there
is no ``carry over'' effect like there was in~(\Burgestep): the value of~$a$
has no effect on the bottom-right entry in~(\Knuthstep) that determines
$\kappa_{i+1}-\nu_{i+1}$. Therefore we can describe the parts of~$\kappa$ by
the following explicit formulae rather than by an algorithm:
$$
\eqalign
{\kappa_0&=A_{k,l}+\max(\mu_0,\nu_0)\cr
 \kappa_{i+1}&=\min(\nu_i,\mu_i)-\\_i+\max(\mu_{i+1},\nu_{i+1})
 \ttex{for all~$i\in\N$.}
}\label(\Knuthgrowtheq)
$$
Comparison with \ref{spin Knuth, equation~(15)} shows that this is the shape
datum for the \caps{RSK}~correspondence, which reflects row-insertion on
semistandard tableaux.

One can give the following interpretation of~(\Knuthgrowtheq) in traditional
terms. If one imagines, during \caps{RSK} insertion, a semistandard Young
tableau~$S$ in whose display the entries~$<l$ fill the diagram~$\set\\$ while
the entries~$l$ fill $\set{\mu/\\}$, and if some insertion of letters~$<l$
into~$S$ makes the shape filled by entries~$<l$ grow to~$\set\nu$, then upon
subsequent row-insertion of $M_{k,l}$ letters~$l$, these will be placed in
row~$0$ starting at position~$\max(\mu_0,\nu_0)$ (i.e., after the unbumped
entries~$l$ of~$S$, if any of them remain, or else after the entries~$<l$
present at that point) making the length of that row
$\max(\mu_0,\nu_0)+A_{k,l}$. That insertion will not bump any letter~$l$, but
the previous insertion of letters~$<l$ will bump $\min(\nu_i,\mu_i)-\\_i$
letters~$l$ from row~$i$ of~$S$, and they will be placed in row~$i+1$ starting
at position $\max(\mu_{i+1},\nu_{i+1})$, making its length
$\max(\mu_{i+1},\nu_{i+1})+\min(\nu_i,\mu_i)-\\_i$. Note in particular that it
is the units that are \emph{not} transferred upwards from the entry at
$(i+k+1,l)$ in~(\bumptransf) that represent letters~$l$ being bumped from
row~$i$ to row~$i+1$.

Now that we have determined the shape datum, the family of shapes
$(M_{(\N-\set{i})\times\set{j}})_{i,j\in\N}$ can be computed using a
(vertically reversed) growth diagram, just like we have done before for other
families of implicit shapes of submatrices. This time in fact we obtain a
classical \caps{RSK} growth diagram, growing from the bottom left corner of
the matrix. We illustrate this in diagram~4 for our usual matrix~$M$.

\midinsert
$$\abovedisplayskip=0pt
\vcenter{
\catcode`@=11
\def~{\unhcopy\QEDbox}
\setbox\z@=\hbox{~}\dimen@=.5\ht\z@ \advance\dimen@.5\dp\z@
\let\f=\hfil
\def\d#1 {\vtop\bgroup\smallskip
          \vbox to\z@{\vss\llap{$#1$\kern6pt}\kern 4pt}%
          \futurelet\n\E}
\def\E{\?\ifx\n\f \then\smallskip\egroup \else\afterassignment\e\count@ \fi}
\def\e{\hbox{\loop~\advance\count@\m@ne \ifnum\count@>\z@ \repeat}%
       \nointerlineskip\futurelet\n\E}
\def\g#1{\omit#1\hfill\ignorespaces}
\def\h#1{\lower8pt\hbox{#1} \vrule height 9pt }
\def\zz{\vtop{\smallskip\hbox{$\circ$}}\hfill}
\def\z{\omit\zz}
\offinterlineskip
\everycr{\noalign{\nobreak}}
\tabskip=\z@
\halign{\h# \tabskip=17 pt &&\d#\f\cr
\omit\hfil\strut\vrule
    & \g0 & \g1 & \g2 & \g3 & \g4 &  \g5  &  \g6  &  \g7 \cr
\noalign{\hrule}
0 &\z&{} 2 &{} 3 2  &{} 4 2 2&{} 4 3 2 1&{} 6 4 3 1&{} 8 5 3 1 1
  &{} 8 8 5 3 1\cr
1 &\z&1 1    &0 3 1   &1 3 2 1 &0 4 2 2 &1 5 4 2 &2 5 5 2 1 &0 8 5 5 2 \cr
2 &\z&1 \f\zz&1 2     &0 3 1   &1 3 2   &1 4 3   &0 5 3 1   &3 5 5 3   \cr
3 &\z&0 \f\zz&2 \f\zz &1 1     &0 2     &1 3     &1 3 1     &0 5 3     \cr
4 &\z&0 \f\zz&0 \f\zz &1 \f\zz &1 \f\zz &1 \f\zz &0 1       & 2 3      \cr
5 &\z&0 \f\zz&0 \f\zz &0 \f\zz &0 \f\zz &0 \f\zz &1 \f\zz   & 2 \f\zz  \cr
}
}
$$
\centerpar{{\bf Diagram~4.} Vertically reversed growth diagram
  with integral matrix and shapes
  $(\is{M_{(\N-\set{i})\times\set{j}}})_{i,j\in\N}$, computed using the
  \caps{RSK} shape datum.}
\endinsert

\proclaim Proposition. \RSKprop
If for any $M\in\Mat$ one associates to each grid point $(i,j)\in\N^2$ the
shape $\is{M_{(\N-\set{i})\times\set{j}}}$, then one obtains a vertically
reversed growth diagram for the \caps{RSK} shape datum.
\QED

Along the top border of diagram~4 we find the semistandard Young tableau $S$
of (\Sdef); this is no surprise since we already knew that
$S=(\is{M_{\N\times\set{j}}})_{j\in\N}$, which was computed along the bottom
border in diagram~1. Along the right border we find a reverse semistandard
Young tableau $(\is{M_{(\N-\set{i})\times\N}})_{i\in\N}$ that we did not
encounter before. However, the way in which~$S$, whose integral encoding is
the matrix~$P$ of~(\PQintdiagram), was found here, corresponds to the
description of theorem~\RSKcompareprop; according to comments following
that theorem, the recording tableau of that insertion is the Sch\"utzenberger
dual of the semistandard Young tableau $\bar{L}$ encoded by the matrix $Q\tr$.
Indeed the Sch\"utzenberger dual of
$$
  \bar{L}=  \Young(0,0,0,0,0,1,1,2|1,1,1,1,1,2,2,3|2,2,3,3,4|3,3,4|4)
\ttext{is}
  \bar{L}^*=\Young(4,4,4,3,3,1,1,1|3,3,3,2,2,0,0,0|2,2,2,1,1|1,1,0|0)
,\label(\Lbarstardef)
$$
which is the reverse semistandard Young tableau that can be read off along the
right border of diagram~4.

\subsection The Sch\"utzenberger involution.

We have encountered the Sch\"utzenberger involution in two similar situations,
for binary and for integral matrices. In the former case it defined a
bijection mapping the reverse transpose semistandard Young
tableau~$R=(\is{M_{\N\times(\N-\set{j}}})_{j\in\N}$ to the transpose
semistandard Young tableau~$R^*=(\is{M_{\N\times\set{j}}})_{j\in\N}$; in the
latter case it gave a bijection transforming the semistandard Young tableau
$\bar{L}=(\is{M_{\set{i}\times\N}})_{i\in\N}$ into the reverse semistandard
Young tableau $\bar{L}^*=(\is{M_{(\N-\set{i})\times\N}})_{i\in\N}$. The fact
that the members of these pairs of tableaux mutually determine each other,
independently of~$M$, may seem surprising at first glance, since at
corresponding positions in their sequence they contain implicit shapes of
complementary submatrices of~$M$, which shapes are completely unrelated. On
the other hand, these pairs of tableaux correspond to matrices obtained
from~$M$ by exhausting crystal operations in opposite directions. We have seen
that $R$~and~$R^*$ correspond to the matrices $P$~and~$\widetilde{P}$ obtained
from the binary matrix~$M$ by exhausting upward respectively downward
operations (cf.\ propositions \Knuthbinaryprop\ and~\bingdprop). In the
integral case $\bar{L}$ corresponds to the matrix~$Q$ obtained from~$M$ by
exhausting leftward operations, and one easily sees that $\bar{L}^*$ similarly
corresponds to a matrix~$\widetilde{Q}$ obtained from~$M$ by exhausting
rightward operations (up to some sufficiently large limit); in the case of the
example this is
\bigdisplay
  \widetilde{Q}=
\pmatrix
{0&0&1&1&0&3&0\cr
 0&0&0&2&2&0&3\cr
 0&0&0&0&3&2&0\cr
 0&0&0&0&0&3&2\cr
 0&0&0&0&0&0&3\cr
},\nn
$$
whose columns from right to left give the weights of the rows of the display
of $\bar{L}^*$ from top to bottom. Since for each of the pairs of matrices
$(P,\widetilde{P})$ and $(Q,\widetilde{Q})$ the members mutually determine
each other (by exhaustion of appropriate crystal operations), it is clear that
the same must be true for the pairs of tableaux $(R,R^*)$ and
$(\bar{L},\bar{L}^*)$. Exhaustion in the remaining pairs of opposite
directions also give pairs of matrices that correspond to tableaux related by
the Sch\"utzenberger involution.

We shall not study the Sch\"utzenberger involution in detail here, and just
observe some facts without formal proof. Its original definition
in~\ref{Schutzenberger QR} for standard Young tableaux can be generalised,
using appropriate notions of standardisation, to semistandard or to transpose
semistandard tableaux; doing so one must map tableaux of weight~$\alpha$
either to tableaux of weight~$\revert\alpha{l}$ for some sufficiently
large~$l$, or to reversed tableaux of weight~$\alpha$ (while keeping the
semistandard or transpose semistandard attribute). Above we have implicitly
assumed the latter, more elegant, form of these generalisations.

The original definition of the Sch\"utzenberger involution is by an algorithm
that repeatedly applies a deflation operation (also called evacuation) to a
standard Young tableau~$S$, which removes a square from the diagram of its
shape at every step, until the shape becomes~$(0)$; the resulting decreasing
sequence of shapes defines the Sch\"utzenberger dual~$S^*$ of~$S$. It is
natural to index the sequence defining the tableau~$S^{(i)}$ deflated
$i$~times by $\N-\set{i}$. If one then writes
$S^{(i)}=(\\^{(i,j)})_{j\in\N-\set{i}}$, then the definition of the deflation
procedure amounts to a simple rule for determining $\\^{(i+1,j+1)}$ once
$\\^{(i+1,j)}$, $\\^{(i,j)}$, and $\\^{(i,j+1)}$ are known. Thus
$S=(\\^{(0,j)})_{j=0}^n$ determines $S^*=(\\^{(i,n)})_{i=0}^n$ via a
two-parameter family of shapes governed by a local rule, in much the same way
as a growth diagram, but without matrix entries.

The generalisations to semistandard or to transpose semistandard tableaux can
be formulated in terms of deflation by an entire horizontal respectively
vertical strip at a time, and give rise to similar families of shapes
$(\\^{(i,j)})_{i,j;i\leq{j}}$, in which neighbouring shapes now differ by
horizontal strips respectively by vertical strips. Here the rule for
determining $\\^{(i+1,j+1)}$ given $\\^{(i+1,j)}$, $\\^{(i,j)}$, and
$\\^{(i,j+1)}$ can be formulated as follows, using the notion of tableau
switching for semistandard tableaux given in \ref{L-R intro, definition~2.2.4}
(the transpose semistandard case is similar): applying tableau switching to
the single-strip skew semistandard tableaux $T_0=(\\^{(i+1,j)}\leh\\^{(i,j)})$
and $T_1=(\\^{(i,j)}\leh\\^{(i,j+1)})$ gives similar tableaux
$(T'_1,T'_0)=X(T_0,T_1)$ where $T'_0=(\\^{(i+1,j)}\leh\\^{(i+1,j+1)})$ and
$T'_1=(\\^{(i+1,j+1)}\leh\\^{(i,j+1)})$. An analysis of what tableau switching
means in this particular case leads to the following alternative formulation:
the two skew semistandard tableaux of two horizontal strips each,
$\\^{(i+1,j)}\leh\\^{(i,j)}\leh\\^{(i,j+1)}$ and
$\\^{(i+1,j)}\leh\\^{(i+1,j+1)}\leh\\^{(i,j+1)}$, are related to one another
by a Bender-Knuth involution (the unique applicable one). Interestingly
Bender-Knuth involutions, which we preferred to replace by ladder reversal
involutions for establishing the symmetry of skew Schur functions to get more
pleasant properties, thus still have their proper place in the theory.

The way in which we have encountered the Sch\"utzenberger involution in fact
relates more naturally to a different way to compute it. The matrix~$P$
obtained from $M\in\Mat$ by exhausting upward crystal operations as in
(\PQintdiagram) encodes a unique semistandard Young tableau~$S$, of
shape~$\pi/(0)$ where $\pi=\is{M}$. If $k>0$ is so large that
$M=M_{\set{k}\times\N}$, then the matrix~$\widetilde{P}$ obtained from~$M$ by
exhausting downward crystal operations $\dnm_i$ for $i\in\set{k-1}$ determines
the Sch\"utzenberger dual $S^*=(\mu^{(j)})_{j\in\N}$ of~$S$, namely by
$\mu^{(j)} =\is{\widetilde{P}_{\N\times(\N-\set{j})}}
=\revert{\rsum{\widetilde{P}_{\N\times(\N-\set{j})}}}k$. But if one fixes
$l\geq\pi_0$, and $\rho\in\Part$ is such that
$\set\rho=\set{k}\times\set{l}$, then $\smash{\widetilde{P}}$ is also the
integral encoding of a unique semistandard tableau of shape
$\rho/\pi^{\diamond}$, where $\pi^{\diamond}=\rho-\revert{\pi}k$, namely of
$S^{*\diamond} = (\mu^{(j)\diamond})_{j\in\N}=
(\rho-\revert{\mu^{(j)}}k)_{j\in\N}$. This is the tableau whose display is
obtained from that of~$S^*$ by the rotation by a half turn that maps the
rectangle~$\set\rho$ to itself.

The semistandard tableaux $S$ and~$S^{*\diamond}$ are encoded by integral
matrices $P$ and~$\widetilde{P}$ that can be obtained from one another by
exhausting upward respectively downward crystal operations, so it follows from
proposition~\jdtrelprop\ that $S$ is the rectification by jeu de taquin
of~$S^{*\diamond}$, and since rotation by a half turn transforms inward jeu de
taquin slides into outward slides, $S^{*\diamond}$ is similarly the
``bottom-right rectification'' of~$S$ within the rectangular shape~$\rho$. In
our example one may check, for $(k,l)=(5,8)$, that
\bigdisplay
  \widetilde{P}=\pmatrix
{1&0&0&0&0&0&0\cr
 1&1&1&0&0&0&0\cr
 0&2&1&1&1&0&0\cr
 0&0&1&1&3&3&0\cr
 0&0&0&0&0&1&7\cr
}
\ttext{and}
  S^{*\diamond}=
\Skew(7:0|5:0,1,2|3:1,1,2,3,4|0:2,3,4,4,4,5,5,5|0:5,6,6,6,6,6,6,6)
\,,\nn
$$
and that the rectification of $S^{*\diamond}$ is the tableau~$S$ of~(\Sdef);
thus $S\leftrightarrow{S}^*$ can be computed via rectification. The same
relation exists between $\bar{L}$ and $\bar{L}^{*\diamond}$ in~(\Lbarstardef),
and between $R$ in~(\Rdef) and $R^{*\diamond}$ obtained from~$R^*$
in~(\Rstardef); in the latter case a version of jeu de taquin for transpose
semistandard tableaux is used.

This alternative method of computing the involution was first presented by
Sch\"utzenberger in \ref{Schutzenberger CRGS}. It is not so easy to deduce it
from the original definition of the Sch\"utzenberger involution, but the
opposite is straightforward, once it is established that the rectification of
any semistandard tableau~$T$ of shape $\rho/\\^\diamond$ is of shape~$\\$.
Then the shape $\mu^{(j)}$ in $S^*=(\mu^{(j)})_{j\in\N}$ can be determined as
the shape of the rectification of the skew tableau $S^{(j)}$ obtained by
removing $j$ initial horizontal strips from~$S$, and these shapes can be found
by successively deflating~$S$ by one horizontal strip at a time. By
construction that rectification is given by the sequence of shapes
$(\\^{(j,l)})_{l=j}^n$, where $(\\^{(j,l)})_{j,l; j\leq{l}\leq{n}}$ is the
semistandard tableau switching family used to compute
$S^*=(\\^{(j,n)})_{j=0}^n$ starting from $S=(\\^{(0,l)})_{l=0}^n$.

Now if $P$ is the integral encoding of~$S$, then $P_{\N\times(\N-\set{j})}$ is
the integral encoding of~$S^{(j)}$, for all $j\in\N$, and the integral
encoding of the rectificationof~$S^{(j)}$ can be found by exhausting upwards
crystal operations applied to $P_{\N\times(\N-\set{j})}$. This gives have the
following interpretation for the shapes in the semistandard tableau switching
family: $\\^{(j,l)}=\is{M_{\N\times(\set{l}-\set{j})}}$, which agrees with
$S^*=(\is{P_{\N\times(\N-\set{j})}})_{j\in\N}$. Since vertical crystal
operations do not change the implicit shape,
$(\is{M_{\N\times(\set{l}-\set{j})}})_{j,l; j\leq{l}}$ is a semistandard
tableau switching family for any integral matrix~$M$. The same is true for
$(\is{M_{(\set{k}-\set{i})\times\N}})_{i,k; i\leq{k}}$; these statements hold
in the binary case as well, with the adjustment that the former family is
transpose semistandard. However, in contrast with the situation for the
various growth diagrams we considered, we have not been able to find a simple
direct argument why a local rule for these families of shapes should exist,
without referring to the entire family; for instance one showing that
knowledge of the shapes $\is{M_{\N\times(\set{l}-\set{j+1})}}$,
$\is{M_{\N\times(\set{l}-\set{j})}}$, and
$\is{M_{\N\times(\set{l+1}-\set{j})}}$ suffices to determine
$\is{M_{\N\times(\set{l+1}-\set{j+1})}}$. In principle one could establish the
local rule for that family using those for
$(M_{\set{k}\times\set{l}})_{k,l\in\N}$ and for
$(M_{\set{k}\times(\N-\set{j})})_{k,j\in\N}$, and induction on~$k$, similarly
to the proof of \ref{RSS elementary, theorem~4.1.1}, but this is not very
attractive.

 \par

\subsecno=-1

\newsection Characterisation of the implicit shape.
\labelsec\chainsec

The implicit shape of a matrix~$M$ characterises by definition the double
crystal in which $M$ occurs. We have by now seen many different ways to
compute this shape, but these involve either transforming $M$ into a more
special matrix within its double crystal (by exhausting some operations), or
constructing a growth diagram from~$M$. It would be convenient to have a
characterisation of~$\is{M}$ in terms of~$M$ itself, which would provide a
quantity invariant on double crystals and distinguishing any two distinct
double crystals. In this section we shall present such a characterisation
of~$\is{M}$, which will however usually not provide a practical way to compute
it; the methods given above can therefore be interpreted as algorithms to
compute the invariant that we shall describe below.

It is well known that the length of the longest increasing and decreasing
sequences extracted from a permutation~$\sigma$ (represented as a word) are
given by the lengths of the top row and leftmost column of the tableaux
associated to~$\sigma$ under the Schensted correspondence; extending this
result, it was shown in~\ref{Greene Schensted} that the entire common shape of
those tableaux can be characterised in terms of the permutation matrix (which
is viewed as describing a partially ordered set). A version of Greene's poset
invariant can also be defined for integral and binary matrices, with a similar
interpretation as the shape of the tableaux associated under the \caps{RSK}
correspondence and its dual. Therefore we shall not waste any time trying to
derive the invariant we are looking for, and simply present the proper
formulation. (Note that the word invariant is used in two different senses
here: the partition defined by Greene depends only on the abstract finite
poset defined by the permutation, but its extension to matrices is invariant
under crystal operations, even though these do in general alter the associated
poset.)

\proclaim Definitions.
Let $M\in\bMat$ be a binary matrix and $N\in\Mat$ an integral matrix, and
$m\in\N$.
\defitem $M$ is called a chain if $M_{k,l}=1$ implies $M_{i,j}=0$  for all
         $i<k$ and $j\geq{l}$, and it is called an antichain if $M_{k,l}=1$
         implies $M_{i,j}=0$ for all $i\leq{k}$ and $j<l$.
\defitem $N$ is called a chain if $N_{k,l}>0$ implies $N_{i,j}=0$  for all
         $i<k$ and $j<l$, and it is called an antichain if $N_{k,l}\leq1$ for
         all $k,l\in\N$, while $N_{k,l}=1$ implies $N_{i,j}=0$ for all
         $(i,j)\neq(k,l)$ with $i\leq{k}$ and $j\geq{l}$.
\defitem An $m$-chain or an $m$-antichain in $M$ is a sequence
         $(M^{(k)})_{k\in\set{m}}$ of $m$ binary matrices, all chains
         respectively antichains, such that $M=\sum_{k\in\set{m}}M^{(k)}+M'$
         for some $M'\in\bMat$.
\defitem An $m$-chain or an $m$-antichain in $N$ is a sequence
         $(N^{(k)})_{k\in\set{m}}$ of $m$ integral matrices, all
         chains respectively antichains, such that
         $N=\sum_{k\in\set{m}}N^{(k)}+N'$ for some $N'\in\Mat$.
\defitem The $m$-length of a binary or integral matrix~$A$ is
         the maximum of $\sum_{k\in\set{m}}|A^{(k)}|$, taken over all
         $m$-chains $(A^{(k)})_{k\in\set{m}}$ in $A$, where
         $|A|=\sum_{i,j\in\N}A_{i,j}$. The $m$-breadth of~$A$ is the maximum
         of $\sum_{k\in\set{m}}|A^{(k)}|$, taken over all $m$-antichains
         $(A^{(k)})_{k\in\set{m}}$ in $A$.

The terminology is derived from the partially ordered sets associated as
follows to the matrices. For $M\in\bMat$ an element $x_{i,j}$ is associated to
every bit $M_{i,j}=1$, and $x_{i,j}<x_{k,l}$ holds if and only if $i\leq{k}$
and~$j<l$. For $N\in\Mat$ a chain of $N_{i,j}$ elements is associated to every
entry $N_{i,j}$, and for every element~$x$ of the chain for~$N_{i,j}$ and
every element~$y$ of the chain for another entry~$N_{k,l}$ one has $x<y$ if
and only if $i\geq{k}$ and $j\leq{l}$. In this setting it is convenient to use
some more related terminology. For $m$-chain or $m$-antichain
$(A^{(k)})_{k\in\set{m}}$, either binary or integral, we shall call
$\sum_{k\in\set{m}}|A^{(k)}|$ its size. If $A,B$ are binary or integral
matrices, the matrix $(\min(A_{i,j},B_{i,j}))_{i,j\in\N}$ will be called their
intersection $A\thru{B}$; in case $B=A\thru{B}$ we shall write $A\supset{B}$
and say that $A$ contains~$B$. In particular $A$ contains the sum of any
$m$-(anti)chain in~$A$. Then the intersection of a chain and an antichain has
size at most~$1$, and it follows that the sum of an $m$-chain cannot contain
any antichain of size~$>m$, and that the sum of an $m$-antichain cannot
contain any chain of size~$>m$. The converse is also true, namely that any
matrix that does not contain any antichains of size~$>m$ is the sum of an
$m$-chain, and one that does not contain any chains of size~$>m$ is the sum of
an $m$-antichain, but this is somewhat less obvious.

The definition is formulated so as to emphasise the similarities between the
different cases, but there are notable differences as well. In the binary case
chains and antichains can be interpreted as sets of squares where the matrix
has bits~`$1$', and the sums involved in $m$-(anti)chains obviously must
correspond to disjoint unions of such sets. In the integral case antichains
can also be seen as sets of squares, but the antichains in an $m$-antichain
need not be disjoint, or even distinct. On the other hand one could require in
the definition of $m$-length for the integral case that the nonzero entries in
the chains of the $m$-chain involved occur in disjoint squares: the indicated
maximum would not decrease due to this restriction.

It is not hard to check both in the binary and in the integral case that a
matrix~$M$ is a chain if and only if $\is{M}$ has at most one non-zero part,
in other words if $\is{M}_1=0$, and that $M$ is an antichain if and only if
all parts of $\is{M}$ are $\leq1$, in other words if $\is{M}\tr_1=0$. In
general it is quite hard to determine, for a given $M$ and $m$, an
$m$-(anti)chain in~$M$ of maximal size. In the special cases $M=\diagram\\$
and $M=\diagon\\$ however, this is easy. For $M=\diagram\\$, the rows of~$M$
define chains and the columns antichains, and the $m$-chain or $m$-antichain
formed by the $m$ longest rows respectively columns has the maximal possible
size, since its intersection with each individual column respectively row has
maximal size. For $M=\diagon\\$, each chain contained in~$M$ has at most one
nonzero (diagonal) entry, so an $m$-chain of maximal size contained in~$M$ is
formed by extracting the $m$ largest diagonal entries:
$(M_{\set{i}\times\set{i}})_{i\in\set{m}}$. In this case the size of an
$m$-antichain in~$M$ is at most $\sum_{i\in\N}\min(m,M_{i,i})$, and this size
is obtained by taking the $m$-antichain $(\1(\\\tr_j))_{j\in\set{m}}$, where
$\1(l)$ denotes the $l\times{l}$ identity matrix
$(\Kr{i=j\in\set{l}})_{i,j\in\N}$. Note that in both cases the $m$-length
of~$M$ is $\sum_{i\in\set{m}}\\_i$, and its $m$-breadth is
$\sum_{i\in\set{m}}\\\tr_i$. The main result of this subsection is that this
is true in general, where $\\=\is{M}$.

\proclaim Theorem. \isinvarthm
For any binary or integral matrix $M$, and any $m\in\N$, the $m$-length of~$M$
is equal to $\sum_{i\in\set{m}}\\_i$, and the $m$-breadth of~$M$
is equal to $\sum_{i\in\set{m}}\\\tr_i$, where $\\=\is{M}$.

Since we know the theorem holds for normal forms, it is equivalent to the
statement that for every~$m$ the $m$-length and $m$-breadth of a matrix are
invariant under all crystal operations, and this is the way in which we shall
prove it. Note that one consequence of the theorem is that the sequence of all
$m$-lengths of~$M$ determines that of all $m$-breadths and vice versa. It is
known that this statement is true in general for finite posets, see
\ref{Greene poset, Theorem~1.6}, but a proof for that case is rather more
difficult than for the situation of our theorem.

There are numerous ways to approach a proof of the theorem, many of which
start by reducing to a simplified situation. It is fairly easy to see that the
argument by which the theorem holds for normal forms can be adapted to show
that is holds whenever either upward or leftward moves are exhausted, so that
a proof for operations in one direction would suffice (but this facilitates
only the binary case). The result is also known to hold for permutation
matrices, see \ref{Greene Schensted, Theorem~3.1}, and it is possible to
enlarge matrices, pulling their units apart so as to obtain a permutation
matrix, while not changing their implicit shape or $m$-length and $m$-breadth;
this would provide a proof. One can define a notion of chains in the display
of a semistandard tableau~$T$, and \ref{Schutzenberger CRGS, th\'eor\`eme~3.1}
states a variant of Greene's result for that case, where the shape is that of
the rectification of~$T$; it can be seen that the $m$-length and $m$-breadth
so associated to~$T$ are the same as those of its integral encoding, so this
would provide a more direct proof for the integral case (but since the chains
of~$T$ do not correspond to those of its binary encoding, more work would be
needed for that case). Finally one could try, as is done in the proofs of the
results cited, to decompose the changes to the poset structure into minimal
modifications involving only three elements, and prove that these do not
effect $m$-length or $m$-breadth.

All of these methods however require arguments to justify the reduction, which
are not related directly to the definitions of $m$-length or $m$-breadth, or
of crystal operations; they also depend on nontrivial results about the
Robinson-Schensted correspondence. We believe that for understanding why the
theorem is true, it is most instructive to use a head-on approach, showing
directly that the definition of crystal operations prevents them from changing
any $m$-length or $m$-breadth; in fact this is not too hard to do, although
the situation is more complicated than in the $3$-unit case. By using only the
definition of crystal operations, this approach could even serve as a starting
point for developing their theory (but we believe this would not would make
that development more transparent). We also choose a ``positive'' point of
view, characterising matrices that are the sum of an $m$-(anti)chain
systematically by the existence of such a decomposition, rather than by the
non-existence of any antichain (respectively chain) of size~$m+1$. The proof
we give will consist of a series of lemmas that each prove a part of the
mentioned invariance.

The number of cases to be treated in principle is considerable, since we need
to show both in the binary and the integral case that neither the $m$-length
nor the $m$-breadth can increase nor decrease under crystal operations in any
of the four possible directions. We reduce this number using the usual
symmetry considerations to four distinct cases, which can moreover be regarded
as two pairs of analogous cases. Both in the binary and integral cases,
rotation by a half turn leaves both $m$-length and $m$-breadth invariant, so
for each of these quantities and any pair of opposite directions it will
suffice to show that the quantity will not decrease under a crystal operation
in one direction: by rotation symmetry the inverse crystal operation will then
not decrease the quantity either, so that it is invariant under both
operations. In the binary case, rotation of matrices by a quarter turn
transforms $m$-chains into $m$-antichains, which permits us to consider only
the $m$-length; in the integral case we must consider both $m$-length and
$m$-breadth, but transposition symmetry will allow us to consider transfers in
one direction only.

The binary and the integral case are sufficiently different to merit separate
treatment. We shall start with the somewhat simpler binary case, where as said
we shall consider only the $m$-chains. Much if our discussion will be about
locating bits $M_{i,j}=1$; we shall call this a unit at position~$(i,j)$, and
we shall view the individual chains of an $m$-chain in~$M$ as sets of such
units.

\proclaim Lemma. \binupinvlem
Let a binary matrix $M\in\bMat$ contain an $m$-chain of size~$n$ for some
$m,n\in\N$, and let $M'=\upm_k(M)$ for some $k\in\N$. Then $M'$ also contains
an $m$-chain of size~$n$.

\proof
Let $M_{k,l}=0$ and $M_{k+1,l}=1$ be the bits interchanged by the application
of~$\upm_k$. One may assume that the unit at $(k+1,l)$ occurs in one of the
chains~$C$ of the $m$-chain, for otherwise one can simply keep the same
$m$-chain. If, after removing the unit at~$(k+1,l)$ from~$C$, any of the
chains of the $m$-chain would remain a chain when a unit at~$(k,l)$ is added
to it, then an $m$-chain of size~$n$ in~$M'$ can be formed in that way; we
shall show that after possibly changing the choice of the $m$-chain of
size~$n$ in~$M$, we can always arrange to be in that situation. Let the units
of~$C$ to the left of~$(k+1,l)$ in row~$k+1$ be at positions
$\setof(k+1,j_i):i\in\set{r}\endset$ where $r\in\N$ and $\fold<(j0..r-1)<l$.
If $r=0$ (i.e., if there are no such units) then after removing the unit
at~$(k+1,l)$ we can add a unit at~$(k,l)$ to~$C$ itself, and we are done.

Otherwise definition~\binmovedef\defitemnr1 gives us
$\sum_{j=j_0}^{l-1}M_{k,j}\geq\sum_{j=j_0}^{l-1}M_{k+1,j}\geq{r}$, so there
are at least $r$ different indices $j<l$ such that $M_{k,j}=1$; let the $r$
largest such indices be $\fold<({j'}0..r-1)$. From the fact that
$\sum_{j=j_i}^{l-1}M_{k,j}\geq\sum_{j=j_i}^{l-1}M_{k+1,j}$ for $i\in\set{r}$,
it follows that $j_i\leq{j'_i}$ for all~$i$. None of the units $M_{k,j'_i}=1$
can occur in~$C$ because of its unit at $(k+1,j_0)$. If none of them occur in
any of the other chains of the $m$-chain either, then one can modify~$C$ by
removing the units at $(k+1,j_i)$ and adding units at $(k,j'_i)$, for
$i\in\set{r}$; this reduces us to the case $r=0$ already treated. If at least
one of the units $M_{k,j'_i}=1$ does occur in a chain of the $m$-chain, let
$j'_s$ be the maximal index for which this happens, and let the unit
at~$(k,j'_s)$ occur in the chain~$C'$. If $C'$ contains any units in row~$k$
to the right of the one at~$(k,j'_s)$ then this must be in a column~$>l$, and
a unit at~$(k+1,l)$ can be added to~$C'$. If not, then we replace $C,C'$ as
follows by a new pair of chains $D,D'$ of the same total size: the chain $D$
contains the units of~$C'$ in columns $\leq{j'_s}$, the units at $(k,j'_i)$
for $s<i<r$ and the units of~$C$ in columns~$\geq{l}$, while $D$ contains the
units of~$C$ in columns~$\leq{j_s}$ and the units of~$C'$ in columns~$>j'_s$.
After this change $D$ takes over the role of~$C$, and we are in the case
$r=0$, so we are done.
\QED

This situation of the final case of the proof is schematically represented in
figure~1, in which black circles represent units, the white circle is the bit
$M_{k,l}=0$, the red and blue lines indicate the chains $C$~and~$C'$, with
dotted parts that do not occur in $D$ with the unit at~$(k+1,l)$ replaced by
one at~$(k,l)$, nor in~$D'$; the green lines indicate the new parts of the
latter chains and the yellow dotted lines show the correspondence
$(k+1,j_i)\leftrightarrow(k,j'_i)$.

\midinsert
\centerline{\epsfysize=28mm \epsfbox{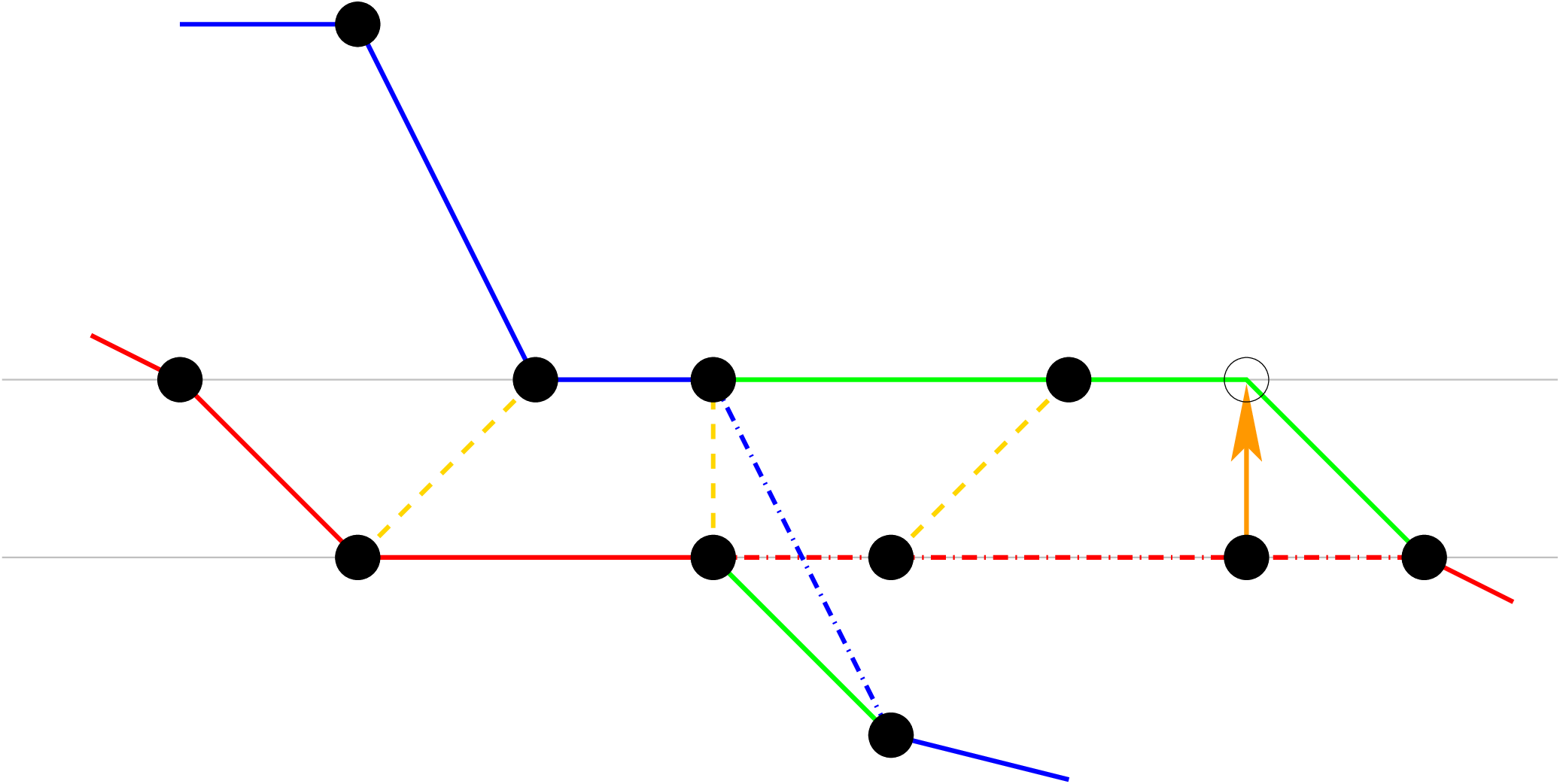}}
\centerpar{{\bf Figure~1.} Schematic view of an partial exchange of chains for
 a vertical move in a binary matrix.}
\endinsert

\proclaim Lemma. \binltinvlem
Let a binary matrix $M\in\bMat$ contain an $m$-chain of size~$n$ for some
$m,n\in\N$, and let $M'=\ltm_l(M)$ for some $l\in\N$. Then $M'$ also contains
an $m$-chain of size~$n$.

Before stating the somewhat technical formal proof, we shall give an informal
reasoning leading to it, ordering the cases that can occur by increasing
complication. Let $M_{k,l}=0$ and $M_{k,l+1}=1$ be the bits interchanged
by~$\ltm_l$. First consider the relatively simple case $m=1$. If the unit at
$(k,l+1)$ does not occur in the given ($1$-)chain~$C$, then $C$ will still be
a chain in~$M'$. If the unit does occur in~$C$, we can just replace it in~$C$
by the unit that was moved to the left, unless $C$ already contains a unit in
column~$l$, necessarily at some $(k',l)$ with $k'<k$. The inequality
$\sum_{i=k'}^{k-1}M_{i,l}\leq\sum_{i=k'}^{k-1}M_{i,l+1}$ from
definition~\binmovedef\defitemnr2 shows that there exists a unit at
$(k'',l+1)$ with $k'\leq{k''}<k$; then by replacing in~$C$ the unit at
$(k,l+1)$ by the one at $(k'',l+1)$, one obtains a chain of size~$n$ in~$M$
that is also in~$M'$.

In the case $m>1$ there is at most one chain~$C$ in the $m$-chain with a unit
at $(k,l+1)$, and if so one can apply the reasoning above to~$C$, in order to
find another chain to replace it in the $m$-chain. The reasoning proceeds as
before, but at the end there is a complication in case the unit at~$(k'',l+1)$
occurs in another chain~$C'$ of the $m$-chain. If $C'$ contains no unit in
column~$l$, then one can resolve the difficulty by ``genetic recombination''
of $C$ and~$C'$ (a crossover of chains somewhat simpler than the one in
figure~1): the unit at~$(k',l)$ and any other units of~$C$ in columns $<l$
combine with the unit at~$(k'',l+1)$ and any other units of~$C'$ in columns
$>l+1$ to form a new chain~$D$, while the remaining units of~$C'$, which are
in columns $<l$, combine with the unit at~$(k,l+1)$ and any other units of~$C$
in columns $>l+1$, to form another new chain~$D'$, which does not contain a
unit in column~$l$. This situation is schematically represented in figure~2,
using the same conventions as in figure~1. If on the other hand $C'$ does
contain a unit in column~$l$, then one can repeat the same kind of argument to
find another unit in column~$l+1$, and so on; eventually there must be such a
unit that is not in a chain with another unit in column~$l$, and all the
pieces of the chains involved can be pieced together to form an $m$-chain
in~$M'$ (unlike for the previous lemma, any number of chains can be involved
in the modifications here). In order to proceed in an orderly fashion the
formal proof shall start with untangling any chains that might be involved.
Figure~3 schematically shows the two types of modification that may be
required at the end of the proof below to reorganise the chains in the most
difficult cases: a shift by several units in column~$l+1$ between chains, or a
cyclic recombination of parts of several chains.

\midinsert
\centerline{\epsfysize=3cm \epsfbox{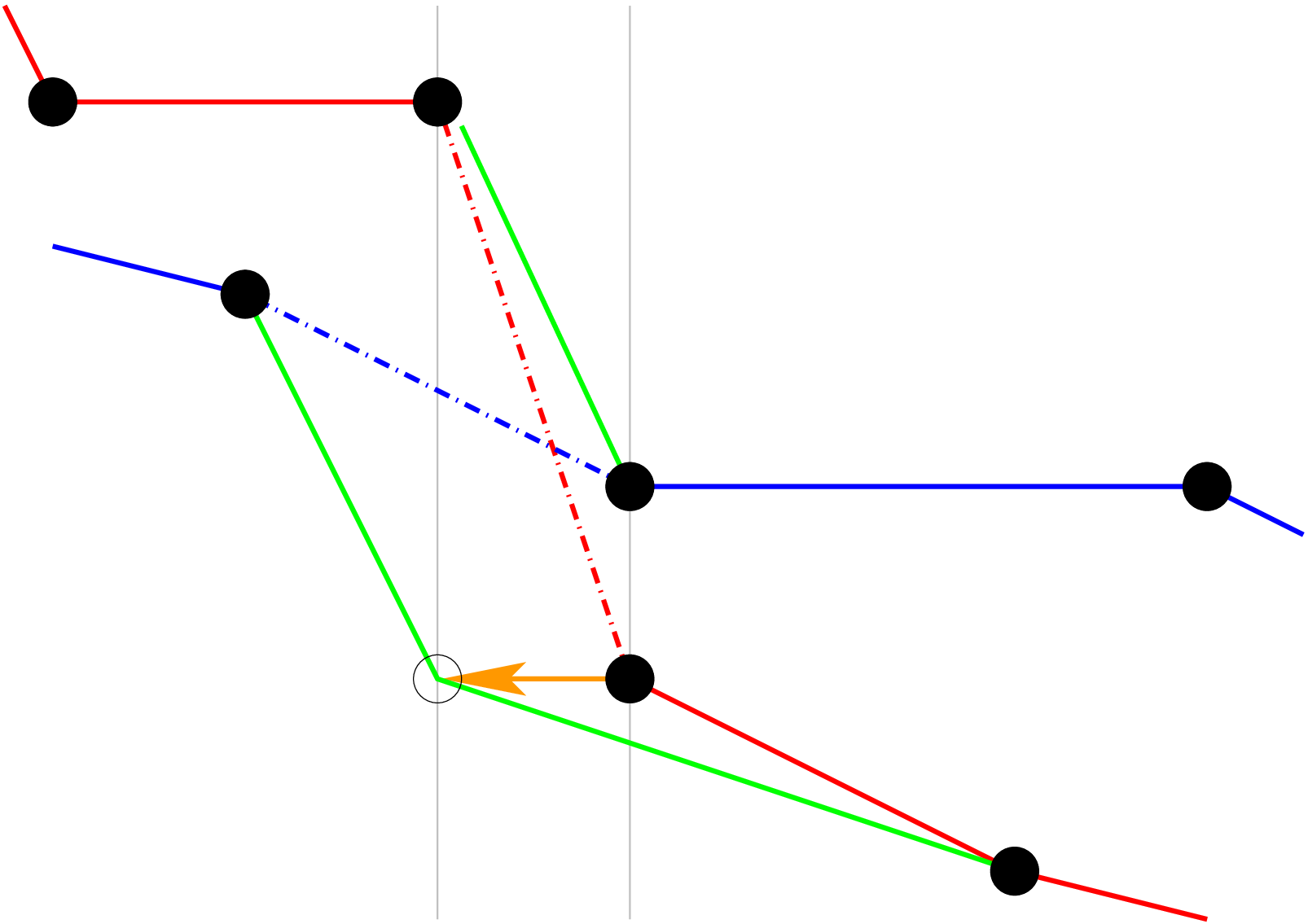}}
\centerpar{{\bf Figure~2.}
Schematic view of a simple recombination of chains for a horizontal move in a
binary matrix.}
\endinsert

\proof
Let the move effected by $\ltm_j$ interchange bits $M_{k,l}=0$ and
$M_{k,l+1}=1$. Consider the subset~$S$ of chains in the $m$-chain that have a
unit both in column $l$ and in column~$l+1$, and put $r=\Card{S}\geq0$.
Let $I$ be the set of rows of those units in column $l$, and $I'$ the set of
rows of those units in column $l+1$. Since each chain of~$S$ contributes an
$i\in{I}$ and an $i'\in{I'}$ with $i\leq{i'}$, it follows after sorting that
if $I=\sing{\li(i0..r-1)}$ and $I'=\sing{\li({i'}0..r-1)}$ with
$\fold<(i0..r-1)$ and $\fold<({i'}0..r-1)$, then $i_j\leq{i'_j}$ for all
$j\in\set{r}$. Then we may assume (after recombining parts of the chains if
necessary) that the units $M_{i_j,l}=1$ and $M_{i'_j,l+1}=1$ occur in the same
chain $C^{(j)}\in{S}$, for all $t\in\set{r}$.

If $k\notin{I'}$, then either none of the chains of the $m$-chain has a unit
at~$(k,l+1)$ so that one can keep the $m$-chain, or that unit occurs in a
chain without unit in column~$l$ so that it can be replaced in that chain by a
unit at~$(k,l)$. Therefore assume that $k=i'_t$ for some $t\in\set{r}$, so
that the chain $C^{(t)}$ contains units at~$(k,l+1)$ and at~$(i_t,l)$; since
$M_{k,l}=0$ one must have $i_t<k$. Now for any $j\leq{t}$,
definition~\binmovedef\defitemnr2 gives
$\sum_{i=i_j}^{k-1}M_{i,l}\leq\sum_{i=i_j}^{k-1}M_{i,l+1}$, and each of the
chains $\foldwith{C^{(#1)}},\ldots(j..t-1)$ contributes equally to both sides,
while $\smash{C^{(t)}}$ contributes only to the left hand side. Then the right
hand side either has a contribution not from any chain in~$S$, or one
from~$C^{(j-1)}$. The latter does not happen for~$j=0$, so there exists
some~$i<k$ with $M_{i,l+1}=1$ and $i\notin{I'}$; now put
$k^-=\max\setof{i<k}:M_{i,l+1}=1,i\notin{I'}\endset$ and
$s=\min\setof{j\leq{t}}:\smash{i'_j}>k^-\endset\leq{t}$. The argument given
shows that $i_s\leq{k^-}$, and that $i_j\leq{i'_{j-1}}$ for $s<j\leq{t}$.

Now if the unit at~$(k^-,l+1)$ does not occur in any chain of the $m$-chain,
then one can replace in~$C^{(s)}$ the unit at~$(i'_s,l+1)$ by one at
$(k^-,l+1)$, and for $s<j\leq{t}$ replace in $C^{(j)}$ the unit
at~$(i'_j,l+1)$ by one at $(i'_{j-1},l+1)$; the $m$-chain of size~$n$ in~$M$
so obtained is also an $m$-chain in~$M'$. Finally if some chain of the
$m$-chain does contain a unit at~$(k^-,l+1)$, the this chain does not lie
in~$S$, and one can ``cyclicly recombine'' that chain~$C$ with the chains
$C^{(j)}$ for $s\leq{j}\leq{t}$: any units of~$C$ in columns~$<l$, are
combined with the unit at~$(k,l+1)$ and any other units of~$C^{(t)}$ in
columns~$>l$ to form a new chain~$D$; the units in columns~$\leq{l}$
of~$C^{(s)}$ are combined with the unit at~$(k^-,l+1)$ and any other units
of~$C$ in columns~$>l$; and for $s<j\leq{t}$ the units in columns~$\leq{l}$
of~$C^{(j)}$ are combined with the units in columns~$>l+1$ of~$C^{(j-1)}$. We
have obtained an $m$-chain in~$M$ of size~$n$, in which the chain~$D$ contains
a unit at~$(k,l+1)$ that can be replaced by a unit at~$(k,l)$ to form an
$m$-chain in~$M'$ of size~$n$.
\QED

\midinsert
\centerline{\epsfysize=3cm \epsfbox{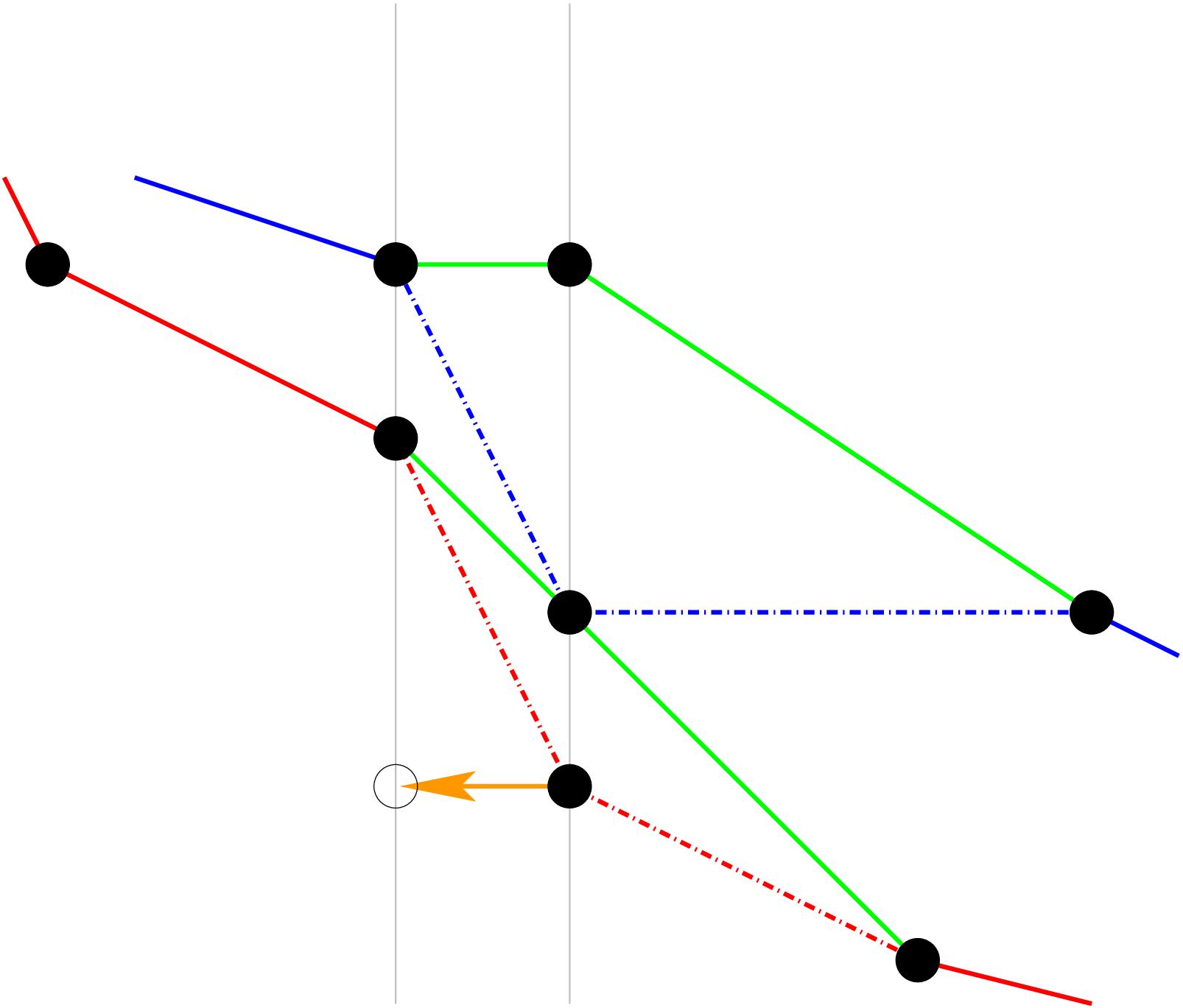}\hss 
            \epsfysize=3cm \epsfbox{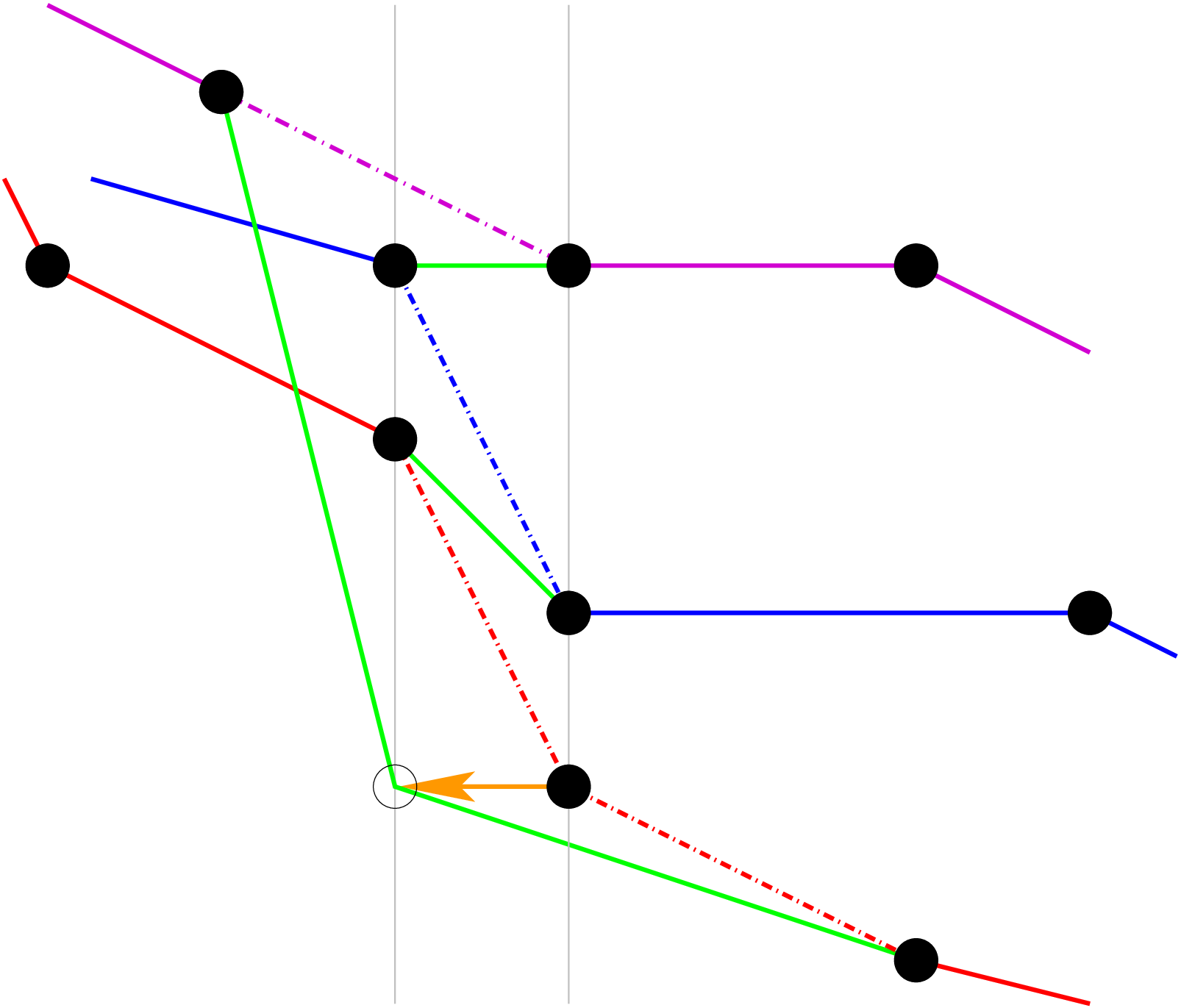}}
\centerpar{{\bf Figure~3.}
Schematic view of a shift between chains by several units,\goodbreak and of a
cyclic recombination of chains, for a horizontal move in a binary matrix.}
\endinsert

We now consider the integral case, where the cases of $m$-chains and
$m$-antichains must be treated separately, but moves in only one direction
need to be considered. As we shall see, these two cases are like ``thick''
versions of those treated in lemmas \binupinvlem~and~\binltinvlem,
respectively. The thickness refers to the fact that many units can share a
same position in the matrix; this seems to be a slight advantage in the case
of chains (which are now more compact), but a disadvantage in the case of
antichains (more effort is needed to tell them apart). To emphasise the
indicated analogy, we shall deal with upward transfers for $m$-chains, and
with leftward transfers for $m$-antichains.

\proclaim Lemma. \intchaininvlem
Let an integral matrix $M\in\Mat$ contain an $m$-chain of size~$n$ for some
$m,n\in\N$, and let $M'=\upm_k(M)$ for some $k\in\N$. Then $M'$ contains an
$m$-chain of size at least~$n$.

\proof
Let the transfer effected by~$\upm_k$ be from $M_{k+1,l}$ to~$M_{k,l}$. One
may assume that one of the chains~$C$ of the $m$-chain contains a unit at
$(k+1,l)$, for otherwise one can simply keep the same $m$-chain. If any chain
of the $m$-chain can accommodate a (new) unit at~$(k,l)$ while remaining a
chain, then an $m$-chain of size~$n$ in~$M'$ can simply be constructed by
removing a unit from~$C$ at~$(k,l+1)$ and adding a unit at~$(k,l)$ to the
chain in question. We shall show that one can always reach this situation by
replacing if necessary the given $m$-chain in~$M$ by another one of size at
least~$n$.

Let $l'=\max\setof{j}:C_{k+1,j}>0\endset\geq{l}$ be the column of the
rightmost unit of~$C$ in row~$k+1$. If $l'=l$ then $C$ itself can accommodate
a unit at~$(k,l)$. Otherwise, if none of the chains in the $m$-chain have any
units at~$(k,j)$ for $l\leq{j}<l'$ then we may replace $C$ by a chain obtained
by removing all its units at positions $(k+1,j+1)$ with $l\leq{j}<l'$, while
adding for those same values of~$j$ all $M_{k,j}$ units of~$M$ at~$(k,j)$.
The inequality $\sum_{j=l}^{\smash{l'}-1}(M_{k,j}-M_{k+1,j+1})\geq0$ from
definition~\intmovedef\defitemnr1 shows that the size of the new chain is at
least that of~$C$, so we have reduced to the case $l'=l$. If some of the
chains does have a unit at some $(k,j)$ with $l\leq{j}<l'$, then let $j_0$ be
the minimal such index~$j$, and $C'$ a chain of our $m$-chain with a unit
at~$(k,j_0)$. If $C'$ has another unit in row~$k$ to the left of the one
at~$(k,j_0)$, necessarily in a column~$<l$, then it can accommodate an
additional unit at~$(k,l)$. If $C'$ has no such unit, then we can replace the
chains $C,C'$ by a pair of new chains~$D,D'$ formed as follows: the chain $D$
contains the units of~$C$ in columns~$\leq{l}$, all $M_{k,j}$ units available
in~$M$ at $(k,j)$ for $l\leq{j}<j_0$, and the units of~$C'$ in rows~$\leq{k}$,
which are in columns~$\geq{j_0}$; the chain $D'$ contains the units of~$C'$ in
rows~$>k$, which are in columns~$\leq{j_0}$, and the units of~$C$ in
columns~$>j_0$. The inequality $\sum_{j=l}^{j_0-1}(M_{k,j}-M_{k+1,j+1})\geq0$
shows that replacing $C,C'$ by $D,D'$ gives an $m$-chain of size at least~$n$,
and $D$ can accommodate an additional unit at~$(k,l)$.
\QED

We note that in the final case we could have included any units of~$C$ at
$(k+1,j_0)$ into~$D'$, but the $m$-chain will be large enough without them; in
fact we could have done without the units of~$C'$ at that position as well.
But there need not be any such units, and there will not be if $n$ is maximal.

\proclaim Lemma.
Let an integral matrix $M\in\Mat$ contain an $m$-antichain of size~$n$ for some
$m,n\in\N$, and let $M'=\ltm_l(M)$ for some $l\in\N$. Then $M'$ also contains
an $m$-antichain of size~$n$.

\proof
Let the application of $\ltm_j$ result in a transfer from $M_{k,l+1}$
to~$M_{k,l}$. Consider the subset~$S$ of antichains in the $m$-antichain the
have a unit both in column $l$ and in column~$l+1$, and put $r=\Card{S}\geq0$.
Let $I,I'$ be the $r$-multisets on~$\N$ such that the units in column~$l$ of
the chains in~$S$ form the multiset $\multisetof(i,l):i\in{I}\endset$, and
that the units in column~$l+1$ of the chains in~$S$ form the multiset
$\multisetof(i,l+1):i\in{I'}\endset$; write $I=\multiset{\li(i0..r-1)}$ and
$I'=\multiset{\li({i'}0..r-1)}$ with $\fold\leq(i0..r-1)$ and
$\fold\leq({i'}0..r-1)$. Since each chain of~$S$ contributes an index $i$
to~$I$ and an index $i'$ to~$I'$ with $i\leq{i'}$, it follows that
$i_j\leq{i'_j}$ for all $j\in\set{r}$. We can then cut each chain in~$S$ into
a part in columns~$\leq{l}$ and a part in columns~$>l$, and recombine all the
parts to form an $r$-antichain $(A^{(j)})_{j\in\set{r}}$ with the same sum
as~$S$, in which the antichain~$A^{(j)}$ has units both at $(i_j,l)$ and at
$(i'_j,l+1)$. We replace $S$ by this $r$-antichain, and set
$A=\sum_{j\in\set{r}}A^{(j)}$.

If $A_{k,l+1}<M_{k,l+1}$ then either our $m$-antichain is already one in~$M'$
(if it does not use all units of~$M$ at~$(k,l+1)$), or it contains an
antichain with a unit at~$(k,l+1)$ and no units in column~$l$, in which we can
then replace its unit at~$(k,l+1)$ by one at $(k,l)$; in both cases we have
found an $m$-antichain of size~$n$ in~$M'$. Now assume that
$A_{k,l+1}=M_{k,l+1}$ and let $t=\max\setof{j}:\smash{i'_j}=k\endset$ be the
index of the last contribution to $A_{k,l+1}$.
Definition~\intmovedef\defitemnr2 gives $\smash{\sum_{i=i_j}^{k-1}}
(M_{i+1,l+1}-M_{i,l})\geq1$ for all $j\leq{t}$, and the antichains
$\foldwith{A^{(#1)}},\ldots(j..t)$ give null contributions to the left hand
side; therefore either $M_{i,l+1}>A_{i,l+1}$ holds for some~$i$ with
$i_j<i<k$, or $j>0$ and $i'_{j-1}>i_j$. Then we can set
$k^-=\max\setof{i<k}:A_{i,l+1}<M_{i,l+1}\endset$ and
$s=\min\setof{j}:i'_j\geq{k^-}\endset\leq{t}$, and we have $i_s<k^-$ and
$i_j<\smash{i'_{j-1}}$ whenever $s<j\leq{t}$.

If our $m$-antichain does not use all $M_{k^-,l+1}$ available units
at~$(k^-,l+1)$, then one can add such a unit to the antichain~$A^{(s)}$ to
replace the one at $(i'_s,l+1)$, and then replace for all $s<j\leq{t}$ the
unit at $(i'_j,l+1)$ in~$A^{(j)}$ by one at~$(i'_{j-1},l+1)$ (some such
replacements may give back the same antichain); thus our $m$-antichain is
replaced by another $m$-antichain of size~$n$ in~$M$ that does not use all
units of~$M$ at~$(k,l+1)$, and is therefore also an $m$-antichain in~$M'$. If
our $m$-antichain does use all units at~$(k^-,l+1)$, then it contains at least
one antichain~$A$ with such a unit and without any unit in column~$l$. Any
units of~$A$ in columns~$<l$ can be combined with the units of~$A^{(t)}$ in
columns~$>l$ to from a new antichain~$B$, while the units of~$A$ in
columns~$>l$ can be combined with the units of~$A^{(s)}$ in columns~$\leq{l}$;
furthermore for $s<j\leq{t}$ the units of~$A^{(j)}$ in columns~$\leq{l}$ can
be combined with the units of~$A^{(j-1)}$ in columns~$>l$. Thus we have
reorganised our $m$-antichain into another one of the same size in~$M$,
containing an antichain~$B$ with a unit at~$(k,l+1)$ that can be replaced by
one at~$(k,l)$.
\QED

We are now set to prove theorem \isinvarthm; the proof is just a
formalisation of the argument given before our sequence of lemmas.

\proofof{theorem \isinvarthm}
If $M,M'$ are a binary or integral $k\times{l}$ matrices, and
$M^\diamond,M^{\prime\diamond}$ are obtained from them by a rotation by a half
turn (mapping the origin to $(k-1,l-1)$), then the relation $M'=\upm_i(M)$
with $i\in\set{k-1}$ is equivalent to
$M^{\prime\diamond}=\dnm_{k-2-i}(M^\diamond)$, while $M'=\ltm_j(M)$ with
$j\in\set{l-1}$ is equivalent to
$M^{\prime\diamond}=\rtm_{l-2-j}(M^\diamond)$. Moreover, the rotation maps
chains to chains and antichains to antichains. Then for each of the four
preceding lemmas, one can derive the corresponding statement for the crystal
operation in the opposite direction by applying the lemma to~$M^\diamond$.
Therefore each lemma may be strengthened to ``$M$ contains an $m$-(anti)chain
of size~$n$ if and only if $M'$ contains one'' (lemma~\intchaininvlem\ speaks
of an $m$-chain of size at least~$n$, but this is only to make its proof
easier to read; such an $m$-chain exists if and only if an $m$-chain of
size~$n$ exists). Furthermore, applying rotation by a quarter turn to the two
lemmas for binary matrices gives the corresponding statements with $m$-chains
replaced by $m$-antichains, and with horizontal and vertical moves
interchanged, while the two lemmas for integral matrices give the
corresponding statements with horizontal and vertical moves interchanged by
transposition. We may conclude both in the binary and the integral case that
if $M$~and~$M'$ are related by any crystal operation, then $M$ contains an
$m$-chain of size~$n$ if and only if $M'$ does so, and the same for
$m$-antichains. The theorem now follows from the fact that it holds for normal
forms, as we had already observed.
\QED

As the proofs above indicate, whenever $M$~and~$M'$ are related by a crystal
operation, one can obtain a $m$-(anti)chain in $M'$ of maximal size by
choosing an appropriate such $m$-(anti)chain in~$M$ and either leaving its sum
unchanged or moving one unit between the same squares as in the application of
the crystal operation. Nevertheless, due to the fact that the moving unit may
migrate from one (anti)chain to another and to the necessity of switching from
one $m$-(anti)chain to another with a different sum between successive crystal
operations, the composition of the maximal $m$-(anti)chains in~$M$ and the
relations between them may be different from those in the normal form
$\diagram\\$ or $\diagon\\$ where $\\=\is{M}$. Notably the distribution of
sizes of individual (anti)chains may be necessarily different: while in normal
forms such $m$-chains always decompose into chains of sizes $\li(\\0..m-1)$,
an $m$-chain of maximal size for $m>1$ need in general not even contain any
individual chain of maximal size~$\\_0$. Here are examples of such a binary
and integral matrix, both with $\\=(5,3,1)$, whose unique $2$-chain of
size~$8$ only decomposes into two chains of size~$4$, and fail to contain the
unique chain of size~$5$:
$$
  \pmatrix{1&1&0&1&1\cr0&0&1&0&0\cr1&1&0&1&1\cr}\in\bMat
\ttext{and}
  \pmatrix{2&0&0&2&0\cr0&0&1&0&0\cr0&2&0&0&2\cr}\in\Mat
\nn
$$
(both examples were directly inspired by the one given in~\ref{Greene
Schensted, \Sec3}). If one would trace the evolution of maximal $m$-chains
from their normal form back to these matrices for~$m=1$ and $m=2$, using the
appropriate lemmas above, then it is clear that at some point the two cases
must diverge; there is therefore no hope of improving those lemmas so that
they keep track of complete decompositions into (anti)chains from which
$m$-(anti)chains of appropriate sizes could be extracted for all $m\in\N$. The
only thing one can say in general about the sizes of individual chains in a
$m$-chain of maximal size is that they lie in the range from $\\_{m-1}$ to
$\\_0$ inclusive: they obviously cannot exceed the maximal size~$\\_0$ of
individual chains, and if any had size strictly less than~$\\_{m-1}$, then
could remove that chain and be left with an $(m-1)$-chain of size exceeding
$\sum_{i\in\set{m-1}}\\_i$, which is impossible.

The property exhibited by the example allows another negative conclusion to be
drawn: if $\bMat$ and $\Mat$ are considered as partially ordered sets
by~`$\subset$' (componentwise comparison) and $\Part$ is considered so by
inclusion of diagrams, then the implicit shape \emph{does not define a poset
morphism}. After all, the unit in the middle of the matrix is contained in any
$1$-chain of maximal size but not in all $2$-chains of maximal size (in fact
there is only one such $m$-chain in either case), so that suppression of this
unit will decrease the initial part of $\is{M}$ without decreasing the sum of
its first two parts; consequently the diagram of $\is{M}=(5,3,1)$ does not
contain that of $\is{M'}=(4,4)$ where $M'$ is obtained from~$M$ by replacing
the middle entry by~$0$.

To prove that a given $m$-chain has maximal size it is not necessary to find
all $m$-chains, but it suffices to produce any $\\_m$-antichain of maximal
size. This is so because of the following general property, in which we define
the ``union'' of two matrices as $C\union{A} =C+A-C\thru{A}
=(\max(C_{i,j},A_{i,j}))_{i,j\in\N}$:

\proclaim Proposition.
Let $M$ be a binary or integral matrix, and $\\=\is{M}$. If for some $k,l>0$
and some sum~$C$ of a $k$-chain in~$M$ and a sum~$A$ of an $l$-antichain
in~$M$ one has $|C\thru{A}|=kl$ and $M=C\union{A}$, then the $k$-chain and the
$l$-antichain are both of maximal size, and $\\_k\leq{l}\leq\\_{k-1}$, in
other words $(k-1,l-1)\in\set\\$ and $(k,l)\notin\set\\$. Conversely for any
$k,l$ satisfying $\\_k\leq{l}\leq\\_{k-1}$, the sum~$C$ of any $k$-chain of
maximal size in~$M$ and the sum~$A$ of any $l$-antichain of maximal size
in~$M$ satisfy $|C\thru{A}|=kl$ and $M=C\union{A}$.

\proof
Since for any sum~$C'$ of a $k$-chain in~$M$ and a sum~$A'$ of an
$l$-antichain in~$M$ one has $|C'\thru{A'}|\leq{kl}$, and therefore
$|C'\union{A'}|\geq|C'|+|A'|-kl$, while on the other hand
$M\supset{C'}\union{A'}$ implies $|C'\union{A'}|\leq|M|=|C|+|A|-kl$, it is
clear that in the first part $|C|$ and $|A|$ are maximal among the sizes of
$k$-chains and of $l$-antichains in~$M$, respectively. Now let $C$ be the sum
of any $k$-chain in~$M$ of maximal size, so $|C|=\sum_{i\in\set{m}}\\_i$, and
let $A$ be the sum of any $l$-antichain in~$M$ of maximal size, so
$|A|=\sum_{j\in\set{l}}\\\tr_j$; put $d=kl-|C\thru{A}|$ and
$d'=|M-C\union{A}|$ (so in the case just considered one has $d=d'=0$). The
expression giving $|C|$ can be interpreted as $\Card\set\\|_{\set{k}\times\N}$
(where $\set\\|_X$ abbreviates $\set\\\thru{X}$), the number of squares in the
topmost $k$ rows of~$\set\\$. Similarly the expression for~$|A|$ can be
interpreted as $\Card\set\\|_{\N\times\set{l}}$. Therefore the equation
$|C|+|A|=|C\thru{A}|+|C\union{A}|$ can be written as
$$
  \Card\set\\|_{\set{k}\times\N}+\Card\set\\|_{\N\times\set{l}}
 =\Card(\set{k}\times\set{l})-d+\Card\set\\-d'
.\nn
$$
Now the only squares in~$\N^2$ with different contributions to the sets
counted on the left and on the right are those in
$D=\set{k}\times\set{l}\setminus\set\\$ or in
$D'=\set\\|_{(\N-\set{k})\times(\N-\set{l})}$, and in both cases they
contribute~$1$ to the right hand side and~$0$ to the left hand side. Therefore
$d+d'=\Card{D}+\Card{D'}$, and we may conclude that $d=d'=0$ if and only if
$D=D'=\emptyset$, which is equivalent to $(k-1,l-1)\in\set\\$ and
$(k,l)\notin\set\\$. The proposition follows.
\QED

Note that although in the case of normal forms one even has $d=\Card{D}$ and
$d'=\Card{D'}$, this is not true in general: in the examples given the sum~$C$
of the $2$-chain of maximal size is also the sum of a $4$-antichain, so in
this case with $(k,l)=(4,2)$ and $A=C$ one has $d=0$ but $D=\sing{(3,1)}$.

While the proposition allows an easy certification of the maximality of a
$k$-chain and an $l$-antichain if an appropriate pair is given, it does not
provide a procedure for obtaining such a pair. The best algorithmic method we
know of to determine such witnesses for (partial information about) the
implicit shape of a matrix~$M$ is to apply crystal operations to~$M$ until
reaching a form in which such witnesses are evident (for instance a matrix
encoding a straight tableau), and then trace back one's steps using the
descriptions in the proofs of the lemmas of this section to find corresponding
witnesses in~$M$. This is somewhat laborious to do by hand, but still feasible
for the matrices of our running examples. Here it is easiest to reduce to the
form that encodes the straight tableau~$S$ of~(\Sdef), namely is the
matrix~$Q$ of~(\PQbindiagram) in the binary case, which is reached after $12$
leftward moves, or the matrix~$P$ of~(\PQintdiagram) in the integral case,
which is reached after $13$ upward unit transfers. After this reduction we can
find maximal $k$-chains and $l$-antichains for every $k$ and~$l$ by taking in
the binary case the $k$ topmost units (or all units if there are less
than~$k$) in every column of~$Q$ respectively the units in the $l$ leftmost
rows of~$Q$, or in the integral case the $k$ topmost rows of~$P$ respectively
the $l$ leftmost units (or all available units) in every row of~$P$ (in both
cases these correspond to the $k$ topmost rows and the $l$ leftmost columns of
the display of~$S$). As an example, tracing that maximal $3$-chain of
size~$13$ and that maximal $4$-antichain of size~$16$ back to~$M$, we find
in~$M$ the $3$-chain coloured red, green and blue in the following displays
\bigdisplay
\let\r=\red \let\g=\green \let\b=\blue
\pmatrix{  0&  1&  0&  0&  1&  0&  0&  0&  0\cr
         \r1&\r1&\r1&  0&  0&  0&  0&  0&  0\cr
	 \g1&  0&  1&  0&  0&  1&  0&  0&  0\cr
	   0&\g1&  0&\r1&  0&  0&  0&  0&  0\cr
	   0&  0&\g1&\g1&\r1&  0&\r1&  0&  0\cr
	 \b1&  0&  0&  0&\g1&  0&  0&\r1&\r1\cr
	   0&\b1&\b1&\b1&\b1&\g1&\g1&\g1&  0\cr}
\qquad
\pmatrix{  1&  0&  1&  0&\r1&\r2&  0\cr
           1&\r1&  0&\r1&\r1&  0&\g3\cr
	   0&\r2&  1&  0&\g1&\g1&  0\cr
	   0&  0&\g1&\g1&\g1&  0&\b2\cr
	   0&  0&  0&  0&  0&\b1&\b2\cr}
,\nn
$$
while the two $4$-antichains one obtains in both cases can be described as
follows by the sets of positions where their individual antichains have units:
in the binary case these are $\{(0,1),(1,0),(2,0),(5,0)\}$,
$\{(0,3),(1,1),(3,1),(6,1)\}$, $\{(1,2),(2,2),(4,2),(6,2)\}$ and
$\{(2,5),(3,3),(4,3),(6,3)\}$, and in the integral case these sets are
$\{(1,0),(2,1),(3,3),(4,6)\}$, $\{(0,0),(2,1),(3,5),(4,6)\}$,
$\{(1,1),(2,2),(3,3),(4,5)\}$ and $\{(0,2),(1,3),(2,4),(3,6)\}$. The most
complicated of these computations is for the $4$-antichain in the integral
case, where cases that require recombination of chains occur frequently; this
explains the fact that the result for this case are rather tangled. There is
also some liberty in interpreting the transitions prescribed in the proofs of
the lemmas, since more than one antichain can contain a given unit; one could
also perform some intermediate untangling of antichains, which would affect
the final result as well.

It is interesting to compare the results above with the original approach
of~\ref{Greene Schensted}. There the roles of chains and antichains are played
by increasing and decreasing subsequences extracted from a given finite
sequence of numbers. In fact only sequences representing permutations are
considered, but that restriction can be lifted if for chains one allows weakly
increasing sequences while for antichains one requires strictly decreasing
sequences; everything generalises immediately by treating equal entries as if
they were increasing from left to right (a price to pay is that the theory
becomes asymmetric, and the opposite convention would have been equally
feasible, but we shall stick to this more usual convention). One can associate
words to matrices in a way that is compatible with this interpretation of
chains and antichains: for a binary matrix one takes concatenation of the
strictly decreasing words whose weights are given by the columns of the matrix
taken from left to right, while for an integral matrix one takes the
concatenation of the weakly increasing words whose weights are given by the
rows of the matrix from bottom to top. While this correspondence seems biased
towards words that have this particular composition, one can in fact model
words bijectively by matrices, by restricting the latter to those whose column
sums (in the binary case) or row sums (in the integral case) are all~$1$ up to
some point beyond which they are all~$0$: the mentioned concatenations then
combine words consisting of a single letter.

The invariance of the maximal size of $m$-chains or of $m$-antichains is then
shown for a certain type of transposition of letters in the word, which
generates an equivalence relation on words that is known by \ref{Knuth PJM,
Theorem~6} to characterise words with the same insertion tableau in the
\caps{RSK}~correspondence. To be precise, the interchange of a pair of distinct
adjacent letters is permitted in the presence of a letter~$l$ directly to the
left or to the right of the pair, whose value is either strictly in between
those of the letters being interchanged, or (in keeping with the
interpretation of equal letters) if either $l$ is to the left of the pair and
its value equal to the larger of the letters being interchanged, or $l$ is to
the right of the pair and its value equal to the smaller one. If $a,b,c$
designate letters with $a<b<c$, these transpositions can therefore be
represented as transformations $bac\leftrightarrow{bca}$,
$acb\leftrightarrow{cab}$, $bab\leftrightarrow{bba}$ or
$aba\leftrightarrow{baa}$ of three consecutive letters; we shall call them
elementary Knuth transformations.

Modelling words by matrices, the fact that elementary Knuth transformations
leave unchanged the \caps{RSK}~insertion tableau for the word is reflected by
the fact that they can be realised by horizontal moves within the three
columns (in the binary case) or vertical transfers between the three rows (in
the integral case) corresponding to the letters concerned. We illustrate this
in the binary case for transition $bac\to{bca}$:
\bigdisplay
\def\mm(#1,#2,#3,#4){{\def~{&\kern-.4em}\pmatrix{#1\cr
\noalign{\smallskip}#2\cr\noalign{\medskip}#3\cr
#4\cr\noalign{\kern-\baselineskip}}}}
\mm(~1,1,~~1,b~a~c)\buildrel\ltm_i\over\longrightarrow
\mm(1,1,~~1,ba~~c)\buildrel\ltm_{i+1}\over\longrightarrow
\mm(1,1,~1~,ba~c)\buildrel\rtm_i\over\longrightarrow
\mm(~1,1,~1~,b~ca)\buildrel\rtm_{i+1}\over\longrightarrow
\mm(~~1,1,~1,b~c~a)
\advance \belowdisplayskip by .4\baselineskip
,\nn
$$
and similarly for a transition $aba\to{baa}$:
\bigdisplay
\def\mm(#1,#2,#3){{\def~{&\kern-.4em}\pmatrix{#1\cr
\noalign{\medskip}#2\cr
#3\cr\noalign{\kern-\baselineskip}}}}
\mm(1~~1,~1,a~b~a)\buildrel\ltm_{i+1}\over\longrightarrow
\mm(1~1~,~1,a~ba)\buildrel\ltm_i\over\longrightarrow
\mm(1~1~,1,ba~a)\buildrel\rtm_{i+1}\over\longrightarrow
\mm(1~~1,1,ba~~a)\buildrel\rtm_i\over\longrightarrow
\mm(~1~1,1,b~a~a)
\advance \belowdisplayskip by .4\baselineskip
.\nn
$$
The word corresponding to the matrix changes only during the application
of~$\rtm_i$ in the first case, and that of~$\ltm_i$ in the second. The
integral case is similar; we show just these essential transitions for
$bac\to{bca}$ and for $bab\to{bba}$ (keeping in mind that we must read from
bottom to top):
\bigdisplay
\def\mm(#1:#2,#3:#4){{\def~{&}\kern1.4 em \pmatrix{\llap{$#2$\quad}#1\cr
\llap{$#4$\quad}#3\cr\noalign{\smallskip}}}}
\mm(1~~1:ac,~1~:b\ )\buildrel\dnm_j\over\longrightarrow\mm(1~~:a,~1~1:bc\ )
\ttext{and}
\mm(1~1:ab,~1:b\ )\buildrel\dnm_j\over\longrightarrow\mm(1~:a,~2:bb\ )
.\nn
$$
It is clear that if these transformations are made inside a larger matrix~$M$,
then the one of the matrices $P,Q$ associated to~$M$ in theorem~\decompthm\
that is obtained by exhausting crystal operations in a direction parallel to
those used here does not change; this matrix is $Q$ in the binary case and $P$
in the integral case. If $M$ encodes a skew tableau~$T$ with at most one
square in any row or column, then that matrix encodes both in the binary and
the integral case the rectification~$S$ of~$T$, which is also the
\caps{RSK}~insertion tableau for the word read off from the display of~$T$
from bottom left to top right.

It is less obvious that by repeating elementary Knuth transformations, one can
transform the word corresponding to~$M$ into the word corresponding to the
matrix that encodes~$S$, which word can be read off from the display of~$S$ in
the reverse Kanji (binary case) of reverse Semitic (integral case) order. The
proof of \ref{Knuth PJM, Theorem~6} is based on the fact that Schensted
insertion can be emulated using elementary Knuth transformations: if
$T'=T\leftarrow{x}$ is the result of row-inserting an entry~$x$ into a
semistandard Young tableau~$T$, and $t$~and~$t'$ are the words obtained by
reading the displays of $T$~and~$T'$ in reverse Semitic order, then $t'$ can
be obtained from $tx$ by a sequence of elementary Knuth transformations (see
also \ref{Fulton, \Sec2.1}; the transpositions proceed in a single sweep from
right to left). A similar result can easily be shown for column insertion
$x\to{T}$, reading the tableaux in reverse Kanji order.

Since we have seen that elementary Knuth transformations can be realised by
crystal operations, one may ask if one can conversely derive, from the
exhaustion of crystal operations that transforms~$M$ into the matrix
encoding~$S$, a sequence of elementary Knuth transformations that realises the
transformation of the corresponding words. Of course, since we remain inside
one Knuth equivalence class, successive words obtained during the exhaustion
of the crystal operations can be linked by elementary Knuth transformations,
but the number of them required at each step turns out to be surprisingly
large. It is shown in \ref{Schutzenberger CRGS, 2.3.bis} how the change to the
reverse Semitic reading of a tableau (with a hole) caused by the upward move
of one entry occurring during an inward jeu de taquin slide can be realised by
a sequence of elementary Knuth transformations (the goal of that statement was
in fact to show that only the directions $bca\to{bac}$, $bba\to{bab}$,
$acb\to{cab}$, or $aba\to{baa}$ are needed here; it should be noted that some
transpositions might also allow an interpretation as a transformation in the
opposite sense within another triplet). That case corresponds precisely to an
upward unit transfer in an integral matrix. It is not quite clear how many
transpositions are produced in the proof of the cited result, although it
seems to be larger than necessary. We shall sharpen that result by providing a
count of the transpositions used; experimental evidence strongly suggests that
the number we give is in fact minimal.

\proclaim Proposition.
Let $M\in\Mat$ allow an upward unit transfer from $(i_0+1,j_0)$ to $(i_0,j_0)$
resulting in $M'=\upm_{i_0}(M)$, and let the words $w,w'$ be the
concatenations of the weakly increasing words with weights given by the rows
of $M$ and of~$M'$, respectively, ordered from bottom to top. Then $w$ can be
transformed into $w'$ by a succession of $n^2$ elementary Knuth
transformations of type $bca\to{bac}$, $bba\to{bab}$, $acb\to{cab}$, or
$aba\to{baa}$, where $n=\sum_{j<j_0}M_{i_0,j}+\sum_{j>j_0}M_{i_0+1,j}$.

\proof
For brevity we take the designation $bca\to{bac}$ to include type
$bba\to{bab}$ as well, and $acb\to{cab}$ is similarly taken to include type
$aba\to{baa}$. Put $k=\sum_{j<j_0}M_{i_0,j}$ and $l=\sum_{j>j_0}M_{i_0+1,j}$,
so that $k+l=n$. Write the rightmost $n+1$ letters of the subword
corresponding to row~$i_0+1$ of~$M$ as $\li(c1..k),j_0,\li(d1..l)$ and the
leftmost $n$ letters of the subword corresponding to row~$i_0$ of~$M$ as
$\li(a1..k),\li(b1..l)$. The inequalities of definition~\intmovedef\defitemnr1
guarantee that the letters $\li(c1..k)$ and $\li(b1..l)$ exist in those rows;
the concatenation of these weekly increasing words forms a
subword~$v=\fold(c1..k)j_0\fold(d1..l)\fold(a1..k)\fold(b1..l)$ of~$w$ in
which all transpositions will take place, and we shall ignore all other
letters. The same inequalities also imply that $a_i<c_i$ for $1\leq{i}\leq{k}$
and $b_i<d_i$ for $1\leq{i}\leq{l}$, and if $l>0$ one has $j_0\leq{b_1}$; thus
$v$ corresponds to the reverse Semitic reading of a tableau at an intermediate
stage of a jeu de taquin slide, when moving~$j_0$ upward in the situation
$$
  \Young(a_1,\cdot,\cdot,a_k,   ,b_1,\cdot,\cdot,b_k|
         c_1,\cdot,\cdot,c_k,j_0,d_1,\cdot,\cdot,d_k)
\,.\nn
$$
We must transformation $v$ into
$v'=\fold(c1..k)\fold(d1..l)\fold(a1..k)j_0\fold(b1..l)$, which we shall do
explicitly by the following choreography. We proceed in $n$ sweeps of $n$
elementary Knuth transformations each, starting the successive sweeps with the
triplets whose rightmost letters are successively $\li(a1..k),\li(b1..l)$, and
proceeding from right to left within each sweep. During sweep~$i$ with
$1\leq{i}\leq{k}$, the letter $a_i$ is first shifted to the left $n-i+1$
places by transformations of type $bca\to{bac}$, after which it ends up to the
right of the letter~$c_i$; that letter is then shifted to the left $i-1$
places across the letters $\li(a1..i-1)$ by transformations of type
$acb\to{cab}$, which results in a word $\fold(c1..i)\fold(a1..i)\fold(ci+1..k)
j_0\fold(d1..l)\fold(ai+1..k)\fold(b1..l)$. After these $k$ sweeps one
therefore has the word $\fold(c1..k)\fold(a1..k) j_0\fold(d1..l)\fold(b1..l)$.
During sweep~$k+i$ with $1\leq{i}\leq{l}$, the letter $b_i$ is first shifted
to the left $l-i$ places across the letters $\li(di+1..l)$ by transformations
of type $bca\to{bac}$ in which it has the role of~$a$, after which it finds
itself to the right of the letter~$d_i$; that letter is then shifted to the
left $k+i$ places by transformations of type $acb\to{cab}$, which results in a
word $\fold(c1..k)\fold(d1..i)\fold(a1..k)
j_0\fold(b1..i)\fold(di+1..l)\fold(bi+1..l)$. After these sweeps we have the
word~$v'$.
\QED

For a leftward move in a binary matrix from $(i_0,j_0+1)$ to $(i_0,j_0)$,
there is a similar realisation of the change in the corresponding words, which
is the change in the reverse Kanji reading of a tableau that occurs during a
leftward step in a jeu de taquin slide, using $n^2$ elementary Knuth
transformations of the same types, where
$n=\sum_{i<i_0}M_{i,j_0}+\sum_{i>i_0}M_{i,j_0+1}$. The proof is the same,
interchanging left and right in the word, and with the implied interchange of
weak and strict inequalities.

The proposition indicates that the number of elementary Knuth transformations
corresponding to the exhaustion of upward crystal operations can be
considerably larger than the number of such operations; in our example we
needed $13$ upward unit transfers, which correspond to $150$ elementary Knuth
transformations. This can be improved a bit when there are successive
transfers between the same pair of rows: if there are $k$ such transfers then
all involve the same value of~$n$, and it appears that by an adaptation of the
transposition procedure one can achieve the combined transformation of the
word in $n(n+k-1)$ elementary Knuth transformations instead of $kn^2$; in the
example this allows the total number of steps to be reduced to~$122$. This is
still worse than the number needed for the traditional emulation of the
\caps{RSK}~insertion, although it is of the same order of magnitude: a
straightforward application of that procedure to the word corresponding to~$M$
gives a sequence of $120$ transpositions, which can be reduced to~$88$ by
observing that the first $32$ transpositions bring back the word to its
original form, due to the circumstance that the word corresponding to the
three bottom rows of~$M$ is already the reverse Semitic reading of a straight
tableau. It should be noted that our matrix~$M$ was chosen as the integral
encoding of a skew semistandard tableau that is fairly close to its
rectification, whence the numbers above are relatively small. These data
suggest that for the purpose of finding maximal $m$-(anti)chains in a binary
or integral matrix, there is some advantage in using the transformations
described in our lemmas rather than the (somewhat simpler) ones corresponding
to elementary Knuth transformations.

Another approach to finding maximal $m$-(anti)chains would be to try to
augment the inductive computation of $\is{M}$ via growth diagrams with
additional data, keeping track of maximal $m$-(anti)chains as witnesses for
the implicit shape of submatrices. It appears however that this can only be
done for~$m=1$.

 \par

\subsecno=-1 

\newsection Remaining proofs.
\labelsec\proofsec

We shall finally provide, for completeness, proofs of statements for which
they were so far omitted. The first such statement was the fact that crystal
operations can be defined in such a way for any semistandard tableau~$T$ that
they correspond both to vertical crystal operations on its binary
encoding~$M$, and to horizontal crystal operations on it integral
encoding~$N$.

\proofof{proposition~\bininteqprop}
We know that if either $\upm_m(M)$ or~$\ltm_m(N)$ is defined, then like
$M$~or~$N$ it satisfies the tableau condition for~$\\/\kappa$ (propositions
\bininvarprop~and~\intinvarprop), and is therefore the binary or integral
encoding of some tableau~$T'$. Moreover $T'$ can easily be seen to differ
from~$T$ only in the partition~$\\^{(m+1)}$, for which some square~$(k,l)$ has
been added to its diagram to form the diagram of the corresponding partition
in~$T'$ (in terms of the display of~$T$, an entry~$m+1$ in that square has
been changed to~$m$). We have to show that if the other type of encoding is
taken of $T$~and of~$T'$, then the matrices obtained are also related by the
appropriate raising operation. By interchanging the roles of~$T$ and of~$T'$,
this will also prove the statement about $\dnm_m$ and $\rtm_m$.

Suppose first that the binary encodings of $T$~and~$T'$ are
$M$~and~$\upm_m(M)$, and let their integral encodings be $N,N'$, respectively.
Then with $k,l$ as above, we know that the bits $M_{m,l}=0$ and $M_{m+1,l}=1$
are interchangeable in~$M$, and we must show that a unit transfer is possible
from $N_{k,m+1}$ to~$N_{k,m}$ in~$N$. Consider the left member
$\sum_{i=k}^{\smash{k'}-1}(N_{i,m}-N_{i+1,m+1})$ of the second inequality in
definition~\intmovedef\defitemnr2, which we must show to be nonnegative for
any~$k'>k$. It counts the number of squares of the skew diagram
$\set{\\^{(m+1)}/\\^{(m)}}$ that lie in rows~$i\in\set{k'}-\set{k}$, minus the
number of squares of $\set{\\^{(m+2)}/\\^{(m+1)}}$ in rows
$i+1\in\set{k'+1}-\set{k+1}$. Now these squares can also be described as those
that lie in columns $j\in\set{l}-\set{l'}$ where
$l'=\\^{(m+1)}_{k'}\leq{l}=\\^{(m+1)}_k$, since squares of
$\set{\\^{(m+1)}/\\^{(m)}}$ in row~$k$ are in columns~$<l$ while squares of
$\set{\\^{(m+2)}/\\^{(m+1)}}$ in row~$k$ are in columns~$\geq{l}$, and
similarly for $k'$ and~$l'$. It follows that
$\sum_{i\in\set{k'}-\set{k}}(N_{i,m}-N_{i+1,m+1})
=\sum_{j\in\set{l}-\set{l'}}(M_{m,j}-M_{m+1,j})$ which is nonnegative by
hypothesis. For the first inequality in definition~\intmovedef\defitemnr2 the
argument is similar. Here it is easiest to take the version of that inequality
for $N'$, which is equivalent to the one for~$N$: for any $k'<k$ one puts
$l'=\\^{(m+1)}_{k'}$ (which satisfies $l'\geq{l}+1$), and one shows that
$\sum_{i\in\set{k}-\set{k'}}(N'_{i+1,m+1}-N'_{i,m})
=\sum_{j\in\set{l'}-\set{l+1}}(M_{m+1,j}-M_{m,j})$, which is nonnegative. This
concludes the proof that $N'=\ltm_m(N)$.

Now suppose that the integral encodings of $T$~and~$T'$ are
$N$~and~$\ltm_m(N)$, and let their binary encodings be $M,M'$, respectively.
The argument we have just given shows that, among the inequalities in
definition~\binmovedef\defitemnr1 that must be shown to hold in order to
justify $M'=\upm_m(M)$, those for which $l'\in\setof\\^{(m+1)}_{k'}:k'>k\endset
\union\sing{l} \union\setof\\^{(m+1)}_{k'}-1:k'<k\endset$
are satisfied. This does not cover all inequalities, but the values in rows
$m,m+1$ of~$M$ that are incorporated into the summations on~$j$ between such
values of~$l'$ have a simple structure, as they account for squares of
$\set{\\^{(m+2)}/\\^{(m+1)}}$ in row~$k'$, and squares of
$\set{\\^{(m+1)}/\\^{(m)}}$ in row~$k'-1$, for some fixed~$k'$ depending on
the interval considered. In such an interval one has from left to right: a
possibly empty sequence of columns~$0\choose1$, a possibly empty sequence
either of columns~$0\choose0$ or~$1\choose1$, and a possibly empty sequence of
columns~$1\choose0$; similarly the columns incorporated into the partial sums
after the largest such value of~$l'$ can only be $0\choose1$ or~$0\choose0$.
It follows easily from this description that the required inequalities also
hold at the remaining values of~$l'$.
\QED

The next statement left to prove is that the involutions defined by crystal
operations on matrices satisfy the Coxeter relations of type~$A$. We did not
use this fact anywhere, but it does represent an interesting aspect of crystal
operations, so we present the proof, even if this turns out to require quite a
bit of work.

\proofof{lemma \actionlem}
This is a statement about the structure of crystals, so by
proposition~\bininteqprop\ it suffices to consider the involutions~$\hr_l$ in
the case of integral matrices. We have already seen that the $\hr_l$ are
involutions (the first Coxeter relation), and the fact that $\hr_j$
and~$\hr_l$ commute when $|l-j|\geq2$ (the second Coxeter relation) follows
from the fact that they operate on disjoint pairs of rows. Hence the only
non-obvious point is the third Coxeter relation, but here a direct proof would
be difficult, because for instance the effect of applying $\upm_k$ to~$M$ on
the operation of $\upm_{k+1}$ depends on~$M$ in a complicated way. Instead we
shall proceed by reducing the general case to one where a detailed analysis is
feasible. In general, if $M'$ is obtained from $M$ by applying an invertible
operation that commutes with $\hr_l$ and with~$\hr_{l+1}$, then the equation
$\hr_l\hr_{l+1}\hr_l(M)=\hr_{l+1}\hr_l\hr_{l+1}(M)$ holds if and only if
$\hr_l\hr_{l+1}\hr_l(M')=\hr_{l+1}\hr_l\hr_{l+1}(M')$ does; we shall use such
operations to reduce $M$ to a simpler form.

First of all we replace $M$ by $M_{\N\times(\set{l+3}-\set{l})}$ (the inverse
operation of adding back the columns outside the range $\sing{l,l+1,l+2}$
commutes with $\hr_l$ and~$\hr_{l+1}$, which only look at the columns in that
range), and we may then as well assume that~$k=0$. Secondly we may assume that
upward crystal operations applied to~$M$, which commute with $\hr_l$
and~$\hr_{l+1}$, are exhausted. Thus we have reduced $M$ to the form
$$
  M=\pmatrix{a&b&c\cr0&d&e\cr0&0&f}\in\Mat
\ttex{with $a\geq{d}\geq{f}$ and $a+b\geq{d}+e$}
.\label(\Mabcdef)
$$
Although the proposition can be directly proved for this situation, as we
shall indicate below, one can take a shortcut, as in the proof
in~\ref{Lothaire}, by further reducing the situation using the fact that the
operations $\hr_j$ commute with cyclic permutations of the nonzero rows of the
matrix, as stated in lemma~\forward\cycliccommlem\ below. Doing so to the
matrix~$M$ of~(\Mabcdef), moving row~$0$ below the two following ones, one can
again exhaust upwards crystal operations. In fact using~$\upm_1$ all $a+b+c$
units that were cycled to row~$2$ can be transferred up to row~$1$, and then
$\upm_0$ can be applied $a+b+c-\min(d,c+f)$ times to obtain a matrix
$$
  M'=\pmatrix{a&b+d&c'\cr0&0&e'\cr0&0&0}
\ttex{where $e'=\min(d,c+f)$ and $c'=c+f-e'+e$}
\label(\Mabce)
$$
in which upwards transfers are once again exhausted. Thus the case of matrices
of the form~(\Mabcdef) has been reduced to the case of matrices of the
form~(\Mabce). The latter case is contained in the former, and we can repeat
the steps the we applied to further reduce to the
case of a one-line matrix $M''=(a~~b+d~~c+e+f)$. Now for such matrices one has
$\csum{M''}=M''_0$, and since we know that the effect of $\hr_i$ on
$\csum{M}$ is to transpose the parts $i$~and~$i+1$, it is clear that
$\hr_0\hr_1\hr_0(M'')=\hr_1\hr_0\hr_1(M'')=(c+e+f~~b+d~~a)$.
\QED

We shall now state formally the commutation of operations $\hr_j$ with cyclic
permutations of rows that we used. It is valid for binary matrices as well as
for integral ones, and it also has an analogue for~$\vr_i$; since the
statement is of independent interest we state it in its most general
form.

\proclaim Lemma. \cycliccommlem
Let $n,m\in\N$, $d\in\set{n}$ and $d'\in\set{m}$. Denote by $R$ either the set
of matrices $M\in\bMat$ or of matrices $M\in\Mat$, such that
$M=M_{\set{n}\times\set{m}}$, and by $\pi_{\rm{r}}$ and~$\pi_{\rm{c}}$
operations on~$R$ defined by the cyclic permutation of the rows
$\setof{M_i}:i\in\set{n}\endset$ by $d$ places (so
$\pi_{\rm{r}}(M)_{i,j}=M_{i+d\bmod{n},j}$ for all $M\in{R}$,
$(i,j)\in\set{n}\times\set{m}$), respectively by the cyclic permutation of the
columns $\setof{M\tr_j}:j\in\set{m}\endset$ by $d'$ places (so
$\pi_{\rm{c}}(M)_{i,j}=M_{i,j+d'\bmod{m}}$). Then for all $k,l\in\N$ and
$M\in{R}$ one has $\hr_l(\pi_{\rm{r}}(M))=\pi_{\rm{r}}(\hr_l(M))$ and
$\vr_k(\pi_{\rm{c}}(M))=\pi_{\rm{c}}(\vr_k(M))$.

\proof
By symmetry it will suffice to consider the statement concerning $\hr_l$.
Clearly only the columns $l,l+1$ that are affected by $\hr_l$ are of any
concern, so we may assume they are the only ones: $l=0$ and $m=2$. Note that
the effect of $\pi_{\rm{r}}$ and $\hr_0$ on row and column sums is determined
without knowing more details of~$M$, and these effects commute, so our only
worry is that both sides of the equation might give different matrices with
the same row and column sums. We shall start considering the special case of
integral $2\times2$ matrices (so $n=2$ as well) with $d=1$, so
that~$\pi_{\rm{r}}$ is the interchange of the two rows. We put
$\rsum{M}=(s,t)=\rsum{\hr_0(M)}$, then $M$ has the form $M=\Mtt(a,s-a,b,t-b)$.
One has $\delta=\csum{M}_1-\csum{M}_0 =s+t-2(a+b)$ in (\hreq), and the same
value applies to $\pi_{\rm{r}}(M)$. By symmetry between $M$ and
$\pi_{\rm{r}}(M)$ we may assume that $s\leq{t}$; we distinguish the cases with
$s\leq{a}+b\leq{t}$ from those where either $a+b<s$ or $t<a+b$. The condition
$s\leq{a}+b\leq{t}$ means on one hand that $\min(t-b,t-b-\delta)\geq{a}$,
so that in this case the transfer for computing $\hr_0(M)$ takes place
entirely in the bottom row, and on the other hand that
$s-a\leq\min(b,b-\delta)$ so that the transfer for computing
$\hr_0(\pi_{\rm{r}}(M))$ takes place entirely in the top row of
$\pi_{\rm{r}}(M) =\Mtt(b,t-b,a,s-a)$; as a consequence one has
$\hr_0(\pi_{\rm{r}}(M)) =\Mtt(b-\delta,t-b+\delta,a,s-a)
=\pi_{\rm{r}}(\hr_0(M))$. If instead one has $a+b<s\leq{t}$, then it follows
that $\delta>0$ and in fact $t-b>a>t-b-\delta$, whence we may write $\hr_0(M)
=(\ltm_0)^\delta(M) =(\ltm_0)^{s-(a+b)}((\ltm_0)^{t-(a+b)}(M))$, in which
expression $(\ltm_0)^{t-(a+b)}$ operates in the bottom row, after which the
factor $(\ltm_0)^{s-(a+b)}$ operates in the top row. The computation of
$\hr_0(\pi_{\rm{r}}(M))$ and can be performed similarly using
$s-a>b>s-a-\delta$; the two computations can be summarised as follows:
\bigdisplay
  \vcenter{\halign{ \hss$#{}$&\hss$#{}$\hss&\hss$#{}$\hss&\hss$#{}$\hss
                   &${}#$\hss \cr
    M=
  & \pmatrix{a&s-a\cr b&t-b} \buildrel (\ltm_0)^{t-(a+b)}\over\longrightarrow
  & \pmatrix{a&s-a\cr t-a&a} \buildrel (\ltm_0)^{s-(a+b)}\over\longrightarrow
  & \pmatrix{s-b&b\cr t-a&a}
  & =\hr_0(M)
  ,\cr\noalign{\nobreak\vskip\smallskipamount}
   \pi_{\rm{r}}(M)=
  & \pmatrix{b&t-b\cr a&s-a} \buildrel (\ltm_0)^{s-(a+b)}\over\longrightarrow
  & \pmatrix{b&t-b\cr s-b&b} \buildrel (\ltm_0)^{t-(a+b)}\over\longrightarrow
  & \pmatrix{t-a&a\cr s-b&b}
  & =\hr_0(\pi_{\rm{r}}(M))
  .\cr
  }}
\nn
$$
Clearly $\hr_0(\pi_{\rm{r}}(M))=\pi_{\rm{r}}(\hr_0(M))$ in this case. In the
remaining case $t<a+b$ one has, with $a'=t-b$ and $b'=s-a$, that $a'+b'<s$ and
$M=\Mtt(s-b',b',t-a',a')$, so in fact $M=\hr_0(M')$ for a matrix
$M'=\Mtt(a',s-a',b',t-b')$ to which the previous case applies, which settles
this final case.

Now we return to the more general case of $\hr_0$ operating on a $n\times2$
matrix~$M$ (either binary or integral). Write~$M={S\choose{T}}$ where $S$ is
$d\times2$ matrix and $T$ a $(n-d)\times2$ matrix, so that
$\pi_{\rm{r}}(M)={T\choose{S}}$. Any horizontal crystal operation in~$M$
occurring inside $S$ or~$T$ will also be valid in that part considered in
isolation. But definition~\binmovedef\defitemnr2 or \intmovedef\defitemnr2
implies that if conversely $S\to{S'}$ is a horizontal crystal operation that
is valid in isolation, then ${S\choose{T}}\to{S'\choose{T}}$ is valid in~$M$
only if $\min\(\nrm_0(S),\nrm_0(S')\)\geq\nlm_0(T)$, while such an operation
$T\to{T'}$ gives rise to one ${S\choose{T}}\to{S\choose{T'}}$ only if
$\nrm_0(S)\leq\min\(\nlm_0(T),\nlm_0(T')\)$. This means that the moves or
transfers involved in the transition $M\to\hr_0(M)$ will be distributed
between the parts $S$ and~$T$ in the same way as for
$\bar{M}=\Mtt(a,s-a,b,t-b)\in\Mat$ the transfers involved in the transition
$\bar{M}\to\hr_0(\bar{M})$ are distributed between the two rows, where
$a=\nrm_0(S)$, \ $s=a+\nlm_0(S)$, \ $b=\nrm_0(T)$, and $t=b+\nlm_0(T)$. So if
$\hr_0(\bar{M})=\Mtt(a+\delta_0,s-a-\delta_0,b+\delta_1,t-b-\delta_1)$, then
one has $\hr_0(M)=\smash{S'\choose{T'}}$, where either
$S'=(\ltm_0)^{\delta_0}(S)$ and $T'=(\ltm_0)^{\delta_1}(T)$ with
$\delta_0,\delta_1\geq0$, or $S'=(\rtm_0)^{-\delta_0}(S)$ and
$T'=(\rtm_0)^{-\delta_1}(T)$ with $\delta_0,\delta_1\leq0$. We have seen that
in all cases one then also has $\hr_0\Mtt(b,t-b,a,s-a)
=\Mtt(b+\delta_1,t-b-\delta_1,a+\delta_0,s-a-\delta_0)$, and reasoning for
$\pi_{\rm{r}}(M)={T\choose{S}}$ in the same way as for~$M$, one finds
$\hr_0{T\choose{S}}=\smash{T'\choose{S'}}$ with the same values~$S',T'$, which
proves $\hr_0(\pi_{\rm{r}}(M))=\pi_{\rm{r}}(\hr_0(M))$.
\QED

We note that in the basic case of $2\times2$ integral matrices, the identity
$\hr_0(\pi_{\rm{r}}(M))=\pi_{\rm{r}}(\hr_0(M))$ comes about in a somewhat
curious way when $a+b<s$ or $t<a+b$, since the intermediate matrices traversed
during the computation of $\hr_0$ in both members no not correspond to each
other under~$\pi_{\rm{r}}$. The situation can be illustrated by drawing for
fixed row sums~$(s,t)$ individual matrices as dots with coordinates given by
their first column $(a,b)\in\set{s+1}\times\set{t+1}$, and drawing both the
ladders for $\ltm_0$ on one hand, and ladders for the conjugate operation
$\pi_{\rm{r}}\after\ltm_0\after\pi_{\rm{r}}$ on the other hand. This is done
below for $(s,t)=(5,9)$, with the first type of ladders in red on the left and
the conjugated ladders in blue on the right. We have also marked three pairs
of matrices that are interchanged by $\hr_0$ (as well as by
$\pi_{\rm{r}}\after\hr_0\after\pi_{\rm{r}}$). One sees that outside the range
given by $s\leq{a}+b\leq{t}$ (delimited by the diagonal lines) the pairs of
matrices exchanged form the \emph{only} common points of a pair of ladders of
the two types.
\bigdisplay
  \epsfxsize=0.4\hsize  \epsfbox{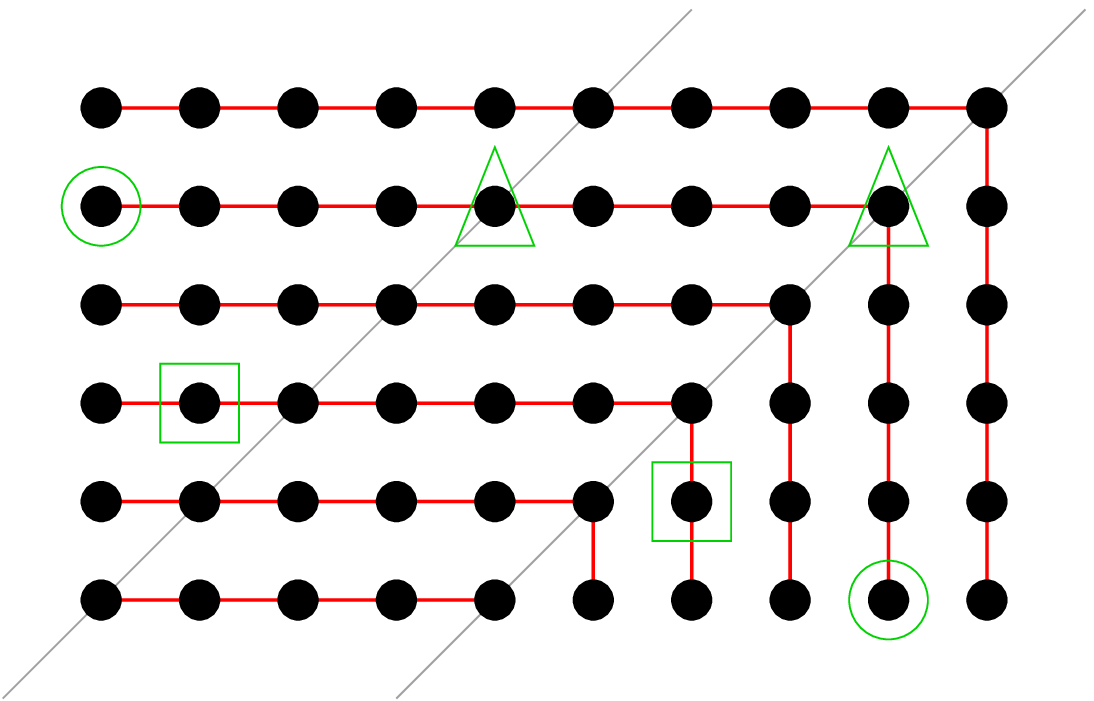}
\quad
  \epsfxsize=0.4\hsize  \epsfbox{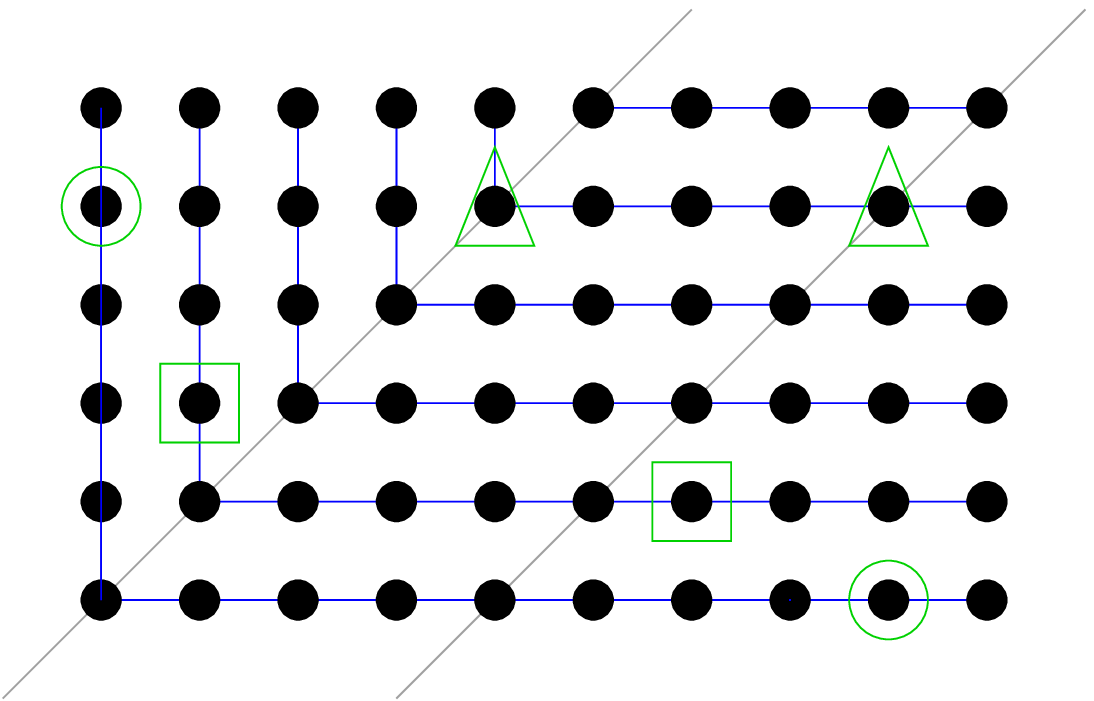}
$$
As a consequence, the simplification in the proof of lemma~\actionlem\ that we
obtained by cyclicly permuting rows, while commuting with the
operations~$\hr_j$, does not preserve the structure of the crystal. Indeed in
the case of one-line matrices ultimately reduced to, horizontal crystal
operations between different pairs of columns commute among each other in the
weak sense that $\ltm_1(\ltm_0(M))=\ltm_0(\ltm_1(M))$ whenever both members
are defined; such a relation does not hold in general crystals.

It is therefore instructive to present as an alternative a direct way to
complete the proof of lemma~\actionlem, without using lemma~\cycliccommlem:
thus one will only have used reductions that preserve the structure of the
crystal. This also provides a less indirect proof for the binary case: we used
the isomorphism of the binary and integral crystals provided by
proposition~\bininteqprop, but the operation of cyclicly permuting rows of an
integral matrix is very hard to translate across this isomorphism (since it
breaks out of the class of matrices corresponding to tableaux of a given
shape). By contrast, the direct computations below, while given for integral
matrices, translate easily into a computation for binary matrices. It suffices
to observe that the matrix~$M$ of (\Mabcdef) corresponds, by the intermediary
of a semistandard Young tableau with three horizontal strips, to a binary
matrix~$B$ with three rows, of which the columns~$B\tr_j$ can be described by
associating to them three bit binary numbers
$c_j=4B_{0,j}+2B_{1,j}+B_{2,j}\in\{0,1,\ldots,7\}$, which then form a weakly
decreasing sequence, and whose frequencies
$m_i=\Card\setof{j\in\N}:c_j=i\endset$ can be deduced from the equations
$m_7+m_6+m_5+m_4=a$, $m_3+m_2=b$, $m_1=c$, $m_7+m_6=d$, $m_5+m_3=e$, $m_7=f$,
and $m_4m_3=0$ (which holds since columns $B\tr_j=(1,0,0)$ and
$B\tr_{j'}=(0,1,1)$ cannot both be present).

Before we do the computation, we make one more simplification: from~$M$ we can
subtract $f$ times the $3\times3$ identity matrix, since $f$ is the smallest
diagonal entry, and the operation is easily seen to commute with the
operations $\ltm_0$, $\rtm_0$ $\ltm_1$, and~$\rtm_1$(in terms of binary
matrices this corresponds to removing the $m_7$ initial columns
$B\tr_j=(1,1,1)$, which play a completely passive role). Thus we can assume
$f=0$, and write $M=\({a\atop0}~{b\atop{d}}~{c\atop{e}}\)\in\Mat$.

Since $a\geq{d}$, all units of the coefficient $M_{1,1}=d$ in~$M$ are fixed
for the operations $\ltm_0$, whence
$$\hr_0\pmatrix{a&b&c\cr0&d&e}=\pmatrix{b+d& a-d& c\cr0&d&e\cr}
.\nn
$$
Computing an expression for $\hr_1(M)$, which will be of the form
$\({a\atop0}~{b'\atop{d'}}~{c'\atop{e'}}\)$, is a bit more complicated.
Assume first that $\min(b,e)=0$, then one will also have $\min(b',e')=0$; from
the fact that moreover $\smash{\Mtt(b',c',d',e')}$ must have rows sums
$(b+c,d+e)$ and column sums $(c+e,b+d)$, one can deduce the expression
$c'=\min(b+c,b+d)=b+\min(c,d)$ for this case. We can reduce the general case
to this one by removing $\min(b,e)$ units both from $b$ and from~$e$, which
should leave $c'$ unchanged; it follows that one must replace $b$ in the
expression for~$c'$ by $b-\min(b,e)$ for the general case. One then finds
$c'=b-\min(b,e)+\min(c,d)$, $b'=c+\min(b,e)-\min(c,d)$,
$d'=e-\min(b,e)+\min(c,d)$ and $e'=d+\min(b,e)-\min(c,d)$, in other words
$$
  \hr_1\pmatrix{a&b&c\cr0&d&e}
  =\pmatrix{a&c+\min(b,e)-\min(c,d)&b-\min(b,e)+\min(c,d)\cr
            0&e-\min(b,e)+\min(c,d)&d+\min(b,e)-\min(c,d)}
.\label(\hroneeq)
$$
Now the verification that $\hr_0\hr_1\hr_0(M)=\hr_1\hr_0\hr_1(M)$ is a
straightforward but laborious computation: one determines the appropriate
expressions in $a,b,c,d,e$, which are piecewise linear ones built using the
operators `$+$', `$-$', and~`$\min$', and one shows that the
corresponding parameters of $\hr_0\hr_1\hr_0(M)$ and $\hr_1\hr_0\hr_1(M)$ are
given by equivalent expressions. The relevant expressions are given in
table~1.
\midinsert
$$\abovedisplayskip=0pt
\vcenter{
\catcode`@=11
\def\rmatrix#1{\left(\null\,\vcenter{\normalbaselines\m@th
    \ialign{\hfil$##$\hfil&&\qquad\hfil$##$\crcr
      \mathstrut\crcr\noalign{\kern-\baselineskip}
      #1\crcr\mathstrut\crcr\noalign{\kern-\baselineskip}}}\,\right)}
\def\\{\cr\noalign{\hrule}height 19pt depth 16 pt&}
\halign{\vrule#& \ \hfil$#$\hfil\ &
                 \vrule#& \ \hfil$#$\hfil\ & \vrule#
\span\\
             M && \pmatrix{a& b& c\cr0&d&e\cr}
&\\
        \hr_0(M)&& \pmatrix{b+d& a-d& c\cr0&d&e\cr}
&\\ \hr_1\hr_0(M)&&
    \rmatrix{b+d&c+\min(a-d,e)-\min(c,d)&a-d-\min(a-d,e)+\min(c,d)\cr
              0 &e-\min(a-d,e)+\min(c,d)&  d+\min(a-d,e)-\min(c,d)\cr}
&\\ \hr_0\hr_1\hr_0(M)&&
    \rmatrix{c+e& b+\min(a-e,d)-\min(c,d)&a-d-\min(a-d,e)+\min(c,d)\cr
              0 & e-\min(a-d,e)+\min(c,d)&  d+\min(a-d,e)-\min(c,d)\cr}
&\\ \hr_1(M)&& \rmatrix{a&c+\min(b,e)-\min(c,d)&b-\min(b,e)+\min(c,d)\cr
                        0&e-\min(b,e)+\min(c,d)&d+\min(b,e)-\min(c,d)\cr}
&\\ \hr_0\hr_1(M)&&
    \rmatrix{c+e& a-e+\min(b,e)-\min(c,d)& b-\min(b,e)+\min(c,d)\cr
              0 &   e-\min(b,e)+\min(c,d)& d+\min(b,e)-\min(c,d)\cr}
&\\ \hr_1\hr_0\hr_1(M)&&
    \rmatrix{c+e& b+\min(a-e,d)-\min(c,d)& a-e-\min(a-e,d)+\min(c,d)\cr
              0 & d-\min(a-e,d)+\min(c,d)&   e+\min(a-e,d)-\min(c,d)\cr}
&\cr\noalign{\hrule}}}
$$
\centerpar{{\bf Table~1.}\quad Verification of
 $\hr_0\hr_1\hr_0(M)=\hr_1\hr_0\hr_1(M)$ for $M=\pmatrix{a&b&c\cr0&d&e\cr}$.}
\endinsert
Since all row and column sums are known beforehand, only one coefficient not
in the leftmost column really needs to be computed at each step, but we have
retained for each coefficient a form close to the one obtained from
application of the defining formula, possibly after some simplification. The
most important simplifications take place when (\hroneeq) is applied to
compute $\hr_1(\hr_0\hr_1(M))$; for instance the coefficient at index~$(0,1)$
is computed as $c'+\min(b',e')-\min(c',d')$ where
$\({a'\atop0}~{b'\atop{d'}}~{c'\atop{e'}}\)=\hr_0\hr_1(M)$, giving
$$ \eqalign
{ \rlap{$b-\min(b,e)+\min(c,d)$}
 \cr\rlap{\quad${}+\min\(a-e+\min(b,e)-\min(c,d)~,~d+\min(b,e)-\min(c,d)\)$}
 \cr\rlap{\quad${}-\min\(b-\min(b,e)+\min(c,d)~,~e-\min(b,e)+\min(c,d)\)$}
\cr&=b-\min(b,e)+\min(c,d)
\cr&\phantom=\quad+\(\min(b,e)-\min(c,d)\)+\min(a-e,d)
           +\(\min(b,e)-\min(c,d)\)-\min(b,e)
\cr&=b+\min\(a-e,d\)-\min(c,d)
.\cr}
$$
Since the result is identical to the corresponding coefficient of
$\hr_0\hr_1\hr_0(M)$, all other coefficients must be equivalent as well, which
is also easily checked explicitly. This computation demonstrates the potential
complications in composing piecewise linear expressions and deciding their
equivalence; here the exercise remained doable because there is effectively
just one component defined in a piecewise manner, with the other components
depending linearly on it. The general situation in the statement of
lemma~\actionlem\ for integral matrices is much more complicated than that
of~(\Mabcdef), so exhausting operations~$\upm_k$ during the initial
simplification, while preserving the structure of the crystal, does greatly
simplify the problem viewed as an equivalence of piecewise linear expressions.

There is one more proof that we have omitted, namely that of
theorem~\binarygrowththm, which states that for binary matrices the
computation of implicit shapes using French normalisation leads to the
row-insertion shape datum for $(\Part,\lev,\leh,\set2)$. We did suggest that
the assiduous reader could check this, but we did not prove her existence.
So here is a proof; those who already found one may skip this.

\proofof{theorem~\binarygrowththm}
We return to the French normalisation process described in the
display~(\Frenchgrowth), and analyse the exhaustion of
$\setof\dnm_i:i\in\set{k}\endset$ applied to~$A$, on the left. To the left of
the vertical line, the submatrix $A_{\set{k+1}\times\set{l}}$, consisting of
$\frdiag\\k$ with below it $\nu\tr-\\\tr$ in row~$k$, is transformed into
$\frdiag\nu{k+1}$, moving down one place every unit in a column~$j$ with
$\nu\tr_j-\\\tr_j=0$. The exhaustion of the~$\dnm_i$ proceeds by weakly
decreasing~$i$, and after operations $\dnm_{i'}$ with $i<i'<k$ have been
exhausted, row~$i+1$ will contain a unit in a column~$j<l$ if and only if
$\nu\tr_j-\\\tr_j=1$ and $\nu\tr_j\geq k-i$. Such a unit could block a
potential downward move in column~$l$ if in addition the bit directly above it
is $A_{i,j}=0$; such a pair $\smash{0\choose1}$ will exist if $\nu\tr_j=k-i$
and $\\\tr_j=k-i-1$, in other words when $(k-i-1,j)\in\set{\nu/\\}$. Indeed
the existence of at least one such pair will rule out any downward move
between rows $i,i+1$ in column~$l$, since there can be no pairs
$\smash{1\choose0}$ to the left of column~$l$ that could be blocked instead.
In the example of~(\Frenchgrowth) we see this phenomenon for $i=1$: since
$\set{\nu/\\}$ has squares in row~$k-i-1=3$ (namely in columns $0$ and~$1$),
the unit at~$(6,1)$ is blocked, although the unit at~$(6,2)$ has been moved
down by~$\dnm_2$ before $\dnm_1$ is applied.

We now consider column~$l$ to see how its initial value of $\revert{\mu-\\}k$
followed by a bit $A_{k,l}=M_{k,l}$ is transformed into
$\revert{\kappa-\nu}{k+1}$, thereby determining~$\kappa$. If after
exhausting~$\dnm_{k-i}$, the bit at~$(k-i,l)$ is~`$1$', then we shall say
there is a carry into row~$i$; initially there is a carry into row~$0$ if and
only if $M_{k,l}=1$. Now if there is a carry into row~$i$, then
$\set{\kappa/\nu}$ will certainly have a square in row~$i$, and if moreover
$\set{\mu/\\}$ already had a square in row~$i$, there will be a carry into
row~$i+1$. If there is no carry into row~$i$ and $\set{\mu/\\}$ has no square
in row~$i$, then neither will $\set{\kappa/\nu}$, nor will there be a carry
into row~$i+1$ (here column~$l$ contains~$\smash{0\choose0}$). If there is no
carry into row~$i$ but $\set{\mu/\\}$ does have square in row~$i$ (there is a
pair $\smash{1\choose0}$ in column~$l$), then that square is also in
$\set{\kappa/\nu}$ if $\nu_i=\\_i$ (the unit in column~$l$ moves down), or
else it gives rise to a carry into row~$i+1$ (if $\set{\nu/\\}$ has any
squares in row~$i$, so that the pair in column~$l$ is blocked).

Now comparing with the description of the row-insertion shape datum for
$(\Part,\lev,\leh,\set2)$ given before theorem~\binarygrowththm, any square
$(i,j)\in\set{\mu/\\}\thru\set{\nu/\\}$ is an element of~$S$ that is absent
from~$\set\\$, and it assures a carry into row~$i+1$; if there was also a
carry into row~$i$, then the square of~$\set{\kappa/\nu}$ in row~$i$ is some
$(i,j')\in{T}$ with $j'>j$. If there is a carry into row~$i$ but no square
of~$\set{\mu/\\}$, then $(i,\nu_i)\in{T}$ is a square of~$\set{\kappa/\nu}$,
and conversely a square~$(i,j)\in{T}$ can only lie in~$\set\kappa$ if there is
a carry into row~$i$, because if there is not, any square of
$\set{\kappa/\nu}$ in row~$i$ is also in~$\set{\mu/\\}$ and therefore not
in~$T$. Thus $\kappa$ is as given by the shape datum.
\QED

Alternatively one could determine $\kappa$ by analysing what happens during
the French normalisation of $A_{\set{k}\times\set{l+1}}$. This is interesting
and proceeds rather differently than in the given proof, but this time we
really leave this for the enthusiasts to discover for themselves.

 \par

\Supplement Conclusion.

After introducing crystal operations on binary and integral matrices, we have
performed constructions with them related to various known facets of the
\caps{RSK}~correspondence. We have found links with all of the facets
considered, and in many cases obvious questions suggested by the results found
have led us to investigate in detail the precise relations with the existing
theory. This required some meticulous work, but the result has shed new light
on many aspects of that theory, and has resulted in a framework in which those
aspects are more tightly and explicitly integrated with each other than they
were before. In fact we have been surprised by the extraordinary versatility
of these crystal operations, whose original humble vocation was to streamline
the proof of a version of the Littlewood-Richardson rule found in
\ref{alternating Schur}, and by the dimensions this paper has taken as a
consequence. To close, we mention one more obvious question that we shall not
answer here, namely whether for the notions of crystals in other types than
$A_n$ there exists a similar notion of double crystals.

\Supplement References.

\input \jobname.ref

\iftoc
  \closeout\tocfile \vfil\eject
  \def\tocitem#1=#2\onpage#3.
    {\par\line{\hbox to 1.5pc{\bf #1\hss}\quad#2\dotfill#3}\ignorespaces}
  \catcode`@=11
  \input \jobname.toc
\fi

\bye